\documentclass{article}
\usepackage{amsfonts,amssymb,enumerate}
\usepackage{pstricks}
\usepackage{pst-grad}
\usepackage{pst-plot}
\usepackage[tiling]{pst-fill}

\usepackage{pstcol}

\usepackage{graphics}
\usepackage[all, knot]{xy}
 \xyoption{arc}
\usepackage{fancyheadings}
\newcommand{\cat}[1]{\ensuremath{\mbox{\bfseries {\upshape {#1}}}}}
\newcommand{\numroman}{\renewcommand{\labelenumi}{\roman{enumi})}}
\newcommand{\numarabic}{\renewcommand{\labelenumi}{\arabic{enumi})}}

\newcommand{\et}{\hspace{-0.08in}{\bf .}\hspace{0.1in}}
\newcommand{\BOX}{\hbox {$\sqcap$ \kern -1em $\sqcup$}}
\newcommand{\qed}{\hskip 2em \hbox{\BOX} \vskip 2ex}

\newcommand{\To}{\Rightarrow}

\newcommand{\Hom}{{\rm Hom}}

\renewcommand{\to}{\rightarrow}
\newcommand{\maps}{\colon}

\newcommand{\Adj}{{\rm Adj}}
\newcommand{\Mon}{{\rm Mon}}
\newcommand{\Comon}{{\rm Comon}}
\newcommand{\Ambi}{{\rm Ambi}}
\newcommand{\Frob}{{\rm Frob}}
\newcommand{\wAmbi}{{\rm pAmbi}}
\newcommand{\wFrob}{{\rm pFrob}}
\newcommand{\wAdj}{{\rm pAdj}}
\newcommand{\wMon}{{\rm pMon}}

\newcommand{\Proof}{\noindent {\bf Proof. }}

\newcommand{\R}{{\mathbb R}}
\newcommand{\Z}{{\mathbb Z}}
\newcommand{\T}{{\mathbb T}}

\newcommand{\F}{{\mathbb F}}

\newcommand{\scs}{\scriptstyle}
\newcommand{\ten}{\otimes }
\newcommand{\INT}{\bigotimes}

\newtheorem{thm}{Theorem}    
\newtheorem{cor}[thm]{Corollary}

\newtheorem{lem}[thm]{Lemma}

\newtheorem{prop}[thm]{Proposition}
\newtheorem{defn}[thm]{Definition}

        \newcommand{\be}{\begin{equation}}
        \newcommand{\ee}{\end{equation}}
        \newcommand{\ba}{\begin{eqnarray}}
        \newcommand{\ea}{\end{eqnarray}}
        \newcommand{\ban}{\begin{eqnarray*}}
        \newcommand{\ean}{\end{eqnarray*}}
        \newcommand{\barr}{\begin{array}}
        \newcommand{\earr}{\end{array}}

\newcommand{\adjunction}[4]{%
  \ensuremath{\xymatrix{%
    {#1} \ar@<4pt>[r]^{{#3}} \ar@{}[r]|{\scriptscriptstyle{\bot}} &%
    {#2} \ar@<4pt>[l]^{{#4}}%
  }}%
}

\newcommand{\radjunction}[4]{%
  \ensuremath{\xymatrix{%
    {#1} \ar@<4pt>[r]^{{#3}} \ar@{}[r]|{\scriptscriptstyle{\top}} &%
    {#2} \ar@<4pt>[l]^{{#4}}%
  }}%
}

\newcommand{\ambijunction}[4]{%
  \ensuremath{\xymatrix{%
    {#1} \ar@<4pt>[r]^{{#3}} \ar@{}[r]|{\scriptscriptstyle{\top \;\bot}} &%
    {#2} \ar@<4pt>[l]^{{#4}}%
  }}%
}

\newcommand{\padjunction}[4]{%
  \ensuremath{\xymatrix{%
    {#1} \ar@<4pt>[r]^{{#3}} \ar@{}[r]|{\scriptscriptstyle{\bot_p}} &%
    {#2} \ar@<4pt>[l]^{{#4}}%
  }}%
}

\newcommand{\pradjunction}[4]{%
  \ensuremath{\xymatrix{%
    {#1} \ar@<4pt>[r]^{{#3}} \ar@{}[r]|{\scriptscriptstyle{\top_p}} &%
    {#2} \ar@<4pt>[l]^{{#4}}%
  }}%
}

\newcommand{\pseudoambijunction}[4]{%
  \ensuremath{\xymatrix{%
    {#1} \ar@<4pt>[r]^{{#3}} \ar@{}[r]|{\scriptscriptstyle{\top_p \;\bot_p}} &%
    {#2} \ar@<4pt>[l]^{{#4}}%
  }}%
}

\SelectTips{cm}{}

\psset{linewidth=0.3pt,dimen=middle}
\psset{xunit=.70cm,yunit=0.70cm}

\newcommand{\mult}{
  \pscustom[fillcolor=lightgray, fillstyle=solid]{
        \psbezier(.9,1.5)(.9,.6)(.25,.7)(.3,0)
        \psline(-0.3,0)
        \psbezier(-0.3,0)(-.25,.7)(-.9,.6)(-.9,1.5)
        \psline(-.3,1.5)
        \psbezier(-.3,1.5)(-.4,.8)(0.4,.8)(.3,1.5)
        \psline(.9,1.5)
    }
}

\newcommand{\comult}{
  \pscustom[fillcolor=lightgray, fillstyle=solid]{
        \psbezier(.9,0)(.9,.9)(.25,.7)(.3,1.5)
        \psline(-0.3,1.5)
        \psbezier(-0.3,1.5)(-.25,.7)(-.9,.9)(-.9,0)
        \psline(-.3,0)
        \psbezier(-.3,0)(-.4,.7)(0.4,.7)(.3,0)
        \psline(.9,0)
    }
}

 \newcommand{\birth}{
 \pscustom[fillcolor=lightgray, fillstyle=solid]{
        \psbezier(-.3,0)(-.3,.7)(0.3,.7)(.3,0)
        \psline(-.3,0)
    }
 }

  \newcommand{\death}{
 \pscustom[fillcolor=lightgray, fillstyle=solid]{
        \psbezier(-.3,0)(-.3,-.7)(0.3,-.7)(.3,0)
        \psline(-.3,0)    }
 }

\newcommand{\zag}{
   \pscustom[fillcolor=lightgray, fillstyle=solid]{
        \psbezier(.9,0)(1,1.4)(-1,1.4)(-.9,0)
        \psline(-.3,0)
        \psbezier(-.3,0)(-.4,.7)(0.4,.7)(.3,0)
        \psline(.9,0)
    }
}

\newcommand{\longzag}{
   \pscustom[fillcolor=lightgray, fillstyle=solid]{
        \psbezier(1.5,0)(1.6,1.6)(-1.6,1.6)(-1.5,0)
        \psline(-.9,0)
        \psbezier(-.9,0)(-1,.9)(1,.9)(.9,0)
        \psline(1.5,0)
    }
}

\newcommand{\zig}{
   \pscustom[fillcolor=lightgray, fillstyle=solid]{
        \psbezier(.9,0)(1,-1.4)(-1,-1.4)(-.9,0)
        \psline(-.3,0)
        \psbezier(-.3,0)(-.4,-.7)(0.4,-.7)(.3,0)
        \psline(.9,0)
    }
}

\newcommand{\smallident}{
    \pspolygon[fillcolor=lightgray,fillstyle=solid](-.3,0)(.3,0)(.3,.55)(-.3,.55)(-.3,0)
}

\newcommand{\medident}{
    \pspolygon[fillcolor=lightgray,fillstyle=solid](-.3,0)(.3,0)(.3,.75)(-.3,.75)(-.3,0)
}

\newcommand{\ident}{
    \pspolygon[fillcolor=lightgray,fillstyle=solid](-.3,0)(.3,0)(.3,1.5)(-.3,1.5)(-.3,0)
}

\newcommand{\longident}{
    \pspolygon[fillcolor=lightgray,fillstyle=solid](-.3,0)(.3,0)(.3,2.05)(-.3,2.05)(-.3,0)
}

\newcommand{\curveright}{
  \pscustom[fillcolor=lightgray, fillstyle=solid]{
        \psbezier(.9,1.5)(.8,.6)(.25,.9)(.3,0)
        \psline(-0.3,0)
        \psbezier(-0.3,0)(-.35,.9)(.2,.6)(.3,1.5)
        \psline(.9,1.5)
    }
}

\newcommand{\curveleft}{
  \pscustom[fillcolor=lightgray, fillstyle=solid]{
        \psbezier(-.9,1.5)(-.8,.6)(-.25,.9)(-.3,0)
        \psline(0.3,0)
        \psbezier(0.3,0)(.35,.9)(-.2,.6)(-.3,1.5)
        \psline(-.9,1.5)
    }
}

\newcommand{\forkQ}{
  \pscustom[fillcolor=lightgray, fillstyle=solid, linecolor=lightgray]{
    \psbezier(1.5,1.5)(1.5,.6)(.25,.6)(.3,0)
    \psline(-.3,0)
    \psbezier(-0.3,0)(-.25,.6)(-1.5,.6)(-1.5,1.5)
    \psline(-.9,1.5)
    \psbezier(-.9,1.5)(-.9,.85)(-.3,.85)(-.3,1.5)
    \psline(.3,1.5)
    \psbezier(.3,1.5)(.3,.85)(.9,.85)(.9,1.5)
    }
    \psbezier(1.5,1.5)(1.5,.6)(.25,.6)(.3,0)
    \psbezier(-0.3,0)(-.25,.6)(-1.5,.6)(-1.5,1.5)
    \psbezier(.3,1.5)(.3,.85)(.9,.85)(.9,1.5)
    \psbezier(-.3,1.5)(-.3,.85)(-.9,.85)(-.9,1.5)
  }

\newcommand{\multQ}{
          \pscustom[fillcolor=lightgray, fillstyle=solid, linecolor=lightgray]{
        \psbezier(.9,1.5)(.9,.6)(.25,.7)(.3,0)
        \psline(-0.3,0)
        \psbezier(-0.3,0)(-.25,.7)(-.9,.6)(-.9,1.5)
        \psline(-.3,1.5)
        \psbezier(-.3,1.5)(-.4,.8)(0.4,.8)(.3,1.5)
        \psline(.9,1.5)
    }
    \psbezier(.9,1.5)(.9,.6)(.25,.7)(.3,0)
    \psbezier(-0.3,0)(-.25,.7)(-.9,.6)(-.9,1.5)
    \psbezier(-.3,1.5)(-.4,.8)(0.4,.8)(.3,1.5)
}

\newcommand{\comultQ}{
  \pscustom[fillcolor=lightgray, fillstyle=solid, linecolor=lightgray]{
        \psbezier(.9,0)(.9,.9)(.25,.7)(.3,1.5)
        \psline(-0.3,1.5)
        \psbezier(-0.3,1.5)(-.25,.7)(-.9,.9)(-.9,0)
        \psline(-.3,0)
        \psbezier(-.3,0)(-.4,.7)(0.4,.7)(.3,0)
        \psline(.9,0)
    }
     \psbezier(.9,0)(.9,.9)(.25,.7)(.3,1.5)
        \psbezier(-0.3,1.5)(-.25,.7)(-.9,.9)(-.9,0)
        \psbezier(-.3,0)(-.4,.7)(0.4,.7)(.3,0)
}

 \newcommand{\birthQ}{
 \pscustom[fillcolor=lightgray, fillstyle=solid, linecolor=lightgray]{
        \psbezier(-.3,0)(-.3,.7)(0.3,.7)(.3,0)
        \psline(-.3,0)
    }
    \psbezier(-.3,0)(-.3,.7)(0.3,.7)(.3,0)
 }

  \newcommand{\deathQ}{
 \pscustom[fillcolor=lightgray, fillstyle=solid, linecolor=lightgray]{
        \psbezier(-.3,0)(-.3,-.5)(0.3,-.5)(.3,0)
        \psline(-.3,0)    }
        \psbezier(-.3,0)(-.3,-.5)(0.3,-.5)(.3,0)
 }

\newcommand{\zagQ}{
   \pscustom[fillcolor=lightgray, fillstyle=solid, linecolor=lightgray]{
        \psbezier(.9,0)(1,1.4)(-1,1.4)(-.9,0)
        \psline(-.3,0)
        \psbezier(-.3,0)(-.4,.7)(0.4,.7)(.3,0)
        \psline(.9,0)
    }
        \psbezier(.9,0)(1,1.4)(-1,1.4)(-.9,0)
        \psbezier(-.3,0)(-.4,.7)(0.4,.7)(.3,0)
}

\newcommand{\smallidentQ}{
    \pspolygon[fillcolor=lightgray, fillstyle=solid, linecolor=lightgray](-.3,0)(.3,0)(.3,.55)(-.3,.55)(-.3,0)
    \psline(-.3,0)(-.3,.55)
    \psline(.3,0)(.3,.55)
}

\newcommand{\medidentQ}{
    \pspolygon[fillcolor=lightgray, fillstyle=solid, linecolor=lightgray](-.3,0)(.3,0)(.3,.75)(-.3,.75)(-.3,0)
    \psline(-.3,0)(-.3,.75)
    \psline(.3,0)(.3,.75)
}

\newcommand{\identQ}{
    \pspolygon[fillcolor=lightgray, fillstyle=solid, linecolor=lightgray](-.3,0)(.3,0)(.3,1.5)(-.3,1.5)(-.3,0)
    \psline(-.3,0)(-.3,1.5)
    \psline(.3,0)(.3,1.5)
}

\newcommand{\curverightQ}{
  \pscustom[fillcolor=lightgray, fillstyle=solid, linecolor=lightgray]{
        \psbezier(.9,1.5)(.8,.6)(.25,.9)(.3,0)
        \psline(-0.3,0)
        \psbezier(-0.3,0)(-.35,.9)(.2,.6)(.3,1.5)
        \psline(.9,1.5)
    }
    \psbezier(.9,1.5)(.8,.6)(.25,.9)(.3,0)
        \psbezier(-0.3,0)(-.35,.9)(.2,.6)(.3,1.5)
}

\newcommand{\curveleftQ}{
  \pscustom[fillcolor=lightgray, fillstyle=solid, linecolor=lightgray]{
        \psbezier(-.9,1.5)(-.8,.6)(-.25,.9)(-.3,0)
        \psline(0.3,0)
        \psbezier(0.3,0)(.35,.9)(-.2,.6)(-.3,1.5)
        \psline(-.9,1.5)
    }
    \psbezier(-.9,1.5)(-.8,.6)(-.25,.9)(-.3,0)
        \psbezier(0.3,0)(.35,.9)(-.2,.6)(-.3,1.5)
}

\newcommand{\multD}{
  \pscustom[fillcolor=gray, fillstyle=solid]{
        \psbezier(.9,1.5)(.9,.6)(.25,.7)(.3,0)
        \psline(-0.3,0)
        \psbezier(-0.3,0)(-.25,.7)(-.9,.6)(-.9,1.5)
        \psline(-.3,1.5)
        \psbezier(-.3,1.5)(-.4,.8)(0.4,.8)(.3,1.5)
        \psline(.9,1.5)
    }
}

  \newcommand{\deathD}{
 \pscustom[fillcolor=gray, fillstyle=solid]{
        \psbezier(-.3,0)(-.3,-.7)(0.3,-.7)(.3,0)
        \psline(-.3,0)    }
 }

\newcommand{\zagD}{
   \pscustom[fillcolor=gray, fillstyle=solid]{
        \psbezier(.9,0)(1,1.4)(-1,1.4)(-.9,0)
        \psline(-.3,0)
        \psbezier(-.3,0)(-.4,.7)(0.4,.7)(.3,0)
        \psline(.9,0)
    }
}

\newcommand{\smallidentD}{
    \pspolygon[fillcolor=gray,fillstyle=solid](-.3,0)(.3,0)(.3,.55)(-.3,.55)(-.3,0)
}

\newcommand{\medidentD}{
    \pspolygon[fillcolor=gray,fillstyle=solid](-.3,0)(.3,0)(.3,.75)(-.3,.75)(-.3,0)
}

\newcommand{\identD}{
    \pspolygon[fillcolor=gray,fillstyle=solid](-.3,0)(.3,0)(.3,1.5)(-.3,1.5)(-.3,0)
}

\newcommand{\twocobb}{
\psset{xunit=.30cm,yunit=0.30cm}
\begin{pspicture}(7.5,7.5)
\pscustom[fillstyle=gradient,
    gradbegin=white, gradend=gray,gradmidpoint=0,gradangle=120]{
            \pscurve(3,7)(3.5,5.5)(4.5,5.5)(5,7)
            \pscurve(5,7)(6,6.7)(6.9,7)
 \psbezier(7,6.9)(6,3.5)(4.9,5)(5,1)
            \psbezier(5,1)(4.6,.35)(3.4,.35)(3,1)
 \psbezier(3,1)(3.1,5)(2,3.5)(1,6.9)
             \pscurve(1,7)(2,6.7)(3,7) }
 \psellipse[fillcolor=lightgray,fillstyle=gradient,
    gradbegin=white, gradend=gray,gradmidpoint=1,gradangle=45](2,7)(1,0.4)
 \psellipse[fillstyle=gradient,
    gradbegin=white, gradend=gray,gradmidpoint=1,gradangle=45](6,7)(1,0.4)
 \begin{psclip}{ \pspolygon[linestyle=none](3,1)(5,1)(5,1.5)(3,1.5)(3,1) }
\psellipse[linestyle=dotted](4,1)(1,0.4)
 \end{psclip}
\end{pspicture}
\quad
 \begin{pspicture}(3,6)
 \pscustom[fillcolor=lightgray,fillstyle=gradient,
    gradbegin=white, gradend=gray,gradmidpoint=0,gradangle=120]{
 \psbezier(1,3.5)(1.4,2.85)(2.6,2.85)(3,3.5)
 \psbezier(3,3.5)(3,5.7)(1,5.7)(1,3.5)
 }
 \begin{psclip}{ \pspolygon[linestyle=none](1,3.5)(1,4.5)(3,4.5)(3,3.5)(1,3.5) }
\psellipse[linestyle=dotted](2,3.5)(1,0.4)
 \end{psclip}
\end{pspicture}
\quad
\begin{pspicture}(8,7.5)
\pscustom[fillcolor=lightgray,fillstyle=gradient,
    gradbegin=white, gradend=gray,gradmidpoint=0,gradangle=120]{
            \pscurve(3,1)(3.5,2.5)(4.5,2.5)(5,1)
            \psbezier(5,1)(5.4,.35)(6.6,.35)(7,1)
 \psbezier(7,1.1)(6,5.5)(4.9,3)(5,7)
            \psline(5,7)(3,7)
 \psbezier(3,7)(3.1,3)(2,5.5)(1,1.1)
             \psbezier(1,1)(1.4,.35)(2.6,.35)(3,1) }
 \psellipse[fillcolor=lightgray,fillstyle=gradient,
    gradbegin=white, gradend=gray,gradmidpoint=1,gradangle=45](4,7)(1,0.4)
 \begin{psclip}{ \pspolygon[linestyle=none](0,1)(0,7)(0,1.5)(8,1.5)(8,1) }
\psellipse[linestyle=dotted](2,1)(1,0.4)\psellipse[linestyle=dotted](6,1)(1,0.4)
 \end{psclip}
\end{pspicture}
 \begin{pspicture}(3,6)
 \pscustom[fillcolor=lightgray,fillstyle=gradient,
    gradbegin=white, gradend=gray,gradmidpoint=0,gradangle=120]{
 \psline(1,5)(3,5)
 \psbezier(3,5)(3,2.8)(1,2.8)(1,5)
 }
 \psellipse[fillcolor=lightgray,fillstyle=gradient,
    gradbegin=white, gradend=gray,gradmidpoint=1,gradangle=45](2,5)(1,0.4)
\end{pspicture}
\quad
\begin{pspicture}(7.5,7.5)
\pscustom[fillcolor=lightgray,fillstyle=gradient,
    gradbegin=white, gradend=gray,gradmidpoint=0,gradangle=125]{
       \psline(1,7)(5,1)
       \psbezier(5,1)(5.4,.35)(6.6,.35)(7,1)
       \psline(3,7)\psline(1,7)
 }
 \pscustom[fillstyle=gradient,
    gradbegin=white, gradend=gray,gradmidpoint=0,gradangle=57]{
       \psline(5,7)(1,1)
       \psbezier(1,1)(1.4,.35)(2.6,.35)(3,1)
       \psline(3,1)\psline(7,7)
 }
 \psellipse[fillcolor=lightgray,fillstyle=gradient,
    gradbegin=white, gradend=gray,gradmidpoint=1,gradangle=45](2,7)(1,0.4)
 \psellipse[fillstyle=gradient,
    gradbegin=white, gradend=gray,gradmidpoint=1,gradangle=45](6,7)(1,0.4)
 \begin{psclip}{ \pspolygon[linestyle=none](1,1)(7,1)(7,1.5)(1,1.5)(1,1) }
\psellipse[linestyle=dotted](2,1)(1,0.4)
\psellipse[linestyle=dotted](6,1)(1,0.4)
 \end{psclip}
 \psline[linestyle=dotted](1,7)(5,1)
 \psline[linestyle=dotted](3,7)(7,1)
\end{pspicture}
}

\newcommand{\twothick}{
 \psset{xunit=.30cm,yunit=0.30cm}
\begin{pspicture}(7.5,7.5)
\pscustom[fillcolor=lightgray,fillstyle=solid]{
            \pscurve(3,7)(3.5,5.5)(4.5,5.5)(5,7)
            \psline(5,7)(7,7)
 \psbezier(7,6.9)(6,3.5)(4.9,5)(5,1)
            \psline(5,1)(3,1)
 \psbezier(3,1)(3.1,5)(2,3.5)(1,6.9)
             \psline(1,7)(3,7) }
\end{pspicture}
\quad
 \begin{pspicture}(3,6)
 \pscustom[fillcolor=lightgray,fillstyle=solid]{
 \psline(1,3.5)(3,3.5)
 \psbezier(3,3.5)(3,5.7)(1,5.7)(1,3.5)
 }
\end{pspicture}
\quad
\begin{pspicture}(8,7.5)
\pscustom[fillcolor=lightgray,fillstyle=solid]{
            \pscurve(3,1)(3.5,2.5)(4.5,2.5)(5,1)
            \psline(5,1)(7,1)
 \psbezier(7,1.1)(6,5.5)(4.9,3)(5,7)
            \psline(5,7)(3,7)
 \psbezier(3,7)(3.1,3)(2,5.5)(1,1.1)
             \psline(1,1)(3,1) }
\end{pspicture}
 \begin{pspicture}(3,6)
 \pscustom[fillcolor=lightgray,fillstyle=solid]{
 \psline(1,5)(3,5)
 \psbezier(3,5)(3,2.8)(1,2.8)(1,5)
 }
\end{pspicture}
}

\newcommand{\threethick}{
\begin{center}
\makebox[0pt]{ $
 \xy
 (0,0)*{
   \begin{pspicture}(4,5.5)
   \pspolygon[fillcolor=lightgray,fillstyle=gradient,
    gradbegin=white, gradend=darkgray,gradmidpoint=.32,gradangle=88](1.9,4.4)(2.4,3.6)(2.4,3)(1.9,4.4)
\pspolygon(1,2)(1,5)(4,5)(4,2)(1,2) \psline(0,0)(1,2)
      \pscustom[fillcolor=lightgray,fillstyle=solid]{
  \psline(.4,3)(.8,3.8)
  \psbezier(.8,3.8)(1.2,4.8)(2.2,4.7)(2.2,3.8)
  \psbezier(2.2,3.8)(2.2,3.6)(2.5,3.6)(2.7,4)
  \psline(3.2,5)
  \psline(3.6,5)
  \psline(3.6,3)
  }
   \pscustom[fillcolor=lightgray,fillstyle=gradient,
    gradbegin=white, gradend=darkgray,gradmidpoint=.32,gradangle=88]{
  \psline(.8,3)(1.2,3.8)
  \psbezier(1.2,3.8)(1.35,4.2)(1.8,4.3)(1.8,3.8)
  \psbezier(1.8,3.8)(1.8,3.2)(2.7,3.1)(3,3.8)
  \psline(3,3.8)(3.6,5)
  \psline(3.6,2)
  \psline(.8,0)
  \psline(.8,3)
  }
  \pspolygon[fillcolor=gray,fillstyle=solid](.8,0)(.4,0)(.4,3)(.8,3)(.8,0)
    \pspolygon(0,0)(0,3)(3,3)(3,0)(0,0)
  \pspolygon(0,3)(3,3)(4,5.0)(1,5.0)(0,3)
  \pspolygon(3,0)(3,3)(4,5.0)(4,2)(3,0)
  \psline[linestyle=dashed](.4,0)(3.2,2)
  \psline[linestyle=dashed](3.2,4)(3.2,2)
  \psline[linestyle=dotted](4,2)(1,2)
  \psline[linestyle=dotted](1,5)(1,2)
  \psline[linestyle=dotted](0,0)(1,2)
 \end{pspicture}};
 \endxy
\qquad \quad \xy
 (0,0)*{
    \begin{pspicture}(4,5.5) 
\pspolygon(1,2)(1,5)(4,5)(4,2)(1,2) \psline(0,0)(1,2)
  \pscustom[fillcolor=lightgray,fillstyle=solid]{
  \psbezier(.8,4)(.9,4.4)(1.5,4.4)(1.4,4)
  \psbezier(1.4,4)(1.3,3.6)(.7,3.6)(.8,4)}
  \pscustom[fillcolor=lightgray,fillstyle=gradient,
    gradbegin=white, gradend=darkgray,gradmidpoint=.2,gradangle=90]{
    \psbezier(1.4,4)(1.3,3.6)(.7,3.6)(.8,4)
    \psline(.8,1)
  \psbezier(.8,1)(.7,.6)(1.3,.6)(1.4,1)
    \psline(1.4,4)}
  \pscustom[fillcolor=lightgray,fillstyle=gradient,
    gradbegin=lightgray, gradend=darkgray,gradmidpoint=.4,gradangle=88]{
  \psbezier(2.3,0)(2.5,.6)(3.1,.4)(3.2,1)
  \psbezier(3.2,1)(3.4,1.5)(3.1,1.5)(3.3,2)
  \psline(3.3,5)
  \psbezier(3.3,5)(3.1,4.5)(3.4,4.5)(3.2,4)
  \psbezier(3.2,4)(3.1,3.4)(2.5,3.6)(2.3,3)
  \psline(2.3,5)
  }
  \pscustom[fillcolor=lightgray,fillstyle=solid]{
  \psbezier(1.9,3)(2.2,3.4)(1.9,3.8)(2,4)
  \psbezier(2,4)(2.1,4.6)(2.8,4.5)(2.9,5)
  \psline(3.3,5)
    \psbezier(3.3,5)(3.1,4.5)(3.4,4.5)(3.2,4)
  \psbezier(3.2,4)(3.1,3.4)(2.5,3.6)(2.3,3)
  \psline(1.9,3)}
  \pscustom[fillcolor=lightgray,fillstyle=gradient,
    gradbegin=white,gradend=black,gradmidpoint=1,gradangle=90]{
  \psbezier(2.4,4)(2.2,3.6)(2.6,3.6)(2.8,4)
  \psbezier(2.8,4)(3,4.4)(2.6,4.4)(2.4,4)}
  \pspolygon[fillcolor=gray,
  fillstyle=solid](1.9,0)(2.3,0)(2.3,3)(1.9,3)(1.9,0)
  \psbezier[linestyle=dashed](.8,1)(.9,1.4)(1.5,1.4)(1.4,1)
  \psbezier[linestyle=dashed](2.4,1)(2.2,.6)(2.6,.6)(2.8,1)
  \psbezier[linestyle=dashed](2.8,1)(3,1.4)(2.6,1.4)(2.4,1)
  \psbezier[linestyle=dashed](1.9,0)(2.2,0.4)(1.9,0.8)(2,1)
  \psbezier[linestyle=dashed](2,1)(2.1,1.6)(2.8,1.5)(2.9,2)
  \psline[linestyle=dashed](2.9,2)(2.9,3.4)
  \psline[linestyle=dotted](4,2)(1,2)
  \psline[linestyle=dotted](1,5)(1,2)
  \psline[linestyle=dotted](0,0)(1,2)
    \pspolygon(0,0)(0,3)(3,3)(3,0)(0,0)
  \pspolygon(0,3)(3,3)(4,5.0)(1,5.0)(0,3)
  \pspolygon(3,0)(3,3)(4,5.0)(4,2)(3,0)
 \end{pspicture}};
 \endxy
 \qquad \quad
 \xy (0,0)*{};
 (0,2)*{
   \begin{pspicture}[.5](4,6)
\pspolygon(2,2)(2,5)(5,5)(5,2)(2,2) \psline(0,0)(2,2)
  \pspolygon(0,3)(3,3)(5,5.0)(2,5.0)(0,3)
\pscustom[fillcolor=lightgray,fillstyle=gradient,
    gradbegin=white, gradend=darkgray,gradmidpoint=.5,gradangle=80]{
    \psline(3.4,5)(3.4,4)
    \psline(3,4)
    \psline(3,5)
    \psline(3.4,5)
 }
   \pscustom[fillcolor=lightgray,fillstyle=solid]{
    \psline(.4,3)(1.3,4)
  \psbezier(1.3,4)(2.2,5)(2.6,4.3)(3,5)
  \psline(3.4,5)
      \psbezier(3.4,5)(2.8,3.6)(4,4.4)(4.2,5)
      \psline(4.6,5)
      \psbezier(4.6,5)(3.8,3.6)(3.6,4.4)(2.8,3.4)
      \psline(2.8,3.4)(2.4,3)
      \psline(.4,3)
  }
  \pscustom[fillcolor=lightgray,fillstyle=gradient,
    gradbegin=white, gradend=darkgray,gradmidpoint=.5,gradangle=80]{
    \psbezier(.8,0)(1.2,.8)(2.6,.8)(2,0)
    \psline(2,3)
      \psline(2,3)(2.8,3.8)
      \psbezier(2.8,3.8)(3.2,4.3)(2.4,4.9)(1.2,3.4)
      \psline(1.2,3.4)(.8,3)
      \psline(.8,0)
  }
  \pscustom[fillcolor=lightgray,fillstyle=gradient,
    gradbegin=white, gradend=darkgray,gradmidpoint=.32,gradangle=88]{
  \psline(2.4,0)(2.6,.2)
  \psbezier(2.6,.2)(2.8,.6)(2.2,.8)(2.8,1.2)
  \psbezier(2.8,1.2)(3.4,1.9)(3.9,1.1)(4.6,2)
  \psline(4.6,5)
  \psbezier(4.6,5)(3.8,3.6)(3.6,4.4)(2.8,3.4)
    \psline(2.8,3.4)(2.4,3)
    \psline(2.4,0)
  }
    \pspolygon[fillcolor=gray,fillstyle=solid](.8,0)(.4,0)(.4,3)(.8,3)(.8,0)
  \pspolygon[fillcolor=gray,fillstyle=solid](2.4,0)(2,0)(2,3)(2.4,3)(2.4,0)
  \psbezier[linestyle=dashed](4.2,2)(3.9,1.6)(2.6,1.6)(3.4,2)
  \psline[linestyle=dashed](.4,0)(1,.6)
  \psbezier[linestyle=dashed](1,.6)(2,1.2)(2.4,1)(3,2)
  \psline[linestyle=dashed](3,2)(3,3.4)
  \psline[linestyle=dashed](3.4,2)(3.4,3.7)
  \psline[linestyle=dashed](4.2,2)(4.2,4.3)
  \psline[linestyle=dotted](5,2)(2,2)
  \psline[linestyle=dotted](0,0)(2,2)
  \psline[linestyle=dotted](2,5)(2,2)
    \pspolygon(0,0)(0,3)(3,3)(3,0)(0,0)
    \pspolygon(3,0)(3,3)(5,5.0)(5,2)(3,0)
 \end{pspicture}
 };
 \endxy $
 }
\end{center}
}

\newcommand{\semistrictbig}{
\xy 0;/r.21pc/:
    (-25,20)*+{A \ten B}="TL";
    (0,20)*+{A' \ten B}="TM";
    (25,20)*+{A'' \ten B}="TR";
    (-25,0)*+{A \ten B'}="ML";
    (0,0)*+{A' \ten B'}="MM";
    (25,0)*+{A' \ten B''}="MR";
    (-25,-20)*+{A \ten B''}="BL";
    (0,-20)*+{A' \ten B''}="BM";
    (25,-20)*+{A'' \ten B''}="BR";
        {\ar^{ f \ten B} "TL";"TM"};
        {\ar^{ f' \ten B} "TM";"TR"};
        {\ar^{ f \ten B'} "ML";"MM"};
        {\ar^{ f' \ten B'} "MM";"MR"};
        {\ar^{ f \ten B''} "BL";"BM"};
        {\ar^{ f' \ten B''} "BM";"BR"};
        {\ar_{ A \ten g} "TL";"ML"};
        {\ar_{ A \ten g'} "ML";"BL"};
        {\ar^{ A' \ten g} "TM";"MM"};
        {\ar^{ A' \ten g'} "MM";"BM"};
        {\ar^{ A'' \ten g} "TR";"MR"};
        {\ar^{ A'' \ten g'} "MR";"BR"};
                {\ar@{=>}^{ \INT_{f,g}} (-12,15)*{}; (-12,8)*{}};
                {\ar@{=>}^{ \INT_{f',g}} (12,15)*{}; (12,8)*{}};
                {\ar@{=>}^{ \INT_{f,g'}} (-12,-5)*{}; (-12,-12)*{}};
                {\ar@{=>}^{ \INT_{f',g'}} (12,-5)*{}; (12,-12)*{}};
 \endxy
\xy (0,0)*{=}; \endxy
 \xy 0;/r.21pc/:
    (-15,10)*+{A \ten B}="TL";
    (15,10)*+{A'' \ten B}="TR";
    (-15,-10)*+{A \ten B''}="BL";
    (15,-10)*+{A'' \ten B''}="BR";
        {\ar_{ A \ten g \circ g'} "TL";"BL"};
        {\ar^{ f \circ f' \ten B} "TL";"TR"};
        {\ar^{ A'' \ten g \circ g'} "TR";"BR"};
        {\ar_{ f \circ f' \ten B''} "BL";"BR"};
        {\ar@{=>}^{ \INT_{f \circ f',g \circ g'}} (-2,3)*{}; (-2,-3)*{}};
 \endxy }

\newcommand{\BIGproof}{
\begin{center}
\makebox[0pt]{ $ \psset{xunit=.18cm,yunit=0.18cm}
\xy
 (-10,20)*+{
        \begin{pspicture}[0.5](4.6,6.5)
        \rput(2.6,4.1){\comult}
        \rput(2.6,5.6){\birth}
        \rput(.8,4.1){\ident}
        \rput(4.4,4.1){\ident}
        \rput(.8,5.6){\smallident}
        \rput(4.4,5.6){\smallident}
        \rput(.8,2.6){\ident}
        \rput(2,2.6){\ident}
        \rput(3.8,2.6){\mult}
        \rput(1.4,1.1){\mult}
        \rput(3.8,1.1){\ident}
        \rput(1.4,.55){\smallident}
        \rput(1.4,.55){\death}
        \rput(3.8,1.1){\death}
        \end{pspicture}
 }="tl1";
 (-50,30)*+{
        \begin{pspicture}[0.5](4.6,6.5)
        \rput(2.6,4.1){\comult}
        \rput(2.6,5.6){\birth}
        \rput(.8,4.1){\ident}
        \rput(4.4,4.1){\ident}
        \rput(.8,5.6){\smallident}
        \rput(4.4,5.6){\smallident}
        \rput(.8,2.6){\ident}
        \rput(2,2.6){\ident}
        \rput(3.8,2.6){\mult}
        \rput(1.4,.55){\mult}
        \rput(2,2.05){\smallident}
        \rput(.8,2.05){\smallident}
        \rput(1.4,.55){\death}
        \rput(3.8,2.6){\death}
        \end{pspicture}
 }="tl1l1";
 (-10,0)*+{
        \begin{pspicture}[0.5](4.4,6.5)
        \rput(2,1.1){\mult}
        \rput(3.8,1.1){\ident}
        \rput(3.2,2.6){\comult}
        \rput(1.4,2.6){\ident}
        \rput(1.4,4.1){\ident}
        \rput(3.2,4.1){\mult}
         \rput(2.6,5.6){\birth}
        \rput(1.4,5.6){\smallident}
         \rput(3.8,5.6){\smallident}
         \rput(2,.55){\smallident}
        \rput(2,.55){\death}
        \rput(3.8,1.1){\death}
        \end{pspicture}
 }="tl2";
 (-50,10)*+{
        \begin{pspicture}[0.5](4.4,6.5)
        \rput(2,.55){\mult}
        \rput(3.2,2.6){\comult}
        \rput(1.4,2.6){\ident}
        \rput(1.4,4.1){\ident}
        \rput(3.2,4.1){\mult}
         \rput(2.6,5.6){\birth}
        \rput(1.4,5.6){\smallident}
         \rput(3.8,5.6){\smallident}
         \rput(2.6,2.05){\smallident}
         \rput(1.4,2.05){\smallident}
        \rput(2,.55){\death}
        \rput(3.8,2.6){\death}
        \end{pspicture}
 }="tl2l1";
 (-10,-20)*+{
        \begin{pspicture}[0.5](3.3,6.5)
        \rput(2,1.1){\comult}
        \rput(2,2.6){\mult}
        \rput(2.6,4.1){\mult}
        \rput(1.4,4.1){\curveleft}
         \rput(2,5.6){\birth}
        \rput(3.2,5.6){\smallident}
         \rput(.8,5.6){\smallident}
         \rput(1.4,.55){\smallident}
        \rput(1.4,.55){\death}
        \rput(2.6,1.1){\death}
        \end{pspicture}
 }="tl3";
  (-40,-10)*+{
        \begin{pspicture}[0.5](4.4,6.5)
        \rput(2,.55){\mult}
        \rput(2.6,2.05){\curveright}
        \rput(1.4,2.05){\ident}
        \rput(1.4,3.55){\ident}
        \rput(3.2,3.55){\mult}
         \rput(2.6,5.05){\birth}
        \rput(1.4,5.05){\smallident}
         \rput(3.8,5.05){\smallident}
        \rput(2,.55){\death}
        \end{pspicture}
 }="tl3l1";
  (0,-35)*+{
        \begin{pspicture}[0.5](3.4,3.5)
        \rput(2,1.1){\comult}
        \rput(2,2.6){\mult}
         \rput(1.4,.55){\smallident}
        \rput(1.4,.55){\death}
        \rput(2.6,1.1){\death}
        \end{pspicture}
 }="tl3d1";
 (10,20)*+{
        \begin{pspicture}[0.5](4.6,6.5)
        \rput(2.6,4.1){\comult}
        \rput(.8,4.1){\ident}
        \rput(4.4,4.1){\ident}
         \rput(.8,5.6){\smallident}
        \rput(2.6,5.6){\birth}
         \rput(4.4,5.6){\smallident}
        \rput(1.4,2.6){\mult}
        \rput(4.4,2.6){\ident}
        \rput(3.2,2.6){\ident}
        \rput(1.4,1.1){\ident}
        \rput(3.8,1.1){\mult}
          \rput(1.4,.55){\smallident}
        \rput(1.4,.55){\death}
        \rput(3.8,1.1){\death}
        \end{pspicture}
 }="tr1";
  (50,30)*+{
        \begin{pspicture}[0.5](4.6,6.5)
        \rput(2.6,4.1){\comult}
        \rput(.8,4.1){\ident}
        \rput(4.4,4.1){\ident}
         \rput(.8,5.6){\smallident}
        \rput(2.6,5.6){\birth}
         \rput(4.4,5.6){\smallident}
        \rput(1.4,2.6){\mult}
        \rput(4.4,2.6){\ident}
        \rput(3.2,2.6){\ident}
        \rput(3.8,.55){\mult}
        \rput(3.2,2.055){\smallident}
        \rput(4.4,2.055){\smallident}
        \rput(1.4,2.6){\death}
        \rput(3.8,.55){\death}
        \end{pspicture}
 }="tr1r1";
(10,0)*+{
        \begin{pspicture}[0.5](2.8,6.5)
        \rput(2.6,4.1){\ident}
        \rput(.8,4.1){\mult}
        \rput(1.4,5.6){\birth}
        \rput(.2,5.6){\smallident}
        \rput(2.6,5.6){\smallident}
        \rput(.8,2.6){\comult}
        \rput(2.6,2.6){\ident}
        \rput(2,1.1){\mult}
        \rput(.2,1.1){\ident}
          \rput(.2,.55){\smallident}
        \rput(.2,.55){\death}
        \rput(2,1.1){\death}
        \end{pspicture}
 }="tr2";
 (50,10)*+{
        \begin{pspicture}[0.5](2.8,6.5)
        \rput(2.6,4.1){\ident}
        \rput(.8,4.1){\mult}
        \rput(1.4,5.6){\birth}
        \rput(.2,5.6){\smallident}
        \rput(2.6,5.6){\smallident}
        \rput(.8,2.6){\comult}
        \rput(2.6,2.6){\ident}
        \rput(2,.55){\mult}
        \rput(2.6,2.05){\smallident}
        \rput(1.4,2.05){\smallident}
        \rput(.2,2.6){\death}
        \rput(2,.55){\death}
        \end{pspicture}
 }="tr2r1";
 (40,-30)*+{
        \begin{pspicture}[0.5](2.8,6.5)
        \rput(2.6,3.55){\ident}
        \rput(.8,3.55){\mult}
        \rput(1.4,5.05){\birth}
        \rput(.2,5.05){\smallident}
        \rput(2.6,5.05){\smallident}
        \rput(1.4,2.05){\curveleft}
        \rput(2.6,2.05){\ident}
        \rput(2,.55){\mult}
        \rput(2,.55){\death}
        \end{pspicture}
 }="tr2r1dd";
  (10,-20)*+{
        \begin{pspicture}[0.5](3.4,6.5)
        \rput(2,1.1){\comult}
        \rput(2,2.6){\mult}
        \rput(2.6,4.1){\curveright}
        \rput(1.4,4.1){\mult}
        \rput(2,5.6){\birth}
        \rput(3.2,5.6){\smallident}
         \rput(.8,5.6){\smallident}
         \rput(1.4,.55){\smallident}
        \rput(1.4,.55){\death}
        \rput(2.6,1.1){\death}
        \end{pspicture}
 }="tr3";
  (35,-5)*+{
        \begin{pspicture}[0.5](2.8,6.5)
        \rput(2.6,4.1){\ident}
        \rput(.8,4.1){\mult}
        \rput(1.4,5.6){\birth}
        \rput(.2,5.6){\smallident}
        \rput(2.6,5.6){\smallident}
        \rput(.8,2.6){\comult}
        \rput(2.6,2.6){\ident}
        \rput(2,1.1){\mult}
        \rput(.2,1.1){\ident}
          \rput(2,.55){\smallident}
        \rput(2,.55){\death}
        \rput(.2,1.1){\death}
        \end{pspicture}
 }="tr3r1";
 (15,-35)*+{
        \begin{pspicture}[0.5](3.4,3.5)
        \rput(2,1.1){\comult}
        \rput(2,2.6){\mult}
         \rput(2.6,.55){\smallident}
        \rput(1.4,1.1){\death}
        \rput(2.6,.55){\death}
        \end{pspicture}
 }="tr3d1";
 (25,-20)*+{
        \begin{pspicture}[0.5](3.4,6.5)
        \rput(2,1.1){\comult}
        \rput(2,2.6){\mult}
        \rput(2.6,4.1){\curveright}
        \rput(1.4,4.1){\mult}
        \rput(2,5.6){\birth}
        \rput(3.2,5.6){\smallident}
         \rput(.8,5.6){\smallident}
         \rput(2.6,.55){\smallident}
        \rput(2.6,.55){\death}
        \rput(1.4,1.1){\death}
        \end{pspicture}
 }="tr3d2";
  (-50,-40)*+{
        \begin{pspicture}[0.5](4.4,6.5)
        \rput(2,.55){\mult}
        \rput(3.2,2.6){\comult}
        \rput(1.4,2.6){\ident}
        \rput(1.4,4.1){\ident}
        \rput(3.2,4.1){\ident}
         \rput(2.6,2.05){\smallident}
         \rput(1.4,2.05){\smallident}
        \rput(2,.55){\death}
        \rput(3.8,2.6){\death}
        \end{pspicture}
 }="X2";
   (0,-60)*+{
        \begin{pspicture}[0.5](4.4,2.5)
        \rput(2,.55){\mult}
        \rput(2,.55){\death}
        \end{pspicture}
 }="X3";
 (50,-40)*+{
        \begin{pspicture}[0.5](2.8,6.5)
        \rput(2.6,4.1){\ident}
        \rput(.8,4.1){\ident}
        \rput(.8,2.6){\comult}
        \rput(2.6,2.6){\ident}
        \rput(2,.55){\mult}
        \rput(2.6,2.05){\smallident}
        \rput(1.4,2.05){\smallident}
        \rput(.2,2.6){\death}
        \rput(2,.55){\death}
        \end{pspicture}
 }="Y2";
    {\ar@{=>}^{\INT} "tl1";"tr1"};
    {\ar@{=>}_{p} "X2";"X3"};
    {\ar@{=>}^{q} "Y2";"X3"};
    {\ar@{=>}_{p} "tl3d1";"X3"};
    {\ar@{=>}_{\ell} "tl3l1";"X3"};
    {\ar@{=>}_{c} "tl1";"tl2"};
    {\ar@{=>}_{\INT} "tl1l1";"tl1"};
    {\ar@{=>}^{\INT} "tr1r1";"tr1"};
    {\ar@{=>}_{\INT} "tr2r1";"tr2"};
    {\ar@{=>}^{\INT} "tl2l1";"tl2"};
    {\ar@{=>}^{c} "tl1l1";"tl2l1"};
    {\ar@{=>}^{b} "tr1r1";"tr2r1"};
    {\ar@{=>}^{\INT} "tl1l1";"tr1r1"};
    {\ar@{=>}_{p} "tl3l1";"tl3"};
    {\ar@{=>}_{p} "tl2l1";"tl3l1"};
    {\ar@{=>}^{r} "tr3";"tl3d1"};
    {\ar@{=>}_{\ell} "tl3";"tl3d1"};
     {\ar@{=>}_{\ell} "tl2l1";"X2"};
     {\ar@{=>}_{r} "tr2r1";"Y2"};
    {\ar@{=>}^{b} "tr1";"tr2"};
    {\ar@{=>}^{c} "tr2";"tr3"};
    {\ar@{=>}_{b} "tl2";"tl3"};
    {\ar@{=>}_{a} "tl3";"tr3"};
    {\ar@{=>}^{\INT} "tl3d1";"tr3d1"};
    {\ar@{=>}_{q} "tr3d1";"X3"};
    {\ar@{=>}^{r} "tr3d2";"tr3d1"};
     {\ar@{=>}_{\INT} "tr3";"tr3d2"};
      {\ar@{=>}^{\INT} "tr2";"tr3r1"};
      {\ar@{=>}^{c} "tr3r1";"tr3d2"};
      {\ar@{=>}^{\INT} "tr2r1";"tr3r1"};
      {\ar@{=>}^{q} "tr2r1";"tr2r1dd"};
      {\ar@{=>}^{q} "tr2r1dd";"tr3d2"};
       {\ar@{=>}^{r} "tr2r1dd";"X3"};
 \endxy
$ }
\end{center}
}

\newcommand{\HUGEproof}{
\begin{center}
\makebox[0pt]{ $ \psset{xunit=.15cm,yunit=0.15cm}
 \xy
 (-80,-20)*+{

 }="t7";
                {\ar@{=>}_{a} "b3";"b4"};
                {\ar@{=>}^{\INT} "b3u1";"b3"};
                {\ar@{=>}_{a} "b3u1";"b4u1"};
                {\ar@{=>}^{\INT} "b4u1";"b4"};
                {\ar@{=>}_{a} "t1";"l2"};
                {\ar@{=>}_{z} "t1";"t1d1"};
                {\ar@{=>}_{\INT} "t1d1";"b3u3"};
                {\ar@{=>}_{z} "t2";"t2d1"};
                {\ar@{=>}^{\INT} "l2r1";"b3u2"};
                {\ar@{=>}_{z} "t1d1";"t2d1"};
                {\ar@{=>}_{a} "l2r1";"t1d1";};
                {\ar@{=>}_{\INT} "t2d1";"t3d1"};
                {\ar@{=>}_{a} "t3d1";"t4d1"};
                {\ar@{=>}_{\INT} "t2d2";"t3d2"};
                {\ar@{=>}_{\INT} "t3d1";"t3d2"};
                {\ar@{=>}_{z} "t4";"t4d1"};
                {\ar@{=>}_{z} "t3";"t3d1"};
                {\ar@{=>}_{z} "l2";"l2r1";};
                {\ar@{=>}_{\INT} "l2";"b1"};
                {\ar@{=>}_{z} "b1";"b2"};
                {\ar@{=>}_{\INT} "b2";"b3"};
                {\ar@{=>}^{z} "t1";"t2"};
                {\ar@{=>}^{\INT} "t2";"t3"};
                {\ar@{=>}^{a} "t3";"t4"};
                {\ar@{=>}^{\INT} "t4";"t5"};
                {\ar@{=>}^{n} "t5";"t6"};
                {\ar@{=>}^{\INT} "t6";"t7"};
                {\ar@{=>}_{z} "t5";"t5d1"};
                {\ar@{=>}_{\INT} "t4d1";"t5d1"};
                {\ar@{=>}_{z} "b3u3";"t2d2"};
                {\ar@{=>}_{\INT} "t2d1";"t2d2"};
                {\ar@{=>}_{a} "t3d2";"t4d2"};
                {\ar@{=>}_{\INT} "t4d1";"t4d2"};
                {\ar@{=>}_-{\INT} "t4d2";"t5d2"+(-8,0)};
                {\ar@{=>}_{\INT} "t5d1";"t5d2"};
                {\ar@{=>}_{n} "t5d2";"t6d2"};
                {\ar@{=>}_{\INT} "t6d1";"t6d2"};
                {\ar@{=>}_<<<<{\INT} "t6d3";"r4"};
                {\ar@{=>}_{\INT} "t2d3";"t3d3"};
                {\ar@{=>}_{a} "t3d3";"t4d3"};
                {\ar@{=>}_{\INT} "t6d4";"r5"};
                {\ar@{=>}_{\INT} "t6d3";"t6d4"};
                {\ar@{=>}_{n} "t6d6";"b7"};
                {\ar@{=>}_{a} "t6d5";"t6d6"};
                {\ar@{=>}_{\INT} "t2d4";"t3d4"};
                {\ar@{=>}_<<<<{\INT} "t3d3";"t3d4"};
                {\ar@{=>}_{n} "t5d5";"t6d5"};
                {\ar@{=>}_{n} "t6d5";"r6"};
                {\ar@{=>}_{\INT} "t4d3";"t5d3"};
                {\ar@{=>}_{a} "t3d4";"t4d4"};
                {\ar@{=>}_<<<<{\INT} "t4d3";"t4d4"};
                {\ar@{=>}^{\INT} "t4d4";"t5d4"};
                {\ar@{=>}_<<<<{\INT} "t5d4";"t5d5"};
                {\ar@{=>}_{n} "t5d4";"t6d4"};
                {\ar@{=>}_<<<<{\INT} "t5d3";"t5d4"};
                {\ar@{=>}_{n} "t5d3";"t6d3"};
                {\ar@{=>}_{a} "t2d2";"t2d3"};
                {\ar@{=>}_{a} "t3d2";"t3d3"};
                {\ar@{=>}_{a} "t4d2";"t4d3"};
                {\ar@{=>}_{a} "t5d2";"t5d3"};
                {\ar@{=>}_{a} "t6d2";"t6d3"};
                {\ar@{=>}_{\INT} "t6d2";"r3"};
                {\ar@{=>}_{n} "t5d1";"t6d1"};
                {\ar@{=>}_{\INT} "t6d1";"r2"};
                {\ar@{=>}_{z} "t6";"t6d1"};
                {\ar@{=>}_{z} "t7";"r2"};
                 {\ar@{=>}_{z} "b4u3";"t2d3"};
                {\ar@{=>}^{\INT} "r4";"r5"};
                {\ar@{=>}_{\INT} "b4";"b5"};
                {\ar@{=>}_{n} "b5";"b6"};
                {\ar@{=>}_{\INT} "b5";"t6d6"};
                {\ar@{=>}_{\INT} "t2d4";"t2d5"};
                {\ar@{=>}_{z} "b5u3";"t2d4";};
                {\ar@{=>}_{\INT} "t2d3";"t2d4"};
                {\ar@{=>}_{\INT} "b6";"b7"};
                {\ar@{=>}_{\INT} "r2";"r3"};
                {\ar@{=>}_{n} "r5";"r6"};
                {\ar@{=>}_{a} "r3";"r4"};
                {\ar@{=>}_{a} "b7";"r6"};
                {\ar@{=>}_{a} "b3u2";"b4u2"};
                {\ar@{=>}^{a} "b4u2";"b5u2"};
                {\ar@{=>}^{a} "b5u3";"b5u2"};
                {\ar@{=>}_{a} "b3u2";"b3u3"};
                {\ar@{=>}_{a} "b3u3";"b4u3"};
                {\ar@{=>}_{n} "b5u3";"t3d5"};
                {\ar@{=>}_{n} "t4d5";"b5u2"};
                {\ar@{=>}_{n} "b5u3";"t2d5";};
                {\ar@{=>}_{a} "t3d5";"t4d5"};
                {\ar@{=>}_<<<<{\INT} "t3d4";"t3d5"};
                {\ar@{=>}^{\INT} "b4u3";"b5u3"};
                {\ar@{=>}^{\INT} "t2d5";"t3d5"};
                {\ar@{=>}^{\INT} "t4d5";"t5d5"};
                {\ar@{=>}_<<<<{\INT} "t4d4";"t4d5"};
                {\ar@{=>}^{\INT} "b3u2";"b3u1"};
                {\ar@{=>}^{\INT} "b4u2";"b4u1"};
                {\ar@{=>}^{\INT} "b5u2";"t6d5"};
                {\ar@{=>}^{\INT} "b4u2";"t6d6"};
 \endxy
$ }
\end{center}
}

\newcommand{\HUGEproofII}{
\begin{center}
\makebox[0pt]{ $ \psset{xunit=.15cm,yunit=0.15cm}
 \xy
 (-80,-10)*+{

 }="t7";
                {\ar@{=>}_{a} "b3";"b4"};
                {\ar@{=>}^{\INT} "b3u1";"b3"};
                {\ar@{=>}_{a} "b3u1";"b4u1"};
                {\ar@{=>}^{\INT} "b4u1";"b4"};
                {\ar@{=>}_{a} "t1";"l2"};
                {\ar@{=>}_{z} "t1";"t1d1"};
                {\ar@{=>}_{\INT} "t1d1";"b3u3"};
                {\ar@{=>}_{z} "t2";"t2d1"};
                {\ar@{=>}^{\INT} "l2r1";"b3u2"};
                {\ar@{=>}_{z} "t1d1";"t2d1"};
                {\ar@{=>}_{a} "l2r1";"t1d1";};
                {\ar@{=>}_{\INT} "t2d1";"t3d1"};
                {\ar@{=>}_{a} "t3d1";"t4d1"};
                {\ar@{=>}_{\INT} "t2d2";"t3d2"};
                {\ar@{=>}_{\INT} "t3d1";"t3d2"};
                {\ar@{=>}_{z} "t4";"t4d1"};
                {\ar@{=>}_{z} "t3";"t3d1"};
                {\ar@{=>}_{z} "l2";"l2r1";};
                {\ar@{=>}_{\INT} "l2";"b1"};
                {\ar@{=>}_{z} "b1";"b2"};
                {\ar@{=>}_{\INT} "b2";"b3"};
                {\ar@{=>}^{z} "t1";"t2"};
                {\ar@{=>}^{\INT} "t2";"t3"};
                {\ar@{=>}^{a} "t3";"t4"};
                {\ar@{=>}^{\INT} "t4";"t5"};
                {\ar@{=>}^{n} "t5";"t6"};
                {\ar@{=>}^{\INT} "t6";"t7"};
                {\ar@{=>}_{z} "t5";"t5d1"};
                {\ar@{=>}_{\INT} "t4d1";"t5d1"};
                {\ar@{=>}_{z} "b3u3";"t2d2"};
                {\ar@{=>}_{\INT} "t2d1";"t2d2"};
                {\ar@{=>}_{a} "t3d2";"t4d2"};
                {\ar@{=>}_{\INT} "t4d1";"t4d2"};
                {\ar@{=>}_-{\INT} "t4d2";"t5d2"+(-8,0)};
                {\ar@{=>}_{\INT} "t5d1";"t5d2"};
                {\ar@{=>}_{n} "t5d2";"t6d2"};
                {\ar@{=>}_{\INT} "t6d1";"t6d2"};
                {\ar@{=>}_<<<<{\INT} "t6d3";"r4"};
                {\ar@{=>}_{\INT} "t2d3";"t3d3"};
                {\ar@{=>}_{a} "t3d3";"t4d3"};
                {\ar@{=>}_{\INT} "t6d4";"r5"};
                {\ar@{=>}_{\INT} "t6d3";"t6d4"};
                {\ar@{=>}_{n} "t6d6";"b7"};
                {\ar@{=>}_{a} "t6d5";"t6d6"};
                {\ar@{=>}_{\INT} "t2d4";"t3d4"};
                {\ar@{=>}_<<<<{\INT} "t3d3";"t3d4"};
                {\ar@{=>}_{n} "t5d5";"t6d5"};
                {\ar@{=>}_{n} "t6d5";"r6"};
                {\ar@{=>}_{\INT} "t4d3";"t5d3"};
                {\ar@{=>}_{a} "t3d4";"t4d4"};
                {\ar@{=>}_<<<<{\INT} "t4d3";"t4d4"};
                {\ar@{=>}^{\INT} "t4d4";"t5d4"};
                {\ar@{=>}_<<<<{\INT} "t5d4";"t5d5"};
                {\ar@{=>}_{n} "t5d4";"t6d4"};
                {\ar@{=>}_<<<<{\INT} "t5d3";"t5d4"};
                {\ar@{=>}_{n} "t5d3";"t6d3"};
                {\ar@{=>}_{a} "t2d2";"t2d3"};
                {\ar@{=>}_{a} "t3d2";"t3d3"};
                {\ar@{=>}_{a} "t4d2";"t4d3"};
                {\ar@{=>}_{a} "t5d2";"t5d3"};
                {\ar@{=>}_{a} "t6d2";"t6d3"};
                {\ar@{=>}_{\INT} "t6d2";"r3"};
                {\ar@{=>}_{n} "t5d1";"t6d1"};
                {\ar@{=>}_{\INT} "t6d1";"r2"};
                {\ar@{=>}_{z} "t6";"t6d1"};
                {\ar@{=>}_{z} "t7";"r2"};
                 {\ar@{=>}_{z} "b4u3";"t2d3"};
                {\ar@{=>}^{\INT} "r4";"r5"};
                {\ar@{=>}_{\INT} "b4";"b5"};
                {\ar@{=>}_{n} "b5";"b6"};
                {\ar@{=>}_{\INT} "b5";"t6d6"};
                {\ar@{=>}_{\INT} "t2d4";"t2d5"};
                {\ar@{=>}_{z} "b5u3";"t2d4";};
                {\ar@{=>}_{\INT} "t2d3";"t2d4"};
                {\ar@{=>}_{\INT} "b6";"b7"};
                {\ar@{=>}_{\INT} "r2";"r3"};
                {\ar@{=>}_{n} "r5";"r6"};
                {\ar@{=>}_{a} "r3";"r4"};
                {\ar@{=>}_{a} "b7";"r6"};
                {\ar@{=>}_{a} "b3u2";"b4u2"};
                {\ar@{=>}^{a} "b4u2";"b5u2"};
                {\ar@{=>}^{a} "b5u3";"b5u2"};
                {\ar@{=>}_{a} "b3u2";"b3u3"};
                {\ar@{=>}_{a} "b3u3";"b4u3"};
                {\ar@{=>}_{n} "b5u3";"t3d5"};
                {\ar@{=>}_{n} "t4d5";"b5u2"};
                {\ar@{=>}_{n} "b5u3";"t2d5";};
                {\ar@{=>}_{a} "t3d5";"t4d5"};
                {\ar@{=>}_<<<<{\INT} "t3d4";"t3d5"};
                {\ar@{=>}^{\INT} "b4u3";"b5u3"};
                {\ar@{=>}^{\INT} "t2d5";"t3d5"};
                {\ar@{=>}^{\INT} "t4d5";"t5d5"};
                {\ar@{=>}_<<<<{\INT} "t4d4";"t4d5"};
                {\ar@{=>}^{\INT} "b3u2";"b3u1"};
                {\ar@{=>}^{\INT} "b4u2";"b4u1"};
                {\ar@{=>}^{\INT} "b5u2";"t6d5"};
                {\ar@{=>}^{\INT} "b4u2";"t6d6"};
 \endxy
$ }
\end{center} }

   \title{Frobenius algebras and planar open string topological field theories}
      \author{
      Aaron D.\ Lauda \\
      Department of Pure Mathematics and Mathematical Statistics, \\
      University of Cambridge,\\
      Cambridge CB3 0WB\\
      UK \\
         \\
      email: a.lauda@dpmms.cam.ac.uk \\}

\begin{document}
\maketitle

\begin{abstract}
Motivated by the Moore-Segal axioms for an open-closed topological
field theory, we consider planar open string topological field
theories. We rigorously define a category $\cat{2Thick}$ whose
objects and morphisms can be thought of as open strings and
diffeomorphism classes of planar open string worldsheets. Just as
the category of 2-dimensional cobordisms can be described as the
free symmetric monoidal category on a commutative Frobenius
algebra, $\cat{2Thick}$ is shown to be the free monoidal category
on a noncommutative Frobenius algebra, hence justifying this
choice of data in the Moore-Segal axioms. Our formalism is
inherently categorical allowing us to generalize this result.  As
a stepping stone towards topological membrane theory we define a
2-category of open strings, planar open string worldsheets, and
isotopy classes of 3-dimensional membranes defined by
diffeomorphisms of the open string worldsheets.  This 2-category
is shown to be the free (weak monoidal) 2-category on a
`categorified Frobenius algebra', meaning that categorified
Frobenius algebras determine invariants of these 3-dimensional
membranes.
\end{abstract}

\section{Introduction} \label{secIntro}
It is well known that 2-dimensional topological quantum field
theories are equivalent to commutative Frobenius algebras
\cite{Abrams1,Dij1,Kock}.  This simple algebraic characterization
arises directly from the simplicity of $\cat{2Cob}$, the
2-dimensional cobordism category.  In fact, \cat{2Cob} admits a
completely algebraic description as the free symmetric monoidal
category on a commutative Frobenius algebra. This description is
apparent when one examines the generating cobordisms in
\cat{2Cob}:
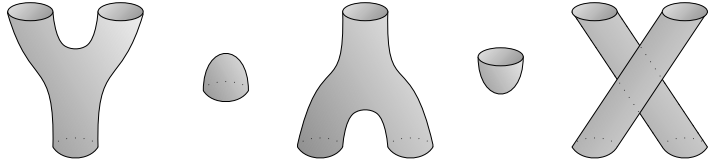
\[
 \twocobb
\]
These cobordisms equip the circle with a commutative algebra and
coalgebra structure, and the relations on the generators of
\cat{2Cob} make the circle into a commutative Frobenius algebra.

This universal property simplifies the construction of functors
from \cat{2Cob} into other symmetric monoidal categories, such as
the category of vector spaces or graded vector spaces. Simply find
an example of a commutative Frobenius algebra in a symmetric
monoidal category $\mathcal{C}$, and the universal property of
\cat{2Cob} determines a functor $\cat{2Cob} \to \mathcal{C}$.
Since a 2-dimensional topological quantum field theory is a
functor from \cat{2Cob} into a symmetric monoidal category, it is
clear that the universal property of $\cat{2Cob}$ greatly
simplifies the study of 2-dimensional topological quantum field
theories.

Topological string theory has recently been a topic of
considerable interest to physicists, and it presents interesting
new problems for mathematicians as well. The celebrated Segal
axioms \cite{Segal} of a conformal field theory yield
2-dimensional topological quantum field theories as the simplest
examples, namely those that do not depend on the conformal
structure of the string backgrounds~\cite{Segal2}.  In this
context, the circle is interpreted as a `closed string', and
2-dimensional cobordisms as `closed string worldsheets'.
Topological string theory for open strings, which inevitably also
contains a closed string sector, can be described as an open
\textit{and} closed topological field theory coupled to
topological gravity. In this case, interactions between the open
and closed strings induce a much richer 2-dimensional topology.

In an open-closed topological field theory, closed strings can
evolve into open strings as demonstrated by the cobordism below:
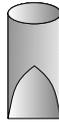
\[ 
 \begin{pspicture}(3,2.5)
  \begin{psclip}{
    \pscustom{
        \psline(1.5,2)(1.5,0)
        \psline(1.5,0)(2.5,0)
        \psline(2.5,0)(2.5,2)
        \psellipse(2,2)(.5,.2)
    }
  }
    \pspolygon[fillcolor=lightgray,fillstyle=gradient,
    gradbegin=white,gradend=gray,gradmidpoint=1,gradangle=90](1.5,0)(1.5,2.4)(2.5,2.4)(2.5,0)(1.5,0)
 \end{psclip}
 \pscustom[fillcolor=lightgray,fillstyle=gradient,
        gradbegin=white, gradend=gray,gradmidpoint=0,gradangle=88]{
    \psline(1.5,2)(1.5,0)
    \psbezier(1.5,0)(1.6,.5)(1.6,.7)(2,1)
    \psbezier(2,1)(2.4,.7)(2.4,.5)(2.5,0)
    \psellipse(2,2)(.5,.2)
 }
 \psellipse[fillcolor=lightgray,fillstyle=gradient,
        gradbegin=white, gradend=gray,gradmidpoint=1,gradangle=90](2,2)(.5,.2)
 \end{pspicture}
\]
The topological sewing conditions for these surfaces were first
analyzed by Cardy and Lewellen~\cite{CL,Lew}, and a category of
such surfaces has also been constructed \cite{BCR}. Our primary
interest will be a subset of the conditions imposed on the open
strings. In particular, the topology generated by worldsheets that
can be embedded in the plane:
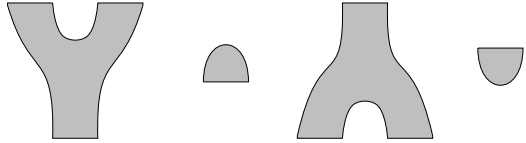
\[
\twothick
\]
A sketch of the axioms for an open-closed field theory was given
by Moore and Segal~\cite{MS} and has been refined in the more
recent work of Lazaroiu~\cite{lar2,lar1}.  In these axioms, the
open strings are represented algebraically by a
\textit{noncommutative} Frobenius algebra. This axiom is intended
to express the topological difference between the circle and the
interval, or closed string and open string. One of the results of
this paper is a justification of this particular axiom.

We define a topological category $\cat{2Thick}$ whose objects are
open strings and whose morphisms are diffeomorphism classes of
planar worldsheets.  More precisely, the objects of this category
are the natural numbers, and the morphisms are diffeomorphism
classes of smooth oriented compact surfaces whose boundary is
equipped with disjoint distinguished intervals. We will sometimes
refer to $\cat{2Thick}$ as the category of `2-dimensional thick
tangles' since it is a `thickened' version of the category of
tangles embedded in the plane. We then define a planar open string
topological field theory as a functor from $\cat{2Thick}$ into a
monoidal category $\mathcal{C}$. We prove that $\cat{2Thick}$ is
the monoidal category freely generated by a (noncommutative)
Frobenius algebra. This implies that a planar open string
topological field theory is equivalent to a (not necessarily
commutative) Frobenius algebra, thereby rigorously establishing
the Moore-Segal axioms for this portion of the open string data.
It also implies that a Frobenius algebra in an arbitrary monoidal
category determines an invariant of 2-dimensional thick tangles.

The idea that all of the known string theories might arise as the
effective limit of a higher-dimensional theory known as $M$-theory
\cite{Town,Witten1} has prompted many to consider `topological
membrane theories'~\cite{Kogan2, Kogan1}.   In topological
membrane theory the 2-dimensional string worldsheets arise as
boundaries of 3-dimensional membranes.  This approach has been
used to show that various types of string theories (e.g.\
heterotic, open, unoriented) arise from the choice of boundary
conditions on a topological membrane~\cite{Kogan3,Castelo}.  These
important results suggest that one should consider a topological
theory that describes open strings and their worldsheets as part
of a larger 3-dimensional topology which reduces to the usual
theory on the boundary.  In this paper we take some first steps
towards realizing such a theory.

Our description of open strings and their worldsheets has the
advantage that it can easily be generalized to distinguish
diffeomorphic worldsheets. Thus, rather than considering
worldsheets only up to diffeomorphism, our formalism can be
extended to describe open strings, \textit{all} planar open string
worldsheets, \textit{and} the 3-dimensional membranes defined by
isotopy classes of diffeomorphisms of planar open string
worldsheets that preserve the initial and final states. Some
examples are given below: \threethick Because of the three levels
of structure present, this type of topology is most naturally
described using 2-categories rather than categories.

Extensions of this sort have already been considered in the
symmetric case of $\cat{2Cob}$.  Using what we will call
$\cat{2Cob}_2$, a 2-categorical extension of the usual
two-dimensional cobordism category $\cat{2Cob}$, Ulrike Tillmann
has considered extended 2-dimensional topological quantum field
theories consisting of 2-functors from $\cat{2Cob}_2$ into a
generalized version of linear categories -- the 2-category of
linear abelian categories \cite{Tillmann}. By encoding
3-dimensional information into a 2-dimensional topological quantum
field theory, Tillmann has established interesting connections
between 2-dimensional TQFT's extended in this way and traditional
3-dimensional topological field theories.

Analogous to Tillmann's construction of extended (closed) string
topological field theories, we consider extended planar open
string topological field theories.  We begin by rigorously
defining the 2-category of open strings, planar open string
worldsheets, and isotopy classes of worldsheets embedded in the
cube, that we denote as $\cat{3Thick}$. We will sometimes refer to
this 2-category as the 2-category of `3-dimensional thick tangles'
since it is a categorification of the category $\cat{2Thick}$.

We will show that $\cat{3Thick}$ can be described by a universal
property that completely characterizes these particular extended
open string field theories.  More precisely, we show that
$\cat{3Thick}$ is the (semistrict monoidal) 2-category freely
generated by a `categorified Frobenius algebra'. A categorified
Frobenius algebra, is just a category with the structure of a
Frobenius algebra where all the usual axioms only hold up to
coherent isomorphism.  This description of $\cat{3Thick}$ implies
that categorified Frobenius algebras, or pseudo Frobenius algebras
as we will refer to them, determine invariants of 3-dimensional
thick tangles.

The fact that our restricted choice of open string worldsheets can
be described as the free monoidal category on a noncommutative
Frobenius algebra is perhaps not that surprising to the expert.
What is interesting is the fact that this result arises quite
naturally from higher-dimensional category theory.  Even more
interesting is that, once seen from a categorical perspective, a
universal property for the 2-category of open string membranes in
the cube is obtained by a straightforward categorification of the
lower dimensional result.

All of the results in this paper are obtained using the
relationship between Frobenius algebras and adjunctions
\cite{Lau1}. An adjunction is just the abstraction of the
definition of adjoint functors between categories.  Indeed, an
adjunction in the 2-category $\cat{Cat}$ whose objects are
categories, morphisms are functors, and whose 2-morphisms are
natural transformations, produces the usual notion of adjoint
functors.  By considering adjunctions in arbitrary 2-categories we
are able to rephrase the definition of a Frobenius algebra in the
language of category theory. Using 2-categorical `string
diagrams', adjunctions provide a bridge that illuminates the
relationship between Frobenius algebras and thick tangles.

This abstract approach has the advantage that, phrased in an
intrinsically categorical way, it can easily be generalized to
define pseudo Frobenius algebras and establish their relationship
with 3-dimensional thick tangles.  Using pseudoadjunctions, the
3-categorical analogs of adjunctions, we are able to define pseudo
Frobenius algebras.  Furthermore, using a version of string
diagrams for 3-categories we are able to show that the 2-category
$\cat{3Thick}$ is the 2-category freely generated by a pseudo
Frobenius algebra. Thus, our results not only provide an algebraic
description of $\cat{2Thick}$ and $\cat{3Thick}$, but also provide
an abstract framework with which these results naturally arise.

\section{The two-dimensional case}

We begin with the 2-dimensional case.  In Section \ref{secFrob} we
review some of the definitions of a Frobenius algebra that will
appear in a categorified form later in the paper.  In Section
\ref{secString} we review the theory of 2-categorical string
diagrams and we use them in Sections \ref{secWA} and \ref{secWAA}
to establish the relationship between Frobenius algebras and
adjunctions.  In Section \ref{sec2Dthick} we define the monoidal
category of 2-dimensional thick tangles and state one of the main
theorems of the paper.  This theorem follows from the categorified
version that will be proven later on.

\subsection{Frobenius Algebras} \label{secFrob}

Frobenius algebras appear throughout many branches of mathematics.
This can be attributed to the robustness of the definition and the
many ways it can be formulated.  Originally, the term Frobenius
algebra referred to an algebra $A$ with the property that $A
\simeq A^*$ as right $A$-modules.  Later, Nakayama provided many
equivalent definitions of Frobenius algebras, including the
characterization of a Frobenius algebra as a finite dimensional
$\Bbbk$-algebra equipped with a linear functional $\varepsilon
\maps A \to \Bbbk$ whose nullspace contains no nontrivial
ideals~\cite{Nak1,Nak2}. Another equivalent definition, motivated
by topological considerations, defines a Frobenius algebra as an
algebra $A$ equipped with a coalgebra structure where the
comultiplication is a map of $A$-modules. This topologically
motivated definition arose in order to rigorously establish the
theorem that a two-dimensional topological quantum field theory is
essentially the same as a commutative Frobenius algebra.   This
fact was first observed by Dijkgraaf~\cite{Dij1}, but it was not
rigorously shown until Frobenius algebras were reformulated in
terms of a coalgebra structure~\cite{Abrams1}.  A modern proof of
this result is given by Kock~\cite{Kock}.

It is clear that there are many equivalent definitions of
Frobenius algebras, each suited for various applications.  We make
no attempt to describe them all nor establish their equivalence.
Rather, we will focus on a few that are relevant to our main
theorem relating Frobenius algebras to thick tangles. Each of
these definitions has a topological nature that we will explain
using a diagrammatic shorthand notation for the relevant maps in
each definition.  These maps are, in a sense, the `topological
building blocks' that each definition provides. Loosely speaking,
our theorem states that this shorthand notation is much more than
a convenient device.  In fact, we will show that the topological
picture is exactly equivalent to the algebraic picture in the
sense that any topological manipulation of the diagrams
corresponds to algebraic manipulations of a Frobenius algebra
compatible with its axioms.

\begin{prop} \label{equivFrob}
Let $A$ be a vector space equipped with morphisms:
\begin{itemize}
    \item $m \maps A \ten A \to A$, and
    \item $\iota\maps \Bbbk \to A$,
\end{itemize}
satisfying the algebra axioms:
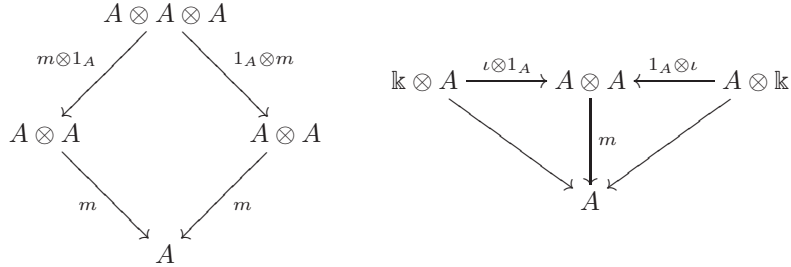
\[
 \xy
 (16,0)*+{A \ten A}="R";
 (0,16)*+{A \ten A \ten A}="T";
 (0,-16)*+{A}="B";
 (-16,0)*+{A \ten A}="L";
    {\ar^{1_A \ten m} "T";"R"};
    {\ar_{m \ten 1_A} "T";"L"};
    {\ar_{ m} "L";"B"};
    {\ar^{ m} "R";"B"};
 \endxy
\quad \quad
 \vcenter{\xy
 (-22,0)*+{\Bbbk \ten A}="L";
 (0,0)*+{A \ten A}="M";
 (22,0)*+{A \ten \Bbbk}="R";
 (0,-16)*+{A}="B";
    {\ar^{\iota \ten 1_A} "L";"M"};
    {\ar_{1_A \ten \iota} "R";"M"};
    {\ar_{ } "L";"B"};
    {\ar^{ } "R";"B"};
    {\ar^{ m} "M";"B"};
 \endxy}
\]
Then the following conditions on a form $\varepsilon \maps A \to
\Bbbk$ are equivalent:

\begin{enumerate}[i.)]
    \item There exists a map $\Delta\maps A \to A \ten A$ that, together
with $\varepsilon \maps A \to \Bbbk$, defines a coalgebra
structure on $A$ satisfying the Frobenius identities. Or more
explicitly, the algebra $A$ is equipped with a map  $\Delta \maps
A \to A \ten A$ such that the following diagrams commute:

\begin{itemize}
 \item the coalgebra axioms:
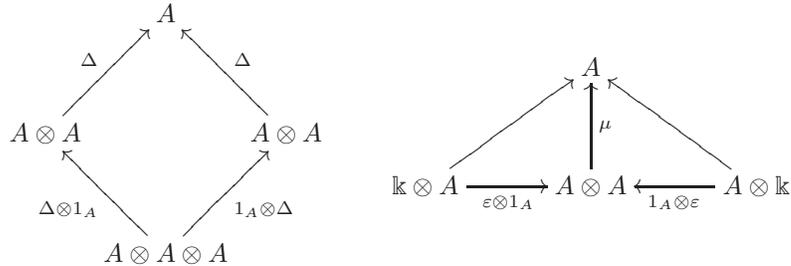
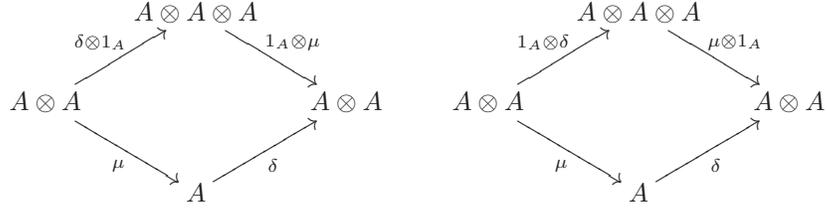
\[
\xy
 (16,0)*+{A \ten A}="R";
 (0,-16)*+{A \ten A \ten A}="T";
 (0,16)*+{A}="B";
 (-16,0)*+{A \ten A}="L";
    {\ar_{1_A \ten \Delta} "T";"R"};
    {\ar^{\Delta \ten 1_A} "T";"L"};
    {\ar^{ \Delta} "L";"B"};
    {\ar_{ \Delta} "R";"B"};
 \endxy
\quad \quad
 \vcenter{\xy
 (-22,0)*+{\Bbbk \ten A}="L";
 (0,0)*+{A \ten A}="M";
 (22,0)*+{A \ten \Bbbk}="R";
 (0,16)*+{A}="B";
    {\ar_{\varepsilon \ten 1_A} "L";"M"};
    {\ar^{1_A \ten \varepsilon} "R";"M"};
    {\ar^{ } "L";"B"};
    {\ar_{ } "R";"B"};
    {\ar_{ \mu} "M";"B"};
 \endxy}
\]
    \item the Frobenius identities:
\[
 \xy
 (20,0)*+{A \ten A}="R";
 (0,12)*+{A \ten A \ten A}="T";
 (0,-12)*+{A}="B";
 (-20,0)*+{A \ten A}="L";
    {\ar^-{1_A \ten \mu} "T";"R"};
    {\ar^-{\delta \ten 1_A} "L";"T"};
    {\ar_-{ \mu} "L";"B"};
    {\ar_-{ \delta} "B";"R"};
 \endxy
\quad \quad
 \xy
 (20,0)*+{A \ten A}="R";
 (0,12)*+{A \ten A \ten A}="T";
 (0,-12)*+{A}="B";
 (-20,0)*+{A \ten A}="L";
    {\ar^-{\mu \ten 1_A} "T";"R"};
    {\ar^-{1_A \ten \delta } "L";"T"};
    {\ar_-{ \mu} "L";"B"};
    {\ar_-{ \delta} "B";"R"};
 \endxy
\]
\end{itemize}

\item There exists a copairing $\rho \maps \Bbbk \to A \ten A$
that equips $A$ with two equivalent comultiplications and counits.
That is to say, the following diagrams commute:

\begin{itemize}
 \item the equivalence of the two coalgebra structures:
 
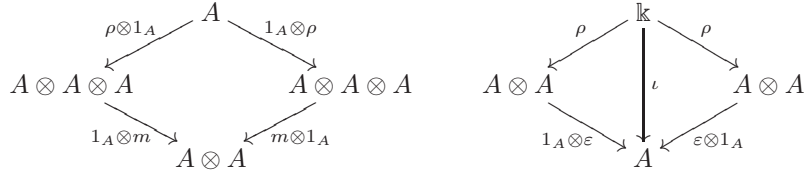
\[
 \xymatrix@R=1.2pc@C=.8pc{
 & A \ar[dr]^-{1_A \ten \rho} \ar[dl]_-{\rho \ten 1_A} \\
  A \ten A \ten A \ar[dr]_-{1_A \ten m} & &
A \ten A \ten A  \ar[dl]^-{m\ten 1_A} \\
 & A \ten A
 }
\qquad
 \xymatrix@R=1.2pc
 {
 & \Bbbk \ar[dr]^-{\rho} \ar[dl]_-{\rho} \ar[dd]^{\iota} \\
 A \ten A \ar[dr]_-{1_A \ten \varepsilon} && A \ten A \ar[dl]^-{\varepsilon \ten
 1_A}\\
 & A
 }
\]
\end{itemize}

\item The form $\varepsilon \maps A \to \Bbbk$ is nondegenerate.
Or more explicitly, the algebra $A$ is equipped with copairing
$\gamma \maps \Bbbk \to A \ten A$ making the following diagrams
commute:
\begin{itemize}
 \item the nondegeneracy of the pairing:
 
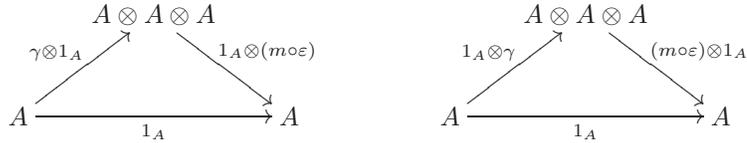
\[
 \xymatrix@C=1.5pc@R=.8pc{
  & A \ten A \ten A \ar[ddr]^-{1_A \ten (m \circ \varepsilon)} & \\
\\
A \ar[uur]^-{\gamma \ten 1_A} \ar[rr]_-{1_A} && A
 }
\qquad \qquad
 \xymatrix@C=1.5pc@R=.8pc{
    & A \ten A \ten A \ar[ddr]^-{(m \circ \varepsilon) \ten 1_A}& \\
\\
  A\ar[uur]^-{1_A \ten \gamma} \ar[rr]_-{1_A}&& A
 }
\]
\end{itemize}

\end{enumerate}
An algebra $A$ equipped with a form $\varepsilon \maps A \to
\Bbbk$ satisfying any of these equivalent conditions is called a
{\em Frobenius algebra}.
\end{prop}

To prove this proposition it is helpful to invoke a short hand
notation for the morphisms and axioms in each of the above
characterizations. Each formulation of Frobenius algebra has a
topological interpretation where the specified maps represent the
topological building blocks and the axioms are topologically
motivated.  We will see later on that these diagrams can be made
completely rigorous using string diagrams from 2-category theory.
Nevertheless, for the reader hesitant in using diagrams for
mathematical proofs, we will explain how to translate this
shorthand notation into the more traditional commutative diagrams.

To begin, notice that the multiplication $m \maps A \ten A \to A$
for the algebra structure takes as input the tensor product of two
copies of $A$ and outputs one copy of $A$.  Similarly, if we
regard the ground field $\Bbbk$ as the tensor product of no copies
of $A$, then the unit for the multiplication $\iota \maps \Bbbk
\to A$ takes no copies of $A$ and produces a single copy of $A$.
We draw these morphisms in our shorthand notation as follows:
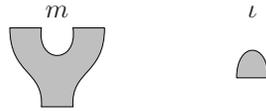
\[
  \begin{pspicture}(2.2,1.8)
  \rput(1,1.8){$m$}
  \rput(1,0){\mult}
 \end{pspicture}
 \qquad \quad
   \begin{pspicture}(2.2,1.8)
  \rput(1,1.8){$\iota$}
  \rput(1,.55){\birth}
 \end{pspicture}
\]
The pictures are read from top to bottom and the horizontal lines
at the top and bottom of each diagram are thought of as a copy of
$A$.  Since the ground field $\Bbbk$ is thought of as no copies of
$A$, it is represented by no horizontal line.  The identity map is
drawn as:
\[
  \begin{pspicture}[0.5](2.2,1.8)
  \rput(1,1.8){$1_A$}
  \rput(1,0){\ident}
 \end{pspicture}
\]
or sometimes as:
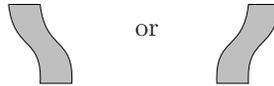
\[
  \begin{pspicture}[0.5](1.5,1.8)
  \rput(1,0){\curveleft}
 \end{pspicture}
 \qquad {\rm or} \qquad
   \begin{pspicture}[0.5](1.5,1.8)
  \rput(.4,0){\curveright}
 \end{pspicture}
\]
for aesthetic purposes.

Similarly, we draw the coalgebra maps as upsidedown versions of
the algebra maps:
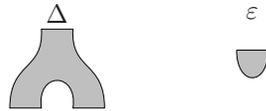
\[
  \begin{pspicture}(2.2,1.8)
  \rput(1,1.8){$\Delta$}
  \rput(1,0){\comult}
 \end{pspicture}
 \qquad \quad
   \begin{pspicture}(2.2,1.8)
  \rput(1,1.8){$\varepsilon$}
  \rput(1,1.1){\death}
 \end{pspicture}
\]
Notice from the diagram for the multiplication and
comultiplication that the tensor product $A \ten A$ is represented
by placing the diagrams side by side.  In a similar manner, we
draw the tensor product of morphisms by placing them side by side
in the diagrams:
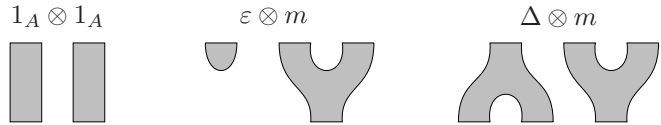
\[
\begin{pspicture}[0.5](2,2.3)
  \rput(.4,0){\ident} \rput(1.6,0){\ident}
  \rput(1,2){$1_A \ten 1_A$}
\end{pspicture}
\qquad \quad
\begin{pspicture}[0.5](3,2.3)
 \rput(2.6,0){\mult} \rput(.6,1.5){\death}
 \rput(1.6,2){$\varepsilon \ten m$}
\end{pspicture}
\qquad \qquad
\begin{pspicture}[0.5](3,2.3)
  \rput(3,0){\mult} \rput(1,0){\comult}
   \rput(2,2){$\Delta \ten m$}
\end{pspicture}
\]

To compose morphisms the diagrams are stacked on top of each
other.  Here are a few illustrative examples:
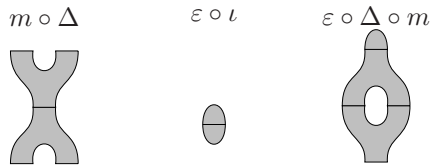
\[ \psset{xunit=.5cm,yunit=0.5cm}
\begin{pspicture}[0.4](2.2,4)
  \rput(1,0){\comult} \rput(1,1.5){\mult}
  \rput(1,3.9){$m \circ \Delta$}
\end{pspicture}
 \qquad \qquad
\begin{pspicture}[0.4](1,3.9)
  \rput(.5,1){\birth} \rput(.5,1){\death}
  \rput(.5,3.9){$\varepsilon \circ \iota$}
\end{pspicture}
\qquad \qquad
\begin{pspicture}[0.4](2.2,3.9)
  \rput(1,0){\mult} \rput(1,1.5){\comult} \rput(1,3){\birth}
  \rput(1,3.9){$\varepsilon \circ \Delta \circ m$}
\end{pspicture}
\]

Using these simple rules we can draw the axioms in the definition
of a Frobenius algebra.  The associativity axiom is depicted as:
\[ \psset{xunit=.5cm,yunit=0.5cm}
\begin{pspicture}[0.5](2,3)
  \rput(1,0){\mult} \rput(1.6,1.5){\curveright} \rput(.4,1.5){\mult}
\end{pspicture}
\qquad = \qquad
 \begin{pspicture}[0.5](2,3)
  \rput(1,0){\mult} \rput(1.6,1.5){\mult} \rput(.4,1.5){\curveleft}
\end{pspicture}
\]
and the unit laws as:
\[ \psset{xunit=.6cm,yunit=0.6cm}
\begin{pspicture}[0.3](2,3)
  \rput(1,0){\mult} \rput(1.6,1.5){\birth} \rput(.4,1.5){\smallident}
\end{pspicture}
\qquad = \quad
\begin{pspicture}[0.3](2,3)
  \rput(1,0){\longident}
\end{pspicture}
\quad = \qquad
\begin{pspicture}[0.3](2,3)
  \rput(1,0){\mult} \rput(1.6,1.5){\smallident} \rput(.4,1.5){\birth}
\end{pspicture}
\]

Figure~\ref{Frobdefs} summarizes the different formulations of
Frobenius algebra presented in Proposition~\ref{equivFrob}. Now we
are ready to prove the proposition.

\begin{figure}
{\bf \caption{Definitions of Frobenius Algebra}}

\begin{center}
{ \makebox[0pt]{
\begin{tabular}{|p{5.5em}|p{10.6em}|p{11em}|p{10em}|}
  \hline
  \bfseries Description & \textit{An algebra equipped with a coalgebra structure satisfying the Frobenius identities} & \textit{An algebra equipped with a copairing inducing two equivalent coalgebra structures}& \textit{An algebra with a nondegenerate form} \\
  \hline
 \bfseries  Morphisms & \[  \psset{xunit=.5cm,yunit=0.5cm}
\begin{pspicture}[0.5](1.2,1.6)
  \rput(.6,0){\mult}
\end{pspicture}
\; \;
\begin{pspicture}[0.5](1.2,1.6)
  \rput(.6,.55){\birth}
\end{pspicture}
\; \;
\begin{pspicture}[0.5](1.2,1.6)
  \rput(.6,0){\comult}
\end{pspicture}
\; \;
\begin{pspicture}[0.5](1,1.6)
  \rput(.6,1.1){\death}
\end{pspicture}
\]&
\[ \psset{xunit=.5cm,yunit=0.5cm}
\begin{pspicture}[0.5](1.2,1.6)
  \rput(.6,0){\mult}
\end{pspicture}
\; \;
\begin{pspicture}[0.5](1.2,1.6)
  \rput(.6,.55){\birth}
\end{pspicture}
\; \;
\begin{pspicture}[0.5](1.2,1.6)
  \rput(.6,.25){\zag}
\end{pspicture}
\;\;
\begin{pspicture}[0.5](1,1.6)
  \rput(.6,1){\death}
\end{pspicture}
\]
 &
  \[  \psset{xunit=.5cm,yunit=0.5cm}
\begin{pspicture}[0.5](1.2,1.6)
  \rput(.6,0){\mult}
\end{pspicture}
\; \;
\begin{pspicture}[0.5](1.2,1.6)
  \rput(.6,.55){\birth}
\end{pspicture}
\; \;
\begin{pspicture}[0.5](1.2,1.6)
  \rput(.6,.25){\zag}
\end{pspicture}
\;\;
\begin{pspicture}[0.5](1,1.6)
  \rput(.6,1){\death}
\end{pspicture}
\]
\\ \hline
 \bfseries Axioms &
 \qquad \textit{algebra axioms and} &
 \qquad \textit{algebra axioms and} &
 \qquad \textit{algebra axioms and} \\
 & \begin{center} $ \psset{xunit=.45cm,yunit=0.45cm}
\begin{pspicture}[0.5](2,3)
  \rput(.8,1.5){\comult} \rput(2,0){\curveleft} \rput(.2,0){\comult}
\end{pspicture}
\; \; = \; \;
\begin{pspicture}[0.5](2,3)
  \rput(1.6,1.5){\comult} \rput(2.2,0){\comult} \rput(.4,0){\curveright}
\end{pspicture}
$ \end{center} &
\begin{center} $ \psset{xunit=.45cm,yunit=0.45cm}
\begin{pspicture}[0.3](3,3)
  \rput(1,0){\mult} \rput(2.2,1.5){\zag}
  \rput(.4,1.5){\ident} \rput(2.8,0){\ident}
\end{pspicture}
\; \; = \; \;
\begin{pspicture}[0.3](3,3)
  \rput(2.2,0){\mult} \rput(1,1.5){\zag}
  \rput(2.8,1.5){\ident} \rput(.4,0){\ident}
\end{pspicture}
$ \end{center}
  &
\begin{center} $ \psset{xunit=.45cm,yunit=0.45cm}
\begin{pspicture}[0.3](3,3)
  \rput(1,.55){\mult} \rput(2.2,2.05){\zag}
  \rput(.4,2.05){\ident} \rput(2.8,.55){\ident}
  \rput(2.8,0){\smallident} \rput(1,.55){\death}
\end{pspicture}
\quad  = \;
\begin{pspicture}[0.3](2,3.6)
  \rput(1,2.05){\ident} \rput(1,.55){\ident} \rput(1,0){\smallident}
\end{pspicture}
$ \end{center}
 \\
  &  \begin{center} $ \psset{xunit=.45cm,yunit=0.45cm}
\begin{pspicture}[0.5](2,2.2)
  \rput(1,.55){\comult} \rput(1.6,.55){\death} \rput(.4,.0){\smallident}
\end{pspicture}
\; = \;
\begin{pspicture}[0.5](.7,2.2)
  \rput(.3,0){\longident}
\end{pspicture}
\; = \;
\begin{pspicture}[0.5](2,2.2)
   \rput(1,.55){\comult} \rput(.4,.55){\death} \rput(1.6,0){\smallident}
\end{pspicture}
$ \end{center}
  &
 \begin{center} $ \psset{xunit=.45cm,yunit=0.45cm}
\begin{pspicture}[0.3](2,2.2)
 \rput(1,.55){\zag}
 \rput(.4,.55){\death}
 \rput(1.6,0){\smallident}
\end{pspicture}
\; = \;
\begin{pspicture}[0.3](1,2.2)
 \rput(.5,0){\smallident}
 \rput(.5,.55){\birth}
\end{pspicture}
\; = \;
\begin{pspicture}[0.3](2,2.2)
 \rput(1,.55){\zag}
 \rput(1.6,.55){\death}
 \rput(.4,0){\smallident}
\end{pspicture}
$ \end{center}
  &
\begin{center} $ \psset{xunit=.45cm,yunit=0.45cm}
\begin{pspicture}[0.3](2,2.2)
  \rput(1,1.05){\ident} \rput(1,-1){\smallident} \rput(1,-.45){\ident}
\end{pspicture}
\; = \quad
\begin{pspicture}[.3](3,2.2)
  \rput(2.2,-.45){\mult} \rput(1,1.05){\zag}
  \rput(2.8,1.05){\ident} \rput(.4,-1){\smallident} \rput(.4,-.45){\ident}
  \rput(2.2,-.45){\death}
\end{pspicture}
$ \end{center}
  \\
& \begin{center} $ \psset{xunit=.45cm,yunit=0.45cm}
\begin{pspicture}[.5](3,3.2)
 \rput(1,1.5){\comult} \rput(2.2,0){\mult}
 \rput(2.8,1.5){\ident} \rput(.4,0){\ident}
\end{pspicture}
\quad = \;
\begin{pspicture}[0.5](2,3.2)
 \rput(1,0){\comult} \rput(1,1.5){\mult}
\end{pspicture}
$ \end{center} &
& \\
& \begin{center} $ \psset{xunit=.45cm,yunit=0.45cm}
\begin{pspicture}[.5](3,3.2)
 \rput(2.2,1.5){\comult} \rput(1,0){\mult}
 \rput(.4,1.5){\ident} \rput(2.8,0){\ident}
\end{pspicture}
\quad = \;
\begin{pspicture}[0.5](2,3.2)
 \rput(1,0){\comult} \rput(1,1.5){\mult}
\end{pspicture}
$ \end{center}& & \\
  \hline
\end{tabular} }}
\end{center}
\label{Frobdefs}
\end{figure}

\paragraph{Proof of Proposition~\ref{equivFrob}. }
(i.\ $\Rightarrow$ iii.) Given an algebra and coalgebra structure
on $A$ define the map $\gamma \maps \Bbbk \to A \ten A$ as
follows:
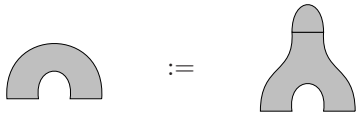
\[
\begin{pspicture}[0.3](2,1.5)
  \rput(1,0){\zag}
\end{pspicture}
\qquad := \qquad
\begin{pspicture}[0.3](2,2.3)
  \rput(1,0){\comult} \rput(1,1.5){\birth}
\end{pspicture}
\]
This copairing makes the form $\varepsilon$ nondegenerate because
the equalities:
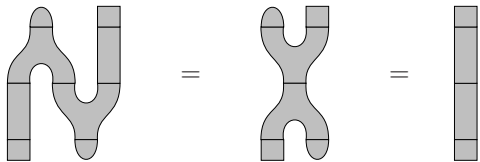
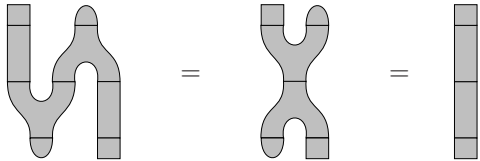
\[ \psset{xunit=.5cm,yunit=0.5cm}
\begin{pspicture}[0.5](3.1,4.2)
  \rput(2.2,.55){\mult}
  \rput(1,2.05){\comult}
  \rput(1,3.55){\birth}
  \rput(2.8,2.05){\ident}
  \rput(2.8,3.55){\smallident}
  \rput(.4,.55){\ident}
  \rput(.4,0){\smallident}
  \rput(2.2,.55){\death}
\end{pspicture}
\qquad = \quad
\begin{pspicture}[0.5](2.5,4.2)
  \rput(1.6,2.05){\mult}
  \rput(1.6,.555){\comult}
  \rput(1,3.55){\birth}
  \rput(2.2,3.55){\smallident}
  \rput(1,0){\smallident}
  \rput(2.2,.55){\death}
\end{pspicture}
\qquad = \quad
\begin{pspicture}[0.5](1.2,4.2)
  \rput(.6,.55){\ident}
  \rput(.6,2.05){\ident}
  \rput(.6,3.55){\smallident}
  \rput(.6,0){\smallident}
\end{pspicture}
\]
and
\[ \psset{xunit=.5cm,yunit=0.5cm}
\begin{pspicture}[0.5](3.1,4.2)
  \rput(1,.55){\mult}
  \rput(2.2,2.05){\comult}
  \rput(2.2,3.55){\birth}
  \rput(.4,2.05){\ident}
  \rput(.4,3.55){\smallident}
  \rput(2.8,.55){\ident}
  \rput(2.8,0){\smallident}
  \rput(1,.55){\death}
\end{pspicture}
\qquad = \quad
\begin{pspicture}[0.5](2.5,4.2)
  \rput(1.6,2.05){\mult}
  \rput(1.6,.555){\comult}
  \rput(2.2,3.55){\birth}
  \rput(1,3.55){\smallident}
  \rput(2.2,0){\smallident}
  \rput(1,.55){\death}
\end{pspicture}
\qquad = \quad
\begin{pspicture}[0.5](1.2,4.2)
  \rput(.6,.55){\ident}
  \rput(.6,2.05){\ident}
  \rput(.6,3.55){\smallident}
  \rput(.6,0){\smallident}
\end{pspicture}
\]
follow from the Frobenius identities.

(iii.\ $\Rightarrow$ ii.) We only need to check that the axioms of
ii.\ are satisfied since ii.\ and iii.\ have the same morphisms.
The first axiom is proved as follows:
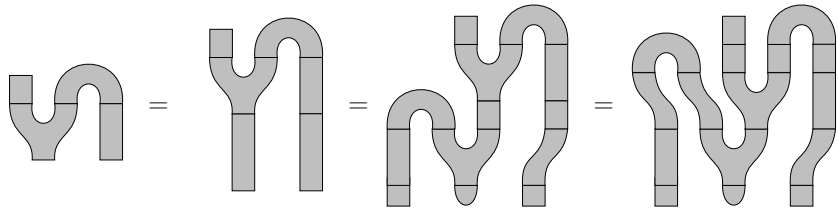
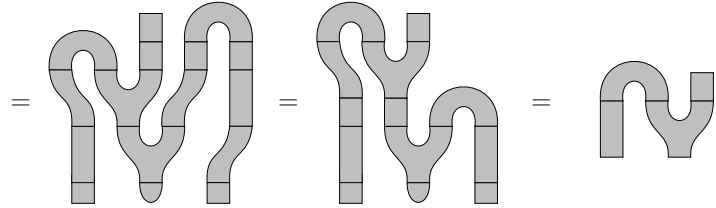
\[ \psset{xunit=.5cm,yunit=0.5cm}
   \begin{pspicture}[0.5](3.6,2.55)
  \rput(1.4,0){\mult}
  \rput(3.2,0){\ident}
  \rput(2.6,1.5){\zag}
  \rput(.8,1.5){\medident}
 \end{pspicture}
  \;\; = \;\;
   \begin{pspicture}[0.5](3.6,4.2)
  \rput(1.4,0){\longident}
  \rput(1.4,2.05){\mult}
  \rput(2.6,3.55){\zag}
  \rput(.8,3.55){\medident}
  \rput(3.2,2.05){\ident}
  \rput(3.2,0){\longident}
 \end{pspicture}
 \;\; = \;\;
   \begin{pspicture}[0.5](4.8,5)
  \rput(2,.55){\death}
  \rput(.2,0){\medident}
  \rput(.2,.55){\ident}
  \rput(2,.55){\mult}
  \rput(2.6,2.05){\medident}
  \rput(2.6,2.8){\mult}
  \rput(.8,2.05){\zag}
  \rput(3.8,4.3){\zag}
  \rput(2,4.3){\medident}
  \rput(4.4,2.8){\ident}
  \rput(4.4,2.05){\medident}
  \rput(3.8,.55){\curveright}
  \rput(3.8,0){\smallident}
 \end{pspicture}
 \;\; = \;\;
   \begin{pspicture}[0.5](5.4,5)
  \rput(2.6,.55){\death}
  \rput(.8,0){\medident}
  \rput(.8,.55){\ident}
  \rput(2.6,.55){\mult}
  \rput(3.8,3.55){\medident}
  \rput(2.6,3.55){\medident}
  \rput(3.2,2.05){\mult}
  \rput(.8,3.55){\zag}
  \rput(.8,2.05){\curveleft}
  \rput(2,2.05){\curveleft}
  \rput(4.4,4.3){\zag}
  \rput(2.6,4.3){\medident}
  \rput(5,3.55){\medident}
  \rput(5,2.05){\ident}
  \rput(4.4,.55){\curveright}
  \rput(4.4,0){\smallident}
 \end{pspicture}
\]

\[ \psset{xunit=.5cm,yunit=0.5cm} \hspace{.8in} = \; \;
   \begin{pspicture}[0.5](5.4,5)
  \rput(2.6,.55){\death}
  \rput(.8,0){\medident}
  \rput(.8,.55){\ident}
  \rput(2.6,.55){\mult}
  \rput(3.8,3.55){\medident}
  \rput(2.6,3.55){\medident}
  \rput(2,2.05){\mult}
  \rput(.8,3.55){\zag}
  \rput(.8,2.05){\curveleft}
  \rput(3.2,2.05){\curveright}
  \rput(4.4,4.3){\zag}
  \rput(2.6,4.3){\medident}
  \rput(5,3.55){\medident}
  \rput(5,2.05){\ident}
  \rput(4.4,.55){\curveright}
  \rput(4.4,0){\smallident}
 \end{pspicture}
   \;\; = \;\;
   \begin{pspicture}[0.5](5.4,5)
  \rput(2.6,.55){\death}
  \rput(.8,0){\medident}
  \rput(.8,.55){\ident}
  \rput(2.6,.55){\mult}
  \rput(2.6,4.3){\medident}
  \rput(2,2.05){\medident}
  \rput(2,2.8){\mult}
  \rput(.8,4.3){\zag}
  \rput(.8,2.8){\curveleft}
  \rput(.8,2.05){\medident}
  \rput(3.8,2.05){\zag}
  \rput(4.4,0){\smallident}
  \rput(4.4,.55){\ident}
 \end{pspicture}
   =
   \begin{pspicture}[0.5](3.6,2.55)
  \rput(3.2,0){\mult}
  \rput(1.4,0){\ident}
  \rput(2,1.5){\zag}
  \rput(3.8,1.5){\medident}
 \end{pspicture}
\]
For the reader who finds these topological manipulations a bit too
cavalier, we will translate this proof into a traditional
commutative diagram.  Although, the equivalence between the
algebraic and topological pictures will be established later on,
an explicit example will make this relationship more apparent.

We want to show that the map
\[
 \xymatrix@C=2.5pc{
 A \ar[r]^-{1_A \ten \rho} & A \ten A \ten A \ar[r]^-{m \ten 1_A} &
 A \ten A
 }
\]
is equal to the map
\[
 \xymatrix@C=2.5pc{
 A \ar[r]^-{\rho \ten 1_A} & A \ten A \ten A \ar[r]^-{1_A \ten m} &
 A \ten A.
 }
\]
In the diagrammatic proof, the first equality, where the bottom of
the diagram is stretched, corresponds to composing with the
identity $1_A \ten 1_A$. The second equality is an application of
the nondegeneracy axiom:
\[
 \xy
   (-40,15)*+{A}="t1";
   (-10,15)*+{A^{\ten 3}}="t2";
   (20,15)*+{A^{\ten 2}}="t3";
   (60,15)*+{A^{\ten 2}}="t4";
   (20,0)*+{A^{\ten 4}}="m1";
   (40,0)*+{A^{\ten 3}}="m2";
        {\ar^{1_A \ten \rho} "t1";"t2" };
        {\ar^{m \ten 1_A} "t2";"t3" };
        {\ar^{1_A \ten 1_A} "t3";"t4" };
        {\ar_{\rho \ten 1_A \ten 1_A} "t3";"m1" };
        {\ar_{ 1_A \ten m \ten 1_A} "m1";"m2" };
        {\ar_{ 1_A \ten \varepsilon \ten 1_A} "m2";"t4" };
 \endxy
\]
where $A^{\ten n}$ is short hand for the $n$-fold tensor product
of $A$.

The next equality holds in any monoidal category as a consequence
of the tensor product being functor. What it amounts to is the
fact that in a monoidal category it does not matter in which order
we apply maps between the tensor product of objects. In
particular, whenever we have morphisms $f \maps A \to A'$ and $g
\maps B \to B'$ in a monoidal category we get a commuting square:
\[
 \xy
    (-14,10)*+{A \ten B}="TL";
    (14,10)*+{A \ten B'}="TR";
    (14,-10)*+{A' \ten B'}="BR";
    (-14,-10)*+{A' \ten B}="BL";
    {\ar^{ A \ten g} "TL";"TR"};
    {\ar^{ f \ten B'} "TR";"BR"};
    {\ar_{ f \ten B} "TL";"BL"};
    {\ar_{ A' \ten g} "BL";"BR"};
 \endxy
\]
In our case we are interested in the special case when $f = \rho$
and $g = m \ten 1_A$.  Hence our commutative diagram becomes:
\[
 \xy
   (-40,15)*+{A}="t1";
   (-10,15)*+{A^{\ten 3}}="t2";
   (20,15)*+{A^{\ten 2}}="t3";
   (60,15)*+{A^{\ten 2}}="t4";
   (20,0)*+{A^{\ten 4}}="m1";
   (40,0)*+{A^{\ten 3}}="m2";
    (-10,0)*+{A^{\ten 5}}="m0";
        {\ar^{1_A \ten \rho} "t1";"t2" };
        {\ar^{m \ten 1_A} "t2";"t3" };
        {\ar^{1_A^{\ten 2}} "t3";"t4" };
        {\ar_{\rho \ten 1_A^{\ten 2}} "t3";"m1" };
        {\ar_{ 1_A \ten m \ten 1_A} "m1";"m2" };
        {\ar_{ 1_A \ten \varepsilon \ten 1_A} "m2";"t4" };
        {\ar_{\rho \ten  1_A^{\ten 3} } "t2";"m0" };
        {\ar_{ 1_A^{\ten 2} \ten m \ten 1_A} "m0";"m1" };
 \endxy
\]
Next we use the associativity of the multiplication:
\[
 \xy
   (-40,15)*+{A}="t1";
   (-10,15)*+{A^{\ten 3}}="t2";
   (20,15)*+{A^{\ten 2}}="t3";
   (60,15)*+{A^{\ten 2}}="t4";
   (20,0)*+{A^{\ten 4}}="m1";
   (40,0)*+{A^{\ten 3}}="m2";
    (-10,0)*+{A^{\ten 5}}="m0";
    (20,-10)*+{A^{\ten 4}}="b1";
        {\ar^{1_A \ten \rho} "t1";"t2" };
        {\ar^{m \ten 1_A} "t2";"t3" };
        {\ar^{1_A^{\ten 2}} "t3";"t4" };
        {\ar_{\rho \ten 1_A^{\ten 2}} "t3";"m1" };
        {\ar^{ 1_A \ten m \ten 1_A} "m1";"m2" };
        {\ar_{ 1_A \ten \varepsilon \ten 1_A} "m2";"t4" };
        {\ar_{\rho \ten  1_A^{\ten 3} } "t2";"m0" };
        {\ar^{ 1_A^{\ten 2} \ten m \ten 1_A} "m0";"m1" };
        {\ar_-{1_A \ten m \ten 1_A^{\ten 2}  \; \; \;} "m0";"b1" };
        {\ar_{ \; \; 1_A \ten m \ten 1_A} "b1";"m2" };
 \endxy
\]
Now we use two applications of the property of monoidal categories
mentioned above:
\[
 \xy
   (-40,15)*+{A}="t1";
   (-10,15)*+{A^{\ten 3}}="t2";
   (20,15)*+{A^{\ten 2}}="t3";
   (60,15)*+{A^{\ten 2}}="t4";
   (20,0)*+{A^{\ten 4}}="m1";
   (40,0)*+{A^{\ten 3}}="m2";
    (-10,0)*+{A^{\ten 5}}="m0";
    (-40,0)*+{A^{\ten 3}}="m00";
    (20,-10)*+{A^{\ten 4}}="b1";
    (0,-20)*+{A^{\ten 2}}="b0";
        {\ar^{1_A \ten \rho} "t1";"t2" };
        {\ar^{m \ten 1_A} "t2";"t3" };
        {\ar^{1_A^{\ten 2}} "t3";"t4" };
        {\ar_{\rho \ten 1_A^{\ten 2}} "t3";"m1" };
        {\ar^{ 1_A \ten m \ten 1_A} "m1";"m2" };
        {\ar^{ 1_A \ten \varepsilon \ten 1_A} "m2";"t4" };
        {\ar_{\rho \ten  1_A^{\ten 3} } "t2";"m0" };
        {\ar^{ 1_A^{\ten 2} \ten m \ten 1_A} "m0";"m1" };
        {\ar_-{1_A \ten m \ten 1_A^{\ten 2}  \; \; \;} "m0";"b1" };
        {\ar_{ \; \; 1_A \ten m \ten 1_A} "b1";"m2" };
        {\ar_{ \rho \ten 1_A} "t1";"m00" };
        {\ar^{ 1_A^{\ten 3} \ten \rho} "m00";"m0" };
        {\ar_{1_A \ten m} "m00";"b0" };
        {\ar_{1_A^{\ten 2} \ten \rho} "b0";"b1" };
 \endxy
\]
Finally, we apply the other nondegeneracy axiom.
\[
 \xy
   (-40,15)*+{A}="t1";
   (-10,15)*+{A^{\ten 3}}="t2";
   (20,15)*+{A^{\ten 2}}="t3";
   (60,15)*+{A^{\ten 2}}="t4";
   (20,0)*+{A^{\ten 4}}="m1";
   (40,0)*+{A^{\ten 3}}="m2";
    (-10,0)*+{A^{\ten 5}}="m0";
    (-40,0)*+{A^{\ten 3}}="m00";
    (20,-10)*+{A^{\ten 4}}="b1";
    (0,-20)*+{A^{\ten 2}}="b0";
        {\ar^{1_A \ten \rho} "t1";"t2" };
        {\ar^{m \ten 1_A} "t2";"t3" };
        {\ar^{1_A^{\ten 2}} "t3";"t4" };
        {\ar_{\rho \ten 1_A^{\ten 2}} "t3";"m1" };
        {\ar^{ 1_A \ten m \ten 1_A} "m1";"m2" };
        {\ar^{ 1_A \ten \varepsilon \ten 1_A} "m2";"t4" };
        {\ar_{\rho \ten  1_A^{\ten 3} } "t2";"m0" };
        {\ar^{ 1_A^{\ten 2} \ten m \ten 1_A} "m0";"m1" };
        {\ar_-{1_A \ten m \ten 1_A^{\ten 2}  \; \; \;} "m0";"b1" };
        {\ar_{ \; \; 1_A \ten m \ten 1_A} "b1";"m2" };
        {\ar_{ \rho \ten 1_A} "t1";"m00" };
        {\ar^{ 1_A^{\ten 3} \ten \rho} "m00";"m0" };
        {\ar_{1_A \ten m} "m00";"b0" };
        {\ar_{1_A^{\ten 2} \ten \rho} "b0";"b1" };
        {\ar@/_4pc/_{1_A^{\ten 2}} "b0";"t4" };
 \endxy
\]

For the second axiom of ii.\ we use the unit laws for the algebra
together with the nondegeneracy axioms.

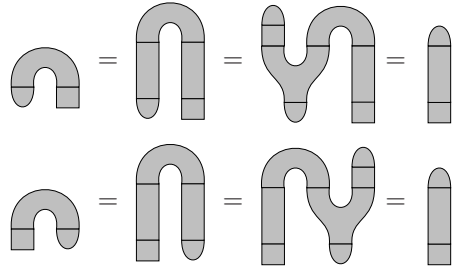
\[ \psset{xunit=.5cm,yunit=0.5cm}
\begin{pspicture}[0.5](2,2.2)
 \rput(1,.55){\zag}
 \rput(.4,.55){\death}
 \rput(1.6,0){\smallident}
\end{pspicture}
\; = \;
\begin{pspicture}[0.5](2,2.8)
 \rput(.4,.55){\ident}
 \rput(1,2.05){\zag}
 \rput(.4,.55){\death}
 \rput(1.6,0){\smallident}
  \rput(1.6,.55){\ident}
\end{pspicture}
\; = \;\begin{pspicture}[0.5](3,3)
  \rput(1,.55){\mult} \rput(2.2,2.05){\zag}
  \rput(.4,2.05){\smallident} \rput(.4,2.6){\birth} \rput(2.8,.55){\ident}
  \rput(2.8,0){\smallident} \rput(1,.55){\death}
\end{pspicture}
\; = \;\begin{pspicture}[0.5](1,3)
  \rput(.5,0){\smallident}  \rput(.5,.55){\ident}  \rput(.5,2.05){\birth}
\end{pspicture}
\]
\[ \psset{xunit=.5cm,yunit=0.5cm}
\begin{pspicture}[0.5](2,2.2)
 \rput(1,.55){\zag}
 \rput(1.6,.55){\death}
 \rput(.4,0){\smallident}
\end{pspicture}
\; = \;
\begin{pspicture}[0.5](2,2.8)
 \rput(.4,.55){\ident}
 \rput(1,2.05){\zag}
 \rput(1.6,.55){\death}
 \rput(.4,0){\smallident}
  \rput(1.6,.55){\ident}
\end{pspicture}
\; = \;\begin{pspicture}[0.5](3,3)
  \rput(2.2,.55){\mult} \rput(2.2,.55){\death} \rput(1,2.05){\zag}
  \rput(2.8,2.05){\smallident} \rput(2.8,2.6){\birth} \rput(.4,.55){\ident}
  \rput(.4,0){\smallident}
\end{pspicture}
\; = \;\begin{pspicture}[0.5](1,3)
  \rput(.5,0){\smallident}  \rput(.5,.55){\ident}  \rput(.5,2.05){\birth}
\end{pspicture}
\]

(ii.\ $\Rightarrow$ i.) To show that ii.\ implies i.\ we must
define the comultiplication $\Delta \maps A \to A \ten A$ map.  We
do this as follows:
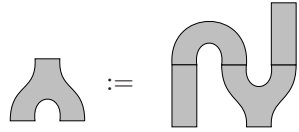
\[ \psset{xunit=.55cm,yunit=0.55cm}
\begin{pspicture}[0.5](1.6,1.5)
  \rput(1,0){\comult}
\end{pspicture}
\quad := \quad
\begin{pspicture}[0.3](3,3)
  \rput(2.2,0){\mult} \rput(1,1.5){\zag}
  \rput(2.8,1.5){\ident} \rput(.4,0){\ident}
\end{pspicture}
\]
We leave it as an exercise to the reader to verify that this
defines a coalgebra structure on $A$. We will verify the Frobenius
identities. The first Frobenius identity is proved by the string
of equalities:
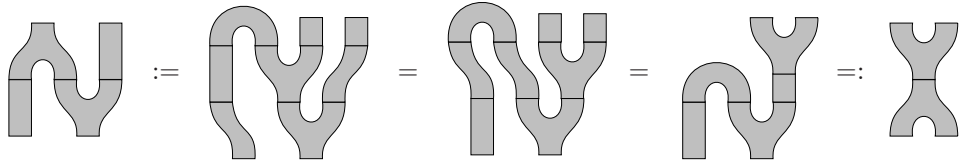
\begin{center}
 \makebox[0em]{
$ \psset{xunit=.5cm,yunit=0.5cm}
\begin{pspicture}[1](3,3)
 \rput(1,3){\comult} \rput(2.2,1.5){\mult}
 \rput(2.8,3){\ident} \rput(.4,1.5){\ident}
\end{pspicture}
\quad := \quad
\begin{pspicture}[0.5](4,4.2)
  \rput(2.6,0){\mult}
  \rput(3.2,1.5){\curveright}
  \rput(2,1.5){\mult}
  \rput(.2,1.5){\ident}
  \rput(.8,0){\curveleft}
  \rput(.8,3){\zag}
  \rput(2.6,3){\medident}
  \rput(3.8,3){\medident}
\end{pspicture}
\quad = \quad
    \begin{pspicture}[0.5](3.6,4)
  \rput(2.6,0){\mult}
  \rput(.8,0){\ident}
  \rput(3.2,1.5){\mult}
  \rput(2,1.5){\curveleft}
  \rput(.8,1.5){\curveleft}
  \rput(3.8,3){\medident}
  \rput(2.6,3){\medident}
  \rput(.8,3){\zag}
 \end{pspicture}
\quad \; = \;
\begin{pspicture}[0.5](4,4.2)
  \rput(2.6,0){\mult}
  \rput(3.2,1.5){\medident}
  \rput(3.2,2.25){\mult}
  \rput(.8,0){\ident}
  \rput(1.4,1.5){\zag}
\end{pspicture}
\; \; =: \; \;
\begin{pspicture}[0.5](2,3)
 \rput(1,0){\comult} \rput(1,1.5){\mult}
\end{pspicture}
$}
\end{center}
which follows from the associativity of the multiplication and the
axioms of a monoidal category.  With just a little more work, the
second Frobenius identity is proved similarly:
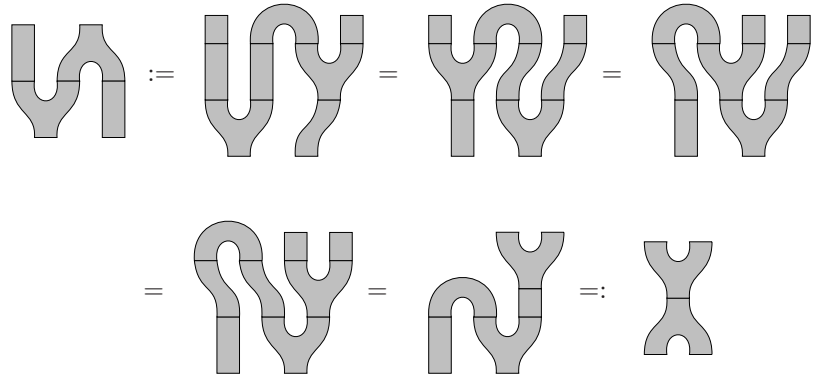
\[ \psset{xunit=.5cm,yunit=0.5cm}
\begin{pspicture}[1](3,3)
 \rput(2.4,3){\comult} \rput(1.2,1.5){\mult}
 \rput(.6,3){\ident} \rput(3,1.5){\ident}
\end{pspicture}
\quad := \quad
   \begin{pspicture}[0.5](3.6,4)
  \rput(.8,0){\mult}
  \rput(2.6,0){\curveright}
  \rput(3.2,1.5){\mult}
  \rput(.2,1.5){\ident}
  \rput(1.4,1.5){\ident}
  \rput(3.8,3){\medident}
  \rput(.2,3){\medident}
  \rput(2,3){\zag}
 \end{pspicture}
\quad = \quad
   \begin{pspicture}[0.5](3.6,4)
  \rput(.8,0){\ident}
  \rput(2.6,0){\mult}
  \rput(2,1.5){\curveright}
  \rput(3.2,1.5){\curveright}
  \rput(.8,1.5){\mult}
  \rput(3.8,3){\medident}
  \rput(.2,3){\medident}
  \rput(2,3){\zag}
 \end{pspicture}
 \quad = \quad
    \begin{pspicture}[0.5](3.6,4)
  \rput(2.6,0){\mult}
  \rput(.8,0){\ident}
  \rput(3.2,1.5){\curveright}
  \rput(2,1.5){\mult}
  \rput(.8,1.5){\curveleft}
  \rput(3.8,3){\medident}
  \rput(2.6,3){\medident}
  \rput(.8,3){\zag}
 \end{pspicture}
 \]
 \[ \psset{xunit=.5cm,yunit=0.5cm}
 \qquad \; \;= \quad
    \begin{pspicture}[0.5](3.6,4)
  \rput(2.6,0){\mult}
  \rput(.8,0){\ident}
  \rput(3.2,1.5){\mult}
  \rput(2,1.5){\curveleft}
  \rput(.8,1.5){\curveleft}
  \rput(3.8,3){\medident}
  \rput(2.6,3){\medident}
  \rput(.8,3){\zag}
 \end{pspicture}
  \quad = \; \;
     \begin{pspicture}[0.5](3.6,4)
  \rput(2.6,0){\mult}
  \rput(.8,0){\ident}
  \rput(3.2,1.5){\medident}
  \rput(1.4,1.5){\zag}
  \rput(3.2,2.25){\mult}
 \end{pspicture}
 \quad =: \quad
\begin{pspicture}[0.5](2,3)
 \rput(1,0){\comult} \rput(1,1.5){\mult}
\end{pspicture}
\]
\qed

The observant reader will have noticed that in proving
Proposition~\ref{equivFrob} we never used the fact that $A$ was a
vector space.  In fact, if we translate all of the pictures above
into commutative diagrams it is clear that the proof of
Proposition~\ref{equivFrob} relied only on the abstract properties
of the maps in each characterization.  This means that if we were
to take all of the diagrams in the discussion above and place them
in some other category $\mathcal{C}$ `sufficiently like
$\cat{Vect}$' then the proof would still be valid, and we can
define Frobenius algebras in the category $\mathcal{C}$ using any
of the above descriptions. This process of writing mathematical
objects using only commutative diagrams and placing them in other
categories where they make sense is what category theorists call
{\em internalization}.  An internalized Frobenius algebra is
sometimes referred to as a {\em Frobenius object}.

Let's consider what kind of additional structure the category
$\mathcal{C}$ must have in order to be `sufficiently like
$\cat{Vect}$', that is, in order to define a Frobenius object in
$\mathcal{C}$.  Since all of our diagrams required the tensor
product $A \ten A$ our category $\mathcal{C}$ should have a
multiplication.  Furthermore, since we used the `unit vector
space', or the ground field $\Bbbk$, our category $\mathcal{C}$
should have a unit for the above multiplication.  But all this
just amounts to the definition of a monoidal category. So the
notion of a Frobenius algebra makes sense in any monoidal category
$\mathcal{C}$, and our diagrammatic proof of
Proposition~\ref{equivFrob} shows that all three characterizations
are equivalent in $\mathcal{C}$. Eventually we will provide yet
another equivalent definition of Frobenius algebra in terms of
adjunctions, but first we will need to develop some categorical
language.

\subsection{String diagrams for 2-categories} \label{secString}

Now we begin the process of rephrasing the definition of a
Frobenius algebra in the language of higher-dimensional category
theory.  Our goal is to define Frobenius algebras using the notion
of an adjunction in a 2-category.  In order to understand the
relationship between adjunctions and Frobenius algebras we will
use string diagrams from 2-category theory.  Once we have defined
Frobenius algebras in terms of adjunctions we can then study the
corresponding string diagrams and the topology of 2-dimensional
thick tangles will begin to emerge.

Before we get too far ahead of ourselves we recall the definition
of a 2-category. Speaking colloquially the idea is that, just as
categories have objects, morphisms between objects, and various
axioms regarding composites and identities, 2-categories have
objects, morphisms between objects, and \textit{2-morphisms}
between morphisms together with axioms for the composites and
identities of both morphisms and 2-morphisms.  What makes
2-categories so interesting is that the axioms for the 1-morphisms
can now hold only up to coherent isomorphism.  Thus, rather than
having composition be associative on the nose, we can instead
require that composition be associative up to isomorphism
satisfying laws of its own. Similarly with the identity
constraints.  When all the axioms of a 2-category hold up to
isomorphism it is called a weak 2-category or
bicategory~\cite{Benabou}.

Below we discuss string diagrams for \textit{strict} 2-categories.
That is, 2-categories where composition is strictly associative
and the identity constraints hold as equations.  We justify this
choice with two reasons. First, there is a coherence theorem for
bicategories which states that every bicategory is biequivalent to
a strict 2-category~\cite{Street2}. Biequivalence is a notion of
equivalence between bicategories in which the usual axioms of an
equivalence only hold up to coherent isomorphism. This notion of
equivalence is often the most natural one to use between
bicategories.  The second reason for considering strict
2-categories rather than bicategories is because this notion
provides the correct framework to understand the topology of
2-dimensional thick tangles.  However, in Section \ref{secThick}
we will see that the topology of 3-dimensional thick tangles does
require the more general notion of a bicategory.

Typically, the objects of a category are represented geometrically
as little bullets, and the morphisms of a category as arrows
between the bullets:
\[
\xy
 (-7,0)*{\bullet}="1";
 (7,0)*{\bullet}="2";
    "1"; "2" **\dir{-} ?(.55)*\dir{>};
(-7,2.5)*{A}; (0,2.5)*{\scs F};
 (7.5,2.5)*{B};
\endxy .
\]
The composite of two morphisms is usually drawn as:
\[
\xy
 (-7,0)*{\bullet}="1";
 (7,0)*{\bullet}="2";
 (21,0)*{\bullet}="3";
    "1"; "2" **\dir{-} ?(.55)*\dir{>};
    "2"; "3" **\dir{-} ?(.55)*\dir{>};
(-7,2.5)*{A}; (0,2.5)*{\scs F};
 (7.5,2.5)*{B}; (14,2.5)*{\scs G};
 (21.5,2.5)*{C};
\endxy .
\]
In a 2-category we have objects and morphisms as before, but now
there are 2-morphisms going between morphisms:
\[
 \xy
  (-15,0)*{\bullet}="1";
  (0,0)*{\bullet}="2";
 "1";"2" **\crv{(-12,9) & (-3,9)};
  "1";"2" **\crv{(-12,-9) & (-3,-9)};
    (-7.5,6.75)*{\scriptstyle >}+(0,3)*{\scriptstyle F};
    (-7.5,-6.75)*{\scriptstyle >}+(0,-3)*{\scriptstyle F'};
    (-17.5,0)*{\scriptstyle A};
  (2.5,0)*{\scriptstyle B};
  {\ar@{=>}_{\scriptstyle \alpha}(-7.5,3)*{};(-7.5,-3)*{}} ;
 \endxy .
 \]
There are two types of composites coming from
 the two ways we can glue these pictures together.  We
 have a {\em vertical composite}:
 \[
 \xy
  (-8,0)*{\bullet}="1";
  (8,0)*{\bullet}="2";
 "1";"2" **\crv{(-5,9) & (5,9)};
  "1";"2" **\crv{(-5,-9) & (5,-9)};
   "1";"2" **\dir{-} ?(.55)*\dir{>};
    (,6.75)*{\scriptstyle >}+(0,3)*{\scriptstyle F};
    (0,-6.75)*{\scriptstyle >}+(0,-3)*{\scriptstyle F'};
    (-10.5,0)*{\scriptstyle A};
  (10.5,0)*{\scriptstyle B};
  {\ar@{=>}_{\scriptstyle \alpha}(0,5)*{};(0,1.5)*{}} ;
  {\ar@{=>}_{\scriptstyle \alpha'}(0,-1.5)*{};(0,-5)*{}} ;
 \endxy
 \]
 and a {\em horizontal composite}:
 \[
 \xy
  (-15,0)*{\bullet}="1";
  (0,0)*{\bullet}="2";
  (15,0)*{\bullet}="3";
 "1";"2" **\crv{(-12,9) & (-3,9)};
  "1";"2" **\crv{(-12,-9) & (-3,-9)};
   "2";"3" **\crv{(3,9) & (12,9)} ;
    "2";"3" **\crv{ (3,-9)& (12,-9) };
    (-7.5,6.75)*{\scriptstyle >}+(0,3)*{\scriptstyle F};
    (7.5,6.75)*{\scriptstyle >}+(0,3)*{\scriptstyle G};
    (-7.5,-6.75)*{\scriptstyle >}+(0,-3)*{\scriptstyle F'};
    (7.5,-6.75)*{\scriptstyle >}+(0,-3)*{\scriptstyle G'};
    (-17.5,0)*{\scriptstyle A};
  (-2.5,0)*{\scriptstyle B};
  (17.5,0)*{\scriptstyle C};
  {\ar@{=>}_{\scriptstyle \beta}(7.5,3)*{};(7.5,-3)*{}} ;
  {\ar@{=>}_{\scriptstyle \alpha}(-7.5,3)*{};(-7.5,-3)*{}} ;
   (-18,0)*{};
  (18,0)*{};
 \endxy .
 \]
 The axioms of a 2-category require that there exist an identity
2-morphism for vertical and horizontal composition, but perhaps
the most interesting axiom of a 2-category is the {\em interchange
law}:
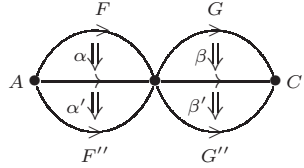
\[
\xy
  (-8,0)*{\bullet}="1";
  (8,0)*{\bullet}="2";
 "1";"2" **\crv{(-5,9) & (5,9)};
  "1";"2" **\crv{(-5,-9) & (5,-9)};
   "1";"2" **\dir{-} ?(.55)*\dir{>};
    (,6.75)*{\scriptstyle >}+(0,3)*{\scriptstyle F};
    (0,-6.75)*{\scriptstyle >}+(0,-3)*{\scriptstyle F''};
    (-10.5,0)*{\scriptstyle A};
  {\ar@{=>}_{\scriptstyle \alpha}(0,5)*{};(0,1.5)*{}} ;
  {\ar@{=>}_{\scriptstyle \alpha'}(0,-1.5)*{};(0,-5)*{}} ;
    (8,0)*{\bullet}="1";
  (24,0)*{\bullet}="2";
 "1";"2" **\crv{(11,9) & (21,9)};
  "1";"2" **\crv{(11,-9) & (21,-9)};
   "1";"2" **\dir{-} ?(.55)*\dir{>};
    (16,6.75)*{\scriptstyle >}+(0,3)*{\scriptstyle G};
    (16,-6.75)*{\scriptstyle >}+(0,-3)*{\scriptstyle G''};
  (26.5,0)*{\scriptstyle C};
  {\ar@{=>}_{\scriptstyle \beta}(16,5)*{};(16,1.5)*{}} ;
  {\ar@{=>}_{\scriptstyle \beta'}(16,-1.5)*{};(16,-5)*{}} ;
 \endxy .
 \]
The interchange law asserts that the above diagram is
unambiguously defined.  That is, the result of vertically
composing and then horizontally composing is the same as first
horizontally composing and then vertically composing.

This sort of notation, often referred to as globular notation, is
common among 2-category theorists. But for our purposes, we will
be interested in a different sort of diagram related to
2-categories --- string diagrams.  String diagrams are just the
Poincar\'e duals of the usual globular diagrams.  We construct
string diagrams from globular diagrams by inverting the dimensions
of the picture.  The objects, typically represented in globular
notation as 0-dimensional points (or bullets), become
2-dimensional surfaces in the new notation.  The morphism
represented by 1-dimensional edges remain 1-dimensional, and the
2-dimensional globes representing 2-morphisms become 0-dimensional
in the string diagram.

To see how this works, let $\mathcal{D}$ be a 2-category. We
depict objects $A$ and $B$ of $\mathcal{D}$ as surfaces:
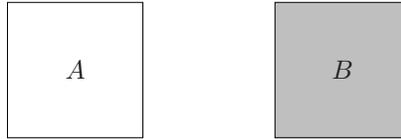
\[ \psset{xunit=.90cm,yunit=0.90cm}
\begin{pspicture}(2,2)
     \pspolygon(0,0)(0,2)(2,2)(2,0)(0,0)
     \rput(1,1){$A$}
 \end{pspicture}
 \qquad \qquad \quad
 \begin{pspicture}(2,2)
\pscustom[fillcolor=lightgray, fillstyle=solid]{
  \pspolygon(0,0)(0,2)(2,2)(2,0)(0,0)
  }
  \rput(1,1){$B$}
 \end{pspicture}
\]
where we have shaded the surface corresponding to $B$ in order to
easily distinguish it from $A$.  Below we show the process of
Poincar\'e dualizing a morphism $F \maps A \to B$ in $\mathcal{D}$
from the globular notation into the string notation:
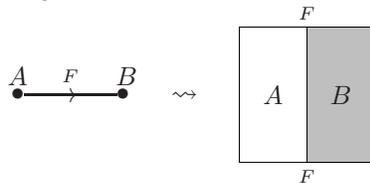
\[ \psset{xunit=.90cm,yunit=0.90cm}
 \xy
 (-7,0)*{\bullet}="1";
 (7,0)*{\bullet}="2";
    "1"; "2" **\dir{-} ?(.55)*\dir{>};
(-7,2.5)*{A}; (0,2.5)*{\scs F};
 (7.5,2.5)*{B};
\endxy
 \quad \rightsquigarrow \quad
  \xy
    (0,0)*+{    \begin{pspicture}(2.2,2)
        \pspolygon(0,0)(0,2)(2,2)(2,0)(0,0)
        \pscustom[fillcolor=lightgray, fillstyle=solid]{
        \pspolygon(1,0)(1,2)(2,2)(2,0)(1,0)
         }
         \rput(.5,1){$A$}
         \rput(1.5,1){$B$}
         \rput(1,2.2){$\scs F$}
         \rput(1,-.2){$\scs F$}
        \end{pspicture}}="";
  \endxy
\]
On the left is a morphism drawn in the usual globular notation,
and on the right the same morphism drawn in string notation.

The composite of morphisms $F \maps A \to B$ and $G \maps B \to C$
in $\mathcal{D}$ is drawn as:
 \ban
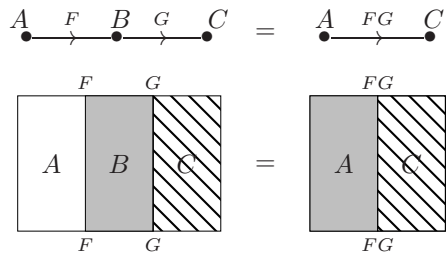
\xy
 (-6,0)*{\bullet}="1";
 (6,0)*{\bullet}="2";
 (18,0)*{\bullet}="3";
    "1"; "2" **\dir{-} ?(.55)*\dir{>};
    "2"; "3" **\dir{-} ?(.55)*\dir{>};
(-7,2.5)*{A}; (0,2.5)*{\scs F};
 (6.5,2.5)*{B}; (12,2.5)*{\scs G};
 (19.5,2.5)*{C};
\endxy
& = &
 \xy
 (-10,0)*{};
 (-7,0)*{\bullet}="1";
 (7,0)*{\bullet}="2";
    "1"; "2" **\dir{-} ?(.55)*\dir{>};
(-7,2.5)*{A}; (0,2.5)*{\scs FG};
 (7.5,2.5)*{C};
\endxy
\\ & & \\ \psset{xunit=.90cm,yunit=0.90cm}
  \xy
    (0,0)*+{    \begin{pspicture}(3,2)
        \pspolygon(0,0)(0,2)(2,2)(2,0)(0,0)
        \pscustom[fillcolor=lightgray, fillstyle=solid]{
        \pspolygon(1,0)(1,2)(2,2)(2,0)(1,0)
         }
         \pscustom[fillcolor=lightgray,fillstyle=vlines,hatchangle=45]{
        \pspolygon(2,0)(2,2)(3,2)(3,0)(2,0)
         }
         \rput(.5,1){$A$}
         \rput(1.5,1){$B$}
         \rput(2.5,1){$C$}
         \rput(1,2.2){$\scs F$}
         \rput(2,2.2){$\scs G$}
         \rput(1,-.2){$\scs F$}
         \rput(2,-.2){$\scs G$}
        \end{pspicture}}="";
  \endxy
 & = &
 \psset{xunit=.90cm,yunit=0.90cm}
 \xy
    (0,0)*+{    \begin{pspicture}(2.2,2)
        \pspolygon(0,0)(0,2)(2,2)(2,0)(0,0)
        \pscustom[fillcolor=lightgray, fillstyle=solid]{
        \pspolygon(1,0)(1,2)(0,2)(0,0)(1,0)
         }
         \pscustom[fillstyle=vlines,hatchangle=45, fillcolor=lightgray]{
        \pspolygon(1,0)(1,2)(2,2)(2,0)(1,0)
         }
         \rput(.5,1){$A$}
         \rput(1.5,1){$C$}
         \rput(1,2.2){$\scs FG$}
         \rput(1,-.2){$\scs FG$}
        \end{pspicture}}="";
  \endxy
 \ean
As a convenient convention, the identity morphism of objects in
$\mathcal{D}$ are not drawn.  This convention allows the
identification:
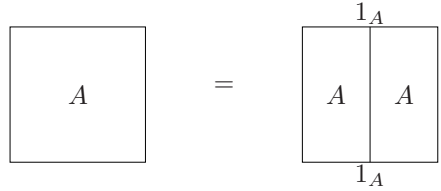
\[ \psset{xunit=.90cm,yunit=0.90cm}
 \xy
   (0,0)*+{\begin{pspicture}(2,2.3)
     \pspolygon(0,0)(0,2)(2,2)(2,0)(0,0)
     \rput(1,1){$A$}
 \end{pspicture}};
 \endxy
 \qquad = \qquad
  \xy
    (0,0)*+{ \begin{pspicture}(2,2.3)
     \pspolygon(0,0)(0,2)(2,2)(2,0)(0,0)
     \psline(1,0)(1,2)
     \rput(.5,1){$A$}
     \rput(1.5,1){$A$}
     \rput(1,2.2){$1_A$}
     \rput(1,-.2){$1_A$}
 \end{pspicture}};
  \endxy
 \]
of string diagrams.

If $F,F' \maps A \to B$ are morphism of $\mathcal{D}$ and $\alpha
\maps F \Rightarrow F'$ is a 2-morphism, then we depict this as:
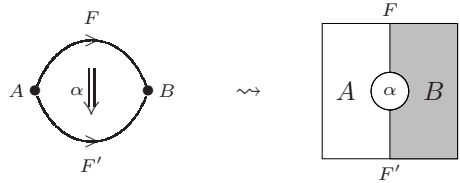
\[
 \xy
  (-15,0)*{\bullet}="1";
  (0,0)*{\bullet}="2";
 "1";"2" **\crv{(-12,9) & (-3,9)};
  "1";"2" **\crv{(-12,-9) & (-3,-9)};
    (-7.5,6.75)*{\scriptstyle >}+(0,3)*{\scriptstyle F};
    (-7.5,-6.75)*{\scriptstyle >}+(0,-3)*{\scriptstyle F'};
    (-17.5,0)*{\scriptstyle A};
  (2.5,0)*{\scriptstyle B};
  {\ar@{=>}_{\scriptstyle \alpha}(-7.5,3)*{};(-7.5,-3)*{}} ;
 \endxy
 \qquad \rightsquigarrow \qquad
 \psset{xunit=.90cm,yunit=0.90cm}
 \xy
 (0,0)*{
\begin{pspicture}(2,2)
        \pspolygon(0,0)(0,2)(2,2)(2,0)(0,0)
        \pscustom[fillcolor=lightgray, fillstyle=solid]{
        \psline(1,0)(1,.75)
        \psarc(1,1){.25}{-90}{90}
        \psline(1,2)
        \psline(2,2) \psline(2,0) \psline(1,0)
         }
         \rput(.35,1){$A$}
         \rput(1.65,1){$B$}
         \rput(1,2.2){$\scs F$}
         \rput(1,-0.2){$\scs F'$}
         \pscircle[fillstyle=solid, fillcolor=white](1,1){.25}
         \rput(1,1){$\scs \alpha$}
         \pscircle(1,1){.25}
 \end{pspicture}};
 \endxy
\]
where the circle surrounding $\alpha$ is thought of as being
0-dimensional. We include it only as a means to label the
2-morphism.  Following the convention that the identity morphisms
are not drawn in string diagrams we omit identity 2-morphisms as
well.  This allows the identification:
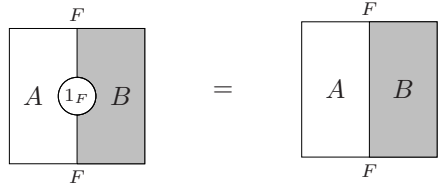
\[
\psset{xunit=.90cm,yunit=0.90cm}
 \xy
    (0,0)*+{
\begin{pspicture}(2,2.2)
        \pspolygon(0,0)(0,2)(2,2)(2,0)(0,0)
        \pscustom[fillcolor=lightgray, fillstyle=solid]{
        \psline(1,0)(1,.75)
        \psarc(1,1){.25}{-90}{90}
        \psline(1,2)
        \psline(2,2) \psline(2,0) \psline(1,0)
         }
         \rput(.35,1){$A$}
         \rput(1.65,1){$B$}
         \rput(1,2.2){$\scs F$}
         \rput(1,-0.2){$\scs F$}
         \pscircle[fillstyle=solid, fillcolor=white](1,1){.25}
         \rput(1,1){$\scs 1_F$}
         \pscircle(1,1){.25}
 \end{pspicture}};
 \endxy
 \qquad = \qquad
  \xy
    (0,0)*+{    \begin{pspicture}(2.2,2)
        \pspolygon(0,0)(0,2)(2,2)(2,0)(0,0)
        \pscustom[fillcolor=lightgray, fillstyle=solid]{
        \pspolygon(1,0)(1,2)(2,2)(2,0)(1,0)
         }
         \rput(.5,1){$A$}
         \rput(1.5,1){$B$}
         \rput(1,2.2){$\scs F$}
         \rput(1,-0.2){$\scs F$}
        \end{pspicture}}="";
  \endxy
 \]
of string diagrams.

Horizontal and vertical composition is achieved in the obvious
way. If $F,F',F'' \maps A \to B$ are morphisms and $\alpha \maps F
\Rightarrow F'$, $\alpha' \maps F' \Rightarrow F''$ are
2-morphisms, then the vertical composite of $\alpha$ and $\alpha'$
is:
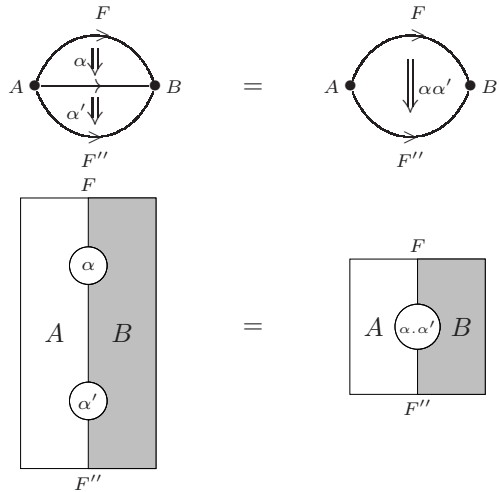
\[
 \xy
  (-8,0)*{\bullet}="1";
  (8,0)*{\bullet}="2";
 "1";"2" **\crv{(-5,9) & (5,9)};
  "1";"2" **\crv{(-5,-9) & (5,-9)};
   "1";"2" **\dir{-} ?(.55)*\dir{>};
    (,6.75)*{\scriptstyle >}+(0,3)*{\scriptstyle F};
    (0,-6.75)*{\scriptstyle >}+(0,-3)*{\scriptstyle F''};
    (-10.5,0)*{\scriptstyle A};
  (10.5,0)*{\scriptstyle B};
  {\ar@{=>}_{\scriptstyle \alpha}(0,5)*{};(0,1.5)*{}} ;
  {\ar@{=>}_{\scriptstyle \alpha'}(0,-1.5)*{};(0,-5)*{}} ;
 \endxy
 \qquad = \qquad
 \xy
  (-8,0)*{\bullet}="1";
  (8,0)*{\bullet}="2";
 "1";"2" **\crv{(-5,9) & (5,9)};
  "1";"2" **\crv{(-5,-9) & (5,-9)};
    (,6.75)*{\scriptstyle >}+(0,3)*{\scriptstyle F};
    (0,-6.75)*{\scriptstyle >}+(0,-3)*{\scriptstyle F''};
    (-10.5,0)*{\scriptstyle A};
  (10.5,0)*{\scriptstyle B};
  {\ar@{=>}^{\scriptstyle \alpha \alpha'}(0,3.5)*{};(0,-3.5)*{}} ;
 \endxy
\]

\[
  \psset{xunit=.90cm,yunit=0.90cm}
\begin{pspicture}(2,2)
        \pspolygon(0,-2)(0,2)(2,2)(2,-2)(0,-2)
        \pscustom[fillcolor=lightgray, fillstyle=solid]{
        \psline(1,0)(1,.75)
        \psarc(1,1){.25}{-90}{90}
        \psline(1,2)
        \psline(2,2) \psline(2,-2) \psline(1,-2) \psline(1,-1.25)
        \psarc(1,-1){.25}{-90}{90}
        \psline(1,0)
         }
         \rput(.5,0){$A$}
         \rput(1.5,0){$B$}
         \rput(1,2.2){$\scs F$}
         \rput(1,-2.2){$\scs F''$}
         \pscircle[fillstyle=solid, fillcolor=white](1,1){.25}
         \pscircle[fillstyle=solid, fillcolor=white](1,-1){.25}
         \rput(1,1){$\scs \alpha$}
         \rput(1,-1){$\scs \alpha'$}
         \pscircle(1,1){.25}
         \pscircle(1,-1){.25}
 \end{pspicture}
 \qquad \quad= \qquad \quad
  \psset{xunit=.90cm,yunit=0.90cm}
 \xy
 (0,0)*{
        \begin{pspicture}(2,2)
        \pspolygon(0,0)(0,2)(2,2)(2,0)(0,0)
        \pscustom[fillcolor=lightgray, fillstyle=solid]{
        \psline(1,0)(1,.7)
        \psarc(1,1){.3}{-90}{90}
        \psline(1,2)
        \psline(2,2) \psline(2,0) \psline(1,0)
         }
         \rput(.35,1){$A$}
         \rput(1.65,1){$B$}
         \rput(1,2.2){$\scs F$}
         \rput(1,-0.2){$\scs F''$}
         \pscircle[fillstyle=solid, fillcolor=white](1,1){.3}
         \rput(1,1){$ \scriptscriptstyle \alpha.\alpha'$}
         \pscircle(1,1){.3}
    \end{pspicture}};
    (0,-20)*{}; 
 \endxy
\]
If $F,F' \maps A \to B$ and $G,G' \maps B \to C$ with $\alpha
\maps F \Rightarrow F'$ and $\beta \maps G \To G'$, then the
horizontal composite is depicted as follows:
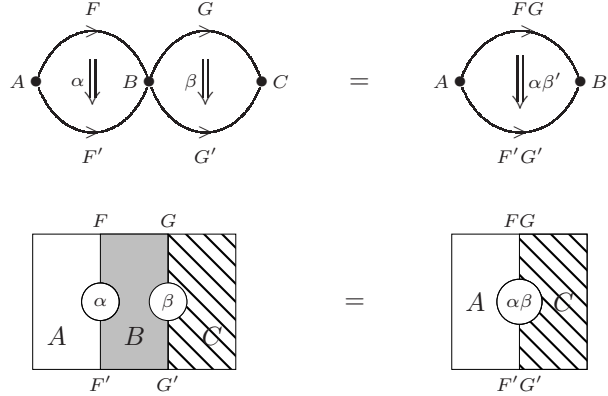
\[
\xy
  (-15,0)*{\bullet}="1";
  (0,0)*{\bullet}="2";
  (15,0)*{\bullet}="3";
 "1";"2" **\crv{(-12,9) & (-3,9)};
  "1";"2" **\crv{(-12,-9) & (-3,-9)};
   "2";"3" **\crv{(3,9) & (12,9)} ;
    "2";"3" **\crv{ (3,-9)& (12,-9) };
    (-7.5,6.75)*{\scriptstyle >}+(0,3)*{\scriptstyle F};
    (7.5,6.75)*{\scriptstyle >}+(0,3)*{\scriptstyle G};
    (-7.5,-6.75)*{\scriptstyle >}+(0,-3)*{\scriptstyle F'};
    (7.5,-6.75)*{\scriptstyle >}+(0,-3)*{\scriptstyle G'};
    (-17.5,0)*{\scriptstyle A};
  (-2.5,0)*{\scriptstyle B};
  (17.5,0)*{\scriptstyle C};
  {\ar@{=>}_{\scriptstyle \beta}(7.5,3)*{};(7.5,-3)*{}} ;
  {\ar@{=>}_{\scriptstyle \alpha}(-7.5,3)*{};(-7.5,-3)*{}} ;
   (-18,0)*{};
  (18,0)*{};
 \endxy
\qquad = \qquad
 \xy
    (12,0)*{};
    (-12,0)*{};
  (-8,0)*{\bullet}="1";
  (8,0)*{\bullet}="2";
 "1";"2" **\crv{(-5,9) & (5,9)};
  "1";"2" **\crv{(-5,-9) & (5,-9)};
    (,6.75)*{\scriptstyle >}+(0,3)*{\scriptstyle FG};
    (0,-6.75)*{\scriptstyle >}+(0,-3)*{\scriptstyle F'G'};
    (-10.5,0)*{\scriptstyle A};
  (10.5,0)*{\scriptstyle B};
  {\ar@{=>}^{\scriptstyle \alpha \beta'}(0,3.5)*{};(0,-3.5)*{}} ;
 \endxy
\]
\medskip

\[
  \psset{xunit=.90cm,yunit=0.90cm}
 \xy
 (0,0)*{
        \begin{pspicture}(3,2)
        \pspolygon(0,0)(0,2)(2,2)(2,0)(0,0)
        \pscustom[fillcolor=lightgray, fillstyle=solid]{
        \psline(1,0)(1,.75)
        \psarc(1,1){.25}{-90}{90}
        \psline(1,2)
        \psline(2,2) \psline(2,0) \psline(1,0)
         }
        \pscustom[fillstyle=vlines,hatchangle=45, fillcolor=lightgray]{
        \psline(2,0)(2,.75)
        \psarc(2,1){.25}{-90}{90}
        \psline(2,2)
        \psline(3,2) \psline(3,0) \psline(2,0)
         }
         \rput(.35,.5){$A$}
         \rput(1.5,.5){$B$}
         \rput(2.65,.5){$C$}
         \rput(1,2.2){$\scs F$}
         \rput(1,-0.2){$\scs F'$}
         \rput(2,2.2){$\scs G$}
         \rput(2,-0.2){$\scs G'$}
         \pscircle[fillstyle=solid, fillcolor=white](1,1){.25}
         \rput(1,1){$\scs \alpha$}
         \pscircle(1,1){.25}
         \pscircle[fillstyle=solid, fillcolor=white](2,1){.25}
         \rput(2,1){$\scs \beta$}
    \end{pspicture}};
 \endxy
  \qquad \quad \; \; \;= \qquad \quad
   \psset{xunit=.90cm,yunit=0.90cm}
 \xy
 (0,0)*{
        \begin{pspicture}(2,2)
        \pspolygon(0,0)(0,2)(2,2)(2,0)(0,0)
        \pscustom[fillstyle=vlines,hatchangle=45, fillcolor=lightgray]{
        \psline(1,0)(1,.7)
        \psarc(1,1){.3}{-90}{90}
        \psline(1,2)
        \psline(2,2) \psline(2,0) \psline(1,0)
         }
         \rput(.35,1){$A$}
         \rput(1.65,1){$C$}
         \rput(1,2.2){$\scs FG$}
         \rput(1,-0.2){$\scs F'G'$}
         \pscircle[fillstyle=solid, fillcolor=white](1,1){.3}
         \rput(1,1){$ \scriptstyle \alpha\beta$}
         \pscircle(1,1){.3}
    \end{pspicture}};
    (0,-13)*{};
 \endxy
\]

The interchange law in the 2-category $\mathcal{D}$ tells us that
each string diagram such as:
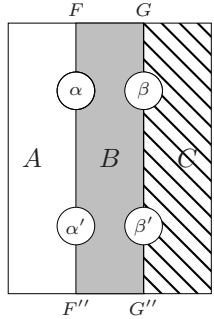
\[
  \psset{xunit=.90cm,yunit=0.90cm}
  \xy
    (0,0)*{
        \begin{pspicture}[.5](2,2.5)
        \pspolygon(0,-2)(0,2)(2,2)(2,-2)(0,-2)
        \pscustom[fillcolor=lightgray, fillstyle=solid]{
        \psline(1,-2)(1,.75)
        \psarc(1,1){.25}{-90}{90}
        \psline(1,2)
        \psline(2,2) \psline(2,-2) \psline(1,-2)
         }
        \pscustom[fillstyle=vlines,hatchangle=45, fillcolor=lightgray]{
        \psline(2,-2)(2,.75)
        \psarc(2,1){.25}{-90}{90}
        \psline(2,2)
        \psline(3,2) \psline(3,-2) \psline(2,-2)
         }
         \rput(.35,0){$A$}
         \rput(1.5,0){$B$}
         \rput(2.65,0){$C$}
         \rput(1,2.2){$\scs F$}
         \rput(1,-2.2){$\scs F''$}
         \rput(2,2.2){$\scs G$}
         \rput(2,-2.2){$\scs G''$}
         \pscircle[fillstyle=solid, fillcolor=white](1,1){.25}
         \rput(1,1){$\scs \alpha$}
         \pscircle[fillstyle=solid, fillcolor=white](1,-1){.25}
         \rput(1,-1){$\scs \alpha'$}
         \pscircle(1,1){.25}
         \pscircle[fillstyle=solid, fillcolor=white](2,1){.25}
         \rput(2,1){$\scs \beta$}
         \pscircle[fillstyle=solid, fillcolor=white](2,-1){.25}
         \rput(2,-1){$\scs \beta'$}
    \end{pspicture}};
    (0,-34)*{}; 
    \endxy
 \]
can be uniquely interpreted as a diagram in $\mathcal{D}$.  Later
we will see that, together with the identity morphisms, the
interchange law justifies vertical and horizontal topological
deformations of these string diagrams. Of course, theorems proving
that any topological deformation of a string diagram produces an
equivalent diagram in $\mathcal{D}$ have been proved by Joyal and
Street~\cite{js1,js2}, but we prefer to demonstrate the needed
deformations directly from the axioms.

Before moving on, as an application of this string diagram
technology, we would like to translate the definition of an
adjunction into string diagrams. An adjunction is a concept that
makes sense in any 2-category, although most of us learn about
adjunctions in the 2-category $\cat{Cat}$ consisting of
categories, functors, and natural transformations.  Adjunctions in
this 2-category are just adjoint functors which are prevalent
throughout mathematics.  We will be interested in adjunctions
because they provide a categorical framework for understanding
Frobenius algebras.  Once we set up this framework we will be able
to categorify it and arrive at a definition of a pseudo Frobenius
algebra.

Let $\mathcal{D}$ be a 2-category.  An {\em adjunction}
$(A,B,L,R,i,e)$ in $\mathcal{D}$ consists of objects $A$ and $B$,
morphism $L \maps A \to B$ and $R \maps B \to A$, and 2-morphisms
$i\maps 1_A \To LR$ and $e \maps RL \To 1_B$ called the {\em unit}
and {\em counit} of the adjunction, such that the following
diagrams:
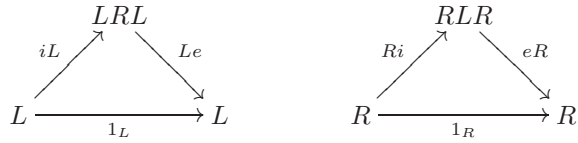
\[
 \xymatrix@C=1.5pc@R=.8pc{
  & LRL \ar[ddr]^-{Le} & \\ \\
 L \ar[uur]^-{iL} \ar[rr]_-{1_L} & & L
 }
\qquad \qquad
 \xymatrix@C=1.5pc@R=.8pc{
  & RLR \ar[ddr]^-{eR} & \\ \\
 R \ar[uur]^-{Ri} \ar[rr]_-{1_R} & & R
 }
\]
commute.  We sometimes refer to these identities as the {\em
zig-zag} identities for reasons that will soon become apparent.

In string notation the morphisms $L$ and $R$ are depicted as:
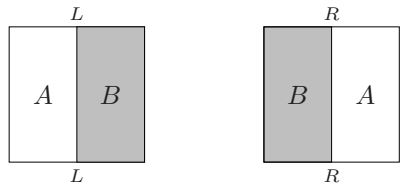
\[ \psset{xunit=.90cm,yunit=0.90cm}
   \begin{pspicture}(2.2,2.2)
  \pspolygon(0,0)(0,2)(2,2)(2,0)(0,0)
\pscustom[fillcolor=lightgray, fillstyle=solid]{
  \pspolygon(1,0)(1,2)(2,2)(2,0)(1,0)
  }
  \rput(.5,1){$A$}
         \rput(1.5,1){$B$}
         \rput(1,2.2){$\scs L$}
         \rput(1,-.2){$\scs L$}
 \end{pspicture}
\qquad \qquad
   \begin{pspicture}(2.2,2.2)
\pscustom[fillcolor=lightgray, fillstyle=solid]{
  \pspolygon(1,0)(1,2)(0,2)(0,0)(1,0)
  }
   \rput(.5,1){$B$}
         \rput(1.5,1){$A$}
         \rput(1,2.2){$\scs R$}
         \rput(1,-.2){$\scs R$}
  \pspolygon(0,0)(0,2)(2,2)(2,0)(0,0)
 \end{pspicture}
\]
and their composites as:
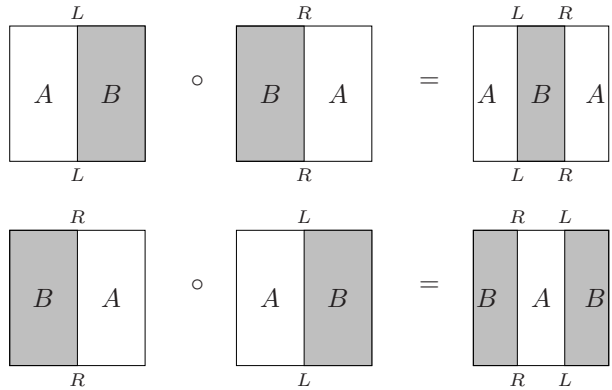
\[ \psset{xunit=.90cm,yunit=0.90cm}
   \begin{pspicture}[.5](2.2,2.2)
  \pspolygon(0,0)(0,2)(2,2)(2,0)(0,0)
\pscustom[fillcolor=lightgray, fillstyle=solid]{
  \pspolygon(1,0)(1,2)(2,2)(2,0)(1,0)
  }
  \rput(.5,1){$A$}
         \rput(1.5,1){$B$}
         \rput(1,2.2){$\scs L$}
         \rput(1,-.2){$\scs L$}
 \end{pspicture}
\quad \circ \quad
   \begin{pspicture}[.5](2.2,2.2)
\pscustom[fillcolor=lightgray, fillstyle=solid]{
  \pspolygon(1,0)(1,2)(0,2)(0,0)(1,0)
  }
   \rput(.5,1){$B$}
         \rput(1.5,1){$A$}
         \rput(1,2.2){$\scs R$}
         \rput(1,-.2){$\scs R$}
  \pspolygon(0,0)(0,2)(2,2)(2,0)(0,0)
 \end{pspicture} \quad = \quad
    \begin{pspicture}[.5](2.2,2.2)
\pscustom[fillcolor=lightgray, fillstyle=solid]{
  \pspolygon(.65,0)(.65,2)(1.35,2)(1.35,0)(.65,0)
  }
         \rput(.2,1){$A$}
         \rput(1,1){$B$}
         \rput(1.8,1){$A$}
         \rput(1.35,2.2){$\scs R$}
         \rput(1.35,-.2){$\scs R$}
         \rput(.65,2.2){$\scs L$}
         \rput(.65,-.2){$\scs L$}
  \pspolygon(0,0)(0,2)(2,2)(2,0)(0,0)
 \end{pspicture}
\]

\[ \psset{xunit=.90cm,yunit=0.90cm}
\begin{pspicture}[.5](2.2,2.2)
\pscustom[fillcolor=lightgray, fillstyle=solid]{
  \pspolygon(1,0)(1,2)(0,2)(0,0)(1,0)
  }
   \rput(.5,1){$B$}
         \rput(1.5,1){$A$}
         \rput(1,2.2){$\scs R$}
         \rput(1,-.2){$\scs R$}
  \pspolygon(0,0)(0,2)(2,2)(2,0)(0,0)
 \end{pspicture}
\quad \circ \quad
   \begin{pspicture}[.5](2.2,2.2)
  \pspolygon(0,0)(0,2)(2,2)(2,0)(0,0)
\pscustom[fillcolor=lightgray, fillstyle=solid]{
  \pspolygon(1,0)(1,2)(2,2)(2,0)(1,0)
  }
  \rput(.5,1){$A$}
         \rput(1.5,1){$B$}
         \rput(1,2.2){$\scs L$}
         \rput(1,-.2){$\scs L$}
 \end{pspicture} \quad = \quad
    \begin{pspicture}[.5](2.2,2.2)
\pscustom[fillcolor=lightgray, fillstyle=solid]{
  \pspolygon(.65,0)(.65,2)(0,2)(0,0)(.65,0)}
  \pscustom[fillcolor=lightgray, fillstyle=solid]{
  \pspolygon(1.35,0)(1.35,2)(2,2)(2,0)(1.35,0)
  }
         \rput(.2,1){$B$}
         \rput(1,1){$A$}
         \rput(1.8,1){$B$}
         \rput(1.35,2.2){$\scs L$}
         \rput(1.35,-.2){$\scs L$}
         \rput(.65,2.2){$\scs R$}
         \rput(.65,-.2){$\scs R$}
  \pspolygon(0,0)(0,2)(2,2)(2,0)(0,0)
 \end{pspicture}
\]

\medskip

With the above diagrams in mind it is easy to see that the unit
and counit of the adjunction can be depicted as:
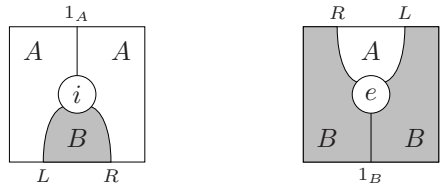
\[
  \psset{xunit=.90cm,yunit=0.90cm}
 \xy
 (0,0)*{
        \begin{pspicture}(2,2.4)
        \pspolygon(0,0)(0,2)(2,2)(2,0)(0,0)
        \pscustom[fillcolor=lightgray, fillstyle=solid]{
        \psbezier(.5,0)(.5,.5)(.65,.8)(.824,.824)
        \psline(1.176,.824)
        \psbezier(1.176,.824)(1.35,.8)(1.5,.5)(1.5,0)
        \psline(.5,0)
         }
         \psline(1,2)(1,1.25)
         \rput(.35,1.65){$A$}
         \rput(1.65,1.65){$A$}
         \rput(1,.35){$B$}
         \rput(1,2.2){$\scs 1_A$}
         \rput(.5,-0.2){$\scs L$}
         \rput(1.5,-0.2){$\scs R$}
         \pscircle[fillstyle=solid, fillcolor=white](1,1){.25}
         \rput(1,1){$ i$}
    \end{pspicture}};
 \endxy
\qquad \qquad \qquad
  \psset{xunit=.90cm,yunit=0.90cm}
 \xy
 (0,0)*{
        \begin{pspicture}(2,2.4)
        \pspolygon(0,0)(0,2)(2,2)(2,0)(0,0)
        \pscustom[fillcolor=lightgray, fillstyle=solid]{
        \psbezier(.5,2)(.5,1.5)(.65,1.2)(.824,1.176)
        \psline(1.176,1.176)
        \psbezier(1.176,1.176)(1.35,1.2)(1.5,1.5)(1.5,2)
        \psline(2,2) \psline(2,0) \psline(0,0) \psline(0,2) \psline(.5,2)
         }
         \psline(1,0)(1,.75)
         \rput(.35,.35){$B$}
         \rput(1.65,.35){$B$}
         \rput(1,1.65){$A$}
         \rput(1,-.2){$\scs 1_B$}
         \rput(.5,2.2){$\scs R$}
         \rput(1.5,2.2){$\scs L$}
         \pscircle[fillstyle=solid, fillcolor=white](1,1){.25}
         \rput(1,1){$ e$}
    \end{pspicture}};
 \endxy
\]

\medskip

\noindent However,  applying the convention that identity
morphisms are not depicted in the string diagrams, these pictures
can be simplified. Further, we can omit the labels of $i$ and $e$
when no confusion is likely to arise, in which case the unit and
counit become:
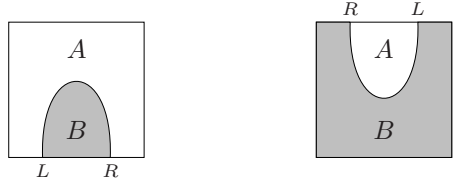
\[\psset{xunit=.90cm,yunit=0.90cm}
   \begin{pspicture}(2.2,2.4)
  \pspolygon(0,0)(0,2)(2,2)(2,0)(0,0)
   \pscustom[fillcolor=lightgray, fillstyle=solid]{
  \psbezier(.5,0)(.45,1.5)(1.55,1.5)(1.5,0)
  \psline(.5,0)
  }
        \rput(1,.4){$B$}
         \rput(1,1.6){$A$}
         \rput(1.5,-.2){$\scs R$}
         \rput(.5,-.2){$\scs L$}
 \end{pspicture}
 \qquad \qquad \qquad
   \begin{pspicture}(2.2,2.4)
  \pspolygon(0,0)(0,2)(2,2)(2,0)(0,0)
   \pscustom[fillcolor=lightgray, fillstyle=solid]{
  \psbezier(.5,2)(.45,.5)(1.55,.5)(1.5,2)
  \psline(2,2)
  \psline(2,0)
   \psline(0,0)
    \psline(0,2)
     \psline(.5,2)
  }
         \rput(1,.4){$B$}
         \rput(1,1.6){$A$}
         \rput(1.5,2.2){$\scs L$}
         \rput(.5,2.2){$\scs R$}
 \end{pspicture}
\]

\medskip

\noindent and in this notation the zig-zag identities become:
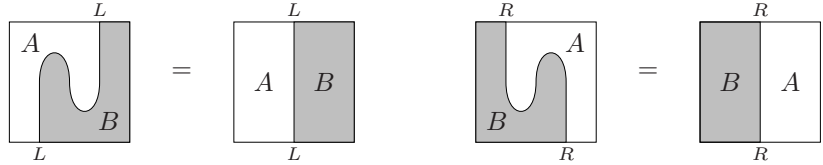
\[ \psset{xunit=.80cm,yunit=0.80cm}
\xy
   (0,0)*+{\begin{pspicture}(2,2.4)
  \pspolygon(0,0)(0,2)(2,2)(2,0)(0,0)
\pscustom[fillcolor=lightgray, fillstyle=solid]{
  \psbezier(.5,1)(.5,1.65)(1,1.65)(1,1)
  \psbezier(1,1)(1,.35)(1.5,.35)(1.5,1)
  \psline(1.5,2)
  \psline(2,2)
  \psline(2,0)
  \psline(.5,0)
  \psline(.5,1)
  }
         \rput(.35,1.65){$A$}
         \rput(1.65,.35){$B$}
         \rput(1.5,2.2){$\scs L$}
         \rput(.5,-.2){$\scs L$}
 \end{pspicture}};
 \endxy
\quad = \quad
  \xy
   (0,0)*+{  \begin{pspicture}(2,2.4)
  \pspolygon(0,0)(0,2)(2,2)(2,0)(0,0)
\pscustom[fillcolor=lightgray, fillstyle=solid]{
  \pspolygon(1,0)(1,2)(2,2)(2,0)(1,0)
  }
  \rput(.5,1){$A$}
         \rput(1.5,1){$B$}
         \rput(1,2.2){$\scs L$}
         \rput(1,-.2){$\scs L$}
 \end{pspicture}};
 \endxy
 \qquad \qquad
 \psset{xunit=.80cm,yunit=0.80cm}
     \xy
    (0,0)*+{   \begin{pspicture}(2,2.4)
  \pspolygon(0,0)(0,2)(2,2)(2,0)(0,0)
\pscustom[fillcolor=lightgray, fillstyle=solid]{
  \psbezier(.5,1)(.5,.35)(1,.35)(1,1)
  \psbezier(1,1)(1,1.65)(1.5,1.65)(1.5,1)
  \psline(1.5,0)
  \psline(0,0)
  \psline(0,2)
  \psline(.5,2)
  \psline(.5,1)
  }
          \rput(1.65,1.65){$A$}
         \rput(.35,.35){$B$}
         \rput(.5,2.2){$\scs R$}
         \rput(1.5,-.2){$\scs R$}
 \end{pspicture}};
  \endxy
 \quad = \quad
 \xy
      (0,0)*+{\begin{pspicture}(2,2.4)
\pscustom[fillcolor=lightgray, fillstyle=solid]{
  \pspolygon(1,0)(1,2)(0,2)(0,0)(1,0)
  }
   \rput(.5,1){$B$}
         \rput(1.5,1){$A$}
         \rput(1,2.2){$\scs R$}
         \rput(1,-.2){$\scs R$}
  \pspolygon(0,0)(0,2)(2,2)(2,0)(0,0)
 \end{pspicture}};
    \endxy
\]
which explains their name.  These identities say that a `zig-zag'
can be `straightened out'.

\subsection{The walking adjunction} \label{secWA}

In the previous section we defined the notion of an adjunction in
a 2-category.  In this section we will study a very special
adjunction --- the `walking adjunction'.  Before defining the
walking adjunction we feel obliged to motivate this seemingly
strange terminology. Imagine you are sitting in a small table in
the back of a crowded pub enjoying a beer with a close friend,
when in walks a fellow with enormous bushy eyebrows.  His eyebrows
are in fact so large that it seems his entire body serves no other
purpose than to provide a frame for these enormous eyebrows to
perch on. In that case, you might be tempted to comment to your
friend: ``Look, there goes the walking pair of eyebrows".  In the
same way, the walking adjunction is the minimal amount of
structure needed in order to have an adjunction; it is the
2-category freely generated by an adjunction.  The 2-category is
merely the frame upon which the adjunction `perches'.  This
`walking' terminology was coined by James Dolan and our
explanation of it is adapted from the expository writings of John
Baez~\cite{ThisWeekII}. In general, we will refer to the free $X$
on generating data $Y$ as the \textit{walking} $Y$.

We begin by explaining what it means for some structure to
`generate a 2-category'.

\begin{defn} \label{generate}
Let $Y$ be a structure that can be defined in an arbitrary
2-category. If $Y$ consists of objects $Y_1,Y_2,\ldots,Y_{n}$,
morphisms $F_1,F_2,\ldots,F_{n'}$, and 2-morphisms
$\alpha_1,\alpha_2,\ldots,\alpha_{n''}$ for $n,n',n'' \in \Z^+$,
then the 2-category $X$ is {\em generated} by $Y$ if:
\begin{enumerate}[(i.)]
    \item Every object of $X$ is some $Y_i$.
    \item Every 1-morphism of $X$ can be obtained by compositions
    from the $F_i$'s and $1_{Y_i}$'s.
    \item Every 2-morphism of $X$ is obtained by horizontal and
    vertical composition from the 2-morphisms $\alpha_i$, and identity 2-morphisms $1_F$
    for arbitrary 1-morphisms $F$.
\end{enumerate}
We say that $X$ is {\em freely generated}\footnote{This can be
stated more elegantly using globular sets~\cite{Bat1}, but for our
purposes this definition suffices.} by $Y$ if the set of $Y$
objects in $\mathcal{C}$ are in bijection with 2-functors from $X$
into $\mathcal{C}$, for every 2-category $\mathcal{C}$.
\end{defn}

\begin{defn} \et \label{DefWalkAdj}
The {\em walking adjunction} $\Adj$ is the 2-category freely
generated by:
\begin{itemize}
    \item objects $A$ and $B$,
    \item morphisms $L \maps A \to B$ and $R \maps B \to A$, and
    \item 2-morphisms $i \maps 1_A \to L \circ R$ and $e \maps R \circ L
\to 1_B$
\end{itemize}
such that the following diagrams:
\[
 \xymatrix@C=1.5pc@R=1pc{
 L \ar[rr]^-{iL} \ar[ddrr]_-{1_L} && LRL \ar[dd]^-{Le} \\
\\
&& L
 }
\qquad \qquad
 \xymatrix@C=1.5pc@R=1pc{
 R \ar[rr]^-{Ri} \ar[ddrr]_-{1_R}
    && RLR \ar[dd]^-{eR} \\
\\
&& R
 }
\]
commute.
\end{defn}

Thus, the walking adjunction $\Adj$, or the 2-category freely
generated by an adjunction, has the property that every adjunction
in a 2-category $\mathcal{C}$ corresponds to a 2-functor $\Adj \to
\mathcal{C}$.  The walking adjunction was first studied under the
name of the free adjunction~\cite{Freeadj}. Another description of
the walking adjunction can be obtained from its categorification
explicitly constructed by Lack~\cite{Lack}.

The walking adjunction turns out to be intimately related to the
walking monoid.  The notion of a monoid makes sense in any
monoidal category, so the walking monoid is just the monoidal
category freely generated by a monoid. Recall that a monoidal
category is a category equipped with a multiplication functor and
a unit object. However, it is sometimes useful to think of
monoidal category as a special kind of 2-category.  More
precisely, a monoidal category is just a one object 2-category:

\begin{center}
\begin{tabular}{|l|l|}
  \hline
  $\mathcal{M}$ -- monoidal category & $\mathcal{C}$ -- 2-category \\
  \hline \hline
  objects                     & morphisms \\
  tensor product of objects   & composition of morphisms \\
  morphisms                    & 2-morphisms \\
  composition                 & vertical composition of 2-morphisms\\
  tensor product of morphisms  & horizontal composition of 2-morphisms \\
  \hline
\end{tabular}
\end{center}

\noindent  The objects of the monoidal category are just the
morphisms of the 2-category $\mathcal{C}$.  The tensor product of
objects comes from the composition of morphisms in $\mathcal{C}$.
Note that every morphism of $\mathcal{C}$ is composable since it
only has one object. The morphisms of the monoidal category
$\mathcal{M}$ are the 2-morphisms of the 2-category $\mathcal{C}$.
The composition of morphisms in $\mathcal{M}$ comes from the
vertical composition of morphisms in $\mathcal{C}$, and the tensor
product of morphisms in $\mathcal{M}$ comes from the horizontal
composition of morphisms in $\mathcal{C}$.

We can also apply this same trick in reverse.  Given a monoidal
category $\mathcal{M}$, we can regard $\mathcal{M}$ as a
2-category $\Sigma(\mathcal{M})$ with one object by applying the
above procedure in reverse.  We sometimes refer to the 2-category
$\Sigma(\mathcal{M})$ as the {\em suspension} of the monoidal
category $\mathcal{M}$. Since a monoidal category is just a one
object 2-category, we can use Definition~\ref{generate} to define
the monoidal category freely generated by a monoid and dually the
monoidal category freely generated on a comonoid.

Below we define the walking monoid $\Mon$ and the walking comonoid
$\Comon$, not to be confused with $\cat{Mon}$ the category whose
objects are monoids, and $\cat{Comon}$ the category whose objects
are comonoids.

\begin{defn}
The {\em walking monoid} $\Mon$ is the monoidal category freely
generated by:
\begin{itemize}
    \item objects $A$ and $I$, and
    \item morphisms $m \maps A \ten A \to A$ and $\iota \maps I \to A$
\end{itemize}
such that
\[
 \xy
 (16,0)*+{A \ten A}="R";
 (0,16)*+{A \ten A \ten A}="T";
 (0,-16)*+{A}="B";
 (-16,0)*+{A \ten A}="L";
    {\ar^{1_A \ten m} "T";"R"};
    {\ar_{m \ten 1_A} "T";"L"};
    {\ar_{ m} "L";"B"};
    {\ar^{ m} "R";"B"};
 \endxy
\quad \quad
 \vcenter{\xy
 (-22,0)*+{I \ten A}="L";
 (0,0)*+{A \ten A}="M";
 (22,0)*+{A \ten I}="R";
 (0,-16)*+{A}="B";
    {\ar^{\iota \ten 1_A} "L";"M"};
    {\ar_{1_A \ten \iota} "R";"M"};
    {\ar_{ } "L";"B"};
    {\ar^{ } "R";"B"};
    {\ar^{ m} "M";"B"};
 \endxy}
\]
commute. Dually, the {\em walking comonoid} $\Comon$ is the
monoidal category freely generated by:
\begin{itemize}
    \item objects $B$ and $I$, and
    \item morphisms $\Delta \maps B \to B\ten B$ and $\varepsilon \maps B\to I$
\end{itemize}
such that
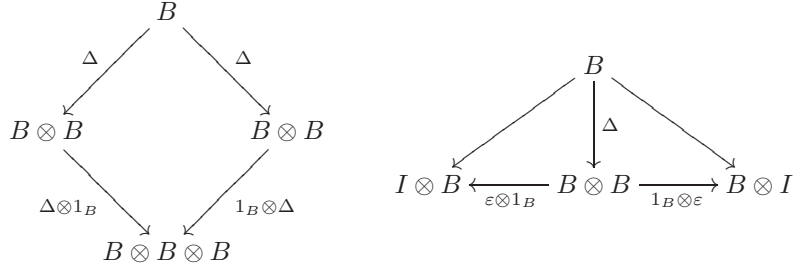
\[
\xy
 (16,0)*+{B \ten B}="R";
 (0,-16)*+{B \ten B \ten B}="T";
 (0,16)*+{B}="B";
 (-16,0)*+{B \ten B}="L";
    {\ar^{1_B \ten \Delta} "T";"R";};
    {\ar_{\Delta \ten 1_B} "T";"L";};
    {\ar_{ \Delta} "L";"B";};
    {\ar^{ \Delta} "R";"B";};
 \endxy
\quad \quad
 \vcenter{\xy
 (-22,0)*+{I \ten B}="L";
 (0,0)*+{B \ten B}="M";
 (22,0)*+{B \ten I}="R";
 (0,16)*+{B}="B";
    {\ar^{\varepsilon \ten 1_B} "M";"L"};
    {\ar_{1_B \ten \varepsilon} "M";"R"};
    {\ar_{ } "L";"B";};
    {\ar^{ } "R";"B";};
    {\ar^{ \Delta} "M";"B";};
 \endxy}
\]
commute.
\end{defn}

The relationship between the walking adjunction and the walking
monoid and comonoid is summed up by:

\begin{thm} \label{WAWM}
The monoidal category $\Hom(A,A)$ in the walking adjunction is the
walking monoid and the monoidal category $\Hom(B,B)$ in the
walking adjunction is the walking comonoid.
\end{thm}

Before proving this theorem we pause briefly to explain why this
theorem is morally true from a topological perspective.  Using
string diagrams it will be clear that the walking monoid is
contained in $\Hom(A,A)$ in the walking adjunction.  The argument
we use was originally developed by M\"{u}ger~\cite{Muger} and was
later elaborated on by Baez~\cite{ThisWeek}.  Unfortunately, the
main difficulty of this proof is showing the opposite inclusion:
that the walking monoid arises as $\Hom(A,A)$ in the walking
adjunction.  Or to put it another way, the difficulty lies in
showing that $\Hom(A,A)$ has \textit{only} the relations of the
walking monoid and no additional relations. We will see that both
inclusions can be proved using a bit of abstract category theory;
but the topological arguments provide the right intuition for
understanding the more abstract result.

In the walking adjunction, all of the objects in the monoidal
category $\Hom(A,A)$ are generated by the morphism $LR$.  We will
show that $LR$ is equipped with the structure of a monoid. The
multiplication on the object $LR$ is defined, using the counit of
the adjunction, as $m:= LeR\maps LRLR \to LR$. We depict this in
string notation as follows:
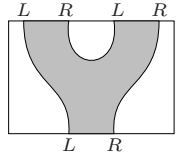
\[ \psset{xunit=1cm,yunit=1cm}
 \xy
  (0,-15)*{};
  (35,0)*{
   \begin{pspicture}(2.2,2.2)
  \pspolygon(-1.1,0)(-1.1,1.5)(1.1,1.5)(1.1,0)(0,0)
  \rput(0,0){\mult}
  \rput(-.9,1.65){$\scs L$}
  \rput(-.35,1.65){$\scs R$}
  \rput(.35,1.65){$\scs L$}
  \rput(.9,1.65){$\scs R$}
  \rput(-.3,-.15){$\scs L$}
  \rput(.3,-.15){$\scs R$}
 \end{pspicture}};
 \endxy
\]
where we have been slightly artistic with the strings representing
the identity morphisms on $L$ and $R$.  This diagram is meant to
be reminiscent of the short hand notation used in
Section~\ref{secFrob}.  The unit $\iota$ for the monoid defined to
be:
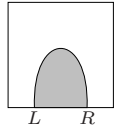
\[
   \begin{pspicture}(2.2,2.2)
  \pspolygon(0,0)(0,2)(2,2)(2,0)(0,0)
   \pscustom[fillcolor=lightgray, fillstyle=solid]{
  \psbezier(.5,0)(.45,1.5)(1.55,1.5)(1.5,0)
  \psline(.5,0)
  }
         \rput(1.5,-.2){$\scs R$}
         \rput(.5,-.2){$\scs L$}
 \end{pspicture}
\]
the unit for the adjunction $\iota:= i \maps 1_A \to LR$.

Now we can use the axioms of a 2-category, together with the
axioms for an adjunction, to show that this multiplication is
associative and that the unit satisfies the unit axioms.  For the
multiplication to be associative we must have an equality of
string diagrams:
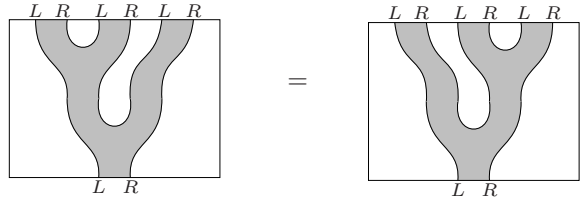
\[ 
\begin{pspicture}[0.5](3,3.4)
  \rput(1,0){\multQ} \rput(1.6,1.5){\curverightQ} \rput(.4,1.5){\multQ}
  \rput(-.5,3.18){$\scs L$}
  \rput(0,3.18){$\scs R$}
  \rput(.75,3.18){$\scs L$}
  \rput(1.3,3.18){$\scs R$}
  \rput(1.95,3.18){$\scs L$}
  \rput(2.5,3.18){$\scs R$}
  \rput(.7,-.18){$\scs L$}
  \rput(1.3,-.18){$\scs R$}
  \pspolygon(-1,0)(-1,3)(3,3)(3,0)(0,0)
\end{pspicture}
\qquad \; = \qquad \qquad
\begin{pspicture}[0.5](3.4,3.5)
  \rput(1,0){\multQ} \rput(1.6,1.5){\multQ} \rput(.4,1.5){\curveleftQ}
    \rput(-.5,3.18){$\scs L$}
  \rput(0,3.18){$\scs R$}
  \rput(.75,3.18){$\scs L$}
  \rput(1.3,3.18){$\scs R$}
  \rput(1.95,3.18){$\scs L$}
  \rput(2.5,3.18){$\scs R$}
  \rput(.7,-.18){$\scs L$}
  \rput(1.3,-.18){$\scs R$}
    \pspolygon(-1,0)(-1,3)(3,3)(3,0)(0,0)
\end{pspicture}
\]

\noindent To prove that these two string diagrams are equal we can
convert them back into the more traditional globular notation.
Notice that nothing interesting occurs with the identity
morphisms, $1_L$ on the far left, and $1_R$ on the far right so
the interesting part is what is happening in the middle.
Translating this into globular notation the proof is as follows:
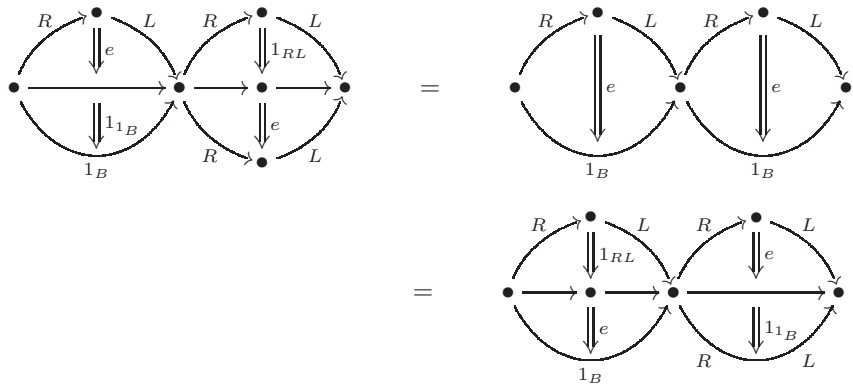
\[
 \xy
 (-23,0)*{};
(23,0)*{}; (-22,0)*+{\bullet}="L"; (0,0)*+{\bullet}="M";
(-11,10)*+{\bullet}="LT";
   (-11,8)="L1";
   (-11,2)="L2";
   (-11,-8)="L4";
   (-11,-2)="L3";
   "L";"LT" **\crv{(-20,7.5)};
               ?(.92)*\dir{>};
   "L";"M" **\crv{(-17,-12) & (-5,-12)};
               ?(.95)*\dir{>};
   "LT";"M" **\crv{(-2,7.5)};
               ?(.92)*\dir{>};
   {\ar "L";"M"};
   {\ar@{=>}^{e} "L1";"L2"};
   {\ar@{=>}^{1_{1_B}} "L3";"L4"};
       (-18,9)*{\scriptstyle R};
       (-4,9)*{\scriptstyle L};
       (-11,-11)*{\scriptstyle 1_B};
   (0,0)*+{\bullet}="L"; (22,0)*+{\bullet}="M";
(11,10)*+{\bullet}="LT"; (11,-10)*+{\bullet}="RB";
(11,0)*+{\bullet}="XX";
   (11,8)="L1";
   (11,2)="L2";
   (11,-8)="L4";
   (11,-2)="L3";
   "L";"LT" **\crv{(2,7.5)};
               ?(.92)*\dir{>};
   "L";"RB" **\crv{(2,-7.5)};
               ?(.92)*\dir{>};
   "LT";"M" **\crv{(20,7.5)};
               ?(.92)*\dir{>};
   "RB";"M" **\crv{(20,-7.5)};
               ?(.92)*\dir{>};
   {\ar "L";"XX"}; {\ar "XX";"M"};
   {\ar@{=>}^{1_{RL}} "L1";"L2"};
   {\ar@{=>}^{e} "L3";"L4"};
       (4,9)*{\scriptstyle R};
       (18,9)*{\scriptstyle L};
       (4,-9)*{\scriptstyle R};
       (18,-9)*{\scriptstyle L};
\endxy
\qquad = \qquad \xy (-23,0)*{}; (23,0)*{}; (-22,0)*+{\bullet}="L";
(0,0)*+{\bullet}="M"; (22,0)*+{\bullet}="R";
(-11,10)*+{\bullet}="LT"; (11,10)*+{\bullet}="RB";
   (-11,7)="L1";
   (-11,-7)="L2";
   (11,7)="R1";
   (11,-7)="R2";
   "L";"M" **\crv{(-17,-12) & (-5,-12)};
               ?(.95)*\dir{>};
   "M";"R" **\crv{ (5,-12)& (17,-12)};
               ?(.95)*\dir{>};
   "L";"LT" **\crv{(-20,7.5)};
               ?(.92)*\dir{>};
   "LT";"M" **\crv{(-2,7.5)};
               ?(.92)*\dir{>};
   "M";"RB" **\crv{(2,7.5)};
               ?(.92)*\dir{>};
   "RB";"R" **\crv{(20,7.5)};
               ?(.92)*\dir{>};
   {\ar@{=>}^{e} "L1";"L2"};
   {\ar@{=>}^{e} "R1";"R2"};
       (-18,9)*{\scriptstyle R};
       (-4,9)*{\scriptstyle L};
       (4,9)*{\scriptstyle R};
       (18,9)*{\scriptstyle L};
       (-11,-11)*{\scriptstyle 1_B};
       (11,-11)*{\scriptstyle 1_B};
\endxy
\]
\[
 \xy
(-23,0)*{}; (23,0)*{};
\endxy
\qquad = \qquad
 \xy
 (-23,0)*{};
(23,0)*{}; (-22,0)*+{\bullet}="L"; (0,0)*+{\bullet}="M";
(-11,0)*+{\bullet}="XX"; (-11,10)*+{\bullet}="LT";
   (-11,8)="L1";
   (-11,2)="L2";
   (-11,-8)="L4";
   (-11,-2)="L3";
      "L";"LT" **\crv{(-20,7.5)};
               ?(.92)*\dir{>};
   "LT";"M" **\crv{(-2,7.5)};
               ?(.92)*\dir{>};
    "L";"M" **\crv{(-17,-12) & (-5,-12)};
               ?(.95)*\dir{>};
   {\ar "L";"XX"}; {\ar "XX";"M"};
   {\ar@{=>}^{1_{RL}} "L1";"L2"};
   {\ar@{=>}^{e} "L3";"L4"};
       (-18,9)*{\scriptstyle R};
       (-4,9)*{\scriptstyle L};
       (-11,-11)*{\scriptstyle 1_B};
   (0,0)*+{\bullet}="L"; (22,0)*+{\bullet}="M";
(11,10)*+{\bullet}="LT";
   (11,8)="L1";
   (11,2)="L2";
   (11,-8)="L4";
   (11,-2)="L3";
   "L";"LT" **\crv{(2,7.5)};
               ?(.92)*\dir{>};
   "L";"M" **\crv{(5,-12) & (17,-12)};
               ?(.95)*\dir{>};
   "LT";"M" **\crv{(20,7.5)};
               ?(.92)*\dir{>};
   {\ar "L";"M"};
   {\ar@{=>}^{e} "L1";"L2"};
   {\ar@{=>}^{1_{1_B}} "L3";"L4"};
       (4,9)*{\scriptstyle R};
       (18,9)*{\scriptstyle L};
       (4,-9)*{\scriptstyle R};
       (18,-9)*{\scriptstyle L};
\endxy
\]
which amounts to nothing more than the interchange law and the
axiom for vertical composition of identities in a 2-category. In
this case, $1_{RL}.e = e = e.1_{1_B}$.

The unit axioms for the monoid require the following equations of
string diagrams:
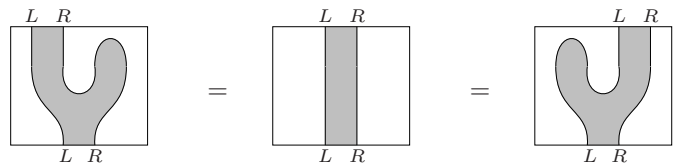
\[ 
\begin{pspicture}[0.3](2.3,3)
  \rput(1,0){\multQ} \rput(1.6,1.5){\birthQ} \rput(.4,1.5){\medidentQ}
  \pspolygon(-.3,0)(-.3,2.25)(2.3,2.25)(2.3,0)(0,0)
  \rput(.1,2.45){$\scs L$}
  \rput(.7,2.45){$\scs R$}
  \rput(.75,-.2){$\scs L$}
  \rput(1.3,-.2){$\scs R$}
\end{pspicture}
\qquad = \qquad
\begin{pspicture}[0.3](2.3,3)
  \rput(1,0){\identQ} \rput(1,1.5){\medidentQ}
  \pspolygon(-.3,0)(-.3,2.25)(2.3,2.25)(2.3,0)(0,0)
  \rput(.7,2.45){$\scs L$}
  \rput(1.3,2.45){$\scs R$}
  \rput(.7,-.2){$\scs L$}
  \rput(1.3,-.2){$\scs R$}
\end{pspicture}
\qquad = \qquad
\begin{pspicture}[0.3](2.3,3)
  \rput(1,0){\multQ} \rput(1.6,1.5){\medidentQ} \rput(.4,1.5){\birthQ}
  \pspolygon(-.3,0)(-.3,2.25)(2.3,2.25)(2.3,0)(0,0)
  \rput(1.3,2.45){$\scs L$}
  \rput(1.9,2.45){$\scs R$}
  \rput(.75,-.2){$\scs L$}
  \rput(1.3,-.2){$\scs R$}
\end{pspicture}
\]
\medskip

\noindent But these axioms follow directly from the zig-zag axioms
in the definition of an adjunction.  Thus, it is clear that the
walking monoid is contained in $\Hom(A,A)$ within the walking
adjunction.

Similarly, we will show that the morphism $RL$ is a comonoid in
the monoidal category $\Hom(B,B)$.  We define a comultiplication
for $RL$ to be the morphism $\Delta := RiL \maps RL \to RLRL$,
drawn diagrammatically as:
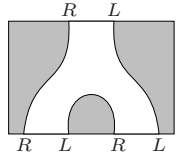
\[ \psset{xunit=1cm,yunit=1cm}
 \xy
  (0,-15)*{};
  (35,0)*{
   \begin{pspicture}[.5](2.2,1.9)
   \pspolygon[fillcolor=lightgray, fillstyle=solid](-1.1,0)(-1.1,1.5)(1.1,1.5)(1.1,0)(0,0)
  \pscustom[fillcolor=white, fillstyle=solid]{
        \psbezier(.9,0)(.8,.9)(.25,.6)(.3,1.5)
        \psline(-0.3,1.5)
        \psbezier(-0.3,1.5)(-.25,.6)(-.8,.9)(-.9,0)
        \psline(-.3,0)
        \psbezier(-.3,0)(-.4,.7)(0.4,.7)(.3,0)
        \psline(.9,0)
    }
  \rput(-.9,-.15){$\scs R$}
  \rput(-.35,-.15){$\scs L$}
  \rput(.35,-.15){$\scs R$}
  \rput(.9,-.15){$\scs L$}
  \rput(-.3,1.65){$\scs R$}
  \rput(.3,1.655){$\scs L$}
\end{pspicture}};
 \endxy
\]
The counit for the comultiplication is:
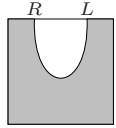
\[
   \begin{pspicture}(2.2,2.2)
  \pspolygon(0,0)(0,2)(2,2)(2,0)(0,0)
   \pscustom[fillcolor=lightgray, fillstyle=solid]{
  \psbezier(.5,2)(.45,.5)(1.55,.5)(1.5,2)
  \psline(2,2)
  \psline(2,0)
   \psline(0,0)
    \psline(0,2)
     \psline(.5,2)
  }
         \rput(1.5,2.2){$\scs L$}
         \rput(.5,2.2){$\scs R$}
 \end{pspicture}
\]
the counit for the adjunction $\varepsilon := e \maps 1_B \to RL$.
By similar arguments as those above, it follows that the walking
comonoid is contained in $\Hom(B,B)$ in the walking adjunction.

As we mentioned above, the main difficulty in proving Theorem
\ref{WAWM} is not in showing that the walking monoid is contained
in $\Hom(A,A)$ within the walking adjunction, but rather, in
showing that the walking monoid actually arises as $\Hom(A,A)$ in
the walking adjunction.  Doing so requires one to prove that the
monoid generating the walking monoid is $LR$, where
$(A,B,L,R,i,e)$ is the adjunction generating the walking
adjunction.  As a consequence, this means that for every monoid
$T$ in a monoidal category $\mathcal{D}$ there exists an
adjunction $(A,B,L,R,i,e)$ in some 2-category $\mathcal{C}$, where
$\mathcal{D}$ is a subcategory of $\Hom(A,A)$ and $T = LR$.
Loosely speaking, we must show that every monoid arises from an
adjunction.

The problem of showing that every monoid arises from an adjunction
turns out to be very related to the problem of showing that every
\textit{monad} arises from an adjunction.  Recall that a monad on
an object $A$ in a 2-category $\mathcal{C}$ is just a monoid in
the monoidal category $\Hom(A,A)$.  Given a monoid in a monoidal
category $\mathcal{D}$, we can regard $\mathcal{D}$ as a
2-category $\Sigma(\mathcal{D})$ with only one object, say
$\bullet$.  Then a monoid $T$ \textit{in} $\mathcal{D}$ becomes a
monad $\T$ \textit{on} the object $\bullet$ in the 2-category
$\Sigma(\mathcal{D})$.  Hence, showing that every monoid arises
from an adjunction can be deduced from showing that every monad
arises from an adjunction.

When monads where first discovered in the 1950's this question as
to whether every monad comes from an adjunction was raised by
Hilton and others~\cite{MacLane}. At this time, people where
mostly interested in monads that arose from adjoint functors, or
adjunctions in the 2-category $\cat{Cat}$. If $\mathcal{A}$ is a
category and $T \maps A \to A$ is a functor defining a monad on
$\mathcal{A}$, then in this context two well known solutions
appeared, the Kleisli construction
$\adjunction{\mathcal{A}}{\mathcal{A}_{\T}}{}{}$~\cite{Kleisli},
and the Eilenberg-Moore construction
$\adjunction{\mathcal{A}}{\mathcal{A}^{\T}}{}{}$~\cite{EM}. These
two solutions are, in a certain sense, the initial and terminal
solution to the problem of constructing such an adjunction.
Similarly, a comonoid in $\Hom(\mathcal{A},\mathcal{A})$ is known
as a comonad, and these constructions work equally well to create
a pair of adjoint functors where the functor $\mathcal{A}^{\T} \to
\mathcal{A}$ is now the left adjoint.

For our purposes we will need to consider monads in 2-categories
other than $\cat{Cat}$.  In particular, we would like to consider
the 2-category $\Sigma(\Mon)$, the suspension of the walking
monoid.  Unfortunately, the Eilenberg-Moore and the Kleisli
construction do not work in the completely general context of an
arbitrary 2-category.  While it is true that every adjunction in a
2-category $\mathcal{C}$ produces a monad, it is not always true
that one can find an adjunction in $\mathcal{C}$ generating a
given monad. The failure of this construction can be attributed to
the lack of an object in $\mathcal{C}$ to play the role of the
Eilenberg-Moore category of algebras (or the lack of a Kleisli
object, but we will focus on Eilenberg-Moore objects in this
paper). When such an object does exist we call it an
Eilenberg-Moore object for the monad $\T$. The existence of
Eilenberg-Moore objects in a 2-category $\mathcal{C}$ is a
completeness property of the 2-category in question. In
particular, $\mathcal{C}$ has Eilenberg-Moore objects if it is
finitely complete as a 2-category~\cite{LS1,Street3}.

Since the 2-category $\Sigma(\Mon)$ has only one object it is
obvious that there will not be an object in $\Sigma(\Mon)$ to play
the role of an Eilenberg-Moore object for the monad $\T$.
Fortunately, there is a categorical construction known as the free
completion under Eilenberg-Moore objects that takes a 2-category
$\mathcal{C}$ and enlarges it into a 2-category
$\cat{EM}(\mathcal{C})$ that has Eilenberg-Moore objects for every
monad $\T$ in $\mathcal{C}$.  This means that every monad $\T$ in
$\mathcal{C}$ arises from an adjunction in
\cat{EM}$(\mathcal{C})$.  Furthermore, there is a fully faithful
embedding $Z \maps \mathcal{C} \to \cat{EM}(\mathcal{C})$ with the
property that for any other 2-category $\mathcal{C}'$ with
Eilenberg-Moore objects, composition with $Z$ induces an
equivalence of categories between the functor category
$[\mathcal{C},\mathcal{C}']$ and the full subcategory of the
functor category $[\cat{EM}(\mathcal{C}), \mathcal{C'}]$
consisting of those 2-functors that preserve Eilenberg-Moore
objects~\cite{LS2}.

This completion is possible because Eilenberg-Moore objects can be
described as a weighted limit~\cite{Street3} whose weight is
finite in the sense of~\cite{kel2}.   This also means that if we
are only interested in an Eilenberg-Moore object for a single
monad $\T$ in $\mathcal{C}$, then we do not have to complete
$\mathcal{C}$ under Eilenberg-Moore objects for \textit{every}
monad $\T$ in $\mathcal{C}$.  We can instead define the 2-category
$\cat{EM}_{\T}(\mathcal{C})$, the free completion of $\mathcal{C}$
under an Eilenberg-Moore object for the monad $\T$.  If $\T$ is a
monad on the object $A$ of $\mathcal{C}$, this 2-category will
contain $A$ (identified with its image under the embedding), and
an Eilenberg-Moore object $A^{\T}$ for the monad $\T$. Hence, in
$\cat{EM}_{\T}(\mathcal{C})$ the monad $\T$ is generated by an
adjunction $\adjunction{A}{A^{\T}}{}{}$.

This construction is particularly well suited for the problem at
hand.  The Eilenberg-Moore completion works just as well, when the
2-category $\mathcal{C}$ has only one object: that is, when
$\mathcal{C} = \Sigma(\mathcal{D})$ is the suspension of a
monoidal category $\mathcal{D}$. As we mentioned above, a monad in
the 2-category $\Sigma(\mathcal{D})$ is just a monoid in the
monoidal category $\mathcal{D}$.  Hence, we have shown that every
monoid $T$ in a monoidal category $\mathcal{D}$ arises from an
adjunction in the 2-category
$\cat{EM}_{T}\big(\Sigma(\mathcal{D})\big)$. For more on the
Eilenberg-Moore completion see~\cite{LS2}, or~\cite{Lau1} for an
explicit description of $\cat{EM}(\cat{Vect})$.

We are now ready to prove Theorem~\ref{WAWM}.  Note that this
theorem follows as a decategorification of a theorem due to Lack
\cite{Lack}.  We present the proof for completeness.

\paragraph{Proof of Theorem~\ref{WAWM}.}   Let $T$ be the monoid
generating the walking monoid.  We will show that the 2-category
$\cat{EM}_{T}\big(\Sigma(\Mon)\big)$ is isomorphic to the walking
adjunction $\Adj$.  By the universal property of the walking
adjunction, adjunctions in a 2-category $\mathcal{C}$ correspond
bijectively to 2-functors $\Adj \to \mathcal{C}$.  Since the
2-category $\cat{EM}_{T}\big(\Sigma(\Mon)\big)$ contains an
adjunction $\adjunction{A}{A^{\T}}{}{}$ we get a 2-functor
$\Lambda \maps \Adj \to \cat{EM}_{T}\big(\Sigma(\Mon)\big)$. If
$\Adj$ is generated by the adjunction $(A,B,L,R,i,e)$, then
$\Lambda$ maps this generating adjunction to the adjunction
$\adjunction{A}{A^{\T}}{}{}$ in
$\cat{EM}_{T}\big(\Sigma(\Mon)\big)$.

We now construct the inverse of the 2-functor $\Lambda$. We have
already shown that given an adjunction $(A,B,L,R,i,e)$, then the
map $LR$ is a monad on $A$ (equivalently $LR$ is a monoid in
$\Hom(A,A)$). Hence, by the universal property of the walking
monoid we get a 2-functor $\Sigma(\Mon) \to \Adj$. In the walking
adjunction this monad has an Eilenberg-Moore object, namely $B$.
Thus, the universal property of the Eilenberg-Moore completion
determines a 2-functor $\bar{\Lambda} \maps
\cat{EM}_{T}\big(\Sigma(\Mon)\big) \to \Adj$ that preserves
Eilenberg-Moore objects.  By their construction, it is clear that
the composites of $\Lambda$ and $\bar{\Lambda}$ are equal to the
identity.  Hence, $\cat{EM}_{T}\big(\Sigma(\Mon)\big) \cong \Adj$.
\qed

From our topological perspective it is clear that we could not
define a comultiplication in $\Hom(A,A)$ since this would require
a map:
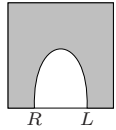
\[
   \begin{pspicture}(2.2,2.2)
  \pspolygon(0,0)(0,2)(2,2)(2,0)(0,0)
   \pscustom[fillcolor=lightgray, fillstyle=solid]{
  \psbezier(.5,0)(.45,1.5)(1.55,1.5)(1.5,0)
  \psline(2,0)
  \psline(2,2)
   \psline(0,2)
    \psline(0,0)
     \psline(.5,0)
  }
  \rput(1.5,-.2){$\scs L$}
         \rput(.5,-.2){$\scs R$}
 \end{pspicture}
\]
which does not exist in the walking adjunction.  In the next
section we will consider the categorical framework where this map
is given, the walking ambidextrous adjunction.

\subsection{The walking ambidextrous adjunction} \label{secWAA}

In this section we further examine adjunctions by looking at
adjunctions that are 2-sided.  This means that in addition to the
2-morphisms $i \maps 1_A \To LR$ and $e \maps RL \To 1_B$, there
are also 2-morphisms $j \maps 1_B \To RL$ and $k \maps LR \To 1_A$
satisfying the triangle identities.  Sometimes category theorists
will specify which of the two possible adjunctions they mean by
referring to one as a left adjunction and the other as a right
adjunction. Since we will consider adjunctions that are both left
and right adjunctions, we will call them \textit{ambidextrous
adjunctions} to indicate their `two-handedness'. Sometimes we call
an ambidextrous adjunction an \textit{ambijunction} for short.

Our primary interest is the walking ambijunction. This has also
been referred to as the free biadjunction. The authors is hesitant
to use this terminology because of possible confusion that may
arise when considering morphisms of bicategories.

\begin{defn} \et \label{DefWalkAmbi}
The {\em walking ambidextrous adjunction} $\Ambi$ is the
2-category freely generated by:
\begin{itemize}
    \item objects $A$ and $B$,
    \item morphisms $L \maps A \to B$ and $R \maps B \to A$, and
    \item 2-morphisms $i \maps 1_A \to L \circ R$, $e \maps R \circ
    L$, $j \maps 1_B \to R \circ L$, and $k \maps L \circ R \to 1_A$
\end{itemize}
such that
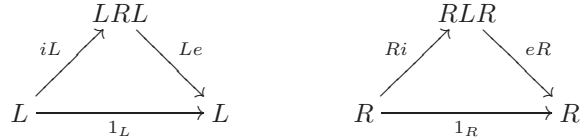
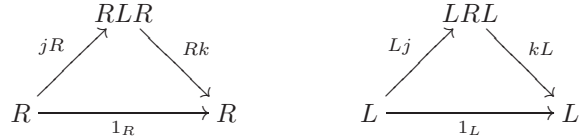
\[
 \xymatrix@C=1.5pc@R=.8pc{
  & LRL \ar[ddr]^-{Le} & \\ \\
 L \ar[uur]^-{iL} \ar[rr]_-{1_L} & & L
 }
\qquad \qquad
 \xymatrix@C=1.5pc@R=.8pc{
  & RLR \ar[ddr]^-{eR} & \\ \\
 R \ar[uur]^-{Ri} \ar[rr]_-{1_R} & & R
 }
\]
and
\[
 \xymatrix@C=1.5pc@R=.8pc{
  & RLR \ar[ddr]^-{Rk} & \\ \\
 R \ar[uur]^-{jR} \ar[rr]_-{1_R} & & R
 }
\qquad \qquad
 \xymatrix@C=1.5pc@R=.8pc{
  & LRL \ar[ddr]^-{kL} & \\ \\
 L \ar[uur]^-{Lj} \ar[rr]_-{1_L} & & L
 }
\]
commute.
\end{defn}

For later convenience we depict the 2-morphisms in the walking
ambidextrous adjunction in string notation:
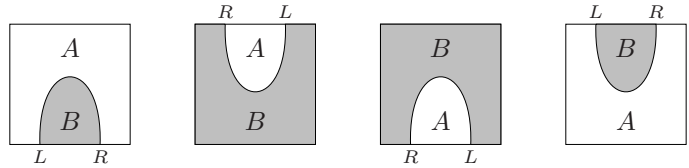
\[\psset{xunit=.80cm,yunit=0.80cm}
   \begin{pspicture}(2.2,2.2)
  \pspolygon(0,0)(0,2)(2,2)(2,0)(0,0)
   \pscustom[fillcolor=lightgray, fillstyle=solid]{
  \psbezier(.5,0)(.45,1.5)(1.55,1.5)(1.5,0)
  \psline(.5,0)
  }
        \rput(1,.4){$B$}
         \rput(1,1.6){$A$}
         \rput(1.5,-.2){$\scs R$}
         \rput(.5,-.2){$\scs L$}
 \end{pspicture}
 \qquad
   \begin{pspicture}(2.2,2.2)
  \pspolygon(0,0)(0,2)(2,2)(2,0)(0,0)
   \pscustom[fillcolor=lightgray, fillstyle=solid]{
  \psbezier(.5,2)(.45,.5)(1.55,.5)(1.5,2)
  \psline(2,2)
  \psline(2,0)
   \psline(0,0)
    \psline(0,2)
     \psline(.5,2)
  }
         \rput(1,.4){$B$}
         \rput(1,1.6){$A$}
         \rput(1.5,2.2){$\scs L$}
         \rput(.5,2.2){$\scs R$}
 \end{pspicture}
 \qquad
   \begin{pspicture}(2.2,2.2)
  \pspolygon(0,0)(0,2)(2,2)(2,0)(0,0)
   \pscustom[fillcolor=lightgray, fillstyle=solid]{
  \psbezier(.5,0)(.45,1.5)(1.55,1.5)(1.5,0)
  \psline(2,0) \psline(2,2) \psline(0,2) \psline(0,0) \psline(.5,0)
  }
        \rput(1,.4){$A$}
         \rput(1,1.6){$B$}
         \rput(1.5,-.2){$\scs L$}
         \rput(.5,-.2){$\scs R$}
 \end{pspicture}
  \qquad
   \begin{pspicture}(2.2,2.2)
  \pspolygon(0,0)(0,2)(2,2)(2,0)(0,0)
   \pscustom[fillcolor=lightgray, fillstyle=solid]{
  \psbezier(.5,2)(.45,.5)(1.55,.5)(1.5,2)
  \psline(.5,2)
  }
         \rput(1,.4){$A$}
         \rput(1,1.6){$B$}
         \rput(1.5,2.2){$\scs R$}
         \rput(.5,2.2){$\scs L$}
 \end{pspicture}
\]
The zig-zag laws for the four maps above are depicted in string
notation as:
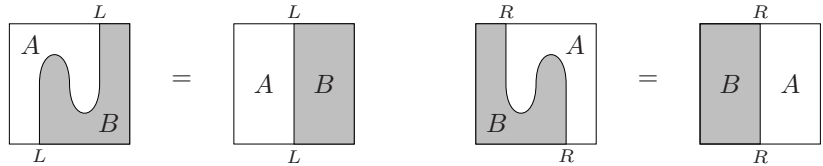
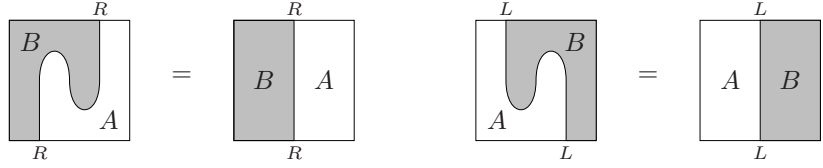
\[ \psset{xunit=.80cm,yunit=0.80cm}
\xy
   (0,0)*+{\begin{pspicture}(2,2.2)
  \pspolygon(0,0)(0,2)(2,2)(2,0)(0,0)
\pscustom[fillcolor=lightgray, fillstyle=solid]{
  \psbezier(.5,1)(.5,1.65)(1,1.65)(1,1)
  \psbezier(1,1)(1,.35)(1.5,.35)(1.5,1)
  \psline(1.5,2)
  \psline(2,2)
  \psline(2,0)
  \psline(.5,0)
  \psline(.5,1)
  }
         \rput(.35,1.65){$A$}
         \rput(1.65,.35){$B$}
         \rput(1.5,2.2){$\scs L$}
         \rput(.5,-.2){$\scs L$}
 \end{pspicture}};
 \endxy
\quad = \quad
  \xy
   (0,0)*+{  \begin{pspicture}(2,2.2)
  \pspolygon(0,0)(0,2)(2,2)(2,0)(0,0)
\pscustom[fillcolor=lightgray, fillstyle=solid]{
  \pspolygon(1,0)(1,2)(2,2)(2,0)(1,0)
  }
  \rput(.5,1){$A$}
         \rput(1.5,1){$B$}
         \rput(1,2.2){$\scs L$}
         \rput(1,-.2){$\scs L$}
 \end{pspicture}};
 \endxy
 \qquad \qquad
 \psset{xunit=.80cm,yunit=0.80cm}
     \xy
    (0,0)*+{   \begin{pspicture}(2,2.2)
  \pspolygon(0,0)(0,2)(2,2)(2,0)(0,0)
\pscustom[fillcolor=lightgray, fillstyle=solid]{
  \psbezier(.5,1)(.5,.35)(1,.35)(1,1)
  \psbezier(1,1)(1,1.65)(1.5,1.65)(1.5,1)
  \psline(1.5,0)
  \psline(0,0)
  \psline(0,2)
  \psline(.5,2)
  \psline(.5,1)
  }
          \rput(1.65,1.65){$A$}
         \rput(.35,.35){$B$}
         \rput(.5,2.2){$\scs R$}
         \rput(1.5,-.2){$\scs R$}
 \end{pspicture}};
  \endxy
 \quad = \quad
 \xy
      (0,0)*+{\begin{pspicture}(2,2.2)
\pscustom[fillcolor=lightgray, fillstyle=solid]{
  \pspolygon(1,0)(1,2)(0,2)(0,0)(1,0)
  }
   \rput(.5,1){$B$}
         \rput(1.5,1){$A$}
         \rput(1,2.2){$\scs R$}
         \rput(1,-.2){$\scs R$}
  \pspolygon(0,0)(0,2)(2,2)(2,0)(0,0)
 \end{pspicture}};
    \endxy
\]

\[ \psset{xunit=.80cm,yunit=0.80cm}
\xy
   (0,0)*+{\begin{pspicture}(2,2.2)
  \pspolygon(0,0)(0,2)(2,2)(2,0)(0,0)
\pscustom[fillcolor=lightgray, fillstyle=solid]{
  \psbezier(.5,1)(.5,1.65)(1,1.65)(1,1)
  \psbezier(1,1)(1,.35)(1.5,.35)(1.5,1)
  \psline(1.5,2)
  \psline(0,2)
  \psline(0,0)
  \psline(.5,0)
  \psline(.5,1)
  }
         \rput(.35,1.65){$B$}
         \rput(1.65,.35){$A$}
         \rput(1.5,2.2){$\scs R$}
         \rput(.5,-.2){$\scs R$}
 \end{pspicture}};
 \endxy
\quad = \quad
 \xy
      (0,0)*+{\begin{pspicture}(2,2.2)
\pscustom[fillcolor=lightgray, fillstyle=solid]{
  \pspolygon(1,0)(1,2)(0,2)(0,0)(1,0)
  }
   \rput(.5,1){$B$}
         \rput(1.5,1){$A$}
         \rput(1,2.2){$\scs R$}
         \rput(1,-.2){$\scs R$}
  \pspolygon(0,0)(0,2)(2,2)(2,0)(0,0)
 \end{pspicture}};
    \endxy
 \qquad \qquad
 \psset{xunit=.80cm,yunit=0.80cm}
     \xy
    (0,0)*+{   \begin{pspicture}(2,2.2)
  \pspolygon(0,0)(0,2)(2,2)(2,0)(0,0)
\pscustom[fillcolor=lightgray, fillstyle=solid]{
  \psbezier(.5,1)(.5,.35)(1,.35)(1,1)
  \psbezier(1,1)(1,1.65)(1.5,1.65)(1.5,1)
  \psline(1.5,0)
  \psline(2,0)
  \psline(2,2)
  \psline(.5,2)
  \psline(.5,1)
  }
          \rput(1.65,1.65){$B$}
         \rput(.35,.35){$A$}
         \rput(.5,2.2){$\scs L$}
         \rput(1.5,-.2){$\scs L$}
 \end{pspicture}};
  \endxy
 \quad = \quad
\xy
   (0,0)*+{  \begin{pspicture}(2,2.2)
  \pspolygon(0,0)(0,2)(2,2)(2,0)(0,0)
\pscustom[fillcolor=lightgray, fillstyle=solid]{
  \pspolygon(1,0)(1,2)(2,2)(2,0)(1,0)
  }
  \rput(.5,1){$A$}
         \rput(1.5,1){$B$}
         \rput(1,2.2){$\scs L$}
         \rput(1,-.2){$\scs L$}
 \end{pspicture}};
 \endxy
\]
Notice that the walking ambidextrous adjunction has the same maps
and axioms as the walking adjunction but with the color inverted
versions as well.

As eluded to in the previous section, the importance of the
walking ambidextrous adjunction is its relationship to the walking
Frobenius algebra.  Understanding the relationship between the two
provides a characterization of Frobenius algebras that, phrased in
an intrinsically categorical way, easily admits categorification
to provide a definition of a pseudo Frobenius algebra. Below we
define the walking Frobenius algebra based on
Proposition~\ref{equivFrob} (i.), although any of the three
definitions would produce equivalent monoidal categories.

\begin{defn}
The {\em walking Frobenius algebra} $\Frob$ is the monoidal
category freely generated by:
\begin{itemize}
    \item objects $A$ and $I$, and
    \item morphisms $m \maps A \ten A \to A$, $\iota \maps I \to
    A$, $\Delta \maps A  \to A\ten A$ and $\varepsilon \maps A\to I$
\end{itemize}
such that
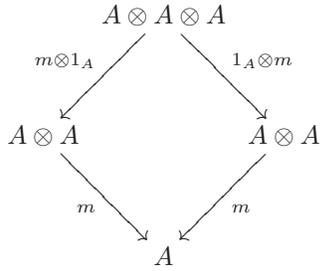
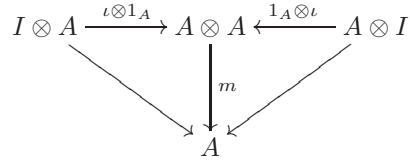
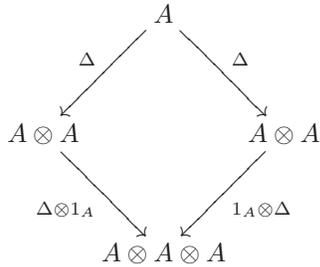
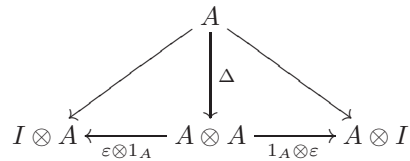
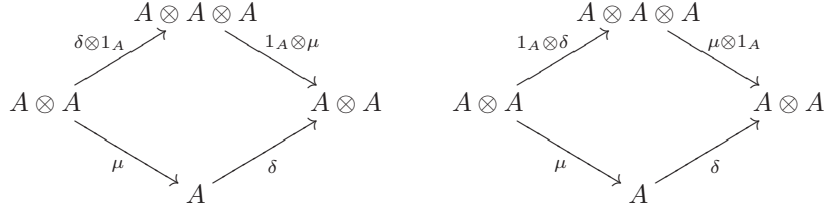
\[
 \xy
 (16,0)*+{A \ten A}="R";
 (0,16)*+{A \ten A \ten A}="T";
 (0,-16)*+{A}="B";
 (-16,0)*+{A \ten A}="L";
    {\ar^{1_A \ten m} "T";"R"};
    {\ar_{m \ten 1_A} "T";"L"};
    {\ar_{ m} "L";"B"};
    {\ar^{ m} "R";"B"};
 \endxy
\quad \quad
 \vcenter{\xy
 (-22,0)*+{I \ten A}="L";
 (0,0)*+{A \ten A}="M";
 (22,0)*+{A \ten I}="R";
 (0,-16)*+{A}="B";
    {\ar^{\iota \ten 1_A} "L";"M"};
    {\ar_{1_A \ten \iota} "R";"M"};
    {\ar_{ } "L";"B"};
    {\ar^{ } "R";"B"};
    {\ar^{ m} "M";"B"};
 \endxy}
\]
\[
\xy
 (16,0)*+{A \ten A}="R";
 (0,-16)*+{A \ten A \ten A}="T";
 (0,16)*+{A}="B";
 (-16,0)*+{A \ten A}="L";
    {\ar^{1_A \ten \Delta} "T";"R";};
    {\ar_{\Delta \ten 1_A} "T";"L";};
    {\ar_{ \Delta} "L";"B";};
    {\ar^{ \Delta} "R";"B";};
 \endxy
\quad \quad
 \vcenter{\xy
 (-22,0)*+{I \ten A}="L";
 (0,0)*+{A \ten A}="M";
 (22,0)*+{A \ten I}="R";
 (0,16)*+{A}="B";
    {\ar^{\varepsilon \ten 1_A} "M";"L"};
    {\ar_{1_A \ten \varepsilon} "M";"R"};
    {\ar_{ } "L";"B";};
    {\ar^{ } "R";"B";};
    {\ar^{ \Delta} "M";"B";};
 \endxy}
\]
and
\[
 \xy
 (20,0)*+{A \ten A}="R";
 (0,12)*+{A \ten A \ten A}="T";
 (0,-12)*+{A}="B";
 (-20,0)*+{A \ten A}="L";
    {\ar^-{1_A \ten \mu} "T";"R"};
    {\ar^-{\delta \ten 1_A} "L";"T"};
    {\ar_-{ \mu} "L";"B"};
    {\ar_-{ \delta} "B";"R"};
 \endxy
\quad \quad
 \xy
 (20,0)*+{A \ten A}="R";
 (0,12)*+{A \ten A \ten A}="T";
 (0,-12)*+{A}="B";
 (-20,0)*+{A \ten A}="L";
    {\ar^-{\mu \ten 1_A} "T";"R"};
    {\ar^-{1_A \ten \delta } "L";"T"};
    {\ar_-{ \mu} "L";"B"};
    {\ar_-{ \delta} "B";"R"};
 \endxy
\]
commute.
\end{defn}
For more on the walking Frobenius algebra see~\cite{ThisWeek,
Kock}.

\begin{thm} \et \label{WalkingBiAjunction}
The monoidal category $\Hom(A,A)$ in the walking ambidextrous
adjunction is the {\em walking Frobenius algebra} (equivalently
$\Hom(B,B)$).
\end{thm}

\paragraph{Proof.} In Theorem~\ref{WAWM} we saw that the
object $LR$ in $\Hom(A,A)$ had a monoidal structure given by $LeR
\maps LRLR \to LR$ with unit $i \maps 1_A \to LR$. Define a
comultiplication on $LR$ by the morphism $LjR \maps  LR \to LRLR$.
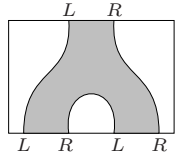
\[ \psset{xunit=1cm,yunit=1cm}
 \xy
  (0,-12)*{};
  (35,0)*{
   \begin{pspicture}(2.2,1.8)
  \pspolygon(-1.1,0)(-1.1,1.5)(1.1,1.5)(1.1,0)(0,0)
  \rput(0,0){\comult}
  \rput(-.9,-.15){$\scs L$}
  \rput(-.35,-.15){$\scs R$}
  \rput(.35,-.15){$\scs L$}
  \rput(.9,-.15){$\scs R$}
  \rput(-.3,1.65){$\scs L$}
  \rput(.3,1.65){$\scs R$}
 \end{pspicture}};
 \endxy
\]
The counit for this comultiplication is $k \maps LR \to 1_A$. 
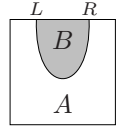
\[
   \begin{pspicture}(2.2,2.2)
  \pspolygon(0,0)(0,2)(2,2)(2,0)(0,0)
   \pscustom[fillcolor=lightgray, fillstyle=solid]{
  \psbezier(.5,2)(.45,.5)(1.55,.5)(1.5,2)
  \psline(.5,2)
  }
         \rput(1,.4){$A$}
         \rput(1,1.6){$B$}
         \rput(1.5,2.2){$\scs R$}
         \rput(.5,2.2){$\scs L$}
 \end{pspicture}
 \]
To show that the comultiplication is coassociative  take the proof
that $\Hom(B,B)$ is a comonoid object and invert the colors of the
shaded regions. All that remains to be shown is that the monoid
and comonoid structures are compatible, that is, they must satisfy
the Frobenius identities:

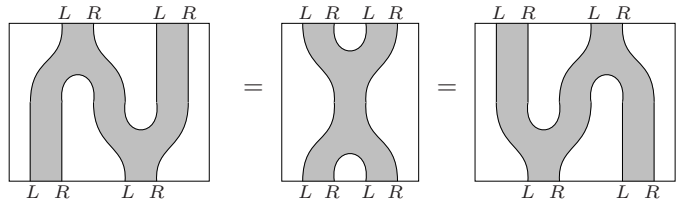
\[ 
\begin{pspicture}[.5](3,3.3)
 \rput(1,1.5){\comultQ} \rput(2.2,0){\multQ}
 \rput(2.8,1.5){\identQ} \rput(.4,0){\identQ}
 \pspolygon(-.3,0)(-.3,3)(3.5,3)(3.5,0)(0,0)
 \rput(.75,3.2){$\scs L$}
  \rput(1.3,3.2){$\scs R$}
  \rput(2.55,3.2){$\scs L$}
  \rput(3.1,3.2){$\scs R$}
  \rput(.15,-.2){$\scs L$}
  \rput(.7,-.2){$\scs R$}
  \rput(1.95,-.2){$\scs L$}
  \rput(2.5,-.2){$\scs R$}
\end{pspicture}
\qquad = \quad
\begin{pspicture}[0.5](2,3.3)
 \rput(1,0){\comultQ} \rput(1,1.5){\multQ}
 \pspolygon(-.3,0)(-.3,3)(2.3,3)(2.3,0)(0,0)
   \rput(.15,-.2){$\scs L$}
  \rput(.65,-.2){$\scs R$}
  \rput(1.35,-.2){$\scs L$}
  \rput(1.85,-.2){$\scs R$}
  \rput(.15,3.2){$\scs L$}
  \rput(.65,3.2){$\scs R$}
  \rput(1.35,3.2){$\scs L$}
  \rput(1.85,3.2){$\scs R$}
\end{pspicture}
\quad = \quad
\begin{pspicture}[.5](3,3.3)
 \rput(2.2,1.5){\comultQ} \rput(1,0){\multQ}
 \rput(.4,1.5){\identQ} \rput(2.8,0){\identQ}
 \pspolygon(-.3,0)(-.3,3)(3.5,3)(3.5,0)(0,0)
   \rput(.15,3.2){$\scs L$}
  \rput(.7,3.2){$\scs R$}
  \rput(1.95,3.2){$\scs L$}
  \rput(2.5,3.2){$\scs R$}
 \rput(.7,-.2){$\scs L$}
  \rput(1.25,-.2){$\scs R$}
  \rput(2.5,-.2){$\scs L$}
  \rput(3.05,-.2){$\scs R$}
\end{pspicture}
\]
\medskip

\noindent This follows directly from the identity axioms and the
interchange law relating vertical and horizontal composition in a
2-category. The first equality is proved in globular notation
below:

\[
\xy (-23,0)*{}; (23,0)*{}; (-22,0)*+{\bullet}="L";
(0,0)*+{\bullet}="M"; (22,0)*+{\bullet}="R";
(-11,10)*+{\bullet}="LT"; (11,-10)*+{\bullet}="RB";
   (-11,7)="L1";
   (-11,-7)="L2";
   (11,7)="R1";
   (11,-7)="R2";
   "L";"M" **\crv{(-17,-12) & (-5,-12)};
               ?(.95)*\dir{>};
   "M";"R" **\crv{ (5,12)& (17,12)};
               ?(.95)*\dir{>};
   "L";"LT" **\crv{(-20,7.5)};
               ?(.92)*\dir{>};
   "LT";"M" **\crv{(-2,7.5)};
               ?(.92)*\dir{>};
   "M";"RB" **\crv{(2,-7.5)};
               ?(.92)*\dir{>};
   "RB";"R" **\crv{(20,-7.5)};
               ?(.92)*\dir{>};
   {\ar@{=>}^{k} "L1";"L2"};
   {\ar@{=>}^{i} "R1";"R2"};
       (-18,9)*{\scriptstyle L};
       (-4,9)*{\scriptstyle R};
       (4,-9)*{\scriptstyle L};
       (18,-9)*{\scriptstyle R};
       (-11,-11)*{\scriptstyle 1_A};
       (11,11)*{\scriptstyle 1_A};
\endxy
\qquad = \qquad
  \xy
  (-23,0)*{};
(23,0)*{}; (-22,0)*+{\bullet}="L"; (0,0)*+{\bullet}="M";
(22,0)*+{\bullet}="R"; (-11,10)*+{\bullet}="LT";
(11,-10)*+{\bullet}="RB";
   (-11,8)="L1";
   (-11,2)="L2";
   (11,-2)="R1";
   (11,-8)="R2";
   (-11,-8)="L4";
   (-11,-2)="L3";
   (11,2)="R4";
   (11,7.5)="R3";
   "L";"M" **\crv{(-17,-12) & (-5,-12)};
               ?(.95)*\dir{>};
   "M";"R" **\crv{ (5,12)& (17,12)};
               ?(.95)*\dir{>};
   "L";"LT" **\crv{(-20,7.5)};
               ?(.92)*\dir{>};
   "LT";"M" **\crv{(-2,7.5)};
               ?(.92)*\dir{>};
   "M";"RB" **\crv{(2,-7.5)};
               ?(.92)*\dir{>};
   "RB";"R" **\crv{(20,-7.5)};
               ?(.92)*\dir{>};
   {\ar "L";"M"};
   {\ar "M";"R"};
   {\ar@{=>}^{k} "L1";"L2"};
   {\ar@{=>}^{i} "R1";"R2"};
   {\ar@{=>}^{1_{1_A}} "L3";"L4"};
   {\ar@{=>}^{1_{1_A}} "R3";"R4"};
       (-18,9)*{\scriptstyle L};
       (-4,9)*{\scriptstyle R};
       (4,-9)*{\scriptstyle L};
       (18,-9)*{\scriptstyle R};
       (-11,-11)*{\scriptstyle 1_A};
       (11,11)*{\scriptstyle 1_A};
\endxy
\]
\[
 \xy
(-23,0)*{}; (23,0)*{};
\endxy
\qquad = \qquad
 \xy
 (-23,0)*{};
(23,0)*{}; (-22,0)*+{\bullet}="L"; (0,0)*+{\bullet}="M";
(-11,10)*+{\bullet}="LT"; (-11,-10)*+{\bullet}="RB";
   (-11,8)="L1";
   (-11,2)="L2";
   (-11,-8)="L4";
   (-11,-2)="L3";
   "L";"LT" **\crv{(-20,7.5)};
               ?(.92)*\dir{>};
   "L";"RB" **\crv{(-20,-7.5)};
               ?(.92)*\dir{>};
   "LT";"M" **\crv{(-2,7.5)};
               ?(.92)*\dir{>};
   "RB";"M" **\crv{(-2,-7.5)};
               ?(.92)*\dir{>};
   {\ar "L";"M"};
   {\ar@{=>}^{k} "L1";"L2"};
   {\ar@{=>}^{i} "L3";"L4"};
       (-18,9)*{\scriptstyle L};
       (-4,9)*{\scriptstyle R};
       (-18,-9)*{\scriptstyle L};
       (-4,-9)*{\scriptstyle R};
\endxy
\]
and the other Frobenius identity follows similarly.  Hence,
$\Hom(A,A)$ in the walking ambidextrous adjunction contains the
walking Frobenius algebra.

To prove the converse we will again borrow some results from monad
theory.  Extending the work of Street~\cite{StreetFrob} and
Eilenberg and Moore~\cite{EM}, the author has shown that every
Frobenius algebra $F$ in an arbitrary monoidal category
$\mathcal{D}$ is generated by an ambidextrous adjunction in the
2-category \cat{EM}$(\Sigma(\mathcal{D}))$.  This construction
uses the fact the a Frobenius algebra $F$ in the monoidal category
$\mathcal{D}$ defines a Frobenius monad
$\F$~\cite{law,StreetFrob}, or a monad and a comonad on the one
object of the 2-category $\Sigma(\mathcal{D})$. This monad and
comonad arising from a Frobenius object has a special property
that makes the Eilenberg-Moore object for the monad isomorphic to
the Eilenberg-Moore object for the comonad. Thus, freely
completing $\mathcal{D}$ under Eilenberg-Moore objects for the
monad suffices to produce an ambidextrous adjunction in
$\cat{EM}\big(\Sigma(\mathcal{D})\big)$ that generates $F$.

We do not have to complete the monoidal category $\mathcal{D}$
under Eilenberg-Moore objects for every Frobenius algebra in
$\mathcal{D}$.  Indeed, if $F$ is a Frobenius algebra in
$\mathcal{D}$ then we can define the free completion of
$\Sigma(\mathcal{D})$ under an Eilenberg-Moore object for the
single Frobenius monad $\F$, denoted
$\cat{EM}_{\F}\big(\Sigma(\mathcal{D})\big)$.  In
$\cat{EM}_{\F}\big(\Sigma(\mathcal{D})\big)$ the monad $\F$ is
generated by an ambidextrous adjunction
$\ambijunction{A}{A^{\T}}{}{}$. Hence, by similar arguments as the
proof of Theorem \ref{WAWM}, we have that the 2-category $\Ambi$
is isomorphic to the 2-category
$\cat{EM}_{\F}\big(\Sigma(\mathcal{\Frob})\big)$, so that $\Frob$
really is $\Hom(A,A)$ in the walking ambidextrous adjunction. \qed

We then have the following corollary:
\begin{cor}
Every 2D topological quantum field theory, in the sense of
Atiyah~\cite{Atiyah}, arises from an ambijunction in the
2-category \cat{EM}$\big(\Sigma(\cat{Vect})\big)$.
\end{cor}

\Proof  A 2D topological quantum field theory is a monoidal
functor from \cat{2Cob} into \cat{Vect}.  It is well known that
such a functor amounts to a commutative Frobenius
algebra~\cite{Kock}. Hence, the result follows. \qed

This result implies that all of the known 2-dimensional
topological quantum field theories that have been constructed in
the axiomatic sense can be understood as arising from an
ambijunction in some 2-category.

\subsection{Two-dimensional thick tangles} \label{sec2Dthick}

Category theory has been used as a language to describe
relationships in topology, especially in algebraic topology;  but
perhaps the most exciting interplay between category theory and
topology comes from understanding various types of topology as
categories with extra structure.  By seeing a category as a
structure rather than a means to describe structure,  progress has
been made in fields that had at first seemed quite mysterious.
Aside from the description of $\cat{2Cob}$ as the free symmetric
monoidal category on a commutative Frobenius algebra, the most
notable instance is the category of tangles in 3-dimensional
space.  This category has a completely algebraic description as
the free braided monoidal category with duals on one object
\cite{tur1,FY,js2,Shum}.   Using this universal property, it is
easy to construct functors from the category of tangles into other
braided monoidal categories with duals, such as the category of
representations of a quantum group. Furthermore, any such functor
determines an invariant of tangles, and in particular, a knot
invariant.  This categorical description plays a vital role in
understanding the Jones polynomial and other `quantum invariants'
of knots~\cite{RT}.

In this section we provide yet another example of this phenomenon.
First we define the topological category $\cat{2Thick}$ of
two-dimensional thick tangles. Analogous to the description of
$\cat{2Cob}$ as the free symmetric monoidal category on a
commutative Frobenius algebra, we will prove that $\cat{2Thick}$
is the free monoidal category on a noncommutative Frobenius
algebra.

\begin{defn}
The monoidal category of {\em two-dimensional thick tangles}
denoted $\cat{2Thick}$ has nonnegative integers as objects. The
1-morphisms from $k$ to $l$ are boundary preserving diffeomorphism
classes of smooth oriented compact surfaces $X$ with boundary
$\partial X$ equipped with disjoint distinguished intervals $i_j^s
\maps I \hookrightarrow
\partial X$, $1 \leq j \leq k$, $i^t_m \maps \hookrightarrow
\partial X$, $1 \leq m \leq l$, equipped with a smooth embedding $d
\maps X \hookrightarrow \R \times [0,1]$ such that
\[
 d^{-1}(\R\times0) = I_1^s \sqcup I^s_2 \sqcup \cdots \sqcup
 I_k^s, \qquad   I^s_j=i^s_j(I) , \qquad d(I^s_j) = [j-\frac{1}{3},j+\frac{1}{3}] \times 0 ,\]
 \[
 d^{-1}(\R\times1) = I_1^t \sqcup I^t_2 \sqcup \cdots \sqcup I_k^t,
 \qquad  I^t_j=i^t_j(I) , \qquad  d(I^t_j) = [j-\frac{1}{3},j+\frac{1}{3}] \times 1.
\]
The image $d(X)$ is called a {\em diagram of two-dimensional thick
tangles}.

Composition $Y \circ X$ of 1-morphisms $\xymatrix@1{k \ar[r]^X & l
\ar[r]^Y & m}$ is defined by sewing of surfaces at boundary
intervals $I^t_j(X)$ and $I^s_j(Y)$.  The identity 1-morphism $1_k
\maps k \to k$ is the union $\coprod_{j=1}^k
[j-\frac{1}{3},j+\frac{1}{3}] \times [0,1]$.  The identity axioms
follow from the isomorphisms $\xymatrix@1{1_l \circ X
\ar[r]^{\sim} & X}$ obtained by taking a neighborhood $(U,I^s_j)
\simeq \left([0,1] \times[0,1],[0,1],[0,1] \times 0\right)$ of the
distinguished interval $I_j^s \subset X$ and by taking any
isomorphism $[0,1] \times [0,1]\bigcup_{[0,1]\times 1 \sim I^s_j}U
\simeq U$.  The tensor product is the disjoint union. The unit
object is 0.
\end{defn}

This monoidal category is actually a decategorified version of the
category of `planar thick tangles' defined by Kerler and
Lyubashenko~\cite{KL}.  Some examples of two-dimensional thick
tangles are shown below\footnote{Note that here we are depicting
the surfaces in a slightly artistic fashion.}:
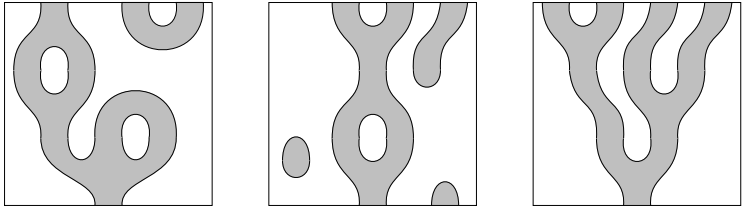
\[ \psset{xunit=.60cm,yunit=0.60cm}
     \begin{pspicture}[.5](4.1,4.8)
    \rput(2.6,1.5){\zagQ}
    \rput(2,0){\forkQ}
    \rput(.8,1.5){\multQ}
    \rput(.8,3){\comultQ}
    \rput(3.2,4.5){\zig}
    \pspolygon(-.3,0)(-.3,4.5)(4.3,4.5)(4.3,0)(-.3,0)
\end{pspicture}
\qquad \quad
     \begin{pspicture}[0.5](4.1,4.8)
  \rput(2,0){\multQ}
  \rput(2,1.5){\comultQ}
  \rput(2,3){\multQ}
  \rput(3.2,3){\curverightQ}
  \rput(3.2,3){\deathQ}
  \rput(3.6,0){\birthQ}
  \rput(.3,1){\deathQ}
  \rput(.3,1){\birthQ}
    \pspolygon(-.3,0)(-.3,4.5)(4.3,4.5)(4.3,0)(-.3,0)
\end{pspicture}
\qquad \quad
     \begin{pspicture}[0.5](4.1,4.8)
  \rput(2,0){\multQ}
  \rput(1.4,1.5){\curveleftQ}
  \rput(2.6,1.5){\multQ}
  \rput(.8,3){\multQ}
  \rput(2,3){\curverightQ}
  \rput(3.2,3){\curverightQ}
    \pspolygon(-.3,0)(-.3,4.5)(4.3,4.5)(4.3,0)(-.3,0)
\end{pspicture}
\]
We mentioned in the introduction that $\cat{2Thick}$ has the
alternative description as the category whose objects are open
strings and whose morphisms are diffeomorphism classes of planar
open string worldsheets.

We now state one of the main theorems of this paper. The proof
follows as a corollary of a result proven in Section
\ref{secPROOFS}.  We will however sketch a proof for the eager
reader.

\begin{thm} \label{decatmainthm}
The category of two-dimensional thick tangles is equivalent to the
free monoidal category on a Frobenius object.
\end{thm}

\noindent {\bf Sketch of Proof. } The height function ${\rm pr}_2
\circ d \maps \partial X \to [ 0,1]$ defines a Morse function on
the surface $X$.  This allows a decomposition of $X$ into the
elementary building blocks defining the data of a Frobenius
algebra. The equality of two surfaces related by the sliding of
handles implies the associativity of the multiplication, the
coassociativity of the comultiplication, and the Frobenius
identities.  The equality of those surfaces obtained by the
cancellation of local maxima with a local minima imply the unit
and counit axioms.  One can check that these are the only allowed
Morse moves for these surfaces so that $\cat{2Thick}$ is freely
generated by a Frobenius algebra. \qed

This result therefore justifies the assignation of a
noncommutative Frobenius algebra to the image of the generators
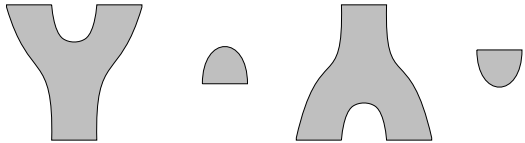
\[
\twothick
\]
in the Moore-Segal axioms of an open-closed topological field
theory.  The idea that the free monoidal category on a
noncommutative Frobenius algebra might be related to a topological
category of this sort has been suggested by Baez~\cite{ThisWeek}
and Kock~\cite{Kock}.

\begin{cor}
A Frobenius algebra in a monoidal category determines an invariant
of 2-dimensional thick tangles.
\end{cor}

\Proof This follows from the universal property of the monoidal
category of thick tangles. \qed

We can define a planar open string topological field theory as a
monoidal functor from $\cat{2Thick}$ into a monoidal category
$\mathcal{C}$.  We then have the following as a simple restatement
of Theorem~\ref{decatmainthm}:

\begin{thm}
A planar open string topological field theory is equivalent to a
(not necessarily commutative) Frobenius algebra in the monoidal
category $\mathcal{C}$.
\end{thm}

\section{Categorification} \label{secCat}

We may now reap the rewards of our characterization of Frobenius
algebras in terms of the walking adjunction.  We will see that, by
assuming a similar relationship holds in higher dimensions, we can
use the existing definition of a categorified adjunction to define
a categorified Frobenius algebra.  We will then show that this
definition of categorified Frobenius algebra is equivalent to what
one might expect by replacing the equations with coherent
isomorphisms in Proposition \ref{equivFrob}.  The advantage of
this approach is that the usual hassle of figuring out coherence
conditions for these isomorphisms can be avoided by utilizing what
is known about categorified adjunctions.

We start off in Section \ref{secPrelim} by defining some
categorical notions that we will need later on.  We then describe
the generalization of string diagrams to the context of
3-categories in Section \ref{secHstring}.  Then we define the
walking pseudo monoid and the walking pseudo Frobenius algebra in
terms of pseudoadjunctions in Sections \ref{secWWA} and
\ref{secWWAA}.  In Section \ref{secExamples} we provide some
examples of pseudo Frobenius algebras before moving on to prove
the main theorem that the walking pseudo Frobenius algebra is
triequivalent to the monoidal bicategory of 3-dimensional thick
tangles.

\subsection{Preliminaries for higher categories} \label{secPrelim}

Adjunctions make sense in a 2-category and even in the more
general context of bicategories or weak 2-categories.  To
categorify the notion of an adjunction we will need to climb one
categorical dimension to the level of 3-categories. Intuitively,
the idea of a 3-category is not hard to grasp. Categories have
objects, morphisms and some axioms; 2-categories have objects,
morphisms, 2-morphisms and some axioms; and 3-categories have
objects, morphisms, 2-morphisms, \textit{and} 3-morphisms together
with some axioms.

Since 3-categories have an extra level of structure, it is
possible for the composition and unit axioms for both 1-morphisms
and 2-morphisms to hold only up to coherent isomorphism.
Similarly, we can require that the the interchange law for
2-morphisms is only satisfied up to isomorphism as well. When all
levels of composites and identities up the level of 3-morphisms
satisfy the usual axioms up to isomorphism, we call this type of
3-category a weak 3-category or tricategory. If the axioms hold as
equations at each level rather than up to isomorphism, then we
call this type of 3-category a strict 3-category.  The most
exciting thing about 3-categories is that this is the first level
in which `weak' structures are no longer equivalent in some sense
to the strict version.  Thus, every tricategory \textit{is not}
triequivalent to a strict 3-category. Triequivalence is the
weakest notion of equivalence that can be defined for
3-categories.  This notion is the one most naturally suited for
relating tricategories.  The standard reference for tricategories
and triequivalences is Gordon, Power and Street~\cite{GPS}.

Although every tricategory is not triequivalent to a strict
3-category, there is a notion of 3-category which every
tricategory is triequivalent to --- a semistrict 3-category.
Semistrict 3-categories are a hybrid notion between strict
3-categories and tricategories.  They represent the strictest
known class of 3-categories that remain triequivalent to
tricategories. Hence,  it is this notion, that of a semistrict
3-category, that will provide a sufficiently general context to
consider categorified adjunctions or pseudoadjunctions.

A semistrict 3-category is defined using enriched category theory
\cite{kel} as a category enriched in $\cat{Gray}$~\cite{GPS}.  For
this reason, semistrict 3-categories are sometimes referred to as
\cat{Gray}-categories.  We will use both terms interchangeably.
$\cat{Gray}$ is the symmetric monoidal closed category whose
underlying category is \cat{2-Cat}; the category whose objects are
2-categories, and whose morphisms are 2-functors. $\cat{2-Cat}$
has a natural monoidal structure given by the cartesian product.
However, enriching over $(\cat{2-Cat}, \times)$ only produces
\textit{strict} 3-categories~\cite{EK}.

Enriching in \cat{Gray} produces a more interesting notion of
3-category because \cat{Gray} has a more interesting monoidal
structure, namely the `Gray'-tensor product. This has the effect
of equipping a \cat{Gray}-category $\mathcal{K}$ with a cubical
functor $M \maps \mathcal{K}(A,B) \times \mathcal{K}(B,C) \to
\mathcal{K}(A,C)$ for all objects $A$,$B$,$C$ in $\mathcal{K}$.
This means that if $f \maps F \To F'$ in $\mathcal{K}(A,B)$, and
$g \maps G \To G'$ in $\mathcal{K}(B,C)$, then, rather than
commuting on the nose, we have an invertible 3-cell in the
following square:
    \[
 \xy
    (-10,8)*+{FG}="TL";
    (10,8)*+{FG'}="TR";
    (10,-8)*+{F'G'}="BR";
    (-10,-8)*+{F'G}="BL";
    {\ar@{=>}^-{Fg} "TL";"TR"};
    {\ar@{=>}^-{fG'} "TR";"BR"};
    {\ar@{=>}_-{fG} "TL";"BL"};
    {\ar@{=>}_-{F'g} "BL";"BR"};
    (0,3)*{}="x";
    (0,-3)*{}="y";
    {\ar@3{->}^{f_g} "x";"y"};
 \endxy
\]
where we write the 3-morphism $M_{f,g}$ as $f_g$ following
Marmolejo~\cite{mar}.  A bit more concretely, in a
\cat{Gray}-category composites are still strictly associative, and
identities behave as identities on the nose, but the interchange
law now holds only up to coherent isomorphism.

\begin{defn}
A semistrict 3-category is a category enriched in \cat{Gray}.
\end{defn}

Using enriched category theory it is also possible to define the
morphisms between two semistrict 3-categories.  A semistrict
3-functor between two semistrict 3-categories is just a
\cat{Gray}-enriched functor, or a \cat{Gray}-functor as it is
often called.  We will also use these two terminologies
interchangeably.

\begin{defn}
A {\em pseudoadjunction} $i,e,I,E \maps L \dashv_p R\maps B \to A$
in a $\cat{Gray}$-category $\mathcal{K}$ consists of:
\begin{itemize}
   \item morphisms $R \maps A \to B$ and $L \maps B \to A$,
   \item 2-morphisms $i \maps 1 \To LR$ and $e \maps RL \To 1$,  and
   \item coherence 3-isomorphisms
   \[ \vcenter{
 \xy
   (-15,0)*+{R}="l";
   (0,15)*+{RLR}="t";
   (15,0)*+{R}="r";
    {\ar@{=>}^{Ri} "l";"t"};
    {\ar@{=>}^{eR} "t";"r"};
    {\ar@{=>}_1 "l";"r"};
        {\ar@3{->}^{I} (0,3);"t"+(0,-6)};
 \endxy}
\qquad {\rm and} \qquad
 \vcenter{\xy
   (-15,0)*+{L}="l";
   (0,15)*+{LRL}="t";
   (15,0)*+{L}="r";
    {\ar@{=>}^{iL} "l";"t"};
    {\ar@{=>}^{Le} "t";"r"};
    {\ar@{=>}_1 "l";"r"};
    {\ar@3{->}_{E} "t"+(0,-6);(0,3)};
 \endxy}
\]
\end{itemize}
such that the following two diagrams are both identities:
\[
 \xy
   (-38,0)*+{RL}="ll";
   (0,14)*+{RL}="t";
   (-14,0)*+{RLRL}="l";
   (14,0)*+{1}="r";
   (0,-14)*+{RL}="b";
        {\ar@{=>}_{eRL} "l";"t"};
        {\ar@{=>}^{RLe} "l";"b"};
        {\ar@{=>}^{e} "t";"r"};
        {\ar@{=>}_{e} "b";"r"};
        {\ar@{=>}^{RiL} "ll";"l"};
        {\ar@/^1pc/@{=>}^1 "ll";"t"};
        {\ar@/_1pc/@{=>}_1 "ll";"b"};
        {\ar@3{->}^{e_e^{-1}} "t"+(2,-11);"b"+(2,11)};
        {\ar@3{->}_{IL} (-16,9);(-16,4)};
        {\ar@3{->}_{RE} (-16,-4);(-16,-9)};
 \endxy
\qquad \qquad
 \xy
   (38,0)*+{LR}="ll";
   (0,14)*+{LR}="t";
   (-14,0)*+{1}="l";
   (14,0)*+{LRLR}="r";
   (0,-14)*+{LR}="b";
        {\ar@{=>}^{i} "l";"t"};
        {\ar@{=>}_{i} "l";"b"};
        {\ar@{=>}_{LRi} "t";"r"};
        {\ar@{=>}^{iLR} "b";"r"};
        {\ar@{=>}^{LeR} "r";"ll"};
        {\ar@/^1pc/@{=>}^1 "t";"ll"};
        {\ar@/_1pc/@{=>}_1 "b";"ll"};
        {\ar@3{->}_{i_i^{-1}} "t"+(-2,-11);"b"+(-2,11)};
        {\ar@3{->}^{LI} (16,9);(16,4)};
        {\ar@3{->}^{ER} (16,-4);(16,-9)};
 \endxy
\]
\end{defn}

Having the aim of defining a pseudo Frobenius algebra to be $LR
\in \Hom(A,A)$ in an ambidextrous pseudoadjunction, it is
important to understand the structure of $\Hom(A,A)$.  Just as
fixing one object in 2-category produces a monoidal category,
fixing one object in a semistrict 3-category produces a semistrict
monoidal 2-category or a \cat{Gray}-monoid as it is sometimes
referred.

\begin{defn}
A semistrict monoidal 2-category is a one object
$\cat{Gray}$-category.
\end{defn}

For later convenience we write this definition in more elementary
terms below.  We employ the convention that the 2-morphism
$M_{f,g}$ will be written as $f_g$ in the context of semistrict
3-categories, and will be written more indicatively as
$\INT_{f,g}$, or simply as $\INT$, in the context of semistrict
monoidal 2-categories.

\begin{prop} \et \label{MonBiCat}
A {\em semistrict monoidal 2-category} consists of a 2-category
$\mathcal{C}$ together with:
\begin{enumerate} \numarabic
 \item An object $I \in \mathcal{C}$.
 \item For any two objects $A,B \in \mathcal{C}$ an object $A \ten B$ in
$\mathcal{C}$.
 \item For any 1-morphism $f \maps A \to A'$ and any object $B \in
\mathcal{C}$ a 1-morphism $f \ten B \maps A \ten B \to A'\ten B$.
 \item For any 1-morphism $g \maps B \to B'$ and any object $A \in
\mathcal{C}$ a 1-morphism $A \ten g \maps A \ten B \to A \ten B'$.
 \item For any object $B \in \mathcal{C}$ and any 2-morphism $\alpha \maps f \to
f'$ a 2-morphism $\alpha \ten B \maps f \ten B \To f' \ten B$.
 \item For any object $A \in \mathcal{C}$ and any 2-morphism $\beta \maps g \To
g'$ a 2-morphism $A \ten \beta \maps A \ten g \To A \ten g'$.
 \item For any two 1-morphisms $f \maps A \to A'$ and $g \maps B \to
B'$ a 2-isomorphism

\[
 \xy
    (-14,10)*+{A \ten B}="TL";
    (14,10)*+{A \ten B'}="TR";
    (14,-10)*+{A' \ten B'}="BR";
    (-14,-10)*+{A' \ten B}="BL";
    {\ar^{ A \ten g} "TL";"TR"};
    {\ar^{ f \ten B'} "TR";"BR"};
    {\ar_{ f \ten B} "TL";"BL"};
    {\ar_{ A' \ten g} "BL";"BR"};
    (0,4)*{}="x";
    (0,-4)*{}="y";
    {\ar@{=>}^{\INT_{f,g}} "x";"y"};
 \endxy
\]
\end{enumerate}
This data is subject to the following conditions.
\begin{enumerate} \numroman
\item For any object $A \in \mathcal{C}$ we have $A \ten - \maps
\mathcal{C} \to \mathcal{C}$ and $- \ten A \maps \mathcal{C} \to
\mathcal{C}$ are 2-functors.

\item For $x$ any object, morphism or 2-morphisms of $\mathcal{C}$
we have $x \ten I = I \ten x = x$.

\item For x any object, morphism or 2-morphism of $\mathcal{C}$,
and for all objects $A,B \in \mathcal{C}$ we have $A \ten (B \ten
x) = (A \ten B) \ten x$, $A \ten (x \ten B) = (A \ten x) \ten B$
and $x \ten (A \ten B) = (x \ten A) \ten B$.

\item For any 1-morphism $f \maps A \to B$, $g \maps B \to B'$ and
$h \maps C \to C'$ in $\mathcal{C}$ we have $\INT_{A \ten g,h} = A
\ten \INT_{g,h}$, $\INT_{f \ten B,h} = \INT_{f, B \ten h}$ and
$\INT_{f, g \ten C} = \INT_{f,g} \ten C$.

\item For any objects $A,B \in \mathcal{C}$ we have $1_A \ten B =
A \ten 1_B = 1_{A \ten B}$, and for any 1-morphism $f \maps A \to
B$, $g \maps B \to B'$ in $\mathcal{C}$ we have $\INT_{1_A,g} =
1_{A \ten g}$ and $\INT_{f,1_B} = 1_{f \ten B}$.

\item For any 1-morphisms $f,h \maps A \to A'$, $g,k \maps B \to
B'$, and any 2-morphisms $\alpha \maps f \To h$, and $\beta \maps
g \To k$,
\[
  \vcenter{\xy (-10,-30)*{}; (30,10)*{};
 (0,0)*++{
 \xymatrix{
 A \ten B \ar[rr]^{f \ten B}
    \ar@/^1.2pc/[dd]^<<<<<{A \ten g}_{}="0"
    \ar@/_1.2pc/[dd]_<<<<<{A \ten k}^{}="1"
    \ar@{=>}"0"; "1" ^{A \ten \beta}
&& A' \ten  B \ar[dd]^{A' \ten g} \\
\\
A \ten B'
    \ar@/^1.2pc/[rr]^>>>>>>{f \ten B'}_{\quad}="0"
    \ar@/_1.2pc/[rr]_{h \ten B'}^{}="1"
    \ar@{=>}"0"; "1" ^>>>>>>{\alpha \ten B'}
&& A' \ten B' }   };
 {\ar@{=>}^{\INT^{-1}_{f,g}} (18,-4);(18,-10)};
\endxy }
 \qquad \qquad \quad = \qquad
  \vcenter{\xy (-8,-30)*{}; (30,10)*{};
 (0,0)*++{
 \xymatrix{
 A \ten B
    \ar[dd]_{A \ten k}
    \ar@/^1.2pc/[rr]^<<<<<{f \ten B}_{}="0"
    \ar@/_1.2pc/[rr]_<<<<<{h \ten B}^{}="1"
    \ar@{=>}"0"; "1" ^{\alpha  \ten B}
&& A' \ten  B
    \ar@/^1.2pc/[dd]^>>>>>>{A'  \ten g}_{}="0"
    \ar@/_1.2pc/[dd]_>>>>>>{A' \ten k}^{}="1"
    \ar@{=>}"0"; "1" ^{A' \ten \beta} \\
\\
A \ten B'
    \ar[rr]_{h \ten B'}^{}
&& A' \ten B' } };
 {\ar@{=>}_{\INT_{h,k}^{-1}} (12,-14);(12,-20)};
\endxy }
\]
 \item For any 1-morphisms $f \maps A \to A'$, $g \maps B \to B'$,
$f' \maps A' \to A''$, $g' \maps B' \to B''$,
\[
  \semistrictbig
\]
\end{enumerate}
\end{prop}

\Proof This is a straightforward verification. \qed

Just as monoids can be defined in any monoidal category,
pseudomonoids can be defined in any semistrict monoidal
2-category.  A pseudomonoid is just a monoid where the equations
describing the associativity and unit constraints are replaced by
coherent isomorphism~\cite{DS}.

\begin{defn}
A pseudomonoid $M$ in the semistrict monoidal 2-category
$\mathcal{C}$ consists of:
\begin{itemize}
    \item An object $M$ of $\mathcal{C}$.
    \item A {\em multiplication} morphism: $m \maps M \ten M \to
    M$.
    \item A {\em unit} for multiplication: $\iota \maps 1 \to M$, and
    \item coherence 2-isomorphisms
\[
 \xy
   (-20,10)*+{M}="L";
   (0,10)*+{M \ten M}="M";
   (20,10)*+{M}="R";
   (0,-10)*+{M}="B";
    {\ar^-{M \ten \iota} "L";"M"};
    {\ar_-{\iota \ten M} "R";"M"};
    {\ar "L";"B"};
    {\ar "R";"B"};
    {\ar^-{m} "M";"B"};
    {\ar@{=>}_{r} "M"+(-4,-4);(-8,2)  };
    {\ar@{=>}^{\ell} "M"+(4,-4) ;(8,2) };
 \endxy
\qquad {\rm and} \qquad
 \xy
   (0,11)*+{M \ten M \ten M}="T";
   (-12,0)*+{M \ten M}="L";
   (12,0)*+{M \ten M}="R";
   (0,-11)*+{M}="B";
   {\ar_{M\ten m} "T";"L"};
   {\ar^{m \ten M} "T";"R"};
   {\ar_{m} "L";"B"};
   {\ar^{m} "R";"B"};
   {\ar@{=>}^{a} "L"+(8,0); "R"+(-8,0) };
 \endxy
\]
\end{itemize}
such that the following two equations are satisfied:
\[
 \vcenter{\xy
   (-24,22)*+{M^4}="tl";
   (-1,22)*+{M^3}="tr";
   (-10,10)*+{M^3}="ml";
   (13,10)*+{M^2}="mr";
   (-24,2)*+{M^3}="o";
   (-10,-10)*+{M^2}="bl";
   (13,-10)*+{M}="br";
        {\ar^{m \ten M^2} "tl";"tr"};
        {\ar_{M^2 \ten m} "tl";"o"};
        {\ar^-{M \ten m \ten M} "tl";"ml"};
        {\ar_{m \ten M} "ml";"mr"};
        {\ar_{M \ten m} "ml";"bl"};
        {\ar^{ m \ten M} "tr";"mr"};
        {\ar_{M \ten m} "o";"bl"};
        {\ar^{m} "mr";"br"};
        {\ar_{m} "bl";"br"};
          {\ar@{=>}_{a M} "tr"+(1,-4); "ml"+(6,2); };
          {\ar@{=>}^{M a } "ml"+(-4.5,-2.5); "o"+(4.5,2.5); };
          {\ar@{=>}^{a} "mr"+(-8,-8); "bl"+(8,8); };
 \endxy}
\qquad = \qquad \vcenter{ \xy
   (10,-10)*+{M}="tl";
   (-13,-10)*+{M^2}="tr";
   (-4,2)*+{M^2}="ml";
   (-27,2)*+{M^3}="mr";
   (10,10)*+{M^2}="o";
   (-4,22)*+{M^3}="bl";
   (-27,22)*+{M^4}="br";
        {\ar_{m} "tl";"tr";};
        {\ar^{m} "tl";"o";};
        {\ar_-{m} "tl";"ml";};
        {\ar^{m \ten M} "ml";"mr";};
        {\ar^{M \ten m} "ml";"bl";};
        {\ar_{M \ten m} "tr";"mr";};
        {\ar^{m \ten M} "o";"bl";};
        {\ar_{ M^2 \ten m} "mr";"br";};
        {\ar^{m \ten M^2} "bl";"br";};
          {\ar@{=>}^{a} "ml"+(-2,-4);  "tr"+(2,4);};
          {\ar@{=>}_{a} "o"+(-4.5,-2.5);  "ml"+(4.5,2.5);};
          {\ar@{=>}^{\INT_{m,m}} "bl"+(-8,-8);  "mr"+(8,8);};
 \endxy}
\]

\[
 \xy
   (-34,0)*+{M^2}="o";
   (0,12)*+{M^2}="t";
   (-12,0)*+{M^3}="l";
   (12,0)*+{M}="r";
   (0,-14)*+{M^2}="b";
        {\ar@/_1pc/_{1} "o";"b" };
        {\ar^{M \ten \iota \ten M} "o";"l" };
        {\ar^{m \ten M} "l";"t" };
        {\ar^{m} "t";"r" };
        {\ar_{m} "b";"r" };
        {\ar_>>>>{M \ten m} "l";"b" };
        {\ar@{=>}_{a^{-1}} "t"+(0,-7);  "b"+(0,7)};
          {\ar@{=>}_<<<{M \ten \ell^{-1} } "l"+(-2,-4);  (-18,-9)};
 \endxy
\qquad = \qquad
 \xy
   (30,0)*+{M}="o";
   (0,12)*+{M^3}="t";
   (-12,0)*+{M^2}="l";
   (12,0)*+{M^2}="r";
        {\ar^{m} "r";"o" };
        {\ar^{M \ten \iota \ten M } "l";"t" };
        {\ar "l";"r" };
        {\ar^{m \ten M} "t";"r" };
         {\ar@{=>}_{ r \ten M} "t"+(2,-5);  "t"+(2,-9)};
 \endxy
\]
\end{defn}

As an example, note that a weak monoidal category is a
pseudomonoid in $\cat{Cat}$. We now move on to describe how some
of these definitions can be understood diagrammatically.

\subsection{String diagrams for higher categories} \label{secHstring}

In this section we will sketch a version of string diagrams for
3-categories.  Unfortunately we will only scratch the surface of
this beautiful subject. We will focus on those aspects needed to
see the relationship between pseudo Frobenius algebras and
3-dimensional thick tangles.  The idea of string diagrams for
3-categories is given in~\cite{Street}.  Carter and Saito also use
3-categorical string diagrams implicitly in their work and explain
its relationship to singularity theory~\cite{CS}.

In the notation that we will use, objects, morphisms and
2-morphisms will be depicted in exactly the same way as before.
What's new is the 3-morphisms which we will simply draw as arrows
\textit{going between string diagrams}.  For example, the morphism
$I \maps 1_R \Rrightarrow Ri.eR$ in the definition of a
pseudoadjunction is depicted as:
\[
\psset{xunit=.70cm,yunit=0.70cm} \xy
      (0,0)*+{\begin{pspicture}(2,2.2)
\pscustom[fillcolor=lightgray, fillstyle=solid]{
  \pspolygon(1,0)(1,2)(0,2)(0,0)(1,0)
  }
   \rput(.5,1){$B$}
         \rput(1.5,1){$A$}
         \rput(1,2.2){$\scs R$}
         \rput(1,-.2){$\scs R$}
  \pspolygon(0,0)(0,2)(2,2)(2,0)(0,0)
 \end{pspicture}};
    \endxy
 \quad \xy {\ar@3{->}^{I} (0,0)*{};(8,0)*{}}; \endxy  \quad
     \xy
    (0,0)*+{   \begin{pspicture}(2,2.2)
  \pspolygon(0,0)(0,2)(2,2)(2,0)(0,0)
\pscustom[fillcolor=lightgray, fillstyle=solid]{
  \psbezier(.5,1)(.5,.35)(1,.35)(1,1)
  \psbezier(1,1)(1,1.65)(1.5,1.65)(1.5,1)
  \psline(1.5,0)
  \psline(0,0)
  \psline(0,2)
  \psline(.5,2)
  \psline(.5,1)
  }
          \rput(1.65,1.65){$A$}
         \rput(.35,.35){$B$}
         \rput(.5,2.2){$\scs R$}
         \rput(1.5,-.2){$\scs R$}
 \end{pspicture}};
  \endxy
\]
Notice that the source and target of 3-morphisms in the diagram
are the 2-morphism represented by the string diagrams.  The
3-morphism $E \maps iL.Le \Rrightarrow 1_L$ is depicted similarly
as:
\[
\psset{xunit=.70cm,yunit=0.70cm} \xy
   (0,0)*+{\begin{pspicture}(2,2.2)
  \pspolygon(0,0)(0,2)(2,2)(2,0)(0,0)
\pscustom[fillcolor=lightgray, fillstyle=solid]{
  \psbezier(.5,1)(.5,1.65)(1,1.65)(1,1)
  \psbezier(1,1)(1,.35)(1.5,.35)(1.5,1)
  \psline(1.5,2)
  \psline(2,2)
  \psline(2,0)
  \psline(.5,0)
  \psline(.5,1)
  }
         \rput(.35,1.65){$A$}
         \rput(1.65,.35){$B$}
         \rput(1.5,2.2){$\scs L$}
         \rput(.5,-.2){$\scs L$}
 \end{pspicture}};
 \endxy
 \quad \xy {\ar@3{->}^{E} (0,0)*{};(8,0)*{}}; \endxy  \quad
  \xy
   (0,0)*+{  \begin{pspicture}(2,2.2)
  \pspolygon(0,0)(0,2)(2,2)(2,0)(0,0)
\pscustom[fillcolor=lightgray, fillstyle=solid]{
  \pspolygon(1,0)(1,2)(2,2)(2,0)(1,0)
  }
  \rput(.5,1){$A$}
         \rput(1.5,1){$B$}
         \rput(1,2.2){$\scs L$}
         \rput(1,-.2){$\scs L$}
 \end{pspicture}};
 \endxy
\]

From now on we will omit the labels on the string diagrams
representing the 2-morphisms.  As usual, the gray colored region
is meant to represent the object $B$ and the white area the object
$A$.  All the other labels can be deduce from the diagrams. Using
this notation we can even depict the pseudoadjunction axioms as
\textit{commutative diagrams of string diagrams}.  For example,
the assertion that
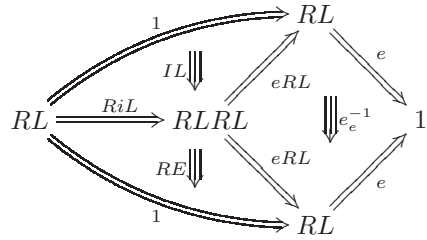
\[
 \xy
   (-38,0)*+{RL}="ll";
   (0,14)*+{RL}="t";
   (-14,0)*+{RLRL}="l";
   (14,0)*+{1}="r";
   (0,-14)*+{RL}="b";
        {\ar@{=>}_{eRL} "l";"t"};
        {\ar@{=>}^{eRL} "l";"b"};
        {\ar@{=>}^{e} "t";"r"};
        {\ar@{=>}_{e} "b";"r"};
        {\ar@{=>}^{RiL} "ll";"l"};
        {\ar@/^1pc/@{=>}^1 "ll";"t"};
        {\ar@/_1pc/@{=>}_1 "ll";"b"};
        {\ar@3{->}^{e_e^{-1}} "t"+(2,-11);"b"+(2,11)};
        {\ar@3{->}_{IL} (-16,9);(-16,4)};
        {\ar@3{->}_{RE} (-16,-4);(-16,-9)};
 \endxy
\]
is equal to the identity can be translated into diagrams. This
axiom says that the composite:
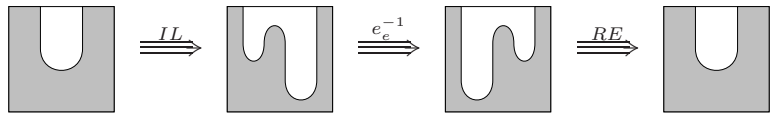
\[
\psset{xunit=.7cm,yunit=.7cm}
\begin{pspicture}[.5](2,2.2)
  \pspolygon(0,0)(0,2)(2,2)(2,0)(0,0)
\pscustom[fillcolor=lightgray, fillstyle=solid]{
  \psbezier(.6,1.2)(.6,.65)(1.4,.65)(1.4,1.2)
  \psline(1.4,2)\psline(2,2) \psline(2,0) \psline(0,0) \psline(0,2)
  \psline(.6,2) \psline(.6,1.2)
  }
 \end{pspicture}
 \quad \xy {\ar@3{->}^{IL} (0,0)*{};(8,0)*{}}; \endxy  \quad
\begin{pspicture}[.5](2,2.2)
  \pspolygon(0,0)(0,2)(2,2)(2,0)(0,0)
\pscustom[fillcolor=lightgray, fillstyle=solid]{
  \psbezier(.3,1.3)(.3,.85)(.7,.85)(.7,1.3)
  \psbezier(.7,1.3)(.7,1.75)(1.1,1.75)(1.1,1.3)
  \psline(1.1,.6)
  \psbezier(1.1,.6)(1.1,.1)(1.7,.1)(1.7,.6)
  \psline(1.7,2)
  \psline(2,2)
  \psline(2,0)
  \psline(0,0)
  \psline(0,2)
  \psline(.3,2)
  \psline(.3,1.3)
  }
 \end{pspicture}
 \quad \xy {\ar@3{->}^{e^{-1}_e} (0,0)*{};(8,0)*{}}; \endxy  \quad
\begin{pspicture}[.5](2,2.2)
  \pspolygon(0,0)(0,2)(2,2)(2,0)(0,0)
\pscustom[fillcolor=lightgray, fillstyle=solid]{
  \psbezier(1.7,1.3)(1.7,.85)(1.3,.85)(1.3,1.3)
  \psbezier(1.3,1.3)(1.3,1.75)(.9,1.75)(.9,1.3)
  \psline(.9,.6)
  \psbezier(.9,.6)(.9,.1)(.3,.1)(.3,.6)
  \psline(.3,2)
  \psline(0,2)
  \psline(0,0)
  \psline(2,0)
  \psline(2,2)
  \psline(1.7,2)
  \psline(1.7,1.3)
  }
 \end{pspicture}
 \quad \xy {\ar@3{->}^{RE} (0,0)*{};(8,0)*{}}; \endxy  \quad
 \begin{pspicture}[.5](2,2.2)
  \pspolygon(0,0)(0,2)(2,2)(2,0)(0,0)
\pscustom[fillcolor=lightgray, fillstyle=solid]{
  \psbezier(.6,1.2)(.6,.65)(1.4,.65)(1.4,1.2)
  \psline(1.4,2)\psline(2,2) \psline(2,0) \psline(0,0) \psline(0,2)
  \psline(.6,2) \psline(.6,1.2)
  }
 \end{pspicture}
\]
is equal to the composite:
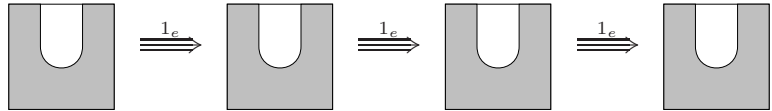
\[
\psset{xunit=.7cm,yunit=.7cm}
\begin{pspicture}[.5](2,2.2)
  \pspolygon(0,0)(0,2)(2,2)(2,0)(0,0)
\pscustom[fillcolor=lightgray, fillstyle=solid]{
  \psbezier(.6,1.2)(.6,.65)(1.4,.65)(1.4,1.2)
  \psline(1.4,2)\psline(2,2) \psline(2,0) \psline(0,0) \psline(0,2)
  \psline(.6,2) \psline(.6,1.2)
  }
 \end{pspicture}
 \quad \xy {\ar@3{->}^{1_e} (0,0)*{};(8,0)*{}}; \endxy  \quad
 \begin{pspicture}[.5](2,2.2)
  \pspolygon(0,0)(0,2)(2,2)(2,0)(0,0)
\pscustom[fillcolor=lightgray, fillstyle=solid]{
  \psbezier(.6,1.2)(.6,.65)(1.4,.65)(1.4,1.2)
  \psline(1.4,2)\psline(2,2) \psline(2,0) \psline(0,0) \psline(0,2)
  \psline(.6,2) \psline(.6,1.2)
  }
 \end{pspicture}
 \quad \xy {\ar@3{->}^{1_e} (0,0)*{};(8,0)*{}}; \endxy  \quad
 \begin{pspicture}[.5](2,2.2)
  \pspolygon(0,0)(0,2)(2,2)(2,0)(0,0)
\pscustom[fillcolor=lightgray, fillstyle=solid]{
  \psbezier(.6,1.2)(.6,.65)(1.4,.65)(1.4,1.2)
  \psline(1.4,2)\psline(2,2) \psline(2,0) \psline(0,0) \psline(0,2)
  \psline(.6,2) \psline(.6,1.2)
  }
 \end{pspicture}
 \quad \xy {\ar@3{->}^{1_e} (0,0)*{};(8,0)*{}}; \endxy  \quad
 \begin{pspicture}[.5](2,2.2)
  \pspolygon(0,0)(0,2)(2,2)(2,0)(0,0)
\pscustom[fillcolor=lightgray, fillstyle=solid]{
  \psbezier(.6,1.2)(.6,.65)(1.4,.65)(1.4,1.2)
  \psline(1.4,2)\psline(2,2) \psline(2,0) \psline(0,0) \psline(0,2)
  \psline(.6,2) \psline(.6,1.2)
  }
 \end{pspicture}
\]
Similarly, the other pseudoadjunction axiom:
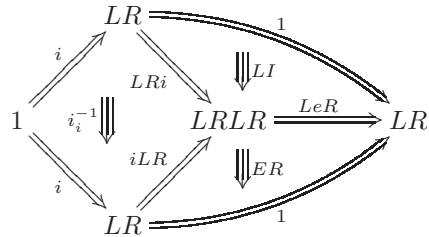
\[
 \xy
   (38,0)*+{LR}="ll";
   (0,14)*+{LR}="t";
   (-14,0)*+{1}="l";
   (14,0)*+{LRLR}="r";
   (0,-14)*+{LR}="b";
        {\ar@{=>}^{i} "l";"t"};
        {\ar@{=>}_{i} "l";"b"};
        {\ar@{=>}_{LRi} "t";"r"};
        {\ar@{=>}^{iLR} "b";"r"};
        {\ar@{=>}^{LeR} "r";"ll"};
        {\ar@/^1pc/@{=>}^1 "t";"ll"};
        {\ar@/_1pc/@{=>}_1 "b";"ll"};
        {\ar@3{->}_{i_i^{-1}} "t"+(-2,-11);"b"+(-2,11)};
        {\ar@3{->}^{LI} (16,9);(16,4)};
        {\ar@3{->}^{ER} (16,-4);(16,-9)};
 \endxy
 \]
 just says that the composite:
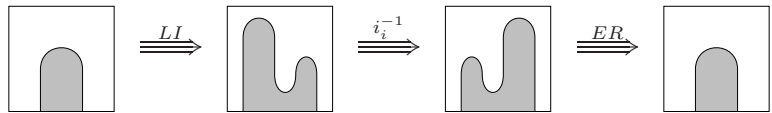
\[
 \begin{pspicture}[.5](2,2.2)
  \pspolygon(0,0)(0,2)(2,2)(2,0)(0,0)
\pscustom[fillcolor=lightgray, fillstyle=solid]{
  \psbezier(.6,.8)(.6,1.35)(1.4,1.35)(1.4,.8)
  \psline(1.4,0)\psline(.6,0) \psline(.6,.8)
  }
 \end{pspicture}
 \quad \xy {\ar@3{->}^{LI} (0,0)*{};(8,0)*{}}; \endxy  \quad
\begin{pspicture}[.5](2,2.2)
  \pspolygon(0,0)(0,2)(2,2)(2,0)(0,0)
\pscustom[fillcolor=lightgray, fillstyle=solid]{
  \psbezier(1.7,.7)(1.7,1.15)(1.3,1.15)(1.3,.7)
  \psbezier(1.3,.7)(1.3,.25)(.9,.25)(.9,.7)
  \psline(.9,1.4)
  \psbezier(.9,1.4)(.9,1.9)(.3,1.9)(.3,1.4)
  \psline(.3,0)
  \psline(1.7,0)
  \psline(1.7,.7)
  }
 \end{pspicture}
  \quad \xy {\ar@3{->}^{i^{-1}_i} (0,0)*{};(8,0)*{}}; \endxy  \quad
\begin{pspicture}[.5](2,2.2)
  \pspolygon(0,0)(0,2)(2,2)(2,0)(0,0)
\pscustom[fillcolor=lightgray, fillstyle=solid]{
  \psbezier(.3,.7)(.3,1.15)(.7,1.15)(.7,.7)
  \psbezier(.7,.7)(.7,.25)(1.1,.25)(1.1,.7)
  \psline(1.1,1.4)
  \psbezier(1.1,1.4)(1.1,1.9)(1.7,1.9)(1.7,1.4)
  \psline(1.7,0)
  \psline(.3,0)
  \psline(.3,.7)
  }
 \end{pspicture}
  \quad \xy {\ar@3{->}^{ER} (0,0)*{};(8,0)*{}}; \endxy  \quad
 \begin{pspicture}[.5](2,2.2)
  \pspolygon(0,0)(0,2)(2,2)(2,0)(0,0)
\pscustom[fillcolor=lightgray, fillstyle=solid]{
  \psbezier(.6,.8)(.6,1.35)(1.4,1.35)(1.4,.8)
  \psline(1.4,0)\psline(.6,0) \psline(.6,.8)
  }
 \end{pspicture}
 \]
is equal to the composite:
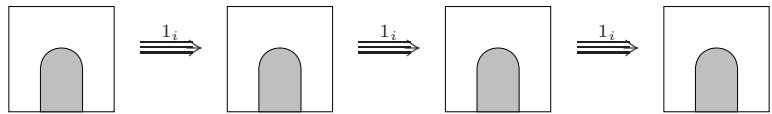
\[
 \begin{pspicture}[.5](2,2.2)
  \pspolygon(0,0)(0,2)(2,2)(2,0)(0,0)
\pscustom[fillcolor=lightgray, fillstyle=solid]{
  \psbezier(.6,.8)(.6,1.35)(1.4,1.35)(1.4,.8)
  \psline(1.4,0)\psline(.6,0) \psline(.6,.8)
  }
 \end{pspicture}
 \quad \xy {\ar@3{->}^{1_i} (0,0)*{};(8,0)*{}}; \endxy  \quad
 \begin{pspicture}[.5](2,2.2)
  \pspolygon(0,0)(0,2)(2,2)(2,0)(0,0)
\pscustom[fillcolor=lightgray, fillstyle=solid]{
  \psbezier(.6,.8)(.6,1.35)(1.4,1.35)(1.4,.8)
  \psline(1.4,0)\psline(.6,0) \psline(.6,.8)
  }
 \end{pspicture}
  \quad \xy {\ar@3{->}^{1_i} (0,0)*{};(8,0)*{}}; \endxy  \quad
 \begin{pspicture}[.5](2,2.2)
  \pspolygon(0,0)(0,2)(2,2)(2,0)(0,0)
\pscustom[fillcolor=lightgray, fillstyle=solid]{
  \psbezier(.6,.8)(.6,1.35)(1.4,1.35)(1.4,.8)
  \psline(1.4,0)\psline(.6,0) \psline(.6,.8)
  }
 \end{pspicture}
  \quad \xy {\ar@3{->}^{1_i} (0,0)*{};(8,0)*{}}; \endxy  \quad
 \begin{pspicture}[.5](2,2.2)
  \pspolygon(0,0)(0,2)(2,2)(2,0)(0,0)
\pscustom[fillcolor=lightgray, fillstyle=solid]{
  \psbezier(.6,.8)(.6,1.35)(1.4,1.35)(1.4,.8)
  \psline(1.4,0)\psline(.6,0) \psline(.6,.8)
  }
 \end{pspicture}
 \]
We sometimes refer to the coherence laws for a pseudoadjunction as
the \textit{triangulator identities}.

Notice that in the 2-categorical context we were able to distort
the diagrams for a given 2-morphisms and get a diagram
\textit{equal} to one that we started with, whereas in this
context the semistrictness makes these types of manipulations
appear as isomorphisms in our diagrams.  For instance the arrows
corresponding to $e^{-1}_e$ and $i^{-1}_i$ above.

It will be helpful to imagine the 3-morphisms as tracing out a
surface starting from the source and extending to the target.  In
this way the source and target are thought of as representing
different time slices of a Morse function on these surfaces. There
is much more to say about 3-categorical string diagrams but we
will stop here.

\subsection{The walking pseudoadjunction} \label{secWWA}

In this section we study a very special pseudoadjunction, the
walking pseudoadjunction.  The existence of the walking
pseudoadjunction was first proven by Lack~\cite{Lack}.  This means
that any pseudoadjunction in a \cat{Gray}-category $\mathcal{K}$
corresponds to a $\cat{Gray}$-functor from the walking
pseudoadjunction into $\mathcal{K}$.  Alternatively, the walking
pseudoadjunction can be thought of as the semistrict 3-category
freely generated by the data of a pseudoadjunction.  Analogous to
Section \ref{secWA} we will show that the semistrict monoidal
category $\Hom(A,A)$ is the walking pseudomonoid.

We begin by defining what it means to freely generate a semistrict
3-category.
\begin{defn} \label{Catgenerate}
Let $Y$ be a structure that can be defined in an arbitrary
semistrict 3-category. If $Y$ consists of objects
$Y_1,Y_2,\ldots,Y_{n}$, morphisms $F_1,F_2,\ldots,F_{n'}$,
2-morphisms $\alpha_1,\alpha_2,\ldots,\alpha_{n''}$, and
3-morphisms $\phi_1,\phi_2, \ldots, \phi_n'''$ for $n,n',n'',n'''
\in \Z^+$, then the semistrict 3-category $X$ is {\em generated}
by $Y$ if:
\begin{enumerate}[(i.)]
    \item Every object of $X$ is some $Y_i$.
    \item Every 1-morphism of $X$ can be obtained by compositions
    from the $F_i$'s and $1_{Y_i}$'s.
    \item Every 2-morphism of $X$ is obtained by taking Gray-tensor products and
    by vertical composition from:
    \begin{itemize}
        \item the 2-morphisms $\alpha_i$, and
        \item the identity 2-morphisms $1_F$
    for arbitrary 1-morphisms $F$.
    \end{itemize}
    \item Every 3-morphism is obtained from the Gray-tensor product of, and
    vertical and horizontal compositions from:
    \begin{itemize}
        \item the 3-morphisms $\phi_i$,
        \item the coherence morphisms
    $f_g$ for arbitrary 2-morphisms $f$ and $g$, and
        \item the identity 3-morphisms $1_F$ for arbitrary 2-morphisms $F$.
    \end{itemize}
\end{enumerate}
We say that $X$ is {\em freely generated} by $Y$ if the set of $Y$
objects in $\mathcal{K}$ are in bijection with
$\cat{Gray}$-functors from $X$ into $\mathcal{K}$, for every
semistrict 3-category $\mathcal{K}$.
\end{defn}

\begin{defn}
The walking pseudoadjunction $\wAdj$ is the semistrict 3-category
freely generated by:
\begin{itemize}
   \item morphisms $R \maps A \to B$ and $L \maps B \to A$,
   \item 2-morphisms $i \maps 1 \To LR$ and $e \maps RL \To 1$,  and
   \item coherence 3-isomorphisms
   \[ \vcenter{
 \xy
   (-15,0)*+{R}="l";
   (0,15)*+{RLR}="t";
   (15,0)*+{R}="r";
    {\ar@{=>}^{Ri} "l";"t"};
    {\ar@{=>}^{eR} "t";"r"};
    {\ar@{=>}_1 "l";"r"};
        {\ar@3{->}^{I} (0,3);"t"+(0,-6)};
 \endxy}
\qquad {\rm and} \qquad
 \vcenter{\xy
   (-15,0)*+{L}="l";
   (0,15)*+{LRL}="t";
   (15,0)*+{L}="r";
    {\ar@{=>}^{iL} "l";"t"};
    {\ar@{=>}^{Le} "t";"r"};
    {\ar@{=>}_1 "l";"r"};
    {\ar@3{->}_{E} "t"+(0,-6);(0,3)};
 \endxy}
\]
\end{itemize}
such that the following two diagrams are both identities:
\[
 \xy
   (-38,0)*+{RL}="ll";
   (0,14)*+{RL}="t";
   (-14,0)*+{RLRL}="l";
   (14,0)*+{1}="r";
   (0,-14)*+{RL}="b";
        {\ar@{=>}_{eRL} "l";"t"};
        {\ar@{=>}^{eRL} "l";"b"};
        {\ar@{=>}^{e} "t";"r"};
        {\ar@{=>}_{e} "b";"r"};
        {\ar@{=>}^{RiL} "ll";"l"};
        {\ar@/^1pc/@{=>}^1 "ll";"t"};
        {\ar@/_1pc/@{=>}_1 "ll";"b"};
        {\ar@3{->}^{e_e^{-1}} "t"+(2,-11);"b"+(2,11)};
        {\ar@3{->}_{IL} (-16,9);(-16,4)};
        {\ar@3{->}_{RE} (-16,-4);(-16,-9)};
 \endxy
\qquad \quad
 \xy
   (38,0)*+{LR}="ll";
   (0,14)*+{LR}="t";
   (-14,0)*+{1}="l";
   (14,0)*+{LRLR}="r";
   (0,-14)*+{LR}="b";
        {\ar@{=>}^{i} "l";"t"};
        {\ar@{=>}_{i} "l";"b"};
        {\ar@{=>}_{LRi} "t";"r"};
        {\ar@{=>}^{iLR} "b";"r"};
        {\ar@{=>}^{LeR} "r";"ll"};
        {\ar@/^1pc/@{=>}^1 "t";"ll"};
        {\ar@/_1pc/@{=>}_1 "b";"ll"};
        {\ar@3{->}_{i_i^{-1}} "t"+(-2,-11);"b"+(-2,11)};
        {\ar@3{->}^{LI} (16,9);(16,4)};
        {\ar@3{->}^{ER} (16,-4);(16,-9)};
 \endxy
\]
\end{defn}

We mentioned in Section \ref{secPrelim} that fixing an object of a
semistrict 3-category produces a semistrict monoidal 2-category.
Thus, using Definition \ref{Catgenerate} it makes sense to talk
about the semistrict monoidal 2-category freely generated by a
pseudomonoid.

\begin{defn}
The walking pseudomonoid is the semistrict monoidal 2-category
freely generated by the data defining a pseudomonoid.
\end{defn}

\noindent Lack has also explicitly construct the walking
pseudomonoid in the context of pseudomonads~\cite{Lack}.

To better understand the walking pseudomonoid and its relationship
to pseudoadjunctions it will be helpful to describe pseudomonoids
using an extension of our shorthand notation from Section
\ref{secFrob}.  In this notation, a pseudomonoid is an object of a
semistrict monoidal 2-category represented by the interval,
equipped with morphisms:
\[
 \begin{pspicture}[.5](1,1.7)
   \rput(.3,0){\mult}
 \end{pspicture}
 \qquad {\rm and} \qquad
 \begin{pspicture}[.5](1,1.7)
   \rput(.3,.75){\birth}
 \end{pspicture}
\]
and coherence 2-isomorphisms:
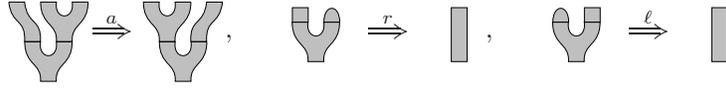
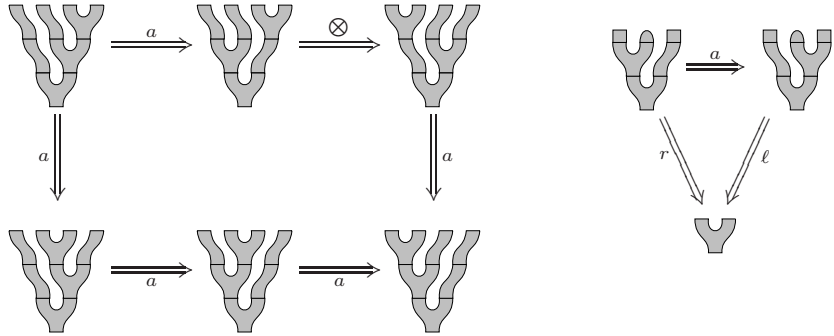
\[ 
\psset{xunit=.35cm,yunit=0.35cm}
\begin{pspicture}[0.5](2.4,3.3)
  \rput(1,0){\mult} \rput(1.6,1.5){\mult} \rput(.4,1.5){\curveleft}
\end{pspicture}
 \; \xy {\ar@2{->}^{a} (0,0)*{};(5,0)*{}}; \endxy  \quad
  \begin{pspicture}[0.5](2.6,3.3)
  \rput(1,0){\mult} \rput(1.6,1.5){\curveright} \rput(.4,1.5){\mult}
 \end{pspicture}
, \qquad  \psset{xunit=.35cm,yunit=0.35cm}
\begin{pspicture}[0.3](2,3)
  \rput(1,0){\mult} \rput(1.6,1.5){\birth} \rput(.4,1.5){\smallident}
\end{pspicture}
 \quad \xy {\ar@2{->}^{r} (0,0)*{};(5,0)*{}}; \endxy  \quad
\begin{pspicture}[0.3](2,3)
  \rput(1,0){\longident}
\end{pspicture},
\qquad  \psset{xunit=.35cm,yunit=0.35cm}
\begin{pspicture}[0.3](2,3)
  \rput(1,0){\mult} \rput(1.6,1.5){\smallident} \rput(.4,1.5){\birth}
\end{pspicture}
\quad \xy {\ar@2{->}^{\ell} (0,0)*{};(5,0)*{}}; \endxy  \quad
\begin{pspicture}[0.3](2,3)
  \rput(1,0){\longident}
\end{pspicture}
\]
such that
\[ \psset{xunit=.30cm,yunit=0.30cm}
\xy
 (-25,15)*+{
    \begin{pspicture}[0.5](4.1,5.5)
  \rput(2,0){\mult}
  \rput(1.4,1.5){\curveleft}
  \rput(2.6,1.5){\mult}
  \rput(.8,3){\curveleft}
  \rput(2,3){\curveleft}
  \rput(3.2,3){\mult}
\end{pspicture}
    }="TL";
 (25,15)*+{
     \begin{pspicture}[0.5](4.1,5.5)
  \rput(2,0){\mult}
  \rput(1.4,1.5){\curveleft}
  \rput(2.6,1.5){\mult}
  \rput(.8,3){\mult}
  \rput(2,3){\curveright}
  \rput(3.2,3){\curveright}
\end{pspicture}
    }="TR";
 (0,15)*+{
\begin{pspicture}[0.5](4.1,5.5)
  \rput(2,0){\mult}
  \rput(1.4,1.5){\mult}
  \rput(2.6,1.5){\curveright}
  \rput(.8,3){\curveleft}
  \rput(2,3){\curveleft}
  \rput(3.2,3){\mult}
\end{pspicture}
    }="TM";
 (-25,-15)*+{
     \begin{pspicture}[0.5](4.1,5.5)
  \rput(2,0){\mult}
  \rput(1.4,1.5){\curveleft}
  \rput(2.6,1.5){\mult}
  \rput(.8,3){\curveleft}
  \rput(2,3){\mult}
  \rput(3.2,3){\curveright}
\end{pspicture}
    }="BL";
  (0,-15)*+{
    \begin{pspicture}[0.5](4.1,5.5)
  \rput(2,0){\mult}
  \rput(1.4,1.5){\mult}
  \rput(2.6,1.5){\curveright}
  \rput(.8,3){\curveleft}
  \rput(2,3){\mult}
  \rput(3.2,3){\curveright}
\end{pspicture}
    }="BM";
  (25,-15)*+{
 \begin{pspicture}[0.5](4.1,5.5)
  \rput(2,0){\mult}
  \rput(1.4,1.5){\mult}
  \rput(2.6,1.5){\curveright}
  \rput(.8,3){\mult}
  \rput(2,3){\curveright}
  \rput(3.2,3){\curveright}
\end{pspicture}
    }="BR";
    {\ar@{=>}^{a} "TL";"TM"};
    {\ar@{=>}_{a} "TL";"BL"};
    {\ar@{=>}^{\INT} "TM";"TR"};
    {\ar@{=>}^{a} "TR";"BR"};
    {\ar@{=>}_{a} "BL";"BM"};
    {\ar@{=>}_{a} "BM";"BR"};
 \endxy
\qquad \qquad \psset{xunit=.30cm,yunit=0.30cm}
 \xy
 (-10,12)*+{
   \begin{pspicture}[0.5](3.4,3.8)
  \rput(2,0){\mult}
  \rput(1.4,1.5){\mult}
  \rput(2.6,1.5){\curveright}
  \rput(3.2,3){\smallident}
  \rput(.8,3){\smallident}
  \rput(2,3){\birth}
 \end{pspicture}
    }="TL";
(10,12)*+{
  \begin{pspicture}[0.5](3.4,3.8)
  \rput(2,0){\mult}
  \rput(2.6,1.5){\mult}
  \rput(1.4,1.5){\curveleft}
  \rput(2,3){\birth}
  \rput(.8,3){\smallident}
  \rput(3.2,3){\smallident}
 \end{pspicture}
    }="TR";
 (0,-10)*+{
   \begin{pspicture}[0.7](2.8,1.8)
  \rput(1.4,0){\mult}
 \end{pspicture}
    }="B";
    {\ar@{=>}^{a} "TL";"TR"};
    {\ar@{=>}_{r} "TL";"B"};
    {\ar@{=>}^{\ell} "TR";"B"};
 \endxy
\]
commute.  With this topological description of pseudomonoids it
will be much easier to prove:

\begin{thm}
The walking pseudomonoid is $\Hom(A,A)$ in the walking
pseudoadjunction.
\end{thm}

\Proof Since all of the morphisms in $\Hom(A,A)$ are generated by
$LR$, we will show that $LR$ is a pseudomonoid.  We will again
appeal to string diagrams for the proof.  We define the
multiplication and unit map as in Section \ref{secWA}.
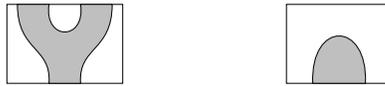
\[
 \xy
  (0,-12)*{};
  (35,0)*{
   \begin{pspicture}(2.2,2.2)
  \pspolygon(-1.1,0)(-1.1,1.5)(1.1,1.5)(1.1,0)(0,0)
  \rput(0,0){\mult}
 \end{pspicture}};
 \endxy
 \qquad \qquad
 \xy
 (0,0)*{
   \begin{pspicture}(2.2,2.2)
  \pspolygon(0,0)(0,1.5)(2,1.5)(2,0)(0,0)
   \pscustom[fillcolor=lightgray, fillstyle=solid]{
  \psbezier(.5,0)(.45,1.2)(1.55,1.2)(1.5,0)
  \psline(.5,0)
  }
 \end{pspicture}}
 \endxy
\]
All that remains is to define the coherence 2-isomorphisms and
show that they satisfy the appropriate coherence axioms.

The map $a \maps (mLR).m \Rrightarrow (LRm).m$ is given by:
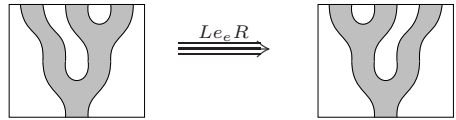
\[ \psset{xunit=.50cm,yunit=0.50cm}  
\begin{pspicture}[0.5](3,3.3)
  \rput(1,0){\multQ} \rput(1.6,1.5){\multQ} \rput(.4,1.5){\curveleftQ}
  \pspolygon(-.8,0)(-.8,3)(2.8,3)(2.8,0)(-.8,0)
\end{pspicture}
 \quad \xy {\ar@3{->}^{Le_e R} (0,0)*{};(12,0)*{}}; \endxy  \qquad
 \quad
\begin{pspicture}[0.5](3,3.3)
  \rput(1,0){\multQ} \rput(1.6,1.5){\curverightQ} \rput(.4,1.5){\multQ}
  \pspolygon(-.8,0)(-.8,3)(2.8,3)(2.8,0)(-.8,0)
\end{pspicture}
\]

In this shorthand notation the first axiom for a pseudomonoid
requires that the diagram of string diagrams:
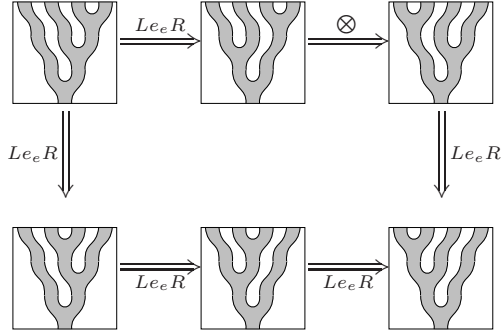
\[ \psset{xunit=.30cm,yunit=0.30cm}
\xy
 (-25,15)*+{
    \begin{pspicture}[0.5](4.1,5.5)
  \rput(2,0){\multQ}
  \rput(1.4,1.5){\curveleftQ}
  \rput(2.6,1.5){\multQ}
  \rput(.8,3){\curveleftQ}
  \rput(2,3){\curveleftQ}
  \rput(3.2,3){\multQ}
    \pspolygon(-.3,0)(-.3,4.5)(4.3,4.5)(4.3,0)(-.3,0)
\end{pspicture}
    }="TL";
 (25,15)*+{
     \begin{pspicture}[0.5](4.1,5.5)
  \rput(2,0){\multQ}
  \rput(1.4,1.5){\curveleftQ}
  \rput(2.6,1.5){\multQ}
  \rput(.8,3){\multQ}
  \rput(2,3){\curverightQ}
  \rput(3.2,3){\curverightQ}
    \pspolygon(-.3,0)(-.3,4.5)(4.3,4.5)(4.3,0)(-.3,0)
\end{pspicture}
    }="TR";
 (0,15)*+{
\begin{pspicture}[0.5](4.1,5.5)
  \rput(2,0){\multQ}
  \rput(1.4,1.5){\multQ}
  \rput(2.6,1.5){\curverightQ}
  \rput(.8,3){\curveleftQ}
  \rput(2,3){\curveleftQ}
  \rput(3.2,3){\multQ}
    \pspolygon(-.3,0)(-.3,4.5)(4.3,4.5)(4.3,0)(-.3,0)
\end{pspicture}
    }="TM";
 (-25,-15)*+{
     \begin{pspicture}[0.5](4.1,5.5)
  \rput(2,0){\multQ}
  \rput(1.4,1.5){\curveleftQ}
  \rput(2.6,1.5){\multQ}
  \rput(.8,3){\curveleftQ}
  \rput(2,3){\multQ}
  \rput(3.2,3){\curverightQ}
    \pspolygon(-.3,0)(-.3,4.5)(4.3,4.5)(4.3,0)(-.3,0)
\end{pspicture}
    }="BL";
  (0,-15)*+{
    \begin{pspicture}[0.5](4.1,5.5)
  \rput(2,0){\multQ}
  \rput(1.4,1.5){\multQ}
  \rput(2.6,1.5){\curverightQ}
  \rput(.8,3){\curveleftQ}
  \rput(2,3){\multQ}
  \rput(3.2,3){\curverightQ}
    \pspolygon(-.3,0)(-.3,4.5)(4.3,4.5)(4.3,0)(-.3,0)
\end{pspicture}
    }="BM";
  (25,-15)*+{
 \begin{pspicture}[0.5](4.1,5.5)
  \rput(2,0){\multQ}
  \rput(1.4,1.5){\multQ}
  \rput(2.6,1.5){\curverightQ}
  \rput(.8,3){\multQ}
  \rput(2,3){\curverightQ}
  \rput(3.2,3){\curverightQ}
  \pspolygon(-.3,0)(-.3,4.5)(4.3,4.5)(4.3,0)(-.3,0)
\end{pspicture}
    }="BR";
    {\ar@{=>}^{Le_eR} "TL";"TM"};
    {\ar@{=>}_{Le_eR} "TL";"BL"};
    {\ar@{=>}^{\INT} "TM";"TR"};
    {\ar@{=>}^{Le_eR} "TR";"BR"};
    {\ar@{=>}_{Le_eR} "BL";"BM"};
    {\ar@{=>}_{Le_eR} "BM";"BR"};
 \endxy
 \]
commutes.  To see that this diagram of string diagrams commutes we
translate it into a traditional commutative diagram:
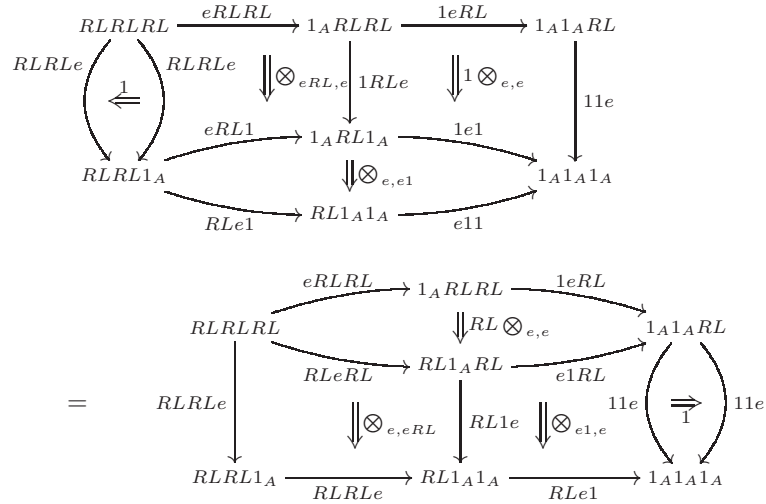
\[
 \def\objectstyle{\scriptstyle}
  \def\labelstyle{\scriptstyle}
 \xy
 (-30,10)*+{RLRLRL}="TL";
 (-30,-10)*+{RLRL1_A}="BL";
 (0,10)*+{1_ARLRL}="TM";
 (0,-5)*+{1_ARL1_A}="M";
 (0,-15)*+{RL1_A1_A}="BM";
 (30,10)*+{1_A1_ARL}="TR";
 (30,-10)*+{1_A1_A1_A}="BR";
    {\ar@/^1.25pc/^>>>>>>>>>>>>{RLRLe} "TL";"BL"};
    {\ar@/_1.25pc/_>>>>>>>>>>>>{RLRLe} "TL";"BL"};
    {\ar@{=>}_{1} (-28,0);(-32,0)};
    {\ar^{eRLRL} "TL";"TM"};
    {\ar^{1RLe} "TM";"M"};
    {\ar@/^.45pc/^{eRL1} "BL";"M"};
    {\ar@/_.45pc/_{RLe1} "BL";"BM"};
    {\ar^{11e} "TR";"BR"};
    {\ar^{1eRL} "TM";"TR"};
    {\ar@/^.45pc/^{1e1} "M";"BR"};
    {\ar@/_.45pc/_{e11} "BM";"BR"};
    {\ar@{=>}^{\INT_{eRL,e}} (-11,6);(-11,0)};
    {\ar@{=>}^{1\INT_{e,e}} (14,6);(14,0)};
    {\ar@{=>}^{\INT_{e,e1}} (0,-8);(0,-12)};
 \endxy \qquad \qquad
\]
\[
\qquad \qquad  =  \qquad
 \def\objectstyle{\scriptstyle}
  \def\labelstyle{\scriptstyle}
 \xy
 (-30,10)*+{RLRLRL}="TL";
 (-30,-10)*+{RLRL1_A}="BL";
 (0,15)*+{1_ARLRL}="TM";
 (0,5)*+{RL1_ARL}="M";
 (0,-10)*+{RL1_A1_A}="BM";
 (30,10)*+{1_A1_ARL}="TR";
 (30,-10)*+{1_A1_A1_A}="BR";
    {\ar@/^1.25pc/^{11e} "TR";"BR"};
    {\ar@/_1.25pc/_{11e} "TR";"BR"};
    {\ar@{=>}_{1} (28,0);(32,0)};
    {\ar_{RLRLe} "BL";"BM"};
    {\ar^{RL1e} "M";"BM"};
    {\ar@/^.45pc/^{eRLRL} "TL";"TM"};
    {\ar@/_.45pc/_{RLeRL} "TL";"M"};
    {\ar_{RLRLe} "TL";"BL"};
    {\ar_{RLe1} "BM";"BR"};
    {\ar@/^.45pc/^{1eRL} "TM";"TR"};
    {\ar@/_.45pc/_{e1RL} "M";"TR"};
    {\ar@{=>}^{\INT_{e,eRL}} (-14,0);(-14,-6)};
    {\ar@{=>}^{\INT_{e1,e}} (11,0);(11,-6)};
    {\ar@{=>}^{RL\INT_{e,e}} (0,12);(0,8)};
 \endxy
\]
which is a consequence of the axioms of a $\cat{Gray}$-category in
Definition~\ref{MonBiCat}, most notably, axioms (vi.) and (vii.).

The unit law isomorphisms $\ell$ and $r$ are given by the diagrams
below:
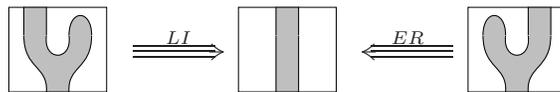
\[ \psset{xunit=.50cm,yunit=0.50cm} 
\begin{pspicture}[0.3](2.3,3)
  \rput(1,0){\multQ} \rput(1.6,1.5){\birthQ} \rput(.4,1.5){\medidentQ}
  \pspolygon(-.3,0)(-.3,2.25)(2.3,2.25)(2.3,0)(0,0)
\end{pspicture}
 \quad \xy {\ar@3{->}^{LI} (0,0)*{};(12,0)*{}}; \endxy  \quad
\begin{pspicture}[0.3](2.3,3)
  \rput(1,0){\identQ} \rput(1,1.5){\medidentQ}
  \pspolygon(-.3,0)(-.3,2.25)(2.3,2.25)(2.3,0)(0,0)
\end{pspicture}
 \quad \xy {\ar@3{->}_{ER} (12,0)*{};(0,0)*{}}; \endxy  \quad
\begin{pspicture}[0.3](2.3,3)
  \rput(1,0){\multQ} \rput(1.6,1.5){\medidentQ} \rput(.4,1.5){\birthQ}
  \pspolygon(-.3,0)(-.3,2.25)(2.3,2.25)(2.3,0)(0,0)
\end{pspicture}
\]
The coherence for these isomorphisms requires that the diagram of
string diagrams:
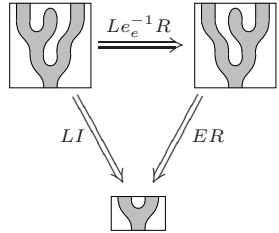
\[ \psset{xunit=.30cm,yunit=0.30cm}
 \xy
 (-12,12)*+{
   \begin{pspicture}[0.5](3.8,3.8)
  \rput(2,0){\multQ}
  \rput(1.4,1.5){\multQ}
  \rput(2.6,1.5){\curverightQ}
  \rput(3.2,3){\medidentQ}
  \rput(.8,3){\medidentQ}
  \rput(2,3){\birthQ}
  \pspolygon(.2,0)(.2,3.75)(3.8,3.75)(3.8,0)(.5,0)
 \end{pspicture}
    }="TL";
(12,12)*+{
  \begin{pspicture}[0.5](3.4,3.8)
  \rput(2,0){\multQ}
  \rput(2.6,1.5){\multQ}
  \rput(1.4,1.5){\curveleftQ}
  \rput(2,3){\birthQ}
  \rput(.8,3){\medidentQ}
  \rput(3.2,3){\medidentQ}
  \pspolygon(.2,0)(.2,3.75)(3.8,3.75)(3.8,0)(.5,0)
 \end{pspicture}
    }="TR";
 (0,-10)*+{
   \begin{pspicture}[0.7](2.8,1.8)
  \rput(1.4,0){\multQ}
  \pspolygon(.2,0)(.2,1.5)(2.6,1.5)(2.6,0)(.2,0)
 \end{pspicture}
    }="B";
    {\ar@{=>}^{\scs Le_{e}^{-1}R} "TL";"TR"};
    {\ar@{=>}_{\scs LI} "TL";"B"};
    {\ar@{=>}^{\scs ER} "TR";"B"};
 \endxy
\]
commutes. This axiom follows from the triangulator identities in
the definition of a pseudoadjunction.  Thus, $\Hom(A,A)$ in the
walking pseudoadjunction contains the walking pseudomonoid.

To prove the converse we again appeal to monad theory.  In this
case, we are interested in pseudomonads.  A pseudomonad on an
object $A$ of a semistrict 3-category is just a pseudomonoid in
the semistrict monoidal 2-category $\Hom(A,A)$ \cite{mar}.  Since
every semistrict monoidal 2-category $\mathcal{C}$ can be regarded
as a one object semistrict 3-category $\Sigma(\mathcal{C})$, a
pseudomonoid in $\mathcal{C}$ amounts to a pseudomonad on the one
object of $\Sigma(\mathcal{C})$. Thus, the problem of showing that
a semistrict monoidal 2-category generated by a pseudomonoid
extends to a semistrict 3-category  generated by a
pseudoadjunction amounts to showing that every pseudomonad arises
from a pseudoadjunction.

Again, the Eilenberg-Moore construction saves the day.  Since
Eilenberg-Moore objects can be defined as a weighted limit, this
notion can be extended using enriched category theory to define
Eilenberg-Moore objects in a $\cat{Gray}$-category. While it is
not true that an arbitrary $\cat{Gray}$-category always posses an
Eilenberg-Moore object for every pseudomonad, using enriched
category theory we can define the free Eilenberg-Moore completion
\cat{EM}$(\mathcal{K})$ of a $\cat{Gray}$-category $\mathcal{K}$
\cite{LS2,Lau1}.  By the theory of such completions we obtain a
\cat{Gray}-functor $Z \maps \mathcal{K} \to \cat{EM}(\mathcal{K})$
with the property that for any \cat{Gray}-category $\mathcal{L}$
with Eilenberg-Moore objects, composition with $Z$ induces an
equivalence of categories between the \cat{Gray}-functor category
$[\mathcal{K},\mathcal{L} ]$ and the full subcategory of the
\cat{Gray}-functor category $[\cat{EM}(\mathcal{K}),\mathcal{L}]$
consisting of those \cat{Gray}-functors which preserve
Eilenberg-Moore objects~\cite{LS2}. Furthermore, $Z$ will be fully
faithful.

In the Eilenberg-Moore completion \cat{EM}$(\mathcal{K})$ every
pseudomonad in the \cat{Gray}-category $\mathcal{K}$ arises from a
pseudoadjunction. If $\T$ is the pseudomonad in $\Sigma(\wMon)$
generating the suspension of the walking pseudomonoid, then we can
freely complete $\Sigma(\wMon)$ under an Eilenberg-Moore object
for just the pseudomonad $\T$. Again, we denote the completion
under an Eilenberg-Moore object for the pseudomonad $\T$ as
$\cat{EM}_{\T}(\mathcal{K})$, and in this $\cat{Gray}$-category
the pseudomonad $\T$ is generated by a pseudoadjunction
$\padjunction{A}{A^{\T}}{}{}$. By the universal property of the
walking pseudomonoid, this determines a $\cat{Gray}$-functor
$\Lambda \maps \wAdj \to \cat{EM}_{\T}(\mathcal{K})$.

To construct an inverse to the $\cat{Gray}$-functor $\Lambda$,
note that we have shown that the morphism $LR$ in $\Hom(A,A)$ of
$\wAdj$ is a pseudomonad. The universal property of the walking
pseudomonoid then determines a $\cat{Gray}$-functor $\wMon \to
\wAdj$. Further, one can check that the pseudomonad $LR$ has an
Eilenberg-Moore object in $\wAdj$, namely $B$.  Hence, by the
universal property of the Eilenberg-Moore completion, we get a
\cat{Gray}-functor $\bar{\Lambda} \maps \cat{EM}_{\T}(\mathcal{K})
\to \wAdj$ that preserves Eilenberg-Moore objects. It is easy to
see that $\Lambda$ and $\bar{\Lambda}$ define an isomorphism of
$\cat{Gray}$-categories. \qed

\subsection{The walking pseudo ambijunction} \label{secWWAA}

An ambidextrous pseudoadjunction, or pseudo ambijunction for
short, is a 2-sided pseudoadjunction.  This means that we have the
additional 2-morphisms $j\maps 1_B \To RL$ and $k \maps LR \to
1_A$ and the additional 3-morphisms
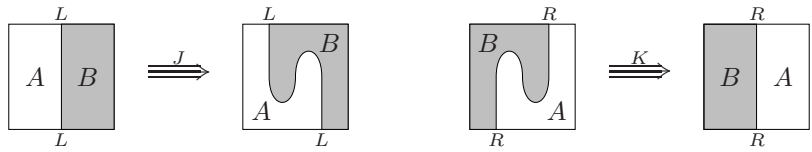
\[ \psset{xunit=.70cm,yunit=0.70cm} \xy
   (0,0)*+{  \begin{pspicture}(2,2.2)
  \pspolygon(0,0)(0,2)(2,2)(2,0)(0,0)
\pscustom[fillcolor=lightgray, fillstyle=solid]{
  \pspolygon(1,0)(1,2)(2,2)(2,0)(1,0)
  }
  \rput(.5,1){$A$}
         \rput(1.5,1){$B$}
         \rput(1,2.2){$\scs L$}
         \rput(1,-.2){$\scs L$}
 \end{pspicture}};
 \endxy
 \quad \xy {\ar@3{->}^{J} (0,0)*{};(8,0)*{}}; \endxy  \quad
     \xy
    (0,0)*+{   \begin{pspicture}(2,2.2)
  \pspolygon(0,0)(0,2)(2,2)(2,0)(0,0)
\pscustom[fillcolor=lightgray, fillstyle=solid]{
  \psbezier(.5,1)(.5,.35)(1,.35)(1,1)
  \psbezier(1,1)(1,1.65)(1.5,1.65)(1.5,1)
  \psline(1.5,0)
  \psline(2,0)
  \psline(2,2)
  \psline(.5,2)
  \psline(.5,1)
  }
          \rput(1.65,1.65){$B$}
         \rput(.35,.35){$A$}
         \rput(.5,2.2){$\scs L$}
         \rput(1.5,-.2){$\scs L$}
 \end{pspicture}};
  \endxy
\qquad \qquad \psset{xunit=.70cm,yunit=0.70cm} \xy
   (0,0)*+{\begin{pspicture}(2,2.2)
  \pspolygon(0,0)(0,2)(2,2)(2,0)(0,0)
\pscustom[fillcolor=lightgray, fillstyle=solid]{
  \psbezier(.5,1)(.5,1.65)(1,1.65)(1,1)
  \psbezier(1,1)(1,.35)(1.5,.35)(1.5,1)
  \psline(1.5,2)
  \psline(0,2)
  \psline(0,0)
  \psline(.5,0)
  \psline(.5,1)
  }
         \rput(.35,1.65){$B$}
         \rput(1.65,.35){$A$}
         \rput(1.5,2.2){$\scs R$}
         \rput(.5,-.2){$\scs R$}
 \end{pspicture}};
 \endxy
 \quad \xy {\ar@3{->}^{K} (0,0)*{};(8,0)*{}}; \endxy  \quad
\xy
      (0,0)*+{\begin{pspicture}(2,2.2)
\pscustom[fillcolor=lightgray, fillstyle=solid]{
  \pspolygon(1,0)(1,2)(0,2)(0,0)(1,0)
  }
   \rput(.5,1){$B$}
         \rput(1.5,1){$A$}
         \rput(1,2.2){$\scs R$}
         \rput(1,-.2){$\scs R$}
  \pspolygon(0,0)(0,2)(2,2)(2,0)(0,0)
 \end{pspicture}};
    \endxy
\]
These 3-morphisms must satisfy the coherence conditions that the
composite:
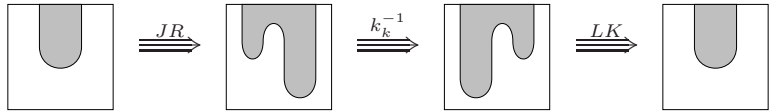
\[
\psset{xunit=.7cm,yunit=.7cm}
\begin{pspicture}[.5](2,2.2)
  \pspolygon(0,0)(0,2)(2,2)(2,0)(0,0)
\pscustom[fillcolor=lightgray, fillstyle=solid]{
  \psbezier(.6,1.2)(.6,.65)(1.4,.65)(1.4,1.2)
  \psline(1.4,2)\psline(.6,2)  \psline(.6,1.2)
  }
 \end{pspicture}
 \quad \xy {\ar@3{->}^{JR} (0,0)*{};(8,0)*{}}; \endxy  \quad
 \begin{pspicture}[.5](2,2.2)
  \pspolygon(0,0)(0,2)(2,2)(2,0)(0,0)
\pscustom[fillcolor=lightgray, fillstyle=solid]{
  \psbezier(.3,1.3)(.3,.85)(.7,.85)(.7,1.3)
  \psbezier(.7,1.3)(.7,1.75)(1.1,1.75)(1.1,1.3)
  \psline(1.1,.6)
  \psbezier(1.1,.6)(1.1,.1)(1.7,.1)(1.7,.6)
  \psline(1.7,2)
  \psline(.3,2)
  \psline(.3,1.3)
  }
 \end{pspicture}
 \quad \xy {\ar@3{->}^{k^{-1}_k} (0,0)*{};(8,0)*{}}; \endxy  \quad
 \begin{pspicture}[.5](2,2.2)
  \pspolygon(0,0)(0,2)(2,2)(2,0)(0,0)
\pscustom[fillcolor=lightgray, fillstyle=solid]{
  \psbezier(1.7,1.3)(1.7,.85)(1.3,.85)(1.3,1.3)
  \psbezier(1.3,1.3)(1.3,1.75)(.9,1.75)(.9,1.3)
  \psline(.9,.6)
  \psbezier(.9,.6)(.9,.1)(.3,.1)(.3,.6)
  \psline(.3,2)
  \psline(1.7,2)
  \psline(1.7,1.3)
  }
 \end{pspicture}
 \quad \xy {\ar@3{->}^{LK} (0,0)*{};(8,0)*{}}; \endxy  \quad
\begin{pspicture}[.5](2,2.2)
  \pspolygon(0,0)(0,2)(2,2)(2,0)(0,0)
\pscustom[fillcolor=lightgray, fillstyle=solid]{
  \psbezier(.6,1.2)(.6,.65)(1.4,.65)(1.4,1.2)
  \psline(1.4,2)\psline(.6,2)  \psline(.6,1.2)
  }
 \end{pspicture}
\]
is equal to the composite:
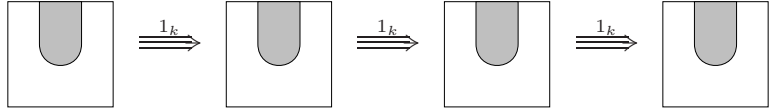
\[
\psset{xunit=.7cm,yunit=.7cm}
\begin{pspicture}[.5](2,2.2)
  \pspolygon(0,0)(0,2)(2,2)(2,0)(0,0)
\pscustom[fillcolor=lightgray, fillstyle=solid]{
  \psbezier(.6,1.2)(.6,.65)(1.4,.65)(1.4,1.2)
  \psline(1.4,2)\psline(.6,2)  \psline(.6,1.2)
  }
 \end{pspicture}
 \quad \xy {\ar@3{->}^{1_k} (0,0)*{};(8,0)*{}}; \endxy  \quad
\begin{pspicture}[.5](2,2.2)
  \pspolygon(0,0)(0,2)(2,2)(2,0)(0,0)
\pscustom[fillcolor=lightgray, fillstyle=solid]{
  \psbezier(.6,1.2)(.6,.65)(1.4,.65)(1.4,1.2)
  \psline(1.4,2)\psline(.6,2)  \psline(.6,1.2)
  }
 \end{pspicture}
 \quad \xy {\ar@3{->}^{1_k} (0,0)*{};(8,0)*{}}; \endxy  \quad
\begin{pspicture}[.5](2,2.2)
  \pspolygon(0,0)(0,2)(2,2)(2,0)(0,0)
\pscustom[fillcolor=lightgray, fillstyle=solid]{
  \psbezier(.6,1.2)(.6,.65)(1.4,.65)(1.4,1.2)
  \psline(1.4,2)\psline(.6,2)  \psline(.6,1.2)
  }
 \end{pspicture}
 \quad \xy {\ar@3{->}^{1_k} (0,0)*{};(8,0)*{}}; \endxy  \quad
\begin{pspicture}[.5](2,2.2)
  \pspolygon(0,0)(0,2)(2,2)(2,0)(0,0)
\pscustom[fillcolor=lightgray, fillstyle=solid]{
  \psbezier(.6,1.2)(.6,.65)(1.4,.65)(1.4,1.2)
  \psline(1.4,2)\psline(.6,2)  \psline(.6,1.2)
  }
 \end{pspicture}
\]
and that the composite:
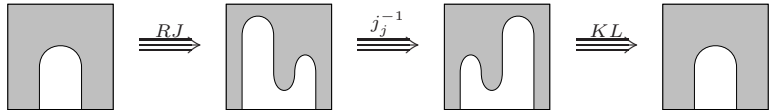
\[
 \begin{pspicture}[.5](2,2.2)
  \pspolygon(0,0)(0,2)(2,2)(2,0)(0,0)
\pscustom[fillcolor=lightgray, fillstyle=solid]{
  \psbezier(.6,.8)(.6,1.35)(1.4,1.35)(1.4,.8)
  \psline(1.4,0)\psline(2,0) \psline(2,2) \psline(0,2)
  \psline(0,0) \psline(.6,0) \psline(.6,.8)
  }
 \end{pspicture}
 \quad \xy {\ar@3{->}^{RJ} (0,0)*{};(8,0)*{}}; \endxy  \quad
\begin{pspicture}[.5](2,2.2)
  \pspolygon(0,0)(0,2)(2,2)(2,0)(0,0)
\pscustom[fillcolor=lightgray, fillstyle=solid]{
  \psbezier(1.7,.7)(1.7,1.15)(1.3,1.15)(1.3,.7)
  \psbezier(1.3,.7)(1.3,.25)(.9,.25)(.9,.7)
  \psline(.9,1.4)
  \psbezier(.9,1.4)(.9,1.9)(.3,1.9)(.3,1.4)
  \psline(.3,0)
  \psline(0,0)
  \psline(0,2)
  \psline(2,2)
  \psline(2,0)
  \psline(1.7,0)
  \psline(1.7,.7)
  }
 \end{pspicture}
  \quad \xy {\ar@3{->}^{j^{-1}_j} (0,0)*{};(8,0)*{}}; \endxy  \quad
\begin{pspicture}[.5](2,2.2)
  \pspolygon(0,0)(0,2)(2,2)(2,0)(0,0)
\pscustom[fillcolor=lightgray, fillstyle=solid]{
  \psbezier(.3,.7)(.3,1.15)(.7,1.15)(.7,.7)
  \psbezier(.7,.7)(.7,.25)(1.1,.25)(1.1,.7)
  \psline(1.1,1.4)
  \psbezier(1.1,1.4)(1.1,1.9)(1.7,1.9)(1.7,1.4)
  \psline(1.7,0)
  \psline(2,0)
  \psline(2,2)
  \psline(0,2)
  \psline(0,0)
  \psline(.3,0)
  \psline(.3,.7)
  }
 \end{pspicture}
  \quad \xy {\ar@3{->}^{KL} (0,0)*{};(8,0)*{}}; \endxy  \quad
 \begin{pspicture}[.5](2,2.2)
  \pspolygon(0,0)(0,2)(2,2)(2,0)(0,0)
\pscustom[fillcolor=lightgray, fillstyle=solid]{
  \psbezier(.6,.8)(.6,1.35)(1.4,1.35)(1.4,.8)
  \psline(1.4,0)\psline(2,0) \psline(2,2) \psline(0,2)
  \psline(0,0) \psline(.6,0) \psline(.6,.8)
  }
 \end{pspicture}
 \]
is equal to the composite:
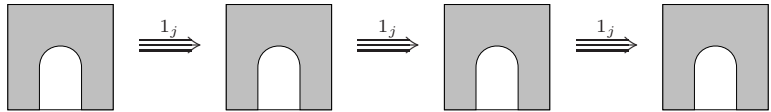
\[
 \begin{pspicture}[.5](2,2.2)
  \pspolygon(0,0)(0,2)(2,2)(2,0)(0,0)
\pscustom[fillcolor=lightgray, fillstyle=solid]{
  \psbezier(.6,.8)(.6,1.35)(1.4,1.35)(1.4,.8)
  \psline(1.4,0)\psline(2,0) \psline(2,2) \psline(0,2)
  \psline(0,0) \psline(.6,0) \psline(.6,.8)
  }
 \end{pspicture}
 \quad \xy {\ar@3{->}^{1_j} (0,0)*{};(8,0)*{}}; \endxy  \quad
 \begin{pspicture}[.5](2,2.2)
  \pspolygon(0,0)(0,2)(2,2)(2,0)(0,0)
\pscustom[fillcolor=lightgray, fillstyle=solid]{
  \psbezier(.6,.8)(.6,1.35)(1.4,1.35)(1.4,.8)
  \psline(1.4,0)\psline(2,0) \psline(2,2) \psline(0,2)
  \psline(0,0) \psline(.6,0) \psline(.6,.8)
  }
 \end{pspicture}
  \quad \xy {\ar@3{->}^{1_j} (0,0)*{};(8,0)*{}}; \endxy  \quad
 \begin{pspicture}[.5](2,2.2)
  \pspolygon(0,0)(0,2)(2,2)(2,0)(0,0)
\pscustom[fillcolor=lightgray, fillstyle=solid]{
  \psbezier(.6,.8)(.6,1.35)(1.4,1.35)(1.4,.8)
  \psline(1.4,0)\psline(2,0) \psline(2,2) \psline(0,2)
  \psline(0,0) \psline(.6,0) \psline(.6,.8)
  }
 \end{pspicture}
  \quad \xy {\ar@3{->}^{1_j} (0,0)*{};(8,0)*{}}; \endxy  \quad
 \begin{pspicture}[.5](2,2.2)
  \pspolygon(0,0)(0,2)(2,2)(2,0)(0,0)
\pscustom[fillcolor=lightgray, fillstyle=solid]{
  \psbezier(.6,.8)(.6,1.35)(1.4,1.35)(1.4,.8)
  \psline(1.4,0)\psline(2,0) \psline(2,2) \psline(0,2)
  \psline(0,0) \psline(.6,0) \psline(.6,.8)
  }
 \end{pspicture}
 \]

\begin{defn}
The walking ambidextrous pseudoadjunction $\wAmbi$ is semistrict
3-category freely generated by a pseudo ambijunction.
\end{defn}

We now define a pseudo Frobenius algebra by categorifying the
relationship between Frobenius algebras and adjunctions.

\begin{defn}
The walking pseudo Frobenius algebra is the semistrict monoidal
2-category $\Hom(A,A)$ in the walking pseudo ambijunction. Hence,
a pseudo Frobenius algebra is $LR \in \Hom(A,A)$ for some pseudo
ambijunction $L \dashv_p R \maps B \to A$.
\end{defn}

Although we obtained this definition by categorifying the
relationship between Frobenius algebras and adjunctions, it is
equivalent to a definition that can be obtained by replacing
equations with isomorphisms and determining the correct coherence
conditions. Below we provide an equivalent definition that can be
viewed as a categorification of Proposition \ref{equivFrob}
(iii.).  In the next section we will see categorifications of the
descriptions given by Proposition \ref{equivFrob} (i.) and (ii.)
and later we will show that these are also equivalent to the
description given above.  The definition we present below is the
simplest and most easily related to pseudo ambijunctions.  It
might be described as a pseudomonoid equipped with a form defining
a `weakly nondegenerate' pairing.  It is also very related to the
notion of a Frobenius pseudomonoid defined by
Street~\cite{StreetFrob}.

\begin{prop} \label{FFR1}
A pseudo Frobenius algebra can be equivalently defined as follows:
A pseudo Frobenius algebra is an object $F$ of a semistrict
monoidal 2-category $\mathcal{C}$ equipped with morphisms:
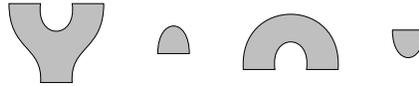
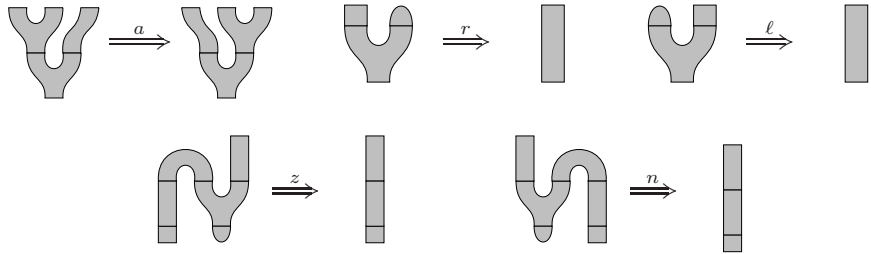
\[ \psset{xunit=.7cm,yunit=0.7cm}
\begin{pspicture}[0.5](1.2,1.6)
  \rput(.6,0){\mult}
\end{pspicture}
\qquad
\begin{pspicture}[0.5](1.2,1.6)
  \rput(.6,.55){\birth}
\end{pspicture}
\qquad
\begin{pspicture}[0.5](1.2,1.6)
  \rput(.6,.25){\zag}
\end{pspicture}
\qquad
\begin{pspicture}[0.5](1,1.6)
  \rput(.6,1){\death}
\end{pspicture}
\]
and 2-isomorphisms:
\[ \psset{xunit=.4cm,yunit=0.4cm}   
\begin{pspicture}[0.5](2.6,3.3)
  \rput(1,0){\mult} \rput(1.6,1.5){\curveright} \rput(.4,1.5){\mult}
\end{pspicture}
 \; \xy {\ar@2{->}^{a} (0,0)*{};(8,0)*{}}; \endxy  \quad
\begin{pspicture}[0.5](3,3.3)
  \rput(1,0){\mult} \rput(1.6,1.5){\mult} \rput(.4,1.5){\curveleft}
\end{pspicture}
\qquad  \psset{xunit=.5cm,yunit=0.5cm}
\begin{pspicture}[0.3](2,3)
  \rput(1,0){\mult} \rput(1.6,1.5){\birth} \rput(.4,1.5){\smallident}
\end{pspicture}
 \quad \xy {\ar@2{->}^{r} (0,0)*{};(6,0)*{}}; \endxy  \quad
\begin{pspicture}[0.3](2,3)
  \rput(1,0){\longident}
\end{pspicture}
\qquad\begin{pspicture}[0.3](2,3)
  \rput(1,0){\mult} \rput(1.6,1.5){\smallident} \rput(.4,1.5){\birth}
\end{pspicture}
\quad \xy {\ar@2{->}^{\ell} (0,0)*{};(6,0)*{}}; \endxy  \quad
\begin{pspicture}[0.3](2,3)
  \rput(1,0){\longident}
\end{pspicture}
\]

\[
 \psset{xunit=.4cm,yunit=0.4cm}
\begin{pspicture}[.5](3,3)
  \rput(2.2,.55){\mult} \rput(1,2.05){\zag}
  \rput(2.8,2.05){\ident} \rput(.4,0){\smallident} \rput(.4,.55){\ident}
  \rput(2.2,.55){\death}
\end{pspicture}
\quad \xy {\ar@2{->}^{z} (0,0)*{};(6,0)*{}}; \endxy  \quad
\begin{pspicture}[0.5](2,3)
  \rput(1,2.05){\ident} \rput(1,0){\smallident} \rput(1,.55){\ident}
\end{pspicture}
\qquad \qquad
\begin{pspicture}[0.5](3,3)
  \rput(1,.55){\mult} \rput(2.2,2.05){\zag}
  \rput(.4,2.05){\ident} \rput(2.8,.55){\ident}
  \rput(2.8,0){\smallident} \rput(1,.55){\death}
\end{pspicture}
\quad \xy {\ar@2{->}^{n} (0,0)*{};(6,0)*{}}; \endxy  \quad
\begin{pspicture}[0.5](2,3.6)
  \rput(1,2.05){\ident} \rput(1,.55){\ident} \rput(1,0){\smallident}
\end{pspicture}
\]
satisfying the pseudomonoid axioms and making the following
diagrams commute:
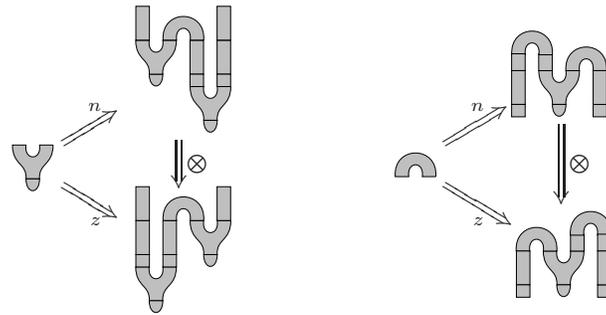
\[ \psset{xunit=.30cm,yunit=0.30cm}
 \xy
 (-10,0)*+{
    \begin{pspicture}[0.5](2.4,2.4)
  \rput(1.4,.55){\mult}
  \rput(1.4,.55){\death}
 \end{pspicture}
    }="TL";
(10,12)*+{
   \begin{pspicture}[0.5](4.8,5.2)
  \rput(1.4,2.6){\death}
  \rput(3.8,.55){\mult}
  \rput(3.8,.55){\death}
  \rput(4.4,2.6){\ident}
  \rput(3.2,2.6){\ident}
  \rput(4.4,2.05){\smallident}
  \rput(3.2,2.05){\smallident}
  \rput(1.4,2.6){\mult}
  \rput(.8,4.1){\ident}
  \rput(4.4,4.1){\ident}
  \rput(2.6,4.1){\zag}
 \end{pspicture}
    }="TR";
 (10,-12)*+{
   \begin{pspicture}[0.5](4.8,5.2)
  \rput(1.4,.55){\mult}
  \rput(1.4,.55){\death}
  \rput(3.8,2.6){\death}
  \rput(3.8,2.6){\mult}
  \rput(2,2.6){\ident}
  \rput(.8,2.6){\ident}
  \rput(2,2.05){\smallident}
  \rput(.8,2.05){\smallident}
  \rput(.8,4.1){\ident}
   \rput(4.4,4.1){\ident}
  \rput(2.6,4.1){\zag}
 \end{pspicture}
    }="BR";
    {\ar@{=>}^{n} "TL";"TR"};
    {\ar@{=>}^{\INT} "TR";"BR"};
    {\ar@{=>}_{z} "TL";"BR"};
 \endxy
 \qquad \qquad \quad
 \xy
 (-10,0)*+{
    \begin{pspicture}[0.5](2.4,1.1)
  \rput(1.4,0){\zag}
 \end{pspicture}
    }="TL";
(10,12)*+{
   \begin{pspicture}[0.5](4.2,3.8)
  \rput(2,.55){\mult}
  \rput(2,.55){\death}
  \rput(3.2,2.05){\zag}
  \rput(1.4,2.05){\medident}
  \rput(.8,2.75){\zag}
  \rput(.2,2.05){\medident}
  \rput(.2,.55){\ident}
  \rput(.2,0){\smallident}
  \rput(3.8,0){\smallident}
  \rput(3.8,.55){\ident}
 \end{pspicture}
    }="TR";
 (10,-12)*+{
   \begin{pspicture}[0.5](3.8,3.8)
  \rput(2,.55){\mult}
  \rput(2,.55){\death}
  \rput(.8,2.05){\zag}
  \rput(2.6,2.05){\medident}
  \rput(3.2,2.75){\zag}
  \rput(3.8,2.05){\medident}
  \rput(.2,.55){\ident}
  \rput(.2,0){\smallident}
  \rput(3.8,0){\smallident}
  \rput(3.8,.55){\ident}
 \end{pspicture}
    }="BR";
    {\ar@{=>}^{n} "TL";"TR"};
    {\ar@{=>}^{\INT} "TR";"BR"};
    {\ar@{=>}_{z} "TL";"BR"};
 \endxy
 \]
Note that in the first diagram above the 2-isomorphisms $\INT$ is
actually two applications of $\INT$.
\end{prop}

\Proof  In the previous section we saw that a pseudomonoid can
always be defined as $LR \in \Hom(A,A)$ for some pseudoadjunction
$\padjunction{A}{B}{L}{R}$.  This follows from the fact that the
walking pseudomonoid is $\Hom(A,A)$ in the walking
pseudoadjunction.  In order to prove the proposition, we need to
show that in an ambidextrous pseudoadjunction the morphisms $LR$
has the additional structure described by the diagrams and axioms
above.  The 2-isomorphisms $z$ and $n$ in $\Hom(A,A)$ are defined
by the pasting composites depicted below:

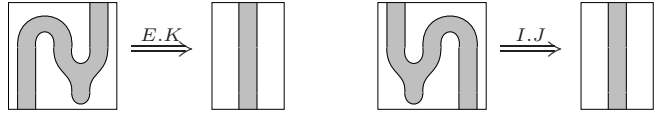
\[
 \psset{xunit=.4cm,yunit=0.4cm}
\begin{pspicture}[.5](3,3.6)
  \rput(2.2,.55){\multQ} \rput(1,2.05){\zagQ}
  \rput(2.8,2.05){\identQ} \rput(.4,0){\smallidentQ} \rput(.4,.55){\identQ}
  \rput(2.2,.55){\deathQ}
  \pspolygon(-.2,0)(-.2,3.55)(3.4,3.55)(3.4,0)(-.2,0)
\end{pspicture}
\quad \xy {\ar@2{->}^{E.K} (0,0)*{};(8,0)*{}}; \endxy  \quad
\begin{pspicture}[0.5](2,3.6)
  \rput(1,2.05){\identQ} \rput(1,0){\smallidentQ} \rput(1,.55){\identQ}
  \pspolygon(-.2,0)(-.2,3.55)(2.2,3.55)(2.2,0)(-.2,0)
\end{pspicture}
\qquad \qquad
\begin{pspicture}[0.5](3,3.6)
  \rput(1,.55){\multQ} \rput(2.2,2.05){\zagQ}
  \rput(.4,2.05){\identQ} \rput(2.8,.55){\identQ}
  \rput(2.8,0){\smallidentQ} \rput(1,.55){\deathQ}
  \pspolygon(-.2,0)(-.2,3.55)(3.4,3.55)(3.4,0)(-.2,0)
\end{pspicture}
\quad \xy {\ar@2{->}^{I.J} (0,0)*{};(8,0)*{}}; \endxy  \quad
\begin{pspicture}[0.5](2,3.6)
  \rput(1,2.05){\identQ} \rput(1,.55){\identQ} \rput(1,0){\smallidentQ}
  \pspolygon(-.2,0)(-.2,3.55)(2.2,3.55)(2.2,0)(-.2,0)
\end{pspicture}
\]
The precise definition of the pasting composites can be deduced
from the diagrams.  These 2-isomorphisms satisfy the specified
axioms since:
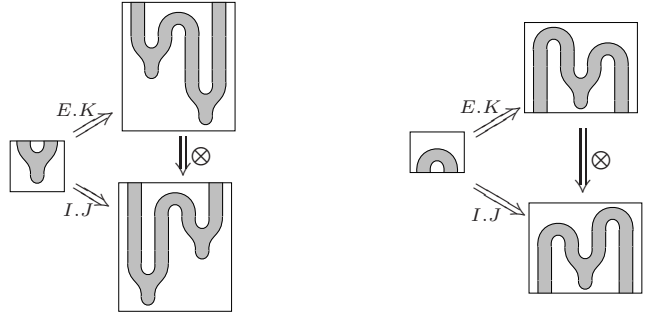
\[ \psset{xunit=.30cm,yunit=0.30cm}
 \xy
 (-10,0)*++{
    \begin{pspicture}[0.5](2.4,2.4)
  \rput(1.4,.55){\multQ}
  \rput(1.4,.55){\deathQ}
    \pspolygon(.2,-.2)(.2,2.05)(2.6,2.05)(2.6,-.2)(.2,-.2)
 \end{pspicture}
    }="TL";
(10,12)*+{
   \begin{pspicture}[0.5](5.6,5.2)
  \rput(1.4,2.6){\deathQ}
  \rput(3.8,.55){\multQ}
  \rput(3.8,.55){\deathQ}
  \rput(4.4,2.6){\identQ}
  \rput(3.2,2.6){\identQ}
  \rput(4.4,2.05){\smallidentQ}
  \rput(3.2,2.05){\smallidentQ}
  \rput(1.4,2.6){\multQ}
  \rput(.8,4.1){\identQ}
  \rput(4.4,4.1){\identQ}
  \rput(2.6,4.1){\zagQ}
  \pspolygon(.1,-.1)(.1,5.6)(5.1,5.6)(5.1,-.1)(.1,-.1)
 \end{pspicture}
    }="TR";
 (10,-12)*+{
   \begin{pspicture}[0.5](5.9,5.2)
  \rput(1.4,.55){\multQ}
  \rput(1.4,.55){\deathQ}
  \rput(3.8,2.6){\deathQ}
  \rput(3.8,2.6){\multQ}
  \rput(2,2.6){\identQ}
  \rput(.8,2.6){\identQ}
  \rput(2,2.05){\smallidentQ}
  \rput(.8,2.05){\smallidentQ}
  \rput(.8,4.1){\identQ}
   \rput(4.4,4.1){\identQ}
  \rput(2.6,4.1){\zagQ}
  \pspolygon(.1,-.1)(.1,5.6)(5.1,5.6)(5.1,-.1)(.1,-.1)
 \end{pspicture}
    }="BR";
    {\ar@{=>}^-{\scs E.K} "TL";"TR"};
    {\ar@{=>}^-{\INT} "TR";(10,-2)};
    {\ar@{=>}_-{\scs I.J} "TL";"BR"};
 \endxy
 \qquad \qquad \quad
 \xy
 (-10,0)*++{
    \begin{pspicture}[0.5](2.4,1.1)
  \rput(1.4,0){\zagQ}
  \pspolygon(.2,0)(.2,1.8)(2.6,1.8)(2.6,0)(.2,0)
 \end{pspicture}
    }="TL";
(10,12)*++{
   \begin{pspicture}[0.5](4.2,3.8)
  \rput(2,.55){\multQ}
  \rput(2,.55){\deathQ}
  \rput(3.2,2.05){\zagQ}
  \rput(1.4,2.05){\medidentQ}
  \rput(.8,2.75){\zagQ}
  \rput(.2,2.05){\medidentQ}
  \rput(.2,.55){\identQ}
  \rput(.2,0){\smallidentQ}
  \rput(3.8,0){\smallidentQ}
  \rput(3.8,.55){\identQ}
    \pspolygon(-.5,0)(-.5,4)(4.5,4)(4.5,0)(-.5,0)
 \end{pspicture}
    }="TR";
 (10,-12)*++{
   \begin{pspicture}[0.5](3.8,3.8)
  \rput(2,.55){\multQ}
  \rput(2,.55){\deathQ}
  \rput(.8,2.05){\zagQ}
  \rput(2.6,2.05){\medidentQ}
  \rput(3.2,2.75){\zagQ}
  \rput(3.8,2.05){\medidentQ}
  \rput(.2,.55){\identQ}
  \rput(.2,0){\smallidentQ}
  \rput(3.8,0){\smallidentQ}
  \rput(3.8,.55){\identQ}
  \pspolygon(-.5,0)(-.5,4)(4.5,4)(4.5,0)(-.5,0)
 \end{pspicture}
    }="BR";
    {\ar@{=>}^-{\scs E.K} "TL";"TR"};
    {\ar@{=>}^{\INT} "TR";"BR"};
    {\ar@{=>}_-{\scs I.J} "TL";"BR"};
 \endxy
 \]
commute from the fact that each set $(I,E)$ and $(J,K)$ satisfy
the triangulator identities.

The converse is proven, again, using monad theory.  Given a pseudo
Frobenius object as defined by the diagrams and axioms in the
proposition, it is clear that the object $F$ of the semistrict
monoidal 2-category $\mathcal{C}$ is a pseudomonoid in
$\mathcal{C}$.  Hence, by definition $\F$ is a pseudomonad on the
one object of the semistrict 3-category $\Sigma(\mathcal{C})$. The
maps $z$ and $n$ provide this pseudomonad with some additional
structure.  In particular, these maps make $\F$ into a {\em
Frobenius pseudomonad}~\cite{Lau1,StreetFrob}.

Using an enriched version of the Eilenberg-Moore completion, the
author has shown that every Frobenius pseudomonad in a semistrict
3-category $\mathcal{K}$ arises from an ambidextrous
pseudoadjunction in $\cat{EM}(\mathcal{K})$, where
$\cat{EM}(\mathcal{K})$ is the free completion of $\mathcal{K}$
under Eilenberg-Moore objects for every pseudomonad in
$\mathcal{K}$~\cite{Lau1}. Since a pseudo Frobenius algebra $F$ in
$\mathcal{C}$ is the same thing as a Frobenius pseudomonad $\F$ on
the one object of $\Sigma(\mathcal{C})$, it follows that $F$
arises from some pseudo ambijunction in
$\cat{EM}\big(\Sigma(\mathcal{C})\big)$.  \qed

We have actually done more than provided an equivalent
characterization of pseudo Frobenius algebras.  It is a simple
extension to show:

\begin{cor}
The walking pseudo Frobenius algebra is the semistrict monoidal
2-category freely generated by a single object, and the morphisms
and 2-isomorphisms of Proposition \ref{FFR1} subject to the axioms
given.
\end{cor}

\Proof  One can check that if $F$ is the pseudo Frobenius algebra
generating the walking pseudo Frobenius algebra, then the
semistrict 3-category $\cat{EM}_{\F}\big(\Sigma(\wFrob)\big)$ is
isomorphic to the walking ambidextrous pseudoadjunction.  Hence,
$\wFrob$ is $\Hom(A,A)$ in $\wAdj$. \qed

\subsection{Examples} \label{secExamples}

We now provide a brief survey of the literature where
higher-dimensional analogs of Frobenius algebras have appeared.

\paragraph{Trivial examples. }
Since every monoidal category can be regarded as a semistrict
monoidal 2-category with only identity 2-morphisms, every
Frobenius algebra is a pseudo Frobenius algebra with only identity
2-morphisms and trivial coherences. \medskip

\paragraph{*-autonomous categories. }
Note that a pseudo Frobenius algebra in the 2-category $\cat{Cat}$
is a weak monoidal category with some extra structure.  Street has
shown that the condition that a monoidal category be Frobenius is
equivalent to the condition that the monoidal category be
$*$-autonomous~\cite{StreetFrob}.  These $*$-autonomous monoidal
categories are known to have an interesting relationship with
quantum groups and quantum groupoids~\cite{DS2}. Combined with our
result relating Frobenius pseudomonoids to pseudo ambijunctions,
the relationship with $*$-autonomous categories may have
implications to quantum groups, as well as the field of linear
logic where $*$-autonomous categories are used extensively.
\medskip

\paragraph{Khovanov's Frobenius functors. }
In {\em A functor-valued invariant of tangles}, Khovanov sketches
a definition of a topological quantum field theory with corners
and suggests that the useful examples of these structures arise
from functors with a 2-sided adjoint, or what he calls `Frobenius
functors'~\cite{khovanov}.  In our language we would say that
these TQFT's with corners arise from ambijunctions in \cat{Cat}.
He then lists and describes in detail some categories with many
Frobenius functors, that is, he lists many examples of
ambijunctions in various 2-subcategories of the 2-category
\cat{Cat}. We repeat his list below although we will not describe
in detail the ambijunctions.  For details the reader is referred
to Khovanov's paper~\cite{khovanov}.

\begin{quote}
Categories with many Frobenius functors:
\end{quote}
\begin{itemize}
    \item Categories of modules over symmetric and Frobenius
    algebras and their derived categories.
    \item Categories of highest weight modules over simple Lie
    algebras and their derived categories.
    \item Derived categories of coherent sheaves on Calabi-Yau
    manifolds.
    \item Fuyaya-Floer categories of lagrangians in a symplectic
    manifold.
\end{itemize}
\medskip

\paragraph{Tillmann's Frobenius categories~\cite{Tillmann}}
Tillmann suggests that in order to encode 3-dimensional
information into a $2$-dimensional topological quantum field
theory one must consider a more interesting version of the
2-dimensional cobordism category, namely \cat{2Cob}$_2$.  The
objects of $\cat{2Cob}_2$ are closed, oriented, compact
$1$-manifolds, and the morphisms are oriented, compact
2-manifolds.  The 2-morphisms of this 2-category are the connected
components of orientation preserving diffeomorphisms of the
2-manifolds . This cobordism 2-category was first studied by
Carmody~\cite{Carmody}.

By extending the category $\cat{2Cob}$ to the 2-category
$\cat{2Cob}_2$, Tillmann defines a modular functor as a monoidal
2-functor from $\cat{2Cob}_2 \to$ $\Bbbk$-\cat{Cat}, where
$\Bbbk$-\cat{Cat} is the 2-category of linear categories, linear
functors and linear natural transformations. Tillmann calls the
image of such a 2-functor a `Frobenius category'.  She goes on to
show that these Frobenius categories are related to 3-dimensional
topological quantum field theory.  In our terminology, these
`Frobenius categories' are a symmetric version of pseudo Frobenius
algebras in the 2-category of $\Bbbk$-linear categories.
\medskip

This example is particularly related to the results of this paper.
In the next section we will discuss the nonsymmetric version of
the cobordism 2-category described above.  The effect of removing
the symmetry requirement amounts to `smashing the cobordisms flat'
into what we call $\cat{3Thick}$, the 2-category of 3-dimensional
thick tangles. This 2-category turns out to be an extension of the
category $\cat{2Thick}$ defined in Section \ref{sec2Dthick}.
Three-dimensional thick tangles are the most important example of
pseudo Frobenius algebras since we will also show that the
monoidal bicategory of 3-dimensional thick tangles is
triequivalent to the walking pseudo Frobenius algebra.  This means
that all of the examples above are the image of a
$\cat{Gray}$-functor from the monoidal bicategory of 3-dimensional
thick tangles into the relevant semistrict monoidal categories.

\section{Thick tangles} \label{secThick}

We now define an extension of the monoidal category of
2-dimensional thick tangles.  It is perhaps not surprising that
this extension will involve an extra level of categorical
structure.  In fact there are numerous examples of higher
categories providing an algebraic description generalizing various
kinds of algebraically defined topological categories. For
instance, just as the category of tangles in 3-dimensions is the
free braided monoidal category with duals on one object, Baez and
Langford have shown that the 2-category of 2-tangles in
4-dimensions is the free semistrict braided monoidal 2-category
with duals on one object~\cite{hdaIV}.  We will give yet another
example of this phenomenon by showing that the monoidal bicategory
of 3-dimensional thick tangles defined below is the walking pseudo
Frobenius algebra.  Note that a monoidal bicategory is just a one
object tricategory.  Hence, by the coherence theorem of
Gordon-Power-Street, every monoidal bicategory is triequivalent to
a semistrict monoidal 2-category.

\begin{defn} \et {\bf (Kerler and Lyubashenko~\cite{KL})}
The monoidal bicategory of {\em 3-dimensional thick tangles}
denoted $\cat{3Thick}$ has nonnegative integers as objects. The
1-morphisms from $k$ to $l$ are smooth oriented compact surfaces
$X$ with boundary $\partial X$ equipped with disjoint
distinguished intervals $i_j^s \maps I \hookrightarrow
\partial X$, $1 \leq j \leq k$, $i^t_m \maps \hookrightarrow
\partial X$, $1 \leq m \leq l$, equipped with a smooth embedding $d
\maps X \hookrightarrow \R \times [0,1]$ such that
\[
 d^{-1}(\R\times0) = I_1^s \sqcup I^s_2 \sqcup \cdots \sqcup
 I_k^s, \qquad   I^s_j=i^s_j(I) , \qquad d(I^s_j) = [j-\frac{1}{3},j+\frac{1}{3}] \times 0 ,\]
 \[
 d^{-1}(\R\times1) = I_1^t \sqcup I^t_2 \sqcup \cdots \sqcup I_k^t,
 \qquad  I^t_j=i^t_j(I) , \qquad  d(I^t_j) = [j-\frac{1}{3},j+\frac{1}{3}] \times 1.
\]
The image $d(X)$ is called a {\em diagram of thick tangles}.

The 2-morphisms $\phi \maps X \Rightarrow Y \maps k \to l$ of
\cat{3Thick} are isotopy classes of oriented homeomorphisms $f
\maps X \to Y$, which preserve the distinguished intervals, i.e.,
$i_{jY}^{s/t} = \big( \xymatrix@1{I \ar[r]^{i_{jY}^{s/t}} & X
\ar[r]^f & Y} \big)$. (Each homeomorphism $f_t \maps X \to Y$, $t
\in [0,1]$ in the isotopy family also preserves the distinguished
intervals.)  Composition $Y \circ X$ of 1-morphisms $\xymatrix@1{k
\ar[r]^X & l \ar[r]^Y & m}$ is defined by sewing of surfaces at
boundary intervals $I^t_j(X)$ and $I^s_j(Y)$.  The unit 1-morphism
$1_k \maps k \to k$ is the union $\coprod_{j=1}^k
[j-\frac{1}{3},j+\frac{1}{3}] \times [0,1]$.  The isomorphism
$\xymatrix@1{1_l \circ X  \ar[r]^{\sim} & X}$ are obtained by
taking a neighborhood $(U,I^s_j) \simeq \left([0,1]
\times[0,1],[0,1],[0,1] \times 0\right)$ of the distinguished
interval $I_j^s \subset X$ and by taking any isomorphism $[0,1]
\times [0,1]\bigcup_{[0,1]\times 1 \sim I^s_j}U \simeq U$.  The
tensor product is the disjoint union.  The unit object is 0.  The
associativity constraints are obvious.
\end{defn}

Here are some examples:
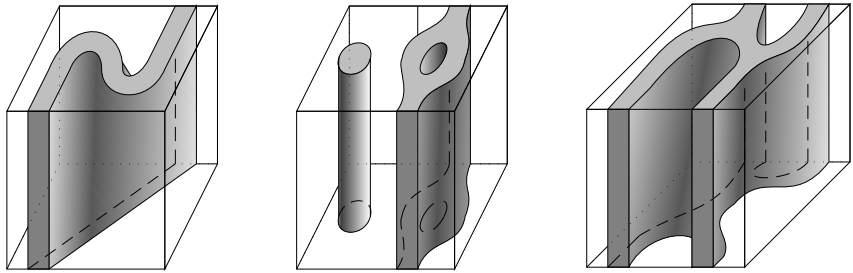
\begin{center} \makebox[0pt]{
$
 \xy
 (0,0)*{
   \begin{pspicture}(4,5.5)
   \pspolygon[fillcolor=lightgray,fillstyle=gradient,
    gradbegin=white, gradend=darkgray,gradmidpoint=.32,gradangle=88](1.9,4.4)(2.4,3.6)(2.4,3)(1.9,4.4)
\pspolygon(1,2)(1,5)(4,5)(4,2)(1,2) \psline(0,0)(1,2)
      \pscustom[fillcolor=lightgray,fillstyle=solid]{
  \psline(.4,3)(.8,3.8)
  \psbezier(.8,3.8)(1.2,4.8)(2.2,4.7)(2.2,3.8)
  \psbezier(2.2,3.8)(2.2,3.6)(2.5,3.6)(2.7,4)
  \psline(3.2,5)
  \psline(3.6,5)
  \psline(3.6,3)
  }
   \pscustom[fillcolor=lightgray,fillstyle=gradient,
    gradbegin=white, gradend=darkgray,gradmidpoint=.32,gradangle=88]{
  \psline(.8,3)(1.2,3.8)
  \psbezier(1.2,3.8)(1.35,4.2)(1.8,4.3)(1.8,3.8)
  \psbezier(1.8,3.8)(1.8,3.2)(2.7,3.1)(3,3.8)
  \psline(3,3.8)(3.6,5)
  \psline(3.6,2)
  \psline(.8,0)
  \psline(.8,3)
  }
  \pspolygon[fillcolor=gray,fillstyle=solid](.8,0)(.4,0)(.4,3)(.8,3)(.8,0)
    \pspolygon(0,0)(0,3)(3,3)(3,0)(0,0)
  \pspolygon(0,3)(3,3)(4,5.0)(1,5.0)(0,3)
  \pspolygon(3,0)(3,3)(4,5.0)(4,2)(3,0)
  \psline[linestyle=dashed](.4,0)(3.2,2)
  \psline[linestyle=dashed](3.2,4)(3.2,2)
  \psline[linestyle=dotted](4,2)(1,2)
  \psline[linestyle=dotted](1,5)(1,2)
  \psline[linestyle=dotted](0,0)(1,2)
 \end{pspicture}};
 \endxy
\qquad \quad \xy
 (0,0)*{
    \begin{pspicture}(4,5.5) 
\pspolygon(1,2)(1,5)(4,5)(4,2)(1,2) \psline(0,0)(1,2)
  \pscustom[fillcolor=lightgray,fillstyle=solid]{
  \psbezier(.8,4)(.9,4.4)(1.5,4.4)(1.4,4)
  \psbezier(1.4,4)(1.3,3.6)(.7,3.6)(.8,4)}
  \pscustom[fillcolor=lightgray,fillstyle=gradient,
    gradbegin=white, gradend=darkgray,gradmidpoint=.2,gradangle=90]{
    \psbezier(1.4,4)(1.3,3.6)(.7,3.6)(.8,4)
    \psline(.8,1)
  \psbezier(.8,1)(.7,.6)(1.3,.6)(1.4,1)
    \psline(1.4,4)}
  \pscustom[fillcolor=lightgray,fillstyle=gradient,
    gradbegin=lightgray, gradend=darkgray,gradmidpoint=.4,gradangle=88]{
  \psbezier(2.3,0)(2.5,.6)(3.1,.4)(3.2,1)
  \psbezier(3.2,1)(3.4,1.5)(3.1,1.5)(3.3,2)
  \psline(3.3,5)
  \psbezier(3.3,5)(3.1,4.5)(3.4,4.5)(3.2,4)
  \psbezier(3.2,4)(3.1,3.4)(2.5,3.6)(2.3,3)
  \psline(2.3,5)
  }
  \pscustom[fillcolor=lightgray,fillstyle=solid]{
  \psbezier(1.9,3)(2.2,3.4)(1.9,3.8)(2,4)
  \psbezier(2,4)(2.1,4.6)(2.8,4.5)(2.9,5)
  \psline(3.3,5)
    \psbezier(3.3,5)(3.1,4.5)(3.4,4.5)(3.2,4)
  \psbezier(3.2,4)(3.1,3.4)(2.5,3.6)(2.3,3)
  \psline(1.9,3)}
  \pscustom[fillcolor=lightgray,fillstyle=gradient,
    gradbegin=white,gradend=black,gradmidpoint=1,gradangle=90]{
  \psbezier(2.4,4)(2.2,3.6)(2.6,3.6)(2.8,4)
  \psbezier(2.8,4)(3,4.4)(2.6,4.4)(2.4,4)}
  \pspolygon[fillcolor=gray,
  fillstyle=solid](1.9,0)(2.3,0)(2.3,3)(1.9,3)(1.9,0)
  \psbezier[linestyle=dashed](.8,1)(.9,1.4)(1.5,1.4)(1.4,1)
  \psbezier[linestyle=dashed](2.4,1)(2.2,.6)(2.6,.6)(2.8,1)
  \psbezier[linestyle=dashed](2.8,1)(3,1.4)(2.6,1.4)(2.4,1)
  \psbezier[linestyle=dashed](1.9,0)(2.2,0.4)(1.9,0.8)(2,1)
  \psbezier[linestyle=dashed](2,1)(2.1,1.6)(2.8,1.5)(2.9,2)
  \psline[linestyle=dashed](2.9,2)(2.9,3.4)
  \psline[linestyle=dotted](4,2)(1,2)
  \psline[linestyle=dotted](1,5)(1,2)
  \psline[linestyle=dotted](0,0)(1,2)
    \pspolygon(0,0)(0,3)(3,3)(3,0)(0,0)
  \pspolygon(0,3)(3,3)(4,5.0)(1,5.0)(0,3)
  \pspolygon(3,0)(3,3)(4,5.0)(4,2)(3,0)
 \end{pspicture}};
 \endxy
 \qquad \quad
 \xy (0,0)*{};
 (0,2)*{
   \begin{pspicture}[.5](4,6)
\pspolygon(2,2)(2,5)(5,5)(5,2)(2,2) \psline(0,0)(2,2)
  \pspolygon(0,3)(3,3)(5,5.0)(2,5.0)(0,3)
\pscustom[fillcolor=lightgray,fillstyle=gradient,
    gradbegin=white, gradend=darkgray,gradmidpoint=.5,gradangle=80]{
    \psline(3.4,5)(3.4,4)
    \psline(3,4)
    \psline(3,5)
    \psline(3.4,5)
 }
   \pscustom[fillcolor=lightgray,fillstyle=solid]{
    \psline(.4,3)(1.3,4)
  \psbezier(1.3,4)(2.2,5)(2.6,4.3)(3,5)
  \psline(3.4,5)
      \psbezier(3.4,5)(2.8,3.6)(4,4.4)(4.2,5)
      \psline(4.6,5)
      \psbezier(4.6,5)(3.8,3.6)(3.6,4.4)(2.8,3.4)
      \psline(2.8,3.4)(2.4,3)
      \psline(.4,3)
  }
  %
  %
  \pscustom[fillcolor=lightgray,fillstyle=gradient,
    gradbegin=white, gradend=darkgray,gradmidpoint=.5,gradangle=80]{
    \psbezier(.8,0)(1.2,.8)(2.6,.8)(2,0)
    \psline(2,3)
      \psline(2,3)(2.8,3.8)
      \psbezier(2.8,3.8)(3.2,4.3)(2.4,4.9)(1.2,3.4)
      \psline(1.2,3.4)(.8,3)
      \psline(.8,0)
  }
  \pscustom[fillcolor=lightgray,fillstyle=gradient,
    gradbegin=white, gradend=darkgray,gradmidpoint=.32,gradangle=88]{
  \psline(2.4,0)(2.6,.2)
  \psbezier(2.6,.2)(2.8,.6)(2.2,.8)(2.8,1.2)
  \psbezier(2.8,1.2)(3.4,1.9)(3.9,1.1)(4.6,2)
  \psline(4.6,5)
  \psbezier(4.6,5)(3.8,3.6)(3.6,4.4)(2.8,3.4)
    \psline(2.8,3.4)(2.4,3)
    \psline(2.4,0)
  }
    \pspolygon[fillcolor=gray,fillstyle=solid](.8,0)(.4,0)(.4,3)(.8,3)(.8,0)
  \pspolygon[fillcolor=gray,fillstyle=solid](2.4,0)(2,0)(2,3)(2.4,3)(2.4,0)
  \psbezier[linestyle=dashed](4.2,2)(3.9,1.6)(2.6,1.6)(3.4,2)
  \psline[linestyle=dashed](.4,0)(1,.6)
  \psbezier[linestyle=dashed](1,.6)(2,1.2)(2.4,1)(3,2)
  \psline[linestyle=dashed](3,2)(3,3.4)
  \psline[linestyle=dashed](3.4,2)(3.4,3.7)
  \psline[linestyle=dashed](4.2,2)(4.2,4.3)
  \psline[linestyle=dotted](5,2)(2,2)
  \psline[linestyle=dotted](0,0)(2,2)
  \psline[linestyle=dotted](2,5)(2,2)
 %
    \pspolygon(0,0)(0,3)(3,3)(3,0)(0,0)
    \pspolygon(3,0)(3,3)(5,5.0)(5,2)(3,0)
 \end{pspicture}
 };
 \endxy$
 }
\end{center}
These surfaces can be interpreted as those topological membranes
that arise from diffeomorphisms of planar open string worldsheets.

We are now ready to state the main theorem, but first, note that
since a monoidal bicategory is just a one object tricategory, and
a semistrict monoidal 2-category can also be regarded as a one
object tricategory, it makes sense to talk about a triequivalence
between them.

\begin{thm} \label{MAINthm}
The monoidal bicategory of 3-dimensional thick tangles is
triequivalent to the walking pseudo Frobenius algebra.
\end{thm}

We will prove this theorem in Section~\ref{secPROOFS}, but first
we give a few corollaries.

\begin{cor}
A pseudo Frobenius object in a semistrict monoidal 2-category
determines an invariant of 3-dimensional thick tangles.
\end{cor}

\Proof This follows from the universal property of the monoidal
bicategory of thick tangles. \qed

\begin{cor}
The monoidal category of two-dimensional thick tangles is
equivalent to the walking Frobenius algebra.
\end{cor}

\Proof  This result follows as a decategorification of
Theorem~\ref{MAINthm}. \qed

If we define an extended planar open string topological field
theory as a monoidal 2-functor from $\cat{3Thick}$ into a monoidal
2-category $\mathcal{C}$, then Theorem~\ref{MAINthm} can be
rephrased as follows:

\begin{thm}
An extended planar open string topological field theory is
equivalent to a pseudo Frobenius algebra in the monoidal
2-category $\mathcal{C}$.
\end{thm}

\subsection{Some lemmas and related definitions} \label{secLemmas}

In this section we state a few lemmas due to Kerler and
Lyubashenko that provide a generators and relations description of
the monoidal bicategory of 3-dimensional thick tangles. Although
some of these definitions are quite long, we will see in the next
section that the semistrict monoidal 2-category with these
generators and relations is monoidally 2-equivalent to the
semistrict monoidal 2-category generated by a pseudo Frobenius
algebra.

In a sense, all of the definitions below can be viewed as
equivalent descriptions of the free semistrict monoidal 2-category
generated by a pseudo Frobenius algebra.

\begin{lem}[Kerler and Lyubashenko~\cite{KL}] \label{lemma1}
The monoidal bicategory of 3-dimensional thick tangles is
triequivalent to the semistrict monoidal 2-category
$\mathcal{F}_1$ generated by one object, morphisms:
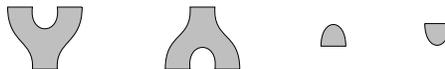
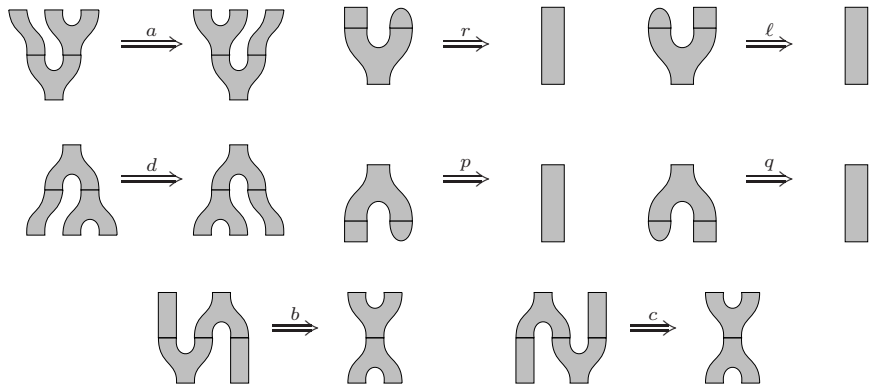
\[  \psset{xunit=.55cm,yunit=0.55cm}

\]
such that the following twenty diagrams:
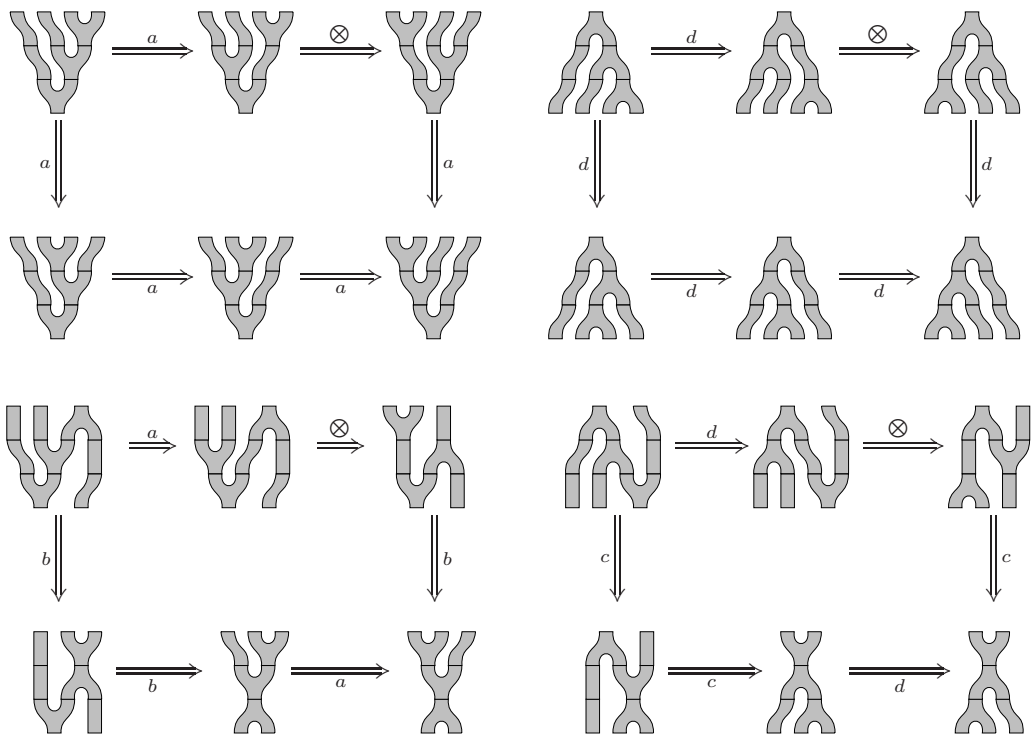
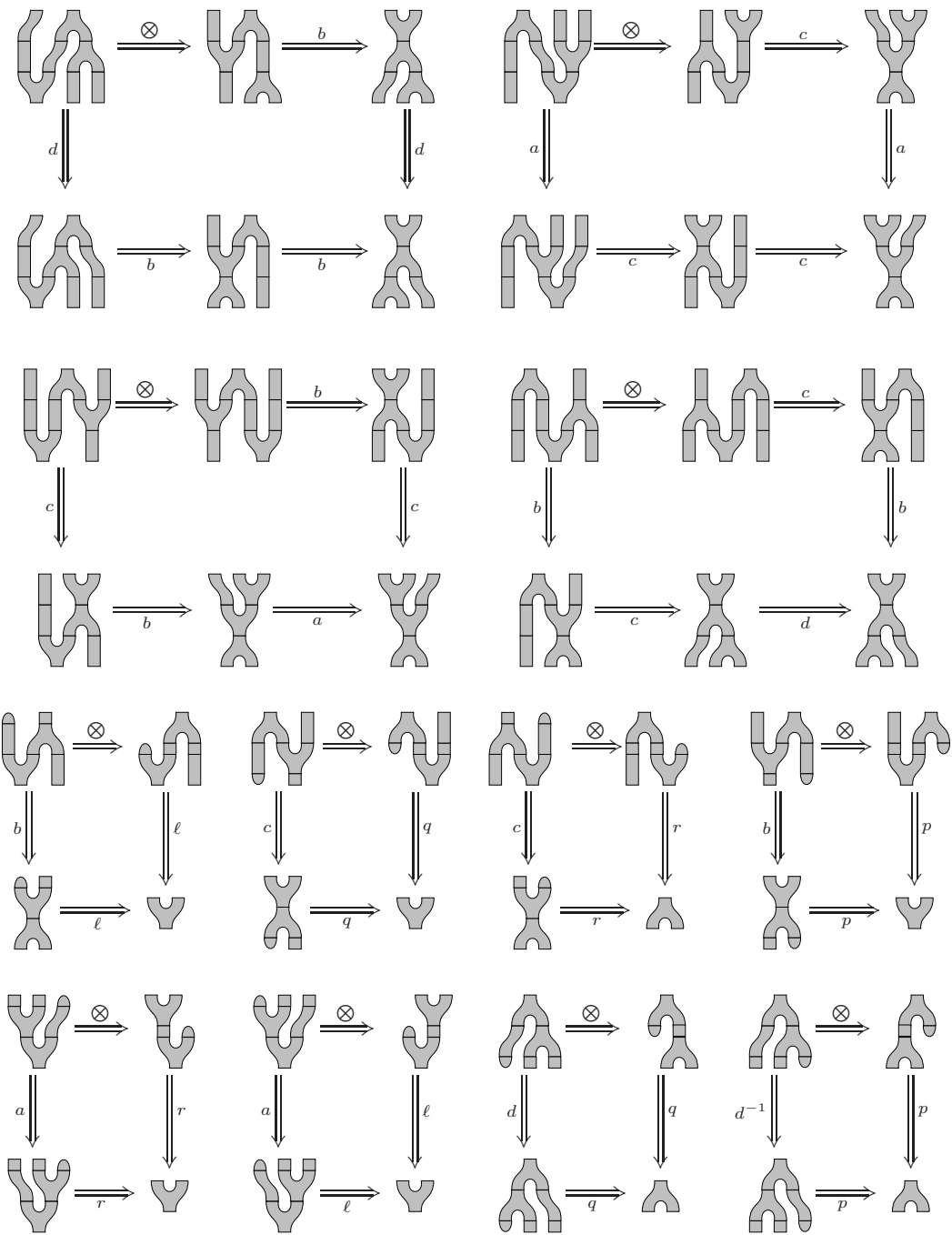
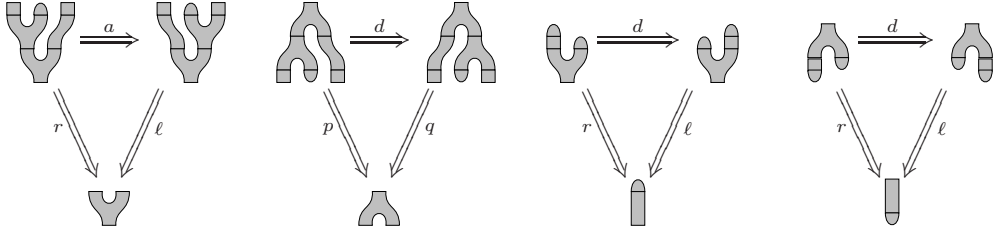
\begin{center}
 \makebox[0pt]{ $\psset{xunit=.30cm,yunit=0.30cm}
\xy
 (-25,15)*+{
    %
    }="B";
    {\ar@{=>}^{d} "TL";"TR"};
    {\ar@{=>}_{r} "TL";"B"};
    {\ar@{=>}^{\ell} "TR";"B"};
 \endxy
 $ }
\end{center}
commute.
\end{lem}

\Proof See Kerler and Lyubashenko~\cite{KL}.\qed

Note that this triequivalence is not a triequivalence in its
weakest form.  First of all, $\cat{3Thick}$ and $\mathcal{F}_1$
have only one object viewed as tricategories.  Hence, the
trifunctors defining this triequivalence are the identity on
objects.  Furthermore, since both 2-categories are freely
generated by one object it is clear that these trifunctors must
also be the identity on 1-morphisms.  Thus, we are not using the
term `triequivalence' in its weakest form.

This definition should remind the reader of the definition of
Frobenius algebra given in Proposition~\ref{equivFrob} (i.).
Notice that all of the axioms that held as equalities are now
replaced by coherent isomorphisms.  It might be described as a
pseudomonoid and pseudocomonoid satisfying the Frobenius
identities up to coherent isomorphism.  Notice also the tremendous
number of coherence laws for this definition. Since the coherence
laws describe how each of generating 2-morphisms behave with
respect to each other, the large number of generating 2-morphisms
means a large number of coherence conditions.

In the next lemma we show that the definition given above is
monoidally 2-equivalent to a definition with a few less coherence
conditions.  To clarify what is meant by monoidal 2-equivalence,
let $\mathcal{D}$ and $\mathcal{D}'$ be monoidally 2-equivalent
categories.  We can regard $\mathcal{D}$ and $\mathcal{D}'$ as one
object $\cat{Gray}$-categories.  In this case, the monoidal
2-equivalence translates to a 3-equivalence of
$\cat{Gray}$-categories.  Hence, we will sometimes say a
3-equivalence of semistrict monoidal 2-categories by which we
simply mean a 3-equivalence of their respective suspensions.

\begin{lem}[Kerler and Lyubashenko~\cite{KL}] \label{lemma2}
The semistrict monoidal 2-category $\mathcal{F}_1$ from Lemma
\ref{lemma1} is 3-equivalent  to the semistrict monoidal
2-category $\mathcal{F}_2$ generated by one object, morphisms:
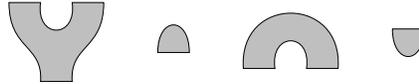
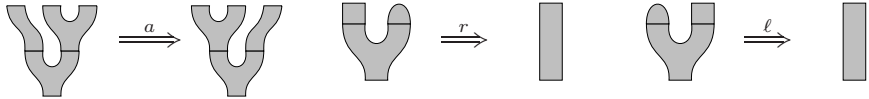
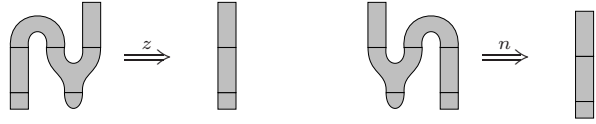
\[ \psset{xunit=.7cm,yunit=0.7cm}
\begin{pspicture}[0.5](1.2,1.6)
  \rput(.6,0){\mult}
\end{pspicture}
\qquad
\begin{pspicture}[0.5](1.2,1.6)
  \rput(.6,.55){\birth}
\end{pspicture}
\qquad
\begin{pspicture}[0.5](1.2,1.6)
  \rput(.6,.25){\zag}
\end{pspicture}
\qquad
\begin{pspicture}[0.5](1,1.6)
  \rput(.6,1){\death}
\end{pspicture}
\]
and 2-isomorphisms:
\[ \psset{xunit=.4cm,yunit=0.4cm}   
\begin{pspicture}[0.5](3,3.3)
  \rput(1,0){\mult} \rput(1.6,1.5){\mult} \rput(.4,1.5){\curveleft}
\end{pspicture}
 \; \xy {\ar@2{->}^{a} (0,0)*{};(8,0)*{}}; \endxy  \quad
\begin{pspicture}[0.5](2.6,3.3)
  \rput(1,0){\mult} \rput(1.6,1.5){\curveright} \rput(.4,1.5){\mult}
\end{pspicture}
\qquad  \psset{xunit=.5cm,yunit=0.5cm}
\begin{pspicture}[0.3](2,3)
  \rput(1,0){\mult} \rput(1.6,1.5){\birth} \rput(.4,1.5){\smallident}
\end{pspicture}
 \quad \xy {\ar@2{->}^{r} (0,0)*{};(6,0)*{}}; \endxy  \quad
\begin{pspicture}[0.3](2,3)
  \rput(1,0){\longident}
\end{pspicture}
\qquad\begin{pspicture}[0.3](2,3)
  \rput(1,0){\mult} \rput(1.6,1.5){\smallident} \rput(.4,1.5){\birth}
\end{pspicture}
\quad \xy {\ar@2{->}^{\ell} (0,0)*{};(6,0)*{}}; \endxy  \quad
\begin{pspicture}[0.3](2,3)
  \rput(1,0){\longident}
\end{pspicture}
\]
\[
 \psset{xunit=.4cm,yunit=0.4cm}
\begin{pspicture}[.5](3,3)
  \rput(2.2,.55){\mult} \rput(1,2.05){\zag}
  \rput(2.8,2.05){\ident} \rput(.4,0){\smallident} \rput(.4,.55){\ident}
  \rput(2.2,.55){\death}
\end{pspicture}
\quad \xy {\ar@2{->}^{z} (0,0)*{};(6,0)*{}}; \endxy  \quad
\begin{pspicture}[0.5](2,3)
  \rput(1,2.05){\ident} \rput(1,0){\smallident} \rput(1,.55){\ident}
\end{pspicture}
\qquad \qquad
\begin{pspicture}[0.5](3,3)
  \rput(1,.55){\mult} \rput(2.2,2.05){\zag}
  \rput(.4,2.05){\ident} \rput(2.8,.55){\ident}
  \rput(2.8,0){\smallident} \rput(1,.55){\death}
\end{pspicture}
\quad \xy {\ar@2{->}^{n} (0,0)*{};(6,0)*{}}; \endxy  \quad
\begin{pspicture}[0.5](2,3.6)
  \rput(1,2.05){\ident} \rput(1,.55){\ident} \rput(1,0){\smallident}
\end{pspicture}
\]
satisfying the pseudomonoid axioms and making the following
diagrams commute:
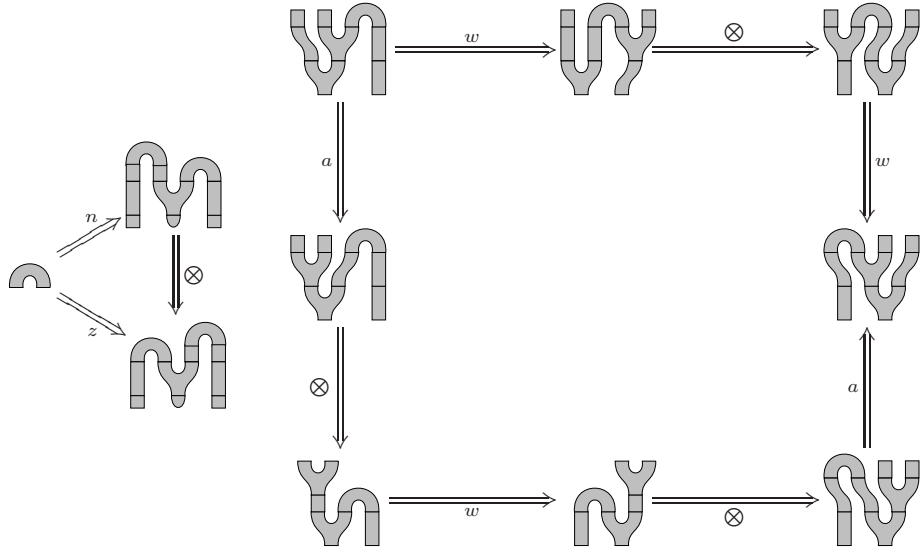
\begin{center}
 \makebox[0pt]{ $\psset{xunit=.30cm,yunit=0.30cm}
 \xy
 (-10,0)*+{
    \begin{pspicture}[0.5](2.4,1.1)
  \rput(1.4,0){\zag}
 \end{pspicture}
    }="TL";
(10,12)*+{
   \begin{pspicture}[0.5](4.2,3.8)
  \rput(2,.55){\mult}
  \rput(2,.55){\death}
  \rput(3.2,2.05){\zag}
  \rput(1.4,2.05){\medident}
  \rput(.8,2.75){\zag}
  \rput(.2,2.05){\medident}
  \rput(.2,.55){\ident}
  \rput(.2,0){\smallident}
  \rput(3.8,0){\smallident}
  \rput(3.8,.55){\ident}
 \end{pspicture}
    }="TR";
 (10,-12)*+{
   \begin{pspicture}[0.5](3.8,3.8)
  \rput(2,.55){\mult}
  \rput(2,.55){\death}
  \rput(.8,2.05){\zag}
  \rput(2.6,2.05){\medident}
  \rput(3.2,2.75){\zag}
  \rput(3.8,2.05){\medident}
  \rput(.2,.55){\ident}
  \rput(.2,0){\smallident}
  \rput(3.8,0){\smallident}
  \rput(3.8,.55){\ident}
 \end{pspicture}
    }="BR";
    {\ar@{=>}^{n} "TL";"TR"};
    {\ar@{=>}^{\INT} "TR";"BR"};
    {\ar@{=>}_{z} "TL";"BR"};
 \endxy
\qquad
\xy
 (-35,-30)*+{
   \begin{pspicture}[0.5](3.6,4)
  \rput(1.4,0){\mult}
  \rput(3.2,0){\ident}
  \rput(.8,1.5){\medident}
  \rput(2.6,1.5){\zag}
  \rput(.8,2.25){\mult}
 \end{pspicture}
    }="LB";
 (-35,0)*+{
   \begin{pspicture}[0.5](4.2,4)
  \rput(1.4,0){\mult}
  \rput(3.8,0){\ident}
  \rput(.8,1.5){\mult}
  \rput(2,1.5){\curveright}
  \rput(3.8,1.5){\ident}
  \rput(.2,3){\medident}
  \rput(1.4,3){\medident}
  \rput(3.2,3){\zag}
 \end{pspicture}
 }="LM";
 (-35,30)*+{
   \begin{pspicture}[0.5](4.2,4)
  \rput(1.4,0){\mult}
  \rput(3.8,0){\ident}
  \rput(2,1.5){\mult}
  \rput(.8,1.5){\curveleft}
  \rput(3.8,1.5){\ident}
  \rput(.2,3){\medident}
  \rput(1.4,3){\medident}
  \rput(3.2,3){\zag}
 \end{pspicture}
 }="LT";
 (0,-30)*+{
   \begin{pspicture}[0.5](3.6,4)
  \rput(2.6,0){\mult}
  \rput(.8,0){\ident}
  \rput(3.2,1.5){\medident}
  \rput(1.4,1.5){\zag}
  \rput(3.2,2.25){\mult}
 \end{pspicture}
 }="MB";
 (0,30)*+{
   \begin{pspicture}[0.5](3.6,4)
  \rput(.8,0){\mult}
  \rput(2.6,0){\curveright}
  \rput(3.2,1.5){\mult}
  \rput(.2,1.5){\ident}
  \rput(1.4,1.5){\ident}
  \rput(3.8,3){\medident}
  \rput(.2,3){\medident}
  \rput(2,3){\zag}
 \end{pspicture}
 }="MT";
 (35,-30)*+{
   \begin{pspicture}[0.5](3.6,4)
  \rput(2.6,0){\mult}
  \rput(.8,0){\ident}
  \rput(3.2,1.5){\mult}
  \rput(2,1.5){\curveleft}
  \rput(.8,1.5){\curveleft}
  \rput(3.8,3){\medident}
  \rput(2.6,3){\medident}
  \rput(.8,3){\zag}
 \end{pspicture}
 }="RB";
 (35,0)*+{
   \begin{pspicture}[0.5](3.6,4)
  \rput(2.6,0){\mult}
  \rput(.8,0){\ident}
  \rput(3.2,1.5){\curveright}
  \rput(2,1.5){\mult}
  \rput(.8,1.5){\curveleft}
  \rput(3.8,3){\medident}
  \rput(2.6,3){\medident}
  \rput(.8,3){\zag}
 \end{pspicture}
  }="RM";
 (35,30)*+{
   \begin{pspicture}[0.5](3.6,4)
  \rput(2.6,0){\mult}
  \rput(.8,0){\ident}
  \rput(3.2,1.5){\curveright}
  \rput(2,1.5){\curveright}
  \rput(.8,1.5){\mult}
  \rput(3.8,3){\medident}
  \rput(.2,3){\medident}
  \rput(2,3){\zag}
 \end{pspicture}
  }="RT";
    {\ar@{=>}_{a} "LT";"LM"};
    {\ar@{=>}_{\INT} "LM";"LB"};
    {\ar@{=>}^{w} "RT";"RM"};
    {\ar@{=>}^{a} "RB";"RM"};
    {\ar@{=>}^{w} "LT";"MT"};
    {\ar@{=>}^{\INT} "MT";"RT"};
    {\ar@{=>}_{w} "LB";"MB"};
    {\ar@{=>}_{\INT} "MB";"RB"};
 \endxy
 $ }
\end{center}
where $w$ is the following composite:

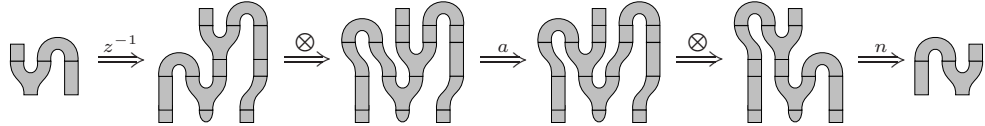
\begin{center}
 \makebox[0pt]{ $\psset{xunit=.30cm,yunit=0.30cm}
   \begin{pspicture}[0.5](3.6,2.55)
  \rput(1.4,0){\mult}
  \rput(3.2,0){\ident}
  \rput(2.6,1.5){\zag}
  \rput(.8,1.5){\medident}
 \end{pspicture}
 \;\; \xy {\ar@{=>}^{\scriptstyle z^{-1}}(-3,0);(3,0)};\endxy \;\;
   \begin{pspicture}[0.5](4.8,5)
  \rput(2,.55){\death}
  \rput(.2,0){\medident}
  \rput(.2,.55){\ident}
  \rput(2,.55){\mult}
  \rput(2.6,2.05){\medident}
  \rput(2.6,2.8){\mult}
  \rput(.8,2.05){\zag}
  \rput(3.8,4.3){\zag}
  \rput(2,4.3){\medident}
  \rput(4.4,2.8){\ident}
  \rput(4.4,2.05){\medident}
  \rput(3.8,.55){\curveright}
  \rput(3.8,0){\smallident}
 \end{pspicture}
 \;\; \xy {\ar@{=>}^{\scriptstyle \INT}(-3,0);(3,0)};\endxy \;\;
   \begin{pspicture}[0.5](5.4,5)
  \rput(2.6,.55){\death}
  \rput(.8,0){\medident}
  \rput(.8,.55){\ident}
  \rput(2.6,.55){\mult}
  \rput(3.8,3.55){\medident}
  \rput(2.6,3.55){\medident}
  \rput(3.2,2.05){\mult}
  \rput(.8,3.55){\zag}
  \rput(.8,2.05){\curveleft}
  \rput(2,2.05){\curveleft}
  \rput(4.4,4.3){\zag}
  \rput(2.6,4.3){\medident}
  \rput(5,3.55){\medident}
  \rput(5,2.05){\ident}
  \rput(4.4,.55){\curveright}
  \rput(4.4,0){\smallident}
 \end{pspicture}
 \;\; \xy {\ar@{=>}^{\scriptstyle a}(-3,0);(3,0)};\endxy \;\;
   \begin{pspicture}[0.5](5.4,5)
  \rput(2.6,.55){\death}
  \rput(.8,0){\medident}
  \rput(.8,.55){\ident}
  \rput(2.6,.55){\mult}
  \rput(3.8,3.55){\medident}
  \rput(2.6,3.55){\medident}
  \rput(2,2.05){\mult}
  \rput(.8,3.55){\zag}
  \rput(.8,2.05){\curveleft}
  \rput(3.2,2.05){\curveright}
  \rput(4.4,4.3){\zag}
  \rput(2.6,4.3){\medident}
  \rput(5,3.55){\medident}
  \rput(5,2.05){\ident}
  \rput(4.4,.55){\curveright}
  \rput(4.4,0){\smallident}
 \end{pspicture}
 \;\; \xy {\ar@{=>}^{\scriptstyle \INT}(-3,0);(3,0)};\endxy \;\;
   \begin{pspicture}[0.5](5.4,5)
  \rput(2.6,.55){\death}
  \rput(.8,0){\medident}
  \rput(.8,.55){\ident}
  \rput(2.6,.55){\mult}
  \rput(2.6,4.3){\medident}
  \rput(2,2.05){\medident}
  \rput(2,2.8){\mult}
  \rput(.8,4.3){\zag}
  \rput(.8,2.8){\curveleft}
  \rput(.8,2.05){\medident}
  \rput(3.8,2.05){\zag}
  \rput(4.4,0){\smallident}
  \rput(4.4,.55){\ident}
 \end{pspicture}
   \xy {\ar@{=>}^{\scriptstyle n}(-3,0);(3,0)};\endxy
   \begin{pspicture}[0.5](3.6,2.55)
  \rput(2.6,0){\mult}
  \rput(.8,0){\ident}
  \rput(1.4,1.5){\zag}
  \rput(3.2,1.5){\medident}
 \end{pspicture}
 $ }
\end{center}
\end{lem}

Notice that this definition is reminiscent of the definition of
Frobenius algebra given in Proposition~\ref{equivFrob} (ii.).
Notice in particular relation of the isomorphism $w$ to the axioms
in Proposition \ref{equivFrob} (ii.).  This definition might be
described as a pseudomonoid equipped with a form and a copairing
defining two coherently isomorphic comultiplications and counits.
In this case, we have been careful to use a minimum number of
generating 2-morphisms in order to minimize the number of
coherence axioms.

To simplify diagrams, we will from now on omit inverses on the
labels for the 2-morphisms.  The correct label should be apparent
from the diagram.

\paragraph{Proof. } We define a strict 2-functor $\Theta \maps \mathcal{F}_1
\to \mathcal{F}_2$ as follows:
\begin{itemize}
 \item On objects $\Theta$ is the identity map.
 \item On the generating morphisms $\Theta$ is given as follows:
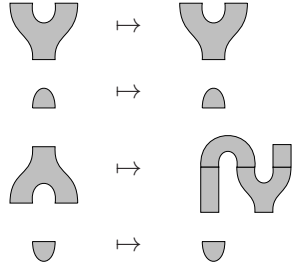
\begin{eqnarray*}
 \psset{xunit=.5cm,yunit=0.5cm}
\begin{pspicture}[.5](1,1.5)
 \rput(.5,0){\mult}
\end{pspicture}
 & \quad \mapsto \quad &
  \psset{xunit=.5cm,yunit=0.5cm}
\begin{pspicture}[.5](1,1.5)
 \rput(.5,0){\mult}
\end{pspicture}
\\
 \psset{xunit=.5cm,yunit=0.5cm}
\begin{pspicture}[.5](1,1.5)
 \rput(.5,.5){\birth}
\end{pspicture}
 & \quad \mapsto \quad &
  \psset{xunit=.5cm,yunit=0.5cm}
\begin{pspicture}[.5](1,1.5)
 \rput(.5,.5){\birth }
\end{pspicture}
\\
 \psset{xunit=.5cm,yunit=0.5cm}
\begin{pspicture}[.5](1,1.5)
 \rput(.5,0){\comult}
\end{pspicture}
 & \quad \mapsto \quad &
  \psset{xunit=.4cm,yunit=0.4cm}
\begin{pspicture}[.5](3.2,2.5)
 \rput(2.3,0){\mult} \rput(2.9,1.5){\medident}\rput(1.1,1.5){\zag} \rput(.5,0){\ident}
\end{pspicture}
\\
\psset{xunit=.5cm,yunit=0.5cm}
\begin{pspicture}[.5](1,1.5)
 \rput(.5,1){\death}
\end{pspicture}
 & \quad \mapsto \quad &
  \psset{xunit=.5cm,yunit=0.5cm}
\begin{pspicture}[.5](1,1.5)
 \rput(.5,1){\death}
\end{pspicture}
\end{eqnarray*}
   \item On the coherence 2-isomorphisms in $\mathcal{F}_1$ given by the commutative
   diagram:
\[
 \xy
    (-14,9)*+{A \ten B}="TL";
    (14,9)*+{A \ten B'}="TR";
    (14,-9)*+{A' \ten B'}="BR";
    (-14,-9)*+{A' \ten B}="BL";
    {\ar^{ A \ten g} "TL";"TR"};
    {\ar^{ f \ten B'} "TR";"BR"};
    {\ar_{ f \ten B} "TL";"BL"};
    {\ar_{ A' \ten g} "BL";"BR"};
    (0,4)*{}="x";
    (0,-4)*{}="y";
    {\ar@{=>}^{\INT_{f,g}} "x";"y"};
 \endxy
\]
$\Theta$ maps $\INT_{f,g}$ to the 2-cell given by the diagram:
\[
 \xy
    (-22,9)*+{\Theta(A) \ten \Theta(B)}="TL";
    (22,9)*+{\Theta(A) \ten \Theta(B')}="TR";
    (22,-9)*+{\Theta(A') \ten \Theta(B')}="BR";
    (-22,-9)*+{\Theta(A') \ten \Theta(B)}="BL";
    {\ar^-{\Theta( A) \ten \Theta(g)} "TL";"TR"};
    {\ar^{\Theta( f) \ten \Theta(B')} "TR";"BR"};
    {\ar_-{\Theta( f) \ten \Theta(B)} "TL";"BL"};
    {\ar_{\Theta( A') \ten \Theta(g)} "BL";"BR"};
    (0,4)*{}="x";
    (0,-4)*{}="y";
    {\ar@{=>}^{\INT_{\Theta(f),\Theta(g)}} "x";"y"};
 \endxy
\]
This choice can consistently be made by the coherence axioms of a
$\cat{Gray}$-category.
    \item
On the generating 2-isomorphisms $\Theta$ is defined as follows: 
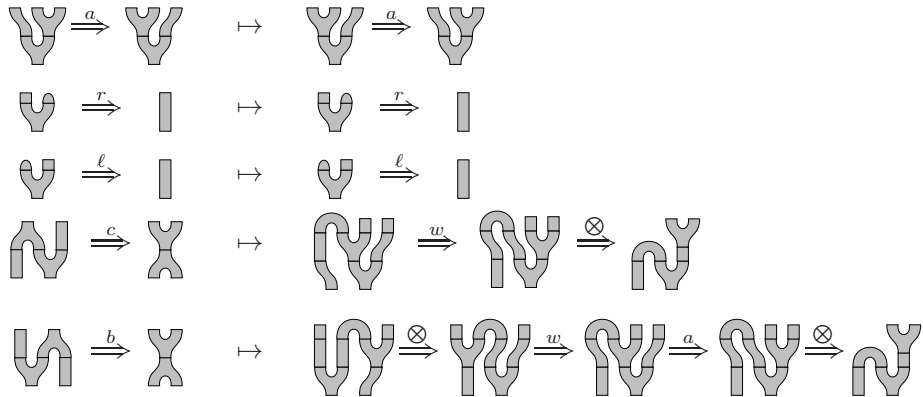
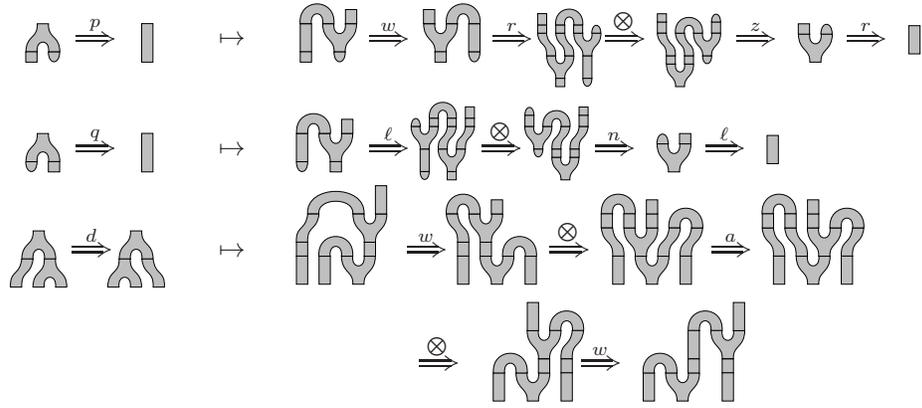
\begin{eqnarray*}
\psset{xunit=.25cm,yunit=0.25cm}
\begin{pspicture}[0.5](2.4,3.3)
  \rput(1,0){\mult} \rput(1.6,1.5){\mult} \rput(.4,1.5){\curveleft}
\end{pspicture}
 \; \xy {\ar@2{->}^{a} (0,0)*{};(5,0)*{}}; \endxy  \quad
\begin{pspicture}[0.5](2.6,3.3)
  \rput(1,0){\mult} \rput(1.6,1.5){\curveright} \rput(.4,1.5){\mult}
\end{pspicture}
&\quad \mapsto \quad&
\psset{xunit=.25cm,yunit=0.25cm}
\begin{pspicture}[0.5](2.6,3.3)
  \rput(1,0){\mult} \rput(1.6,1.5){\curveright} \rput(.4,1.5){\mult}
\end{pspicture}
 \; \xy {\ar@2{->}^{a} (0,0)*{};(5,0)*{}}; \endxy  \quad
\begin{pspicture}[0.5](3,3.3)
  \rput(1,0){\mult} \rput(1.6,1.5){\mult} \rput(.4,1.5){\curveleft}
\end{pspicture}
\\
\psset{xunit=.25cm,yunit=0.25cm}
\begin{pspicture}[0.3](2,3)
  \rput(1,0){\mult} \rput(1.6,1.5){\birth} \rput(.4,1.5){\smallident}
\end{pspicture}
 \quad \xy {\ar@2{->}^{r} (0,0)*{};(5,0)*{}}; \endxy  \quad
\begin{pspicture}[0.3](2,3)
  \rput(1,0){\longident}
\end{pspicture}
&\quad \mapsto \quad&
\psset{xunit=.25cm,yunit=0.25cm}
\begin{pspicture}[0.3](2,3)
  \rput(1,0){\mult} \rput(1.6,1.5){\birth} \rput(.4,1.5){\smallident}
\end{pspicture}
 \quad \xy {\ar@2{->}^{r} (0,0)*{};(5,0)*{}}; \endxy  \quad
\begin{pspicture}[0.3](2,3)
  \rput(1,0){\longident}
\end{pspicture}
\\
\psset{xunit=.25cm,yunit=0.25cm}
\begin{pspicture}[0.3](2,3)
  \rput(1,0){\mult} \rput(1.6,1.5){\smallident} \rput(.4,1.5){\birth}
\end{pspicture}
\quad \xy {\ar@2{->}^{\ell} (0,0)*{};(5,0)*{}}; \endxy  \quad
\begin{pspicture}[0.3](2,3)
  \rput(1,0){\longident}
\end{pspicture}
&\quad \mapsto \quad&
\psset{xunit=.25cm,yunit=0.25cm}
\begin{pspicture}[0.3](2,3)
  \rput(1,0){\mult} \rput(1.6,1.5){\smallident} \rput(.4,1.5){\birth}
\end{pspicture}
\quad \xy {\ar@2{->}^{\ell} (0,0)*{};(5,0)*{}}; \endxy  \quad
\begin{pspicture}[0.3](2,3)
  \rput(1,0){\longident}
\end{pspicture}
\\
\psset{xunit=.25cm,yunit=0.25cm}
\begin{pspicture}[1](3,3)
 \rput(1,3){\comult} \rput(2.2,1.5){\mult}
 \rput(2.8,3){\ident} \rput(.4,1.5){\ident}
\end{pspicture}
  \quad \xy {\ar@{=>}^{c} (0,0);(5,0)}; \endxy \quad
\begin{pspicture}[0.5](1.5,3)
 \rput(.5,0){\comult} \rput(.5,1.5){\mult}
\end{pspicture}
&\quad \mapsto \quad&
 \psset{xunit=.25cm,yunit=0.25cm}
\begin{pspicture}[0.5](4,4.2)
  \rput(2.6,0){\mult}
  \rput(3.2,1.5){\curveright}
  \rput(2,1.5){\mult}
  \rput(.2,1.5){\ident}
  \rput(.8,0){\curveleft}
  \rput(.8,3){\zag}
  \rput(2.6,3){\medident}
  \rput(3.8,3){\medident}
\end{pspicture}
 \quad \xy {\ar@{=>}^{w} (0,0);(5,0)}; \endxy \quad
    \begin{pspicture}[0.5](3.6,4)
  \rput(2.6,0){\mult}
  \rput(.8,0){\ident}
  \rput(3.2,1.5){\mult}
  \rput(2,1.5){\curveleft}
  \rput(.8,1.5){\curveleft}
  \rput(3.8,3){\medident}
  \rput(2.6,3){\medident}
  \rput(.8,3){\zag}
 \end{pspicture}
 \quad \xy {\ar@{=>}^{\INT} (0,0);(5,0)}; \endxy \;
\begin{pspicture}[0.5](4,4.2)
  \rput(2.6,0){\mult}
  \rput(3.2,1.5){\medident}
  \rput(3.2,2.25){\mult}
  \rput(.8,0){\ident}
  \rput(1.4,1.5){\zag}
\end{pspicture}
\\
\psset{xunit=.25cm,yunit=0.25cm}
\begin{pspicture}[1](3,3)
 \rput(2.4,3){\comult} \rput(1.2,1.5){\mult}
 \rput(.6,3){\ident} \rput(3,1.5){\ident}
\end{pspicture}
  \quad \xy {\ar@{=>}^b (0,0);(5,0)}; \endxy \quad
\begin{pspicture}[0.5](1.5,3)
 \rput(.5,0){\comult} \rput(.5,1.5){\mult}
\end{pspicture}
&\quad \mapsto \quad&
\psset{xunit=.25cm,yunit=0.25cm}
   \begin{pspicture}[0.5](3.6,4)
  \rput(.8,0){\mult}
  \rput(2.6,0){\curveright}
  \rput(3.2,1.5){\mult}
  \rput(.2,1.5){\ident}
  \rput(1.4,1.5){\ident}
  \rput(3.8,3){\medident}
  \rput(.2,3){\medident}
  \rput(2,3){\zag}
 \end{pspicture}
\;\; \xy {\ar@{=>}^{\INT} (0,0);(5,0)}; \endxy \;\;
   \begin{pspicture}[0.5](3.6,4)
  \rput(.8,0){\ident}
  \rput(2.6,0){\mult}
  \rput(2,1.5){\curveright}
  \rput(3.2,1.5){\curveright}
  \rput(.8,1.5){\mult}
  \rput(3.8,3){\medident}
  \rput(.2,3){\medident}
  \rput(2,3){\zag}
 \end{pspicture}
\;\;\xy {\ar@{=>}^{w} (0,0);(5,0)}; \endxy \;\;
    \begin{pspicture}[0.5](3.6,4)
  \rput(2.6,0){\mult}
  \rput(.8,0){\ident}
  \rput(3.2,1.5){\curveright}
  \rput(2,1.5){\mult}
  \rput(.8,1.5){\curveleft}
  \rput(3.8,3){\medident}
  \rput(2.6,3){\medident}
  \rput(.8,3){\zag}
 \end{pspicture}
\;\; \xy {\ar@{=>}^{a} (0,0);(5,0)}; \endxy \;\;
    \begin{pspicture}[0.5](3.6,4)
  \rput(2.6,0){\mult}
  \rput(.8,0){\ident}
  \rput(3.2,1.5){\mult}
  \rput(2,1.5){\curveleft}
  \rput(.8,1.5){\curveleft}
  \rput(3.8,3){\medident}
  \rput(2.6,3){\medident}
  \rput(.8,3){\zag}
 \end{pspicture}
 \; \; \xy {\ar@{=>}^{\INT} (0,0);(5,0)}; \endxy
      \begin{pspicture}[0.5](3.6,4)
  \rput(2.6,0){\mult}
  \rput(.8,0){\ident}
  \rput(3.2,1.5){\medident}
  \rput(1.4,1.5){\zag}
  \rput(3.2,2.25){\mult}
 \end{pspicture}
 \\
 \psset{xunit=.25cm,yunit=0.25cm}
\begin{pspicture}[0.3](2,3)
  \rput(1,.55){\comult} \rput(1.6,.55){\death} \rput(.4,.0){\smallident}
\end{pspicture}
\;\; \xy {\ar@2{->}^{p} (0,0)*{};(5,0)*{}}; \endxy  \;\;
\begin{pspicture}[0.3](2,3)
  \rput(1,0){\longident}
\end{pspicture}
&\quad \mapsto \quad&
\psset{xunit=.25cm,yunit=0.25cm}
\begin{pspicture}[0.3](3,3)
   \rput(2.2,.55){\mult}
   \rput(.4,.55){\ident}
   \rput(2.2,.55){\death}
   \rput(1,2.05){\zag}
   \rput(2.8,2.05){\medident}
   \rput(.4,0){\smallident}
\end{pspicture}
\;\; \xy {\ar@2{->}^{w} (0,0)*{};(5,0)*{}}; \endxy  \;\;
\begin{pspicture}[0.3](3,3)
   \rput(1,.55){\mult}
   \rput(2.8,.55){\ident}
   \rput(2.8,.55){\death}
   \rput(2.2,2.05){\zag}
   \rput(.4,2.05){\medident}
   \rput(1,0){\smallident}
\end{pspicture}
\;\; \xy {\ar@2{->}^{r} (0,0)*{};(5,0)*{}}; \endxy  \;\;
\psset{xunit=.2cm,yunit=0.2cm}
\begin{pspicture}[0.5](3,5.5)
   \rput(1,.55){\mult}
   \rput(2.8,.55){\ident}
   \rput(2.8,2.05){\mult}
   \rput(1.6,2.05){\curveleft}
   \rput(.4,2.05){\curveleft}
   \rput(1,3.55){\smallident}
   \rput(2.2,3.55){\smallident}
   \rput(-.2,3.55){\smallident}
   \rput(3.4,3.55){\birth}
   \rput(2.8,.55){\death}
   \rput(1.6,4.1){\zag}
   \rput(-.2,4.1){\medident}
   \rput(1,0){\smallident}
\end{pspicture}
\;\; \xy {\ar@2{->}^{\INT} (0,0)*{};(5,0)*{}}; \endxy  \;\;
\psset{xunit=.2cm,yunit=0.2cm}
\begin{pspicture}[0.5](4.2,5.5)
   \rput(1.4,0){\mult}
   \rput(.8,1.5){\smallident}
   \rput(.2,3.55){\medident}
   \rput(.2,4.3){\smallident}
   \rput(3.8,3.55){\medident}
   \rput(3.8,4.3){\birth}
   \rput(2,1.5){\smallident}
   \rput(.8,2.05){\curveleft}
   \rput(2,2.05){\curveleft}
   \rput(3.2,2.05){\death}
   \rput(3.2,2.05){\mult}
   \rput(2,3.55){\zag}
\end{pspicture}
\;\; \xy {\ar@2{->}^{z} (0,0)*{};(5,0)*{}}; \endxy  \;\;
\psset{xunit=.25cm,yunit=0.25cm}
\begin{pspicture}[0.3](2.3,3)
   \rput(1.4,0){\mult}
   \rput(.8,1.5){\smallident}
   \rput(2,1.5){\birth}
\end{pspicture}
\;\; \xy {\ar@2{->}^{r} (0,0)*{};(5,0)*{}}; \endxy  \;\;
\begin{pspicture}[0.3](2,1.5)
   \rput(.8,0){\ident}
\end{pspicture}
\\
\psset{xunit=.25cm,yunit=0.25cm}
\begin{pspicture}[0.3](2,3)
   \rput(1,.55){\comult} \rput(.4,.55){\death} \rput(1.6,0){\smallident}
\end{pspicture}
\;\; \xy {\ar@2{->}^{q} (0,0)*{};(5,0)*{}}; \endxy  \;\;
\begin{pspicture}[0.3](2,3)
  \rput(1,0){\longident}
\end{pspicture}
&\quad \mapsto \quad&
\psset{xunit=.25cm,yunit=0.25cm}
\begin{pspicture}[0.3](3,3)
   \rput(2,.55){\mult}
   \rput(.2,.55){\ident}
   \rput(.2,.55){\death}
   \rput(.8,2.05){\zag}
   \rput(2.6,2.05){\medident}
   \rput(2,0){\smallident}
\end{pspicture}
\;\; \xy {\ar@2{->}^{\ell} (0,0)*{};(5,0)*{}}; \endxy  \;\;
 \psset{xunit=.2cm,yunit=0.2cm}
\begin{pspicture}[0.3](3,5.5)
   \rput(2,.55){\mult}
   \rput(1.4,2.05){\curveright}
   \rput(2.6,2.05){\curveright}
   \rput(2,3.55){\smallident}
   \rput(3.2,3.55){\smallident}
   \rput(.2,.55){\ident}
   \rput(.2,2.05){\mult}
   \rput(-.4,3.55){\birth}
   \rput(.8,3.55){\smallident}
   \rput(.2,.55){\death}
   \rput(1.4,4.1){\zag}
   \rput(3.2,4.1){\medident}
   \rput(2,0){\smallident}
\end{pspicture}
\;\;\xy {\ar@2{->}^{\INT} (0,0)*{};(5,0)*{}}; \endxy  \;\;
 \psset{xunit=.2cm,yunit=0.2cm}
\begin{pspicture}[0.3](3,5.5)
   \rput(2,0){\mult}
   \rput(1.4,1.5){\smallident}
   \rput(2.6,1.5){\smallident}
   \rput(1.4,2.05){\curveright}
   \rput(2.6,2.05){\curveright}
   \rput(3.2,3.55){\medident}
   \rput(3.2,4.3){\smallident}
   \rput(.2,2.05){\death}
   \rput(.2,2.05){\mult}
   \rput(-.4,3.55){\medident}
   \rput(-.4,4.3){\birth}
   \rput(1.4,3.55){\zag}
\end{pspicture}
\;\; \xy {\ar@2{->}^{n} (0,0)*{};(5,0)*{}}; \endxy  \;\;
\psset{xunit=.25cm,yunit=0.25cm}
\begin{pspicture}[0.3](2.3,3)
   \rput(1.4,0){\mult}
   \rput(2,1.5){\smallident}
   \rput(.8,1.5){\birth}
\end{pspicture}
\;\; \xy {\ar@2{->}^{\ell} (0,0)*{};(5,0)*{}}; \endxy  \;\;
\begin{pspicture}[0.3](2,1.5)
   \rput(.8,0){\ident}
\end{pspicture}
\\
\psset{xunit=.25cm,yunit=0.25cm}
\begin{pspicture}[0.5](3,3.3)
  \rput(1.6,1.5){\comult} \rput(2.2,0){\comult} \rput(.4,0){\curveright}
\end{pspicture}
 \; \xy {\ar@2{->}^{d} (0,0)*{};(5,0)*{}}; \endxy  \;
\begin{pspicture}[0.5](2.6,3.3)
  \rput(1,1.5){\comult} \rput(2.2,0){\curveleft} \rput(.4,0){\comult}
\end{pspicture}
&\quad \mapsto \quad&
\psset{xunit=.25cm,yunit=0.25cm}
\begin{pspicture}[0.3](5,5)
   \rput(3.2,0){\mult}
   \rput(1.4,0){\ident}
   \rput(2,1.5){\zag}
   \rput(2,3.75){\longzag}
   \rput(.2,2.25){\curveright}
   \rput(.2,1.5){\medident}
   \rput(.2,0){\ident}
   \rput(3.8,1.5){\medident}
   \rput(3.8,2.25){\mult}
   \rput(4.4,3.75){\ident}
\end{pspicture}
\;\; \xy {\ar@2{->}^{w} (0,0)*{};(5,0)*{}}; \endxy  \;\;
\begin{pspicture}[0.3](4,5)
   \rput(2,0){\mult}
   \rput(3.8,0){\ident}
   \rput(1.4,1.5){\medident}
   \rput(3.2,1.5){\zag}
   \rput(1.4,2.25){\mult}
   \rput(.2,3.75){\zag}
   \rput(2,3.75){\medident}
   \rput(.2,2.25){\curveleft}
   \rput(.2,1.5){\medident}
   \rput(.2,0){\ident}
\end{pspicture}
\;\; \xy {\ar@2{->}^{\INT} (0,0)*{};(5,0)*{}}; \endxy  \;\;
\begin{pspicture}[0.3](5,5)
   \rput(2.6,0){\mult}
   \rput(2,1.5){\mult}
   \rput(2.6,3){\medident}
   \rput(2.6,3.75){\medident}
   \rput(4.4,1.5){\curveright}
   \rput(4.4,0){\ident}
   \rput(.8,3.75){\zag}
   \rput(.2,3){\medident}
   \rput(1.4,3){\medident}
   \rput(4.4,3){\zag}
   \rput(3.2,1.5){\curveright}
   \rput(.8,1.5){\curveleft}
   \rput(.8,0){\ident}
\end{pspicture}
\;\; \xy {\ar@2{->}^{a} (0,0)*{};(5,0)*{}}; \endxy  \;\;
\begin{pspicture}[0.3](5,5)
   \rput(2.6,0){\mult}
   \rput(3.2,1.5){\mult}
   \rput(2.6,3){\medident}
   \rput(2.6,3.75){\medident}
   \rput(4.4,1.5){\curveright}
   \rput(4.4,0){\ident}
   \rput(.8,3.75){\zag}
   \rput(4.4,3){\zag}
   \rput(2,1.5){\curveleft}
   \rput(.8,1.5){\curveleft}
   \rput(.8,0){\ident}
   \rput(.2,3){\medident}
   \rput(1.4,3){\medident}
\end{pspicture}
\\
& & \qquad \qquad \;\; \xy {\ar@2{->}^{\INT} (0,0)*{};(5,0)*{}};
\endxy \;\;
\psset{xunit=.25cm,yunit=0.25cm}
\begin{pspicture}[0.3](5,5.5)
   \rput(3.2,0){\mult}
   \rput(1.4,0){\ident}
   \rput(2,1.5){\zag}
   \rput(5,3.75){\zag}
   \rput(5,2.25){\curveright}
   \rput(5,1.5){\medident}
   \rput(5,0){\ident}
   \rput(3.8,1.5){\medident}
   \rput(3.8,2.25){\mult}
   \rput(3.2,3.75){\ident}
\end{pspicture}
\;\; \xy {\ar@2{->}^{w} (0,0)*{};(5,0)*{}}; \endxy  \;\;
\begin{pspicture}[0.3](5,5.5)
   \rput(2.6,0){\mult}
   \rput(.8,0){\ident}
   \rput(1.4,1.5){\zag}
   \rput(3.8,3.75){\zag}
   \rput(3.2,2.25){\ident}
   \rput(5,1.5){\medident}
   \rput(5,0){\ident}
   \rput(3.2,1.5){\medident}
   \rput(5,2.25){\mult}
   \rput(5.6,3.75){\ident}
\end{pspicture}
\end{eqnarray*}
\end{itemize}

\noindent One can check that these maps satisfy the required
relations making $\Theta$ into a strict 2-functor.  This amounts
to checking that all 20 of the axioms are satisfied by the maps in
the image of $\Theta$.  Furthermore, it clear from the definition
that $\Theta$ preserves the monoidal structure on the nose. That
is, $\Theta$ is a strictly monoidal strict 2-functor.

We define the other strict 2-functor $\bar{\Theta} \maps
\mathcal{F}_2 \to \mathcal{F}_1$ defining a 2-equivalence of
2-categories as follows:
\begin{itemize}
    \item On objects $\bar{\Theta}$ is the identity map.
    \item On the generating morphisms $\bar{\Theta}$ is defined as
    follows:
    \begin{eqnarray*}
 \psset{xunit=.5cm,yunit=0.5cm}
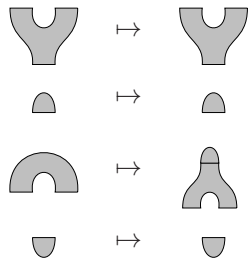
\begin{pspicture}[.5](1,1.5)
 \rput(.5,0){\mult}
\end{pspicture}
 & \quad \mapsto \quad &
  \psset{xunit=.5cm,yunit=0.5cm}
\begin{pspicture}[.5](1,1.5)
 \rput(.5,0){\mult}
\end{pspicture}
\\
 \psset{xunit=.5cm,yunit=0.5cm}
\begin{pspicture}[.5](1,1.5)
 \rput(.5,.5){\birth}
\end{pspicture}
 & \quad \mapsto \quad &
  \psset{xunit=.5cm,yunit=0.5cm}
\begin{pspicture}[.5](1,1.5)
 \rput(.5,.5){\birth }
\end{pspicture}
\\
 \psset{xunit=.5cm,yunit=0.5cm}
\begin{pspicture}[.5](1,1.5)
 \rput(.5,.3){\zag}
\end{pspicture}
 & \quad \mapsto \quad &
  \psset{xunit=.4cm,yunit=0.4cm}
\begin{pspicture}[.5](3.2,2.2)
 \rput(.5,0){\comult}
 \rput(.5,1.5){\birth}
\end{pspicture}
\\
\psset{xunit=.5cm,yunit=0.5cm}
\begin{pspicture}[.5](1,1.5)
 \rput(.5,1){\death}
\end{pspicture}
 & \quad \mapsto \quad &
  \psset{xunit=.5cm,yunit=0.5cm}
\begin{pspicture}[.5](1,1.5)
 \rput(.5,1){\death}
\end{pspicture}
\end{eqnarray*}
    \item On the coherence 2-isomorphisms $\INT$, $\bar{\Theta}$ is
    defined analogously as $\Theta$.  Again, this assignation
    is well-defined by the coherence axioms for a
    $\cat{Gray}$-category.
    \item On the generating 2-isomorphisms $\bar{\Theta}$ is defined as
    follows:
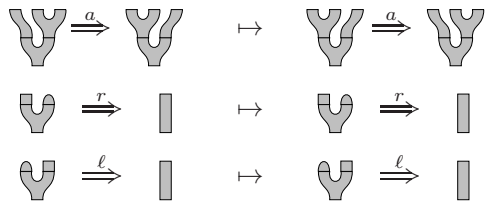
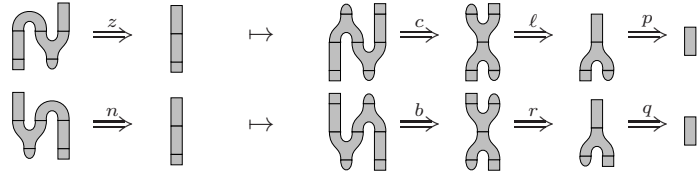
\begin{eqnarray*}
\psset{xunit=.25cm,yunit=0.25cm}
\begin{pspicture}[0.5](2.4,3.3)
  \rput(1,0){\mult} \rput(1.6,1.5){\mult} \rput(.4,1.5){\curveleft}
\end{pspicture}
 \; \xy {\ar@2{->}^{a} (0,0)*{};(5,0)*{}}; \endxy  \quad
\begin{pspicture}[0.5](2.6,3.3)
  \rput(1,0){\mult} \rput(1.6,1.5){\curveright} \rput(.4,1.5){\mult}
\end{pspicture}
&\quad \mapsto \quad&
\psset{xunit=.25cm,yunit=0.25cm}
\begin{pspicture}[0.5](2.6,3.3)
  \rput(1,0){\mult} \rput(1.6,1.5){\curveright} \rput(.4,1.5){\mult}
\end{pspicture}
 \; \xy {\ar@2{->}^{a} (0,0)*{};(5,0)*{}}; \endxy  \quad
\begin{pspicture}[0.5](3,3.3)
  \rput(1,0){\mult} \rput(1.6,1.5){\mult} \rput(.4,1.5){\curveleft}
\end{pspicture}
\\
\psset{xunit=.25cm,yunit=0.25cm}
\begin{pspicture}[0.3](2,3)
  \rput(1,0){\mult} \rput(1.6,1.5){\birth} \rput(.4,1.5){\smallident}
\end{pspicture}
 \quad \xy {\ar@2{->}^{r} (0,0)*{};(5,0)*{}}; \endxy  \quad
\begin{pspicture}[0.3](2,3)
  \rput(1,0){\longident}
\end{pspicture}
&\quad \mapsto \quad&
\psset{xunit=.25cm,yunit=0.25cm}
\begin{pspicture}[0.3](2,3)
  \rput(1,0){\mult} \rput(1.6,1.5){\birth} \rput(.4,1.5){\smallident}
\end{pspicture}
 \quad \xy {\ar@2{->}^{r} (0,0)*{};(5,0)*{}}; \endxy  \quad
\begin{pspicture}[0.3](2,3)
  \rput(1,0){\longident}
\end{pspicture}
\\
\psset{xunit=.25cm,yunit=0.25cm}
\begin{pspicture}[0.3](2,3)
  \rput(1,0){\mult} \rput(1.6,1.5){\smallident} \rput(.4,1.5){\birth}
\end{pspicture}
\quad \xy {\ar@2{->}^{\ell} (0,0)*{};(5,0)*{}}; \endxy  \quad
\begin{pspicture}[0.3](2,3)
  \rput(1,0){\longident}
\end{pspicture}
&\quad \mapsto \quad&
\psset{xunit=.25cm,yunit=0.25cm}
\begin{pspicture}[0.3](2,3)
  \rput(1,0){\mult} \rput(1.6,1.5){\smallident} \rput(.4,1.5){\birth}
\end{pspicture}
\quad \xy {\ar@2{->}^{\ell} (0,0)*{};(5,0)*{}}; \endxy  \quad
\begin{pspicture}[0.3](2,3)
  \rput(1,0){\longident}
\end{pspicture}
\\
\psset{xunit=.25cm,yunit=0.25cm}
\begin{pspicture}[.5](3,3)
  \rput(2.2,.55){\mult} \rput(1,2.05){\zag}
  \rput(2.8,2.05){\ident} \rput(.4,0){\smallident} \rput(.4,.55){\ident}
  \rput(2.2,.55){\death}
\end{pspicture}
\quad \xy {\ar@2{->}^{z} (0,0)*{};(5,0)*{}}; \endxy  \quad
\begin{pspicture}[0.5](2,3.3)
  \rput(1,2.05){\ident} \rput(1,0){\smallident} \rput(1,.55){\ident}
\end{pspicture}
&\quad \mapsto \quad&
\psset{xunit=.25cm,yunit=0.25cm}
\begin{pspicture}[0.5](3.1,4.2)
  \rput(2.2,.55){\mult}
  \rput(1,2.05){\comult}
  \rput(1,3.55){\birth}
  \rput(2.8,2.05){\ident}
  \rput(2.8,3.55){\smallident}
  \rput(.4,.55){\ident}
  \rput(.4,0){\smallident}
  \rput(2.2,.55){\death}
\end{pspicture}
\;\; \xy {\ar@2{->}^{c} (0,0)*{};(5,0)*{}}; \endxy  \;\;
\begin{pspicture}[0.5](2.5,4.2)
  \rput(1.6,2.05){\mult}
  \rput(1.6,.555){\comult}
  \rput(1,3.55){\birth}
  \rput(2.2,3.55){\smallident}
  \rput(1,0){\smallident}
  \rput(2.2,.55){\death}
\end{pspicture}
\;\; \xy {\ar@2{->}^{\ell} (0,0)*{};(5,0)*{}}; \endxy  \;\;
\begin{pspicture}[0.5](2.5,4.2)
  \rput(1.6,2.05){\ident}
  \rput(1.6,.555){\comult}
  \rput(1,0){\smallident}
  \rput(2.2,.55){\death}
\end{pspicture}
\;\; \xy {\ar@2{->}^{p} (0,0)*{};(5,0)*{}}; \endxy  \;\;
\begin{pspicture}[0.5](1.2,1.5)
  \rput(.5,0){\ident}
\end{pspicture}
\\
\psset{xunit=.25cm,yunit=0.25cm}
\begin{pspicture}[0.5](3,3)
  \rput(1,.55){\mult} \rput(2.2,2.05){\zag}
  \rput(.4,2.05){\ident} \rput(2.8,.55){\ident}
  \rput(2.8,0){\smallident} \rput(1,.55){\death}
\end{pspicture}
\quad \xy {\ar@2{->}^{n} (0,0)*{};(5,0)*{}}; \endxy  \quad
\begin{pspicture}[0.5](2,3.6)
  \rput(1,2.05){\ident} \rput(1,.55){\ident} \rput(1,0){\smallident}
\end{pspicture}
&\quad \mapsto \quad&
\psset{xunit=.25cm,yunit=0.25cm}
\begin{pspicture}[0.5](3.1,4.2)
  \rput(1,.55){\mult}
  \rput(2.2,2.05){\comult}
  \rput(2.2,3.55){\birth}
  \rput(.4,2.05){\ident}
  \rput(.4,3.55){\smallident}
  \rput(2.8,.55){\ident}
  \rput(2.8,0){\smallident}
  \rput(1,.55){\death}
\end{pspicture}
\;\; \xy {\ar@2{->}^{b} (0,0)*{};(5,0)*{}}; \endxy  \;\;
\begin{pspicture}[0.5](2.5,4.2)
  \rput(1.6,2.05){\mult}
  \rput(1.6,.555){\comult}
  \rput(2.2,3.55){\birth}
  \rput(1,3.55){\smallident}
  \rput(2.2,0){\smallident}
  \rput(1,.55){\death}
\end{pspicture}
\;\; \xy {\ar@2{->}^{r} (0,0)*{};(5,0)*{}}; \endxy  \;\;
\begin{pspicture}[0.5](2.5,3.8)
  \rput(1.6,2.05){\ident}
  \rput(1.6,.555){\comult}
  \rput(2.2,0){\smallident}
  \rput(1,.55){\death}
\end{pspicture}
\;\; \xy {\ar@2{->}^{q} (0,0)*{};(5,0)*{}}; \endxy  \;\;
\begin{pspicture}[0.5](1.2,1.5)
  \rput(.5,0){\ident}
\end{pspicture}
\end{eqnarray*}
\end{itemize}

\noindent Again, it is a routine and laborious calculation to
check that these maps satisfy the required relations making
$\bar{\Theta}$ into a strict 2-functor.  By construction, it is
clear that $\bar{\Theta}$ strictly preserves the monoidal
structure.

To see that $\Theta$ and $\bar{\Theta}$ define a 2-equivalence of
2-categories we must check that their composites are 2-naturally
isomorphic to the identity.  Since both 2-functors are the
identity on objects it is clear that on objects the composites are
actually equal to the identity map.  One can check that on any
generating morphisms, say $X$, the image of $X$ under the
composites of the above 2-functors is naturally isomorphic to $X$.
Hence, the strict 2-categories $\mathcal{F}_1$ and $\mathcal{F}_2$
are 2-equivalent.  Furthermore, since $\Theta$ and $\bar{\Theta}$
are strict monoidal functors it is clear that $\mathcal{F}_1$ and
$\mathcal{F}_2$ are monoidally 2-equivalent. Regarding
$\mathcal{F}_1$ and $\mathcal{F}_2$ as one object
$\cat{Gray}$-categories, we see then that $\Theta$ and
$\bar{\Theta}$ define a 3-equivalence of $\cat{Gray}$-categories.
\qed

The observant reader will have noticed that the isomorphisms
defining the image of the generators in the above maps are
precisely the proofs of the equivalent definitions of a Frobenius
algebra with the equalities replaced by coherent isomorphisms. In
the next section we will show that these two 2-categories are
actually monoidally 2-equivalent to the walking pseudo Frobenius
algebra.  This means that  an object of a semistrict monoidal
2-category equipped with the morphisms and 2-morphisms satisfying
the axioms of either of the above lemmas serves as an equivalent
definition of pseduo Frobenius algebra. Hence, categorifying any
of the equivalent characterizations from Proposition
\ref{equivFrob} produces monoidally 2-equivalent 2-categories.

\subsection{Proof of main theorem} \label{secPROOFS}

In this section we prove the main result of this paper.  We have
shown in the previous section that the monoidal bicategory of
3-dimensional thick tangles is triequivalent to the semistrict
monoidal category $\mathcal{F}_2$ defined in Lemma \ref{lemma2}.
We now show that the walking pseudo Frobenius algebra is
3-equivalent to $\mathcal{F}_2$ and the main result will follow.

\paragraph{Proof of Theorem \ref{MAINthm}. }
To prove that the walking pseudo Frobenius algebra is
3-equivalent to the semistrict monoidal category $\mathcal{F}_2$
we will use the description of the walking pseudo Frobenius
algebra given in Proposition~\ref{FFR1}.   Notice that the walking
pseudo Frobenius algebra and the semistrict monoidal 2-category
$\mathcal{F}_2$ are both generated by the same objects, morphisms,
and 2-morphisms.  Thus, it suffices to prove that the axioms are
satisfied in both definitions.

\begin{figure}
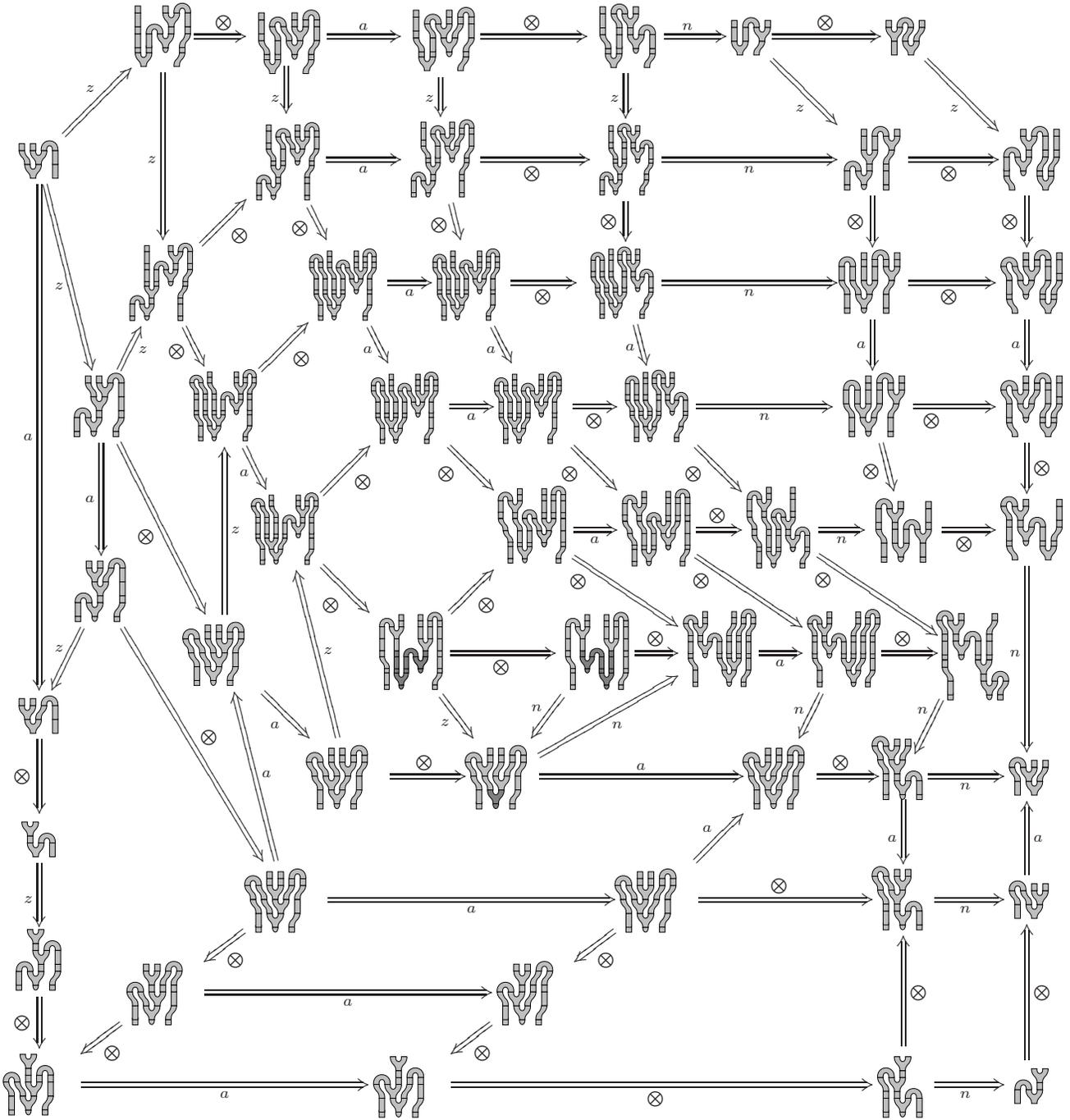
 \label{Huge}
\caption{Proof that $\wFrob$ satisfies the relations of the
$\mathcal{F}_2$.}
\begin{center}
 \HUGEproofII
\end{center}
\end{figure}

To prove that the walking pseudo Frobenius algebra satisfies the
axioms for the generators of the semistrict monoidal 2-category
$\mathcal{F}_2$, we only need to check the axiom for the map $w$
since the other coherence axiom is included in the description of
$\wFrob$ given in Proposition~\ref{FFR1}.  The proof is given in
Figure \ref{Huge}.
The outer rectangle of this diagram is the axiom for the
isomorphism $w$ in the definition of $\mathcal{F}_2$ with $w$
expanded out. All of the outer squares commute by the properties
of a semistrict monoidal 2-category. The innermost triangle with
some of 1-morphisms drawn slightly darker is the first coherence
condition from Proposition \ref{FFR1}. Just below this triangle is
a rather distorted rectangle that is the coherence for
associativity in a pseudomonoid.  Note that some of the arrows
labelled by $\INT$ may mean multiple applications of this
2-morphism.

To show that the semistrict monoidal 2-category $\mathcal{F}_2$
satisfies the axioms of the walking pseudo Frobenius algebra we
must show that the first coherence axiom in Proposition \ref{FFR1}
is satisfied. To simplify the proof, we use the 3-equivalence of
monoidal strict 2-categories $\mathcal{F}_1$ and $\mathcal{F}_2$
from Lemma \ref{lemma2}.  With the twenty axioms from the
generators of $\mathcal{F}_1$ at our disposal, the proof becomes
much simpler.  The proof is shown below:

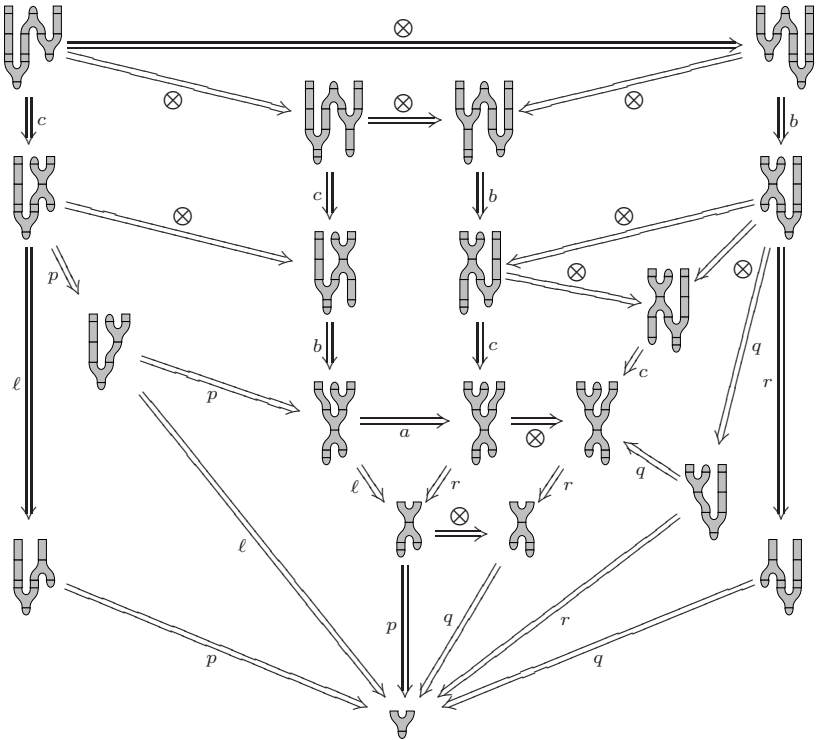
\BIGproof

Since the walking pseudo Frobenius algebra has the same monoidal
structure as $\mathcal{F}_2$ it is clear that the suspension of
the walking pseudo Frobenius algebra and the suspension of
$\mathcal{F}_2$ are 3-equivalent $\cat{Gray}$-categories.

\qed

\section{Conclusion} \label{secConc}

We used the description of the walking adjunction to understand
the walking Frobenius algebra and its relation to 2-dimensional
thick tangles. $\cat{2Thick}$ while not optimal from the
perspective of studying all open string worldsheets has the
advantage that it arises quite naturally from higher-dimensional
category theory. Our success in algebraically characterizing
$\cat{2Thick}$ suggests that the full category of open strings and
their worldsheets might be algebraically characterized by studying
more interesting instances of categorical adjunctions.

An adjunction is an intrinsically categorical concept and since
this notion has already been generalized to the context of
\cat{Gray}-categories we were able to use the relationship between
adjunctions and Frobenius algebras to categorify the notion of a
Frobenius algebra.  The description of the walking pseudo
Frobenius algebra in terms of the walking pseudo ambidextrous
adjunction also allowed us to see the relationship of pseudo
Frobenius algebras to 3-dimensional thick tangles using string
diagrams.  These results suggest that the general machinery of
adjunctions in $n$-categories may prove useful in studying string
membranes and other interesting cobordism categories. In a future
work we will explain how the topology of arbitrary string
membranes (rather than those defined by diffeomorphisms) can be
algebraically described using a generalization of what we have
described in this paper.

\subsubsection*{Acknowledgements}

The author would like to thank John Baez, Alissa Crans, Richard
Garner, Martin Hyland, Steve Lack, Michael Shulman, and Ross
Street for their many discussions and helpful advice.  Thanks also
to David Kagan for helpful correspondence on string theory, and
Alison Waldron for her contribution to the proof of
Theorem~\ref{MAINthm}.  The author is also grateful for the
generosity of the University of Chicago where this work was
completed.

\end{document}